\documentclass[11pt]{article}

\usepackage{enumerate,xspace}
\usepackage{amsmath,amssymb,wasysym}
\usepackage[all]{xy}
\usepackage{proof}
\usepackage[svgnames]{xcolor}
\usepackage{tikz}

\usepackage{stmaryrd} %need for special interpretation bracket
\usepackage{mathtools}
\usepackage{latexsym}
\usepackage{dsfont}
\usepackage{multicol}
\usepackage{cmll}
\usepackage{multirow}
\usepackage{longtable}

\usepackage{lscape}
\usepackage{array}

\usepackage{shuffle}
\DeclareFontFamily{U}{shuffle}{}
\DeclareFontShape{U}{shuffle}{m}{n}{ <-8>shuffle7 <8->shuffle10}{}

\newcommand\abs[1]{\left|#1\right|}

\usepackage{hyperref} %hyperlink package
\hypersetup{
    colorlinks,
    citecolor=red,
    filecolor=red,
    linkcolor=blue,
    urlcolor=red
}

\usetikzlibrary{shapes,snakes}
\pgfdeclarelayer{edgelayer}
\pgfdeclarelayer{nodelayer}
\pgfsetlayers{edgelayer,nodelayer,main}

\definecolor{hexcolor0xff0000}{rgb}{1.000,0.000,0.000}
\definecolor{hexcolor0x000000}{rgb}{0.000,0.000,0.000}
\definecolor{hexcolor0x00ff00}{rgb}{0.000,1.000,0.000}
\definecolor{hexcolor0xffff00}{rgb}{1.000,1.000,0.000}
\definecolor{hexcolor0x000000}{rgb}{0.000,0.000,0.000}
\definecolor{hexcolor0x000000}{rgb}{0.000,0.000,0.000}
\definecolor{white}{rgb}{1.000,1.000,1.000}

\tikzstyle{none}=[inner sep=0pt]
\tikzstyle{port}=[inner sep=0pt]
\tikzstyle{component}=[circle,fill=white,draw=black]
\tikzstyle{integral}=[inner sep=0pt]
\tikzstyle{differential}=[inner sep=0pt]
\tikzstyle{codifferential}=[inner sep=0pt]
\tikzstyle{function}=[regular polygon,regular polygon sides=4,fill=white,draw=black]
\tikzstyle{function2}=[regular polygon,regular polygon sides=4,fill=white,draw=black, inner sep=1pt]
\tikzstyle{duplicate}=[circle,fill=white,draw=black, inner sep=1pt]
\tikzstyle{wire}=[-,draw=black,line width=1.000]
\tikzstyle{object}=[inner sep=2pt]

\newtheorem{observation}{Remark}[section]
\newtheorem{lemma}[observation]{Lemma}  %%share counter with remark
\newtheorem{theorem}[observation]{Theorem}
\newtheorem{definition}[observation]{Definition}
\newtheorem{example}[observation]{Example}

\newtheorem{proposition}[observation]{Proposition} 
\newtheorem{corollary}[observation]{Corollary}

\newcommand{\proof}{\noindent{\sc Proof:}\xspace}
\def\endproof{~\hfill$\Box$\vskip 10pt}

\makeatletter

\newdimen\w@dth

\def\setw@dth#1#2{\setbox\z@\hbox{\scriptsize $#1$}\w@dth=\wd\z@
\setbox\@ne\hbox{\scriptsize $#2$}\ifnum\w@dth<\wd\@ne \w@dth=\wd\@ne \fi
\advance\w@dth by 1.2em}

\def\t@^#1_#2{\allowbreak\def\n@one{#1}\def\n@two{#2}\mathrel
{\setw@dth{#1}{#2}
\mathop{\hbox to \w@dth{\rightarrowfill}}\limits
\ifx\n@one\empty\else ^{\box\z@}\fi
\ifx\n@two\empty\else _{\box\@ne}\fi}}
\def\t@@^#1{\@ifnextchar_ {\t@^{#1}}{\t@^{#1}_{}}}

\def\t@left^#1_#2{\def\n@one{#1}\def\n@two{#2}\mathrel{\setw@dth{#1}{#2}
\mathop{\hbox to \w@dth{\leftarrowfill}}\limits
\ifx\n@one\empty\else ^{\box\z@}\fi
\ifx\n@two\empty\else _{\box\@ne}\fi}}
\def\t@@left^#1{\@ifnextchar_ {\t@left^{#1}}{\t@left^{#1}_{}}}

\def\two@^#1_#2{\def\n@one{#1}\def\n@two{#2}\mathrel{\setw@dth{#1}{#2}
\mathop{\vcenter{\hbox to \w@dth{\rightarrowfill}\kern-1.7ex
                 \hbox to \w@dth{\rightarrowfill}}%
       }\limits
\ifx\n@one\empty\else ^{\box\z@}\fi
\ifx\n@two\empty\else _{\box\@ne}\fi}}
\def\tw@@^#1{\@ifnextchar_ {\two@^{#1}}{\two@^{#1}_{}}}

\def\tofr@^#1_#2{\def\n@one{#1}\def\n@two{#2}\mathrel{\setw@dth{#1}{#2}
\mathop{\vcenter{\hbox to \w@dth{\rightarrowfill}\kern-1.7ex
                 \hbox to \w@dth{\leftarrowfill}}%
       }\limits
\ifx\n@one\empty\else ^{\box\z@}\fi
\ifx\n@two\empty\else _{\box\@ne}\fi}}
\def\t@fr@^#1{\@ifnextchar_ {\tofr@^{#1}}{\tofr@^{#1}_{}}}

\newdimen\W@dth
\def\setW@dth#1#2{\setbox\z@\hbox{$#1$}\W@dth=\wd\z@
\setbox\@ne\hbox{$#2$}\ifnum\W@dth<\wd\@ne \W@dth=\wd\@ne \fi
\advance\W@dth by 1.2em}

\def\T@^#1_#2{\allowbreak\def\N@one{#1}\def\N@two{#2}\mathrel
{\setW@dth{#1}{#2}
\mathop{\hbox to \W@dth{\rightarrowfill}}\limits
\ifx\N@one\empty\else ^{\box\z@}\fi
\ifx\N@two\empty\else _{\box\@ne}\fi}}
\def\T@@^#1{\@ifnextchar_ {\T@^{#1}}{\T@^{#1}_{}}}

\def\T@left^#1_#2{\def\N@one{#1}\def\N@two{#2}\mathrel{\setW@dth{#1}{#2}
\mathop{\hbox to \W@dth{\leftarrowfill}}\limits
\ifx\N@one\empty\else ^{\box\z@}\fi
\ifx\N@two\empty\else _{\box\@ne}\fi}}
\def\T@@left^#1{\@ifnextchar_ {\T@left^{#1}}{\T@left^{#1}_{}}}

\def\Tofr@^#1_#2{\def\N@one{#1}\def\N@two{#2}\mathrel{\setW@dth{#1}{#2}
\mathop{\vcenter{\hbox to \W@dth{\rightarrowfill}\kern-1.7ex
                 \hbox to \W@dth{\leftarrowfill}}%
       }\limits
\ifx\N@one\empty\else ^{\box\z@}\fi
\ifx\N@two\empty\else _{\box\@ne}\fi}}
\def\T@fr@^#1{\@ifnextchar_ {\Tofr@^{#1}}{\Tofr@^{#1}_{}}}

\def\Two@^#1_#2{\def\N@one{#1}\def\N@two{#2}\mathrel{\setW@dth{#1}{#2}
\mathop{\vcenter{\hbox to \W@dth{\rightarrowfill}\kern-1.7ex
                 \hbox to \W@dth{\rightarrowfill}}%
       }\limits
\ifx\N@one\empty\else ^{\box\z@}\fi
\ifx\N@two\empty\else _{\box\@ne}\fi}}
\def\Tw@@^#1{\@ifnextchar_ {\Two@^{#1}}{\Two@^{#1}_{}}}

\def\to{\@ifnextchar^ {\t@@}{\t@@^{}}}
\def\from{\@ifnextchar^ {\t@@left}{\t@@left^{}}}
\def\tofro{\@ifnextchar^ {\t@fr@}{\t@fr@^{}}}
\def\To{\@ifnextchar^ {\T@@}{\T@@^{}}}
\def\From{\@ifnextchar^ {\T@@left}{\T@@left^{}}}
\def\Two{\@ifnextchar^ {\Tw@@}{\Tw@@^{}}}
\def\Tofro{\@ifnextchar^ {\T@fr@}{\T@fr@^{}}}

\makeatother

\title{Differential Categories Revisited}
\author{R.F. Blute \thanks{Department of Mathematics and Statistics, University of Ottawa, rblute@uottawa.ca. Research supported by NSERC.} \quad \quad J.R.B. Cockett \thanks{Department of Computer Science, University of Calgary, robin@ucalgary.ca. Research supported by NSERC.} \quad \quad J-S. P. Lemay \thanks{Department of Computer Science, University of Oxford, jean-simon.lemay@kellogg.ox.ac.uk. Research supported by Kellogg College, the Clarendon Fund, and the Oxford-Google DeepMind Graduate Scholarship.} \quad \quad R.A.G. Seely  \thanks{Department of Mathematics, McGill University, robert.seely@mcgill.ca. Research supported by FRQNT.}} 
\begin{document}
\allowdisplaybreaks

\maketitle
%\section{}
%\subsection{}
\vspace{-1cm}

\begin{abstract}
Differential categories were introduced to provide a minimal categorical doctrine for differential linear logic. Here we revisit the formalism and, in particular, examine the two different approaches to defining differentiation which were introduced. The basic approach used a deriving transformation, while a more refined approach, in the presence of a bialgebra modality, used a codereliction. The latter approach is particularly relevant to linear logic settings, where the coalgebra modality is monoidal and the Seely isomorphisms give rise to a bialgebra modality.  Here, we prove that these apparently distinct notions of differentiation, in the presence of a monoidal coalgebra modality, are completely equivalent. Thus, for linear logic settings, there is only one notion of differentiation.

This paper also presents a number of separating examples for coalgebra modalities including examples which are and are not monoidal, as well as examples which do and do not support differential structure.  Of particular interest is the observation that -- somewhat counter-intuitively -- differential algebras never induce a differential category although they provide a monoidal coalgebra modality.  On the other hand, 
Rota-Baxter algebras -- which are usually associated with integration -- provide an example of a differential category which has a non-monoidal coalgebra modality.  
\end{abstract}

\paragraph{Acknowledgements:} The authors would like to thank the anonymous referee for very helpful and constructive comments in their review, especially regarding the overall structure of this paper.

\newpage 
\tableofcontents

%%%%%%%%%%%%%%%%%%%%%%%%%%%%%%%%%%%%%%%%%%%%%%%%%%%%%%%%%%%%%%%%%%%%%

\section{Introduction}

%%%%%%%%%%%%%%%%%%%%%%%%%%%%%%%%%%%%%%%%%%%%%%%%%%%%%%%%%%%%%%%%%%%%%

Differential linear logic \cite{ehrhard2003differential,ehrhard2006differential} introduced a syntactic proof-theoretic approach to differential calculus and, subsequently, differential categories \cite{blute2006differential} were developed to provide a categorical counterpart for these ideas.   In this categorical approach, two methods for defining the differentiation were introduced based on, respectively, a deriving transformation and a codereliction (Definition \ref{cordef}).  In \cite{fiore2007differential} Fiore proposed an axiomatization for a  deriving transformation, which he called a ``creation operator'', satisfying  additional ``strength laws''.  These laws were very natural to impose in the presence of a monoidal coalgebra modality and finite biproducts.  However, this introduced an apparently stronger notion of differentiation and left open the question of whether a creation operator was a distinct notion of differentiation.  An initial purpose of this paper is to provide a proof that, in the presence of a monoidal coalgebra modality, all these methods of defining differentiation are in fact, equivalent.  Thus, there is only one notion of differentiation in linear logic.

In \cite{fiore2007differential}, Fiore made another interesting observation, namely, that it is much more convenient to work in a setting with finite biproducts.  Furthermore, as one can always complete an additive category to have finite biproducts, Fiore argued that one may as well work in a setting with finite biproducts. Of course, the notion of differentiation in additive categories is tightly coupled to having a coalgebra modality. Thus, when one completes a differential category to have biproducts, one needs to show that this modality can also be extended to the biproduct completion.  It is standard that one can extend a {\em monoidal\/} coalgebra modality to the biproduct completion using the Seely isomorphisms \cite{seely1987linear} and Fiore's work was focused on this case.  However, when the coalgebra modality is {\em not\/} monoidal, there is no reason why the coalgebra modality should extend to the biproduct completion. In fact, the deriving transformations of \cite{blute2006differential} had not assumed that the coalgebra modality was monoidal.  Thus, if one is to entertain these more basic coalgebra modalities, one must be cautious about assuming the presence of biproducts.  A significant aspect of this paper is to point out that there are many natural examples of differential categories in which the coalgebra modality is {\em not\/} monoidal, thus exclusively concentrating on monoidal modalities misses an important part of the geography of differentiation (see the Venn diagram in Section \ref{examples}).

This paper revisits the ``original'' definition of a differential category found in \cite{blute2006differential}. This relies on the idea of a coalgebra modality and a deriving transformation. Familiarity with linear logic may tempt one to think that this is the ``exponential'' modality of linear logic but it is not. It is a strictly more general notion as the modality is not assumed monoidal.  There are many familiar and important examples of differential categories based on a monoidal coalgebra modality -- such as (the opposite of) the free symmetric algebra monad on vector spaces.  However, examples of differential categories based on a modality which is not monoidal are less familiar.  Two compelling examples are smooth functions via the free $\mathcal{C}^\infty$-ring over vector spaces (as mentioned in \cite{blute2006differential}) and the free Rota-Baxter algebra over modules (see \cite{zhang2016monads}), which we prove in Proposition \ref{RBderiving} provides a differential category. 

Here we refer to additive categories with a monoidal coalgebra modality as ``additive linear categories''. The biproduct completion of an additive linear category is then an ``additive monoidal storage category''.  An additive monoidal storage category (which has biproducts) is always an additive linear category (which need not have biproducts) and both always have a monoidal bialgebra modality.  Here, to facilitate the proofs, it is convenient to work with a further intermediate notion we called an ``additive bialgebra modality'': this is a bialgebra modality which has an additional coherence requirement between the additive and bialgebra structure. It is with respect to additive bialgebra modalities that we prove that deriving transformations and coderelictions are equivalent.  

Coderelictions always give deriving transformations, and it was shown, in \cite{blute2006differential} that a deriving transformation for a bialgebra modality is equivalent to a codereliction if and only if the deriving transformation satisfies the $\nabla$-rule. This latter rule was originally thought to be a completely independent requirement. The key new observation is that, for an additive bialgebra modality with a deriving transformation, the $\nabla$-rule is in fact implied. More specifically, while a deriving transformation is assumed to satisfy five rules \textbf{[d.1]} to \textbf{[d.5]} (which include the product rule and the linear rule), for an additive bialgebra modality we prove that, in the presence of the other rules, the $\nabla$-rule is equivalent to the product rule.  Furthermore, when an additive symmetric monoidal category has a monoidal coalgebra modality it is straightforward to show that the modality is an additive bialgebra modality.  Thus for additive linear categories: deriving transformations and coderelictions are equivalent.  

Clearly, additive bialgebra modalities and monoidal coalgebra modalities are closely related. In particular, additive bialgebra modalities can always be extended to the biproduct completion and, furthermore, this biproduct completion has Seely isomorphisms. Thus, additive bialgebra modalities correspond to monoidal storage categories \cite{blute2015cartesian} (also called new 
Seely categories \cite{bierman1995categorical,mellies2003categorical}) having the Seely isomorphisms. However, it is well-known (as the name suggests) that monoidal storage categories have a monoidal coalgebra modality!  Thus, additive bialgebra modalities are, in fact, monoidal coalgebra modalities. This argument provides an abstract proof of the equivalence of the two notions which relies on results dispersed across a number of papers. In order to bring the result into focus for this paper we provide a direct proof (see Appendices \ref{montoadd} and \ref{addtomon}).   We make no claim that the resulting proof is more elegant or shorter: it simply has the merit of collecting a complete demonstration of this equivalence under one roof.

This allows us to complete the first purpose of the paper by observing that, in an additive linear category, coderelictions and deriving transformations always satisfy the ``strength laws''.  Putting all this together one concludes that deriving transformations and creation operators are, in additive linear categories, completely equivalent.
 
 The final section of the paper, Section \ref{examples}, provides separating examples for the categorical structures we have introduced.   Of particular interest is the example of the free differential algebra modality on a module category, which we treat in some detail. It is of particular interest as -- contrary perhaps to expectations -- it is (the dual) of an example of an additive bialgebra modality which {\em does not} admit a deriving transformation.  Furthermore, it does not admit a deriving transformation in the strongest possible sense: assuming that there is a deriving transformation forces the ring to be trivial.  As far as we know these observations are new.  Another interesting example is the free Rota-Baxter algebra modality on a module category: as mentioned above, it is an example of a bialgebra modality which is not additive yet admits a deriving transformation. Again as far as we know this has not been presented in full detail before.
 
 \paragraph{Conventions and the Graphical Calculus:} We shall use diagrammatic order for composition: the composite map $fg$ is the map which first does $f$ then $g$.  Furthermore, to simplify working in symmetric monoidal categories, we will allow ourselves to work in strict symmetric monoidal categories and so will generally suppress the associator and unitor isomorphisms.  For a symmetric monoidal category, we will use $\otimes$ for the tensor product, $K$ for the unit, and $\sigma: A \otimes B \to B \otimes A$ for the symmetry isomorphism. Throughout this paper we shall make extensive use of the graphical calculus \cite{joyal1991geometry} for symmetric monoidal categories as this makes proofs easier to follow. Note that our diagrams are to be read down the page -- from top to bottom -- and we shall often omit labeling wires with objects. We refer the reader to \cite{selinger2010survey} for an introduction to the graphical calculus in monoidal categories and its variations, and to \cite{blute2006differential} for the graphical calculus of a differential category. We will be working with coalgebra modalities which in particular involves an endofunctor $\oc$, and so as in \cite{blute2006differential} we will use functor boxes when dealing with string diagrams involving the endofunctor. A mere map $f: A \to B$ will be encased in a circle while $\oc(f): \oc A \to \oc B$ will be encased in a box: 
 \begin{align*}
 \begin{array}[c]{c}
f
   \end{array}=
 \begin{array}[c]{c}\resizebox{!}{2cm}{%
\begin{tikzpicture}
	\begin{pgfonlayer}{nodelayer}
		\node [style=circle] (0) at (0, 2.25) {$A$};
		\node [style=circle] (1) at (0, -0.25) {$B$};
		\node [style={component}] (2) at (0, 1) {$f$};
	\end{pgfonlayer}
	\begin{pgfonlayer}{edgelayer}
		\draw [style=wire] (0) to (2);
		\draw [style=wire] (2) to (1);
	\end{pgfonlayer}
\end{tikzpicture}
  }% 
\end{array}
 &&  \begin{array}[c]{c}
\oc(f)
   \end{array}=
   \begin{array}[c]{c}\resizebox{!}{2cm}{%
   \begin{tikzpicture}
	\begin{pgfonlayer}{nodelayer}
		\node [style=circle] (0) at (0, 2.25) {$\oc A$};
		\node [style=circle] (1) at (0, -0.25) {$\oc B$};
		\node [style={function}] (2) at (0, 1) {$f$};
	\end{pgfonlayer}
	\begin{pgfonlayer}{edgelayer}
		\draw [style=wire] (0) to (2);
		\draw [style=wire] (2) to (1);
	\end{pgfonlayer}
\end{tikzpicture}
}%
   \end{array}
\end{align*}

%%%%%%%%%%%%%%%%%%%%%%%%%%%%%%%%%%%%%%%%%%%%%%%%%%%%%%%%%%%%%%%%%%%%%%%

\section{Coalgebra Modalities}\label{coalgsec}

%%%%%%%%%%%%%%%%%%%%%%%%%%%%%%%%%%%%%%%%%%%%%%%%%%%%%%%%%%%%%%%%%%%%%%%

In this section we review coalgebra modalities, monoidal coalgebra modalities, and the Seely isomorphisms. Examples can be found in Section \ref{examples}. First, if only to introduce notation and provide a simple graphical calculus example, recall that a \textbf{comonad} on a category $\mathbb{X}$ is a triple $(\oc, \delta, \varepsilon)$ consisting of an endofunctor ${\oc: \mathbb{X} \to \mathbb{X}}$ and two natural transformations $\delta: \oc A \to \oc \oc A$ and ${\varepsilon: \oc A \to A}$ such that the following diagrams commute:
\[\xymatrixcolsep{5pc}\xymatrix{ 
        \oc A  \ar[r]^-{\delta} \ar[d]_-{\delta} \ar@{=}[dr]^-{}& \oc \oc A \ar[d]^-{\varepsilon}  & \oc A  \ar[r]^-{\delta} \ar[d]_-{\delta} & \oc \oc A  \ar[d]^-{\delta}\\
        \oc \oc A \ar[r]_-{\oc(\varepsilon)} & \oc A  & \oc \oc A \ar[r]_-{\oc(\delta)} & \oc \oc \oc A}\]
in the graphical calculus, these equalities are drawn as: 
\begin{equation}\label{comonad}\begin{gathered} 
\resizebox{!}{2cm}{%
\begin{tikzpicture}
	\begin{pgfonlayer}{nodelayer}
		\node [style=component] (0) at (1.25, 0) {$\delta$};
		\node [style=component] (1) at (1.25, -1.25) {$\varepsilon$};
		\node [style=port] (2) at (1.25, -2) {};
		\node [style=port] (3) at (1.25, 0.75) {};
		\node [style=port] (4) at (1.75, -0.5) {$=$};
		\node [style=port] (5) at (2.25, 0.75) {};
		\node [style=port] (6) at (2.25, -2) {};
		\node [style=port] (7) at (2.75, -0.5) {$=$};
		\node [style=component] (8) at (3.25, 0) {$\delta$};
		\node [style=function] (9) at (3.25, -1.25) {$\varepsilon$};
		\node [style=port] (10) at (3.25, -2) {};
		\node [style=port] (11) at (3.25, 0.75) {};
		\node [style=component] (13) at (5.75, 0) {$\delta$};
		\node [style=component] (14) at (5.75, -1.25) {$\delta$};
		\node [style=port] (15) at (5.75, -2) {};
		\node [style=port] (16) at (5.75, 0.75) {};
		\node [style=port] (17) at (6.5, -0.5) {$=$};
		\node [style=component] (18) at (7.25, 0) {$\delta$};
		\node [style=function2] (19) at (7.25, -1.25) {$\delta$};
		\node [style=port] (20) at (7.25, -2) {};
		\node [style=port] (21) at (7.25, 0.75) {};
	\end{pgfonlayer}
	\begin{pgfonlayer}{edgelayer}
		\draw [style=wire] (0) to (1);
		\draw [style=wire] (1) to (2);
		\draw [style=wire] (3) to (0);
		\draw [style=wire] (5) to (6);
		\draw [style=wire] (8) to (9);
		\draw [style=wire] (9) to (10);
		\draw [style=wire] (11) to (8);
		\draw [style=wire] (13) to (14);
		\draw [style=wire] (14) to (15);
		\draw [style=wire] (16) to (13);
		\draw [style=wire] (18) to (19);
		\draw [style=wire] (19) to (20);
		\draw [style=wire] (21) to (18);
	\end{pgfonlayer}
\end{tikzpicture}
}%
\end{gathered}\end{equation}
               
Coalgebra modalities are comonads such that every cofree $\oc$-coalgebra comes equipped with a natural cocommutative comonoid structure. 

\begin{definition}\label{coalgdef} A \textbf{coalgebra modality} \cite{blute2006differential} on a symmetric monoidal category is a quintuple $(\oc, \delta, \varepsilon, \Delta, e)$ consisting of a comonad $(\oc, \delta, \varepsilon)$ and two natural transformations $\Delta: \oc A \to \oc A \otimes \oc A$ and $e: \oc A \to K$ such that $(\oc A, \Delta, e)$ is a cocommutative comonoid, that is, the following diagrams commute: 
\begin{equation}\label{cocomdiagrams}\begin{gathered} \xymatrixcolsep{2.5pc}\xymatrix{\oc A  \ar[r]^-{\Delta} \ar[d]_-{\Delta} & \oc A \otimes \oc A \ar[d]^-{\Delta \otimes 1}& & \oc A \ar[d]^-{\Delta} \ar@{=}[dr]\ar@{=}[dl] &&      \oc A   \ar[r]^-{\Delta}  \ar[dr]_-{\Delta} & \oc A \otimes \oc A \ar[d]^-{\sigma}\\
      \oc A \otimes \oc A \ar[r]_-{1 \otimes \Delta} & \oc A \otimes \oc A \otimes \oc A& \oc A & \oc A \otimes \oc A \ar[l]^-{e \otimes 1} \ar[r]_-{1 \otimes e} & \oc A& &  \oc A \otimes \oc A} \end{gathered}\end{equation}
      and $\delta$ preserves the comultiplication, that is, the following diagram commutes: 
      \begin{equation}\label{}\begin{gathered}\xymatrixcolsep{5pc}\xymatrix{\oc A \ar[r]^-{\Delta}  \ar[d]_-{\delta} & \oc A \otimes \oc A \ar[d]^-{\delta \otimes \delta} \\
      \oc \oc A \ar[r]_-{\Delta} & \oc \oc A \otimes \oc \oc A 
  } \end{gathered}\end{equation}
 \end{definition}
 
In the graphical calculus, the coalgebra modality identities are drawn as follows: 
\begin{equation}\label{coalgeq}\begin{gathered}   
\resizebox{!}{2cm}{%
\begin{tikzpicture}
	\begin{pgfonlayer}{nodelayer}
		\node [style=port] (0) at (1, 3) {};
		\node [style=duplicate] (1) at (1, 2) {$\Delta$};
		\node [style=component] (2) at (0.5, 1) {$e$};
		\node [style=port] (3) at (1.5, 0.5) {};
		\node [style=port] (4) at (2.5, 3) {};
		\node [style=port] (5) at (2.5, 0.5) {};
		\node [style=duplicate] (6) at (4, 2) {$\Delta$};
		\node [style=port] (7) at (3.5, 0.5) {};
		\node [style=component] (8) at (4.5, 1) {$e$};
		\node [style=port] (9) at (4, 3) {};
		\node [style=port] (10) at (1.75, 1.5) {$=$};
		\node [style=port] (11) at (3, 1.5) {$=$};
		\node [style=port] (12) at (-4.25, 3.25) {};
		\node [style=port] (13) at (-4.25, 0.25) {};
		\node [style=port] (14) at (-5.25, 0.25) {};
		\node [style=duplicate] (15) at (-4.75, 1.25) {$\Delta$};
		\node [style=duplicate] (16) at (-4.25, 2.25) {$\Delta$};
		\node [style=port] (17) at (-3.5, 0.25) {};
		\node [style=port] (18) at (-1.75, 3.25) {};
		\node [style=port] (19) at (-1.75, 0.25) {};
		\node [style=port] (20) at (-0.75, 0.25) {};
		\node [style=duplicate] (21) at (-1.25, 1.25) {$\Delta$};
		\node [style=duplicate] (22) at (-1.75, 2.25) {$\Delta$};
		\node [style=port] (23) at (-2.5, 0.25) {};
		\node [style=port] (24) at (-3, 1.75) {$=$};
		\node [style=duplicate] (25) at (6.75, 2) {$\Delta$};
		\node [style=port] (26) at (7.5, 0.5) {};
		\node [style=port] (27) at (6, 0.5) {};
		\node [style=port] (28) at (6.75, 3) {};
		\node [style=duplicate] (29) at (9.25, 2.25) {$\Delta$};
		\node [style=port] (30) at (10, 0.5) {};
		\node [style=port] (31) at (8.5, 0.5) {};
		\node [style=port] (32) at (9.25, 3) {};
		\node [style=port] (33) at (8, 1.5) {$=$};
		\node [style=port] (34) at (11.25, 0.5) {};
		\node [style=port] (35) at (12, 3.5) {};
		\node [style=duplicate] (36) at (12, 1.5) {$\Delta$};
		\node [style=component] (37) at (12, 2.5) {$\delta$};
		\node [style=port] (38) at (12.75, 0.5) {};
		\node [style=port] (39) at (14.25, 3.5) {};
		\node [style=port] (40) at (13.75, 0.5) {};
		\node [style=component] (41) at (13.75, 1.25) {$\delta$};
		\node [style=duplicate] (42) at (14.25, 2.5) {$\Delta$};
		\node [style=port] (43) at (14.75, 0.5) {};
		\node [style=component] (44) at (14.75, 1.25) {$\delta$};
		\node [style=port] (45) at (13, 1.75) {$=$};
	\end{pgfonlayer}
	\begin{pgfonlayer}{edgelayer}
		\draw [style=wire, bend right=15, looseness=1.25] (1) to (2);
		\draw [style=wire, bend left=15] (1) to (3);
		\draw [style=wire] (0) to (1);
		\draw [style=wire] (4) to (5);
		\draw [style=wire, bend left=15, looseness=1.25] (6) to (8);
		\draw [style=wire, bend right=15] (6) to (7);
		\draw [style=wire] (9) to (6);
		\draw [style=wire, bend right=15, looseness=1.25] (16) to (15);
		\draw [style=wire, bend left=15] (16) to (17);
		\draw [style=wire] (12) to (16);
		\draw [style=wire, bend right=15] (15) to (14);
		\draw [style=wire, bend left=15] (15) to (13);
		\draw [style=wire, bend left=15, looseness=1.25] (22) to (21);
		\draw [style=wire, bend right=15] (22) to (23);
		\draw [style=wire] (18) to (22);
		\draw [style=wire, bend left=15] (21) to (20);
		\draw [style=wire, bend right=15] (21) to (19);
		\draw [style=wire, in=90, out=-39] (25) to (26);
		\draw [style=wire] (28) to (25);
		\draw [style=wire, in=90, out=-141] (25) to (27);
		\draw [style=wire] (32) to (29);
		\draw [style=wire, in=90, out=-45, looseness=1.50] (29) to (31);
		\draw [style=wire, in=90, out=-135, looseness=1.50] (29) to (30);
		\draw [style=wire] (35) to (37);
		\draw [style=wire, bend left] (36) to (38);
		\draw [style=wire, bend right] (36) to (34);
		\draw [style=wire] (37) to (36);
		\draw [style=wire, bend right=15, looseness=1.25] (42) to (41);
		\draw [style=wire] (39) to (42);
		\draw [style=wire] (41) to (40);
		\draw [style=wire] (44) to (43);
		\draw [style=wire, in=97, out=-53] (42) to (44);
	\end{pgfonlayer}
\end{tikzpicture}
 }% 
 \end{gathered}\end{equation}

Note that we do not assume that the functor $\oc$ of a coalgebra modality is a monoidal functor -- this will come soon.  And also note that requiring that $\Delta$ and $e$ be natural transformations is equivalent to asking that for each map $f: A \to B$, $\oc(f): \oc A \to \oc B$ is a comonoid morphism. This can be used to show that $\delta$ is in fact also a comonoid morphism. 

\begin{lemma} For any coalgebra modality $(\oc, \delta, \varepsilon, \Delta, e)$, $\delta$ also preserves the counit $e$, that is, the following diagram commutes:
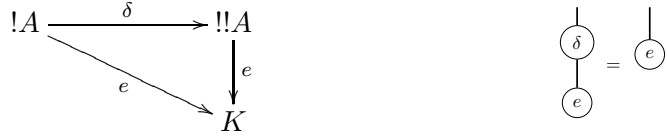
\begin{align*}
\begin{array}[c]{c}
\xymatrixcolsep{3pc}\xymatrixcolsep{5pc}\xymatrix{\oc A \ar[r]^-{\delta} \ar[dr]_-{e} & \oc \oc A \ar[d]^-{e} \\
  & K  
  }
   \end{array}&& 
   \begin{array}[c]{c}
\resizebox{!}{1.5cm}{%
  \begin{tikzpicture}
	\begin{pgfonlayer}{nodelayer}
		\node [style=port] (1) at (0, 3) {};
		\node [style=component] (2) at (0, 1) {$e$};
		\node [style=component] (3) at (0, 2.25) {$\delta$};
		\node [style=port] (5) at (1.5, 3) {};
		\node [style=component] (8) at (1.5, 2) {$e$};
		\node [style=port] (11) at (0.75, 1.75) {$=$};
	\end{pgfonlayer}
	\begin{pgfonlayer}{edgelayer}
		\draw [style=wire] (1) to (3);
		\draw [style=wire] (3) to (2);
		\draw [style=wire] (5) to (8);
	\end{pgfonlayer}
\end{tikzpicture}
}%
   \end{array}
\end{align*}
Therefore, $\delta$ is a comonoid morphism.
\end{lemma}
\begin{proof} By the naturality of $e$ and the comonad identities, we obtain that: 
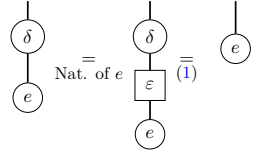
\[\resizebox{!}{2cm}{%
\begin{tikzpicture}
	\begin{pgfonlayer}{nodelayer}
		\node [style=port] (0) at (0.5, 3) {};
		\node [style=component] (1) at (0.5, 1) {$e$};
		\node [style=component] (2) at (0.5, 2.25) {$\delta$};
		\node [style=port] (3) at (4.75, 3) {};
		\node [style=component] (4) at (4.75, 2) {$e$};
		\node [style=port] (5) at (1.75, 1.75) {$=$};
		\node [style=port] (6) at (3, 3) {};
		\node [style=function] (7) at (3, 1.25) {$\varepsilon$};
		\node [style=component] (8) at (3, 2.25) {$\delta$};
		\node [style=component] (9) at (3, 0.25) {$e$};
		\node [style=port] (10) at (3.75, 1.75) {$=$};
		\node [style=port] (11) at (1.75, 1.5) {\small Nat. of $e$};
		\node [style=port] (12) at (3.75, 1.5) {\small (\ref{comonad})};
	\end{pgfonlayer}
	\begin{pgfonlayer}{edgelayer}
		\draw [style=wire] (0) to (2);
		\draw [style=wire] (2) to (1);
		\draw [style=wire] (3) to (4);
		\draw [style=wire] (6) to (8);
		\draw [style=wire] (8) to (7);
		\draw [style=wire] (7) to (9);
	\end{pgfonlayer}
\end{tikzpicture}
}% 
\]
\end{proof} 

We now turn our attention to the case when the coalgebra modality is monoidal. Recall that a symmetric monoidal endofunctor \cite{mac2013categories} on a symmetric monoidal category $\mathbb{X}$ is a triple $(\oc, m_\otimes, m_K)$ consisting of an enfunctor $\oc: \mathbb{X} \to \mathbb{X}$, a natural transformation $m_\otimes: \oc A \otimes \oc B \to \oc (A \otimes B)$, and a map $m_K: K \to \oc K$ such that the following diagrams commute:
\begin{equation}\label{}\begin{gathered} \xymatrixcolsep{1.75pc}\xymatrix{\oc A \otimes \oc B \otimes \oc C  \ar[r]^-{m_\otimes \otimes 1} \ar[d]_-{1 \otimes m_\otimes} & \oc(A \otimes B) \otimes \oc C \ar[d]^-{m_\otimes}& \oc A  \ar@{=}[dr] \ar[d]_-{m_K \otimes 1} \ar[r]^-{1 \otimes m_K} & \oc A \otimes \oc K \ar[d]^-{m_\otimes} & \oc A \otimes \oc B  \ar[r]^-{\sigma}  \ar[d]_-{m_\otimes} & \oc B \otimes \oc A \ar[d]^-{m_\otimes}\\ 
      \oc A \otimes \oc(B \otimes C) \ar[r]_-{m_\otimes} & \oc(A \otimes B \otimes C) & \oc K\otimes \oc A \ar[r]_-{m_\otimes} & \oc A &    \oc(A \otimes B) \ar[r]_-{\oc(\sigma)} &  \oc(B \otimes A)  } \end{gathered}\end{equation}
In the graphical calculus, $m_\otimes$ and $m_K$ are drawn respectively as follows: 
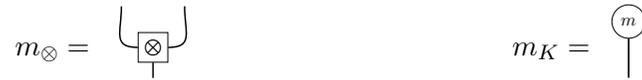
\begin{align*}
m_\otimes = \begin{array}[c]{c}\resizebox{!}{1cm}{%
\begin{tikzpicture}
	\begin{pgfonlayer}{nodelayer}
		\node [style=port] (0) at (2, 2) {};
		\node [style=port] (1) at (0.5, 2) {};
		\node [style={regular polygon,regular polygon sides=4, draw, inner sep=1pt,minimum size=1pt}] (2) at (1.25, 1) {$\bigotimes$};
		\node [style=port] (3) at (1.25, 0.25) {};
	\end{pgfonlayer}
	\begin{pgfonlayer}{edgelayer}
		\draw [style=wire, in=-90, out=180, looseness=1.50] (2) to (1);
		\draw [style=wire, in=0, out=-90, looseness=1.50] (0) to (2);
		\draw [style=wire] (2) to (3);
	\end{pgfonlayer}
\end{tikzpicture}
  }% 
\end{array} 
 && m_K= \begin{array}[c]{c}\resizebox{!}{1cm}{%
\begin{tikzpicture}
	\begin{pgfonlayer}{nodelayer}
		\node [style={circle, draw}] (0) at (1, 1.75) {$m$};
		\node [style=port] (1) at (1, 0.5) {};
	\end{pgfonlayer}
	\begin{pgfonlayer}{edgelayer}
		\draw [style=wire] (0) to (1);
	\end{pgfonlayer}
\end{tikzpicture}}
   \end{array}
\end{align*}
And so the symmetric monoidal endofunctor identities are drawn as follows: 
\begin{equation}\label{smendo}\begin{gathered} \resizebox{!}{2cm}{%
\begin{tikzpicture}
	\begin{pgfonlayer}{nodelayer}
		\node [style=function2] (1) at (3.25, -0.25) {$\bigotimes$};
		\node [style=port] (2) at (3.25, -1.25) {};
		\node [style=component] (3) at (2.5, 0.75) {$m$};
		\node [style=port] (4) at (4, 0.75) {};
		\node [style=port] (5) at (4.75, 0.75) {};
		\node [style=port] (6) at (4.75, -1.25) {};
		\node [style=function2] (8) at (6.25, -0.25) {$\bigotimes$};
		\node [style=port] (9) at (6.25, -1.25) {};
		\node [style=component] (10) at (7, 0.75) {$m$};
		\node [style=port] (11) at (5.5, 0.75) {};
		\node [style=function2] (12) at (10.25, -0.25) {$\bigotimes$};
		\node [style=port] (13) at (10.25, -1.25) {};
		\node [style=port] (14) at (11.25, 1.75) {};
		\node [style=port] (15) at (8.75, 1.75) {};
		\node [style=port] (16) at (10.25, 1.75) {};
		\node [style=function2] (17) at (9.5, 0.75) {$\bigotimes$};
		\node [style=function2] (18) at (13, -0.25) {$\bigotimes$};
		\node [style=port] (19) at (13, -1.25) {};
		\node [style=port] (20) at (12, 1.75) {};
		\node [style=port] (21) at (14.5, 1.75) {};
		\node [style=port] (22) at (13, 1.75) {};
		\node [style=function2] (23) at (13.75, 0.75) {$\bigotimes$};
		\node [style=function2] (24) at (16.5, -0.5) {$\bigotimes$};
		\node [style=port] (25) at (16.5, -1.5) {};
		\node [style=port] (26) at (16, 2.25) {};
		\node [style=port] (27) at (17, 2.25) {};
		\node [style=function2] (28) at (19, 1.25) {$\bigotimes$};
		\node [style=port] (29) at (19, 0.25) {};
		\node [style=port] (30) at (18.25, 2.25) {};
		\node [style=port] (31) at (19.75, 2.25) {};
		\node [style=port] (32) at (18.25, 0.25) {};
		\node [style=port] (33) at (19.75, 0.25) {};
		\node [style=port] (34) at (18.25, -0.75) {};
		\node [style=port] (35) at (19.75, -0.75) {};
		\node [style=port] (36) at (19, -0.75) {};
		\node [style=port] (37) at (19, -1.5) {};
		\node [style=port] (38) at (18.5, 0.25) {};
		\node [style=port] (39) at (19.5, 0.25) {};
		\node [style=port] (40) at (18.5, -0.75) {};
		\node [style=port] (41) at (19.5, -0.75) {};
		\node [style=port] (42) at (4.25, -0.5) {$=$};
		\node [style=port] (43) at (5.25, -0.5) {$=$};
		\node [style=port] (44) at (11.5, -0.25) {$=$};
		\node [style=port] (45) at (17.75, -0.25) {$=$};
	\end{pgfonlayer}
	\begin{pgfonlayer}{edgelayer}
		\draw [style=wire, in=-90, out=180] (1) to (3);
		\draw [style=wire] (1) to (2);
		\draw [style=wire, in=0, out=-90] (4) to (1);
		\draw [style=wire] (5) to (6);
		\draw [style=wire, in=-90, out=0] (8) to (10);
		\draw [style=wire] (8) to (9);
		\draw [style=wire, in=180, out=-90] (11) to (8);
		\draw [style=wire] (12) to (13);
		\draw [style=wire, in=0, out=-90] (14) to (12);
		\draw [style=wire, in=-90, out=180] (17) to (15);
		\draw [style=wire, in=0, out=-90] (16) to (17);
		\draw [style=wire, in=180, out=-90, looseness=1.25] (17) to (12);
		\draw [style=wire] (18) to (19);
		\draw [style=wire, in=180, out=-90] (20) to (18);
		\draw [style=wire, in=-90, out=0] (23) to (21);
		\draw [style=wire, in=180, out=-90] (22) to (23);
		\draw [style=wire, in=0, out=-90, looseness=1.25] (23) to (18);
		\draw [style=wire] (24) to (25);
		\draw [style=wire, in=-90, out=180] (28) to (30);
		\draw [style=wire] (28) to (29);
		\draw [style=wire, in=0, out=-90] (31) to (28);
		\draw [style=wire] (32) to (34);
		\draw [style=wire] (32) to (33);
		\draw [style=wire] (34) to (35);
		\draw [style=wire] (33) to (35);
		\draw [style=wire] (36) to (37);
		\draw [style=wire, in=90, out=-90] (38) to (41);
		\draw [style=wire, in=90, out=-90] (39) to (40);
		\draw [style=wire, in=15, out=-90, looseness=1.50] (26) to (24);
		\draw [style=wire, in=165, out=-90, looseness=1.50] (27) to (24);
	\end{pgfonlayer}
\end{tikzpicture}
}%
 \end{gathered}\end{equation}
A symmetric monoidal comonad \cite{bierman1995categorical} on a symmetric monoidal category is a quintuple $(\oc, \delta, \varepsilon, m_\otimes, m_K)$ consisting of a comonad $(\oc, \delta, \varepsilon)$ and a symmetric monoidal endofunctor $(\oc, m_\otimes, m_K)$ and such that $\delta$ and $\varepsilon$ are monoidal natural transformations, that is, the following diagrams commute:
\begin{equation}\label{}\begin{gathered} \xymatrixcolsep{5pc}\xymatrix{\oc A \otimes \oc B  \ar[d]_-{\delta \otimes \delta} \ar[r]^-{m_\otimes}  & \oc(A \otimes B) \ar[dd]^-{\delta} & K \ar[r]^-{ m_K} \ar[d]_-{m_K} & \oc K \ar[d]^-{\delta}  \\
  \oc\oc A \otimes \oc\oc B  \ar[d]_-{m_\otimes}  & &\oc K \ar[r]_-{\oc(m_K)}& \oc \oc K \\
  \oc (\oc A \otimes \oc B)  \ar[r]_-{\oc(m_\otimes)} & \oc(A \otimes B)
  } \\   \xymatrixcolsep{5pc}\xymatrix{\oc A \otimes \oc B \ar[dr]_-{\varepsilon \otimes \varepsilon} \ar[r]^-{m_\otimes} & \oc(A \otimes B) \ar[d]^-{\varepsilon} & K \ar@{=}[dr]^-{} \ar[r]^-{ m_K} & \oc K \ar[d]^-{\varepsilon} \\
& A \otimes B  & & \oc K 
  } \end{gathered}\end{equation}
  which drawn in the graphical calculus gives: 
  \begin{equation}\label{deltaepmon}\begin{gathered} 
  \resizebox{!}{3.25cm}{%
  \begin{tikzpicture}
	\begin{pgfonlayer}{nodelayer}
		\node [style=function2] (0) at (1.25, 0) {$\bigotimes$};
		\node [style=component] (1) at (1.25, -1.25) {$\delta$};
		\node [style=port] (2) at (0.5, 1) {};
		\node [style=port] (3) at (2, 1) {};
		\node [style=port] (4) at (1.25, -2.25) {};
		\node [style=port] (5) at (2.25, -1) {$=$};
		\node [style=function2] (6) at (4.25, 1) {$\bigotimes$};
		\node [style=component] (7) at (3.5, 2) {$\delta$};
		\node [style=component] (8) at (5, 2) {$\delta$};
		\node [style=port] (9) at (3.5, 3) {};
		\node [style=port] (10) at (5, 3) {};
		\node [style=port] (11) at (3, 0) {};
		\node [style=port] (12) at (5.5, 0) {};
		\node [style=port] (13) at (3, -1.75) {};
		\node [style=port] (14) at (5.5, -1.75) {};
		\node [style=port] (15) at (4.25, 0) {};
		\node [style=port] (16) at (3.5, 0) {};
		\node [style=port] (17) at (5, 0) {};
		\node [style=function2] (18) at (4.25, -1) {$\bigotimes$};
		\node [style=port] (19) at (4.25, -2.5) {};
		\node [style=component] (20) at (7.25, -1) {$\delta$};
		\node [style=component] (21) at (7.25, 0.25) {$m$};
		\node [style=port] (22) at (7.25, -2) {};
		\node [style=function2] (23) at (9, -1) {$m$};
		\node [style=component] (24) at (9, 0.25) {$m$};
		\node [style=port] (25) at (9, -2) {};
		\node [style=function2] (26) at (11.75, 0) {$\bigotimes$};
		\node [style=component] (27) at (11.75, -1) {$\varepsilon$};
		\node [style=port] (28) at (12.5, 1) {};
		\node [style=port] (29) at (11, 1) {};
		\node [style=port] (30) at (11.25, -2) {};
		\node [style=port] (31) at (12.25, -2) {};
		\node [style=component] (32) at (13.5, -0.5) {$\varepsilon$};
		\node [style=port] (33) at (13.5, 0.75) {};
		\node [style=port] (34) at (13.5, -1.75) {};
		\node [style=port] (35) at (14.5, -1.75) {};
		\node [style=port] (36) at (14.5, 0.75) {};
		\node [style=component] (37) at (14.5, -0.5) {$\varepsilon$};
		\node [style=component] (38) at (16.25, 0) {$m$};
		\node [style=component] (39) at (16.25, -1) {$\varepsilon$};
		\node [style=port] (41) at (8, -0.5) {$=$};
		\node [style=port] (42) at (12.75, -0.5) {$=$};
		\node [style=port] (43) at (17, -0.5) {$=$};
	\end{pgfonlayer}
	\begin{pgfonlayer}{edgelayer}
		\draw [style=wire, in=-90, out=180, looseness=1.50] (0) to (2);
		\draw [style=wire] (0) to (1);
		\draw [style=wire, in=0, out=-90, looseness=1.50] (3) to (0);
		\draw [style=wire] (1) to (4);
		\draw [style=wire, in=-90, out=180, looseness=1.50] (6) to (7);
		\draw [style=wire, in=0, out=-90, looseness=1.50] (8) to (6);
		\draw [style=wire] (9) to (7);
		\draw [style=wire] (10) to (8);
		\draw [style=wire] (11) to (12);
		\draw [style=wire] (12) to (14);
		\draw [style=wire] (11) to (13);
		\draw [style=wire] (13) to (14);
		\draw [style=wire] (6) to (15);
		\draw [style=wire, in=-90, out=180, looseness=1.50] (18) to (16);
		\draw [style=wire, in=0, out=-90, looseness=1.50] (17) to (18);
		\draw [style=wire] (18) to (19);
		\draw [style=wire] (21) to (20);
		\draw [style=wire] (20) to (22);
		\draw [style=wire] (24) to (23);
		\draw [style=wire] (23) to (25);
		\draw [style=wire, in=-90, out=0, looseness=1.50] (26) to (28);
		\draw [style=wire, in=180, out=-90, looseness=1.50] (29) to (26);
		\draw [style=wire] (26) to (27);
		\draw [style=wire, in=90, out=-120] (27) to (30);
		\draw [style=wire, in=90, out=-63] (27) to (31);
		\draw [style=wire] (33) to (32);
		\draw [style=wire] (32) to (34);
		\draw [style=wire] (36) to (37);
		\draw [style=wire] (37) to (35);
		\draw [style=wire] (38) to (39);
	\end{pgfonlayer}
\end{tikzpicture}
}%
 \end{gathered}\end{equation}
The symmetric monoidal comonad coherences are precisely what is required so that the co-Eilenberg-Moore category be a symmetric monoidal category such that the forgetful functor preserves the symmetric monoidal structure strictly. 

   \begin{definition}\label{mcoalgdef} A \textbf{monoidal coalgebra modality} \cite{blute2015cartesian} on a symmetric monoidal category is a septuple $(\oc, \delta, \varepsilon, m_\otimes, m_K, \Delta, e)$ consisting of a coalgebra modality $(\oc, \delta, \varepsilon, \Delta, e)$ and a symmetric monoidal comonad $(\oc, m, m_K, \delta, \varepsilon)$, and such that $\Delta$ and $e$ are monoidal transformations (or equivalently $m_\otimes$ and $m_K$ are comonoid morphisms), that is, the following diagrams commute: 
   \begin{equation}\label{}\begin{gathered}\xymatrixcolsep{3.5pc}\xymatrix{\oc A \otimes \oc B  \ar[d]_-{\Delta \otimes \Delta} \ar[r]^-{m_\otimes}  & \oc(A \otimes B) \ar[dd]^-{\Delta} & \oc A \otimes \oc B \ar[d]_-{e \otimes e} \ar[r]^-{m_\otimes} & \oc(A \otimes B) \ar[d]^-{e}  \\
  \oc A \otimes \oc A \otimes \oc B \otimes \oc B \ar[d]_-{1 \otimes \sigma \otimes 1}  && K \otimes K \ar@{=}[r]_-{}& K  \\
  \oc A \otimes \oc B \otimes \oc A \otimes \oc B \ar[r]_-{m_\otimes \otimes m_\otimes} & \oc(A \otimes B) \otimes \oc(A \otimes B) 
  } \\
  \xymatrixcolsep{5pc}\xymatrix{ K \ar[r]^-{m_K} \ar@{=}[d]_-{} & \oc K \ar[d]^-{\Delta} & K \ar@{=}[dr]^-{} \ar[r]^-{m_K} & \oc K \ar[d]^-{e_K} \\
K \otimes K \ar[r]_-{m_K \otimes m_K}  &\oc K \otimes \oc K && K
  }  \end{gathered}\end{equation}
and also that $\Delta$ and $e$ are $\oc$-coalgebra morphisms, that is, the following diagrams commute: 
\begin{equation}\label{}\begin{gathered} \xymatrixcolsep{3.5pc}\xymatrix{\oc A \ar[d]_-{\Delta} \ar[rr]^-{\delta} & & \oc \oc A \ar[d]^-{\oc(\Delta)} & \oc A \ar[d]_-{e} \ar[r]^-{\delta} & \oc \oc A \ar[d]^-{\oc(e)} \\
    \oc A \otimes \oc A \ar[r]_-{\delta \otimes \delta} & \oc \oc A \otimes \oc \oc A \ar[r]_-{m_\otimes} & \oc(\oc A \otimes \oc A) &   K \ar[r]_-{m_K} & \oc(K)
  } \end{gathered}\end{equation}
 A \textbf{linear category}\footnote{Note that here we do not require linear categories to be closed, following \cite{blute2015cartesian}.} \cite{bierman1995categorical,blute2015cartesian} is a symmetric monoidal category with a monoidal coalgebra modality.  
 \end{definition}

In the graphical calculus, that $\Delta$ and $e$ are both monoidal transformations and $\oc$-colagebra morphisms is expressed as follows: 
\begin{equation}\label{mcm1}\begin{gathered} 
\resizebox{!}{2.5cm}{%
\begin{tikzpicture}
	\begin{pgfonlayer}{nodelayer}
		\node [style=port] (0) at (1.75, 1) {};
		\node [style=port] (1) at (3.25, 1) {};
		\node [style=function2] (2) at (2.5, 0) {$\bigotimes$};
		\node [style=duplicate] (3) at (2.5, -1) {$\Delta$};
		\node [style=port] (4) at (3.25, -2) {};
		\node [style=port] (5) at (1.75, -2) {};
		\node [style=port] (6) at (6.75, -2) {};
		\node [style=port] (7) at (6.75, 2) {};
		\node [style=duplicate] (8) at (6.75, 1) {$\Delta$};
		\node [style=port] (9) at (4.75, 2) {};
		\node [style=function2] (10) at (6.75, -1) {$\bigotimes$};
		\node [style=function2] (11) at (4.75, -1) {$\bigotimes$};
		\node [style=duplicate] (12) at (4.75, 1) {$\Delta$};
		\node [style=port] (13) at (4.75, -2) {};
		\node [style=duplicate] (14) at (9.75, 0) {$\Delta$};
		\node [style=port] (15) at (10.5, -1) {};
		\node [style=port] (16) at (9, -1) {};
		\node [style=component] (17) at (9.75, 1) {$m$};
		\node [style=port] (18) at (12.75, -1) {};
		\node [style=port] (19) at (11.75, -1) {};
		\node [style=component] (20) at (11.75, 0.5) {$m$};
		\node [style=component] (21) at (12.75, 0.5) {$m$};
		\node [style=port] (22) at (14.5, 1) {};
		\node [style=port] (23) at (16, 1) {};
		\node [style=function2] (24) at (15.25, 0) {$\bigotimes$};
		\node [style=component] (25) at (15.25, -1) {$e$};
		\node [style=port] (26) at (18, 0.5) {};
		\node [style=port] (27) at (17, 0.5) {};
		\node [style=component] (28) at (17, -1) {$e$};
		\node [style=component] (29) at (18, -1) {$e$};
		\node [style=component] (30) at (19.75, -1) {$e$};
		\node [style=component] (31) at (19.75, 0.5) {$m$};
		\node [style=port] (32) at (3.5, -0.5) {$=$};
		\node [style=port] (33) at (11, 0) {$=$};
		\node [style=port] (34) at (16.5, 0) {$=$};
		\node [style=port] (35) at (20.5, -0.25) {$=$};
	\end{pgfonlayer}
	\begin{pgfonlayer}{edgelayer}
		\draw [style=wire, in=-90, out=180, looseness=1.50] (2) to (0);
		\draw [style=wire, in=0, out=-90, looseness=1.50] (1) to (2);
		\draw [style=wire, bend right] (3) to (5);
		\draw [style=wire, bend left] (3) to (4);
		\draw [style=wire] (2) to (3);
		\draw [style=wire] (9) to (12);
		\draw [style=wire] (11) to (13);
		\draw [style=wire, in=180, out=-135, looseness=1.25] (12) to (11);
		\draw [style=wire] (7) to (8);
		\draw [style=wire] (10) to (6);
		\draw [style=wire, in=0, out=-45, looseness=1.25] (8) to (10);
		\draw [style=wire, in=0, out=-150] (8) to (11);
		\draw [style=wire, in=180, out=-30] (12) to (10);
		\draw [style=wire, bend right] (14) to (16);
		\draw [style=wire, bend left] (14) to (15);
		\draw [style=wire] (17) to (14);
		\draw [style=wire] (20) to (19);
		\draw [style=wire] (21) to (18);
		\draw [style=wire, in=-90, out=180, looseness=1.50] (24) to (22);
		\draw [style=wire, in=0, out=-90, looseness=1.50] (23) to (24);
		\draw [style=wire] (24) to (25);
		\draw [style=wire] (28) to (27);
		\draw [style=wire] (29) to (26);
		\draw [style=wire] (31) to (30);
	\end{pgfonlayer}
\end{tikzpicture}
}%
 \end{gathered}\end{equation}
\begin{equation}\label{mcm2}\begin{gathered} 
\resizebox{!}{2.5cm}{%
\begin{tikzpicture}
	\begin{pgfonlayer}{nodelayer}
		\node [style=duplicate] (0) at (1.75, 0) {$\Delta$};
		\node [style=port] (1) at (2.5, -1) {};
		\node [style=port] (2) at (1, -1) {};
		\node [style=port] (3) at (0.5, -1) {};
		\node [style=port] (4) at (3, -1) {};
		\node [style=port] (5) at (3, 0.5) {};
		\node [style=port] (6) at (0.5, 0.5) {};
		\node [style=component] (7) at (1.75, 1.25) {$\delta$};
		\node [style=port] (8) at (1.75, 2) {};
		\node [style=port] (9) at (1.75, -2) {};
		\node [style=port] (10) at (1.75, -1) {};
		\node [style=duplicate] (11) at (5, 0.75) {$\Delta$};
		\node [style=component] (12) at (5.75, -0.25) {$\delta$};
		\node [style=component] (13) at (4.25, -0.25) {$\delta$};
		\node [style=port] (14) at (5, 1.75) {};
		\node [style=function2] (15) at (5, -1.25) {$\bigotimes$};
		\node [style=port] (16) at (5, -2) {};
		\node [style=function] (17) at (8.75, -0.75) {$e$};
		\node [style=component] (22) at (8.75, 0.25) {$\delta$};
		\node [style=port] (23) at (8.75, 1.25) {};
		\node [style=port] (24) at (8.75, -1.75) {};
		\node [style=component] (27) at (10.25, 0.25) {$e$};
		\node [style=component] (28) at (10.25, -0.75) {$m$};
		\node [style=port] (29) at (10.25, -1.75) {};
		\node [style=port] (30) at (10.25, 1.25) {};
		\node [style=port] (31) at (3.5, -0.25) {$=$};
		\node [style=port] (32) at (9.5, -0.25) {$=$};
	\end{pgfonlayer}
	\begin{pgfonlayer}{edgelayer}
		\draw [style=wire, bend right] (0) to (2);
		\draw [style=wire, bend left] (0) to (1);
		\draw [style=wire] (6) to (5);
		\draw [style=wire] (6) to (3);
		\draw [style=wire] (3) to (4);
		\draw [style=wire] (5) to (4);
		\draw [style=wire] (8) to (7);
		\draw [style=wire] (10) to (9);
		\draw [style=wire] (7) to (0);
		\draw [style=wire, bend right] (11) to (13);
		\draw [style=wire, bend left] (11) to (12);
		\draw [style=wire] (14) to (11);
		\draw [style=wire] (15) to (16);
		\draw [style=wire, in=0, out=-90, looseness=1.75] (12) to (15);
		\draw [style=wire, in=180, out=-90, looseness=1.75] (13) to (15);
		\draw [style=wire] (23) to (22);
		\draw [style=wire] (22) to (17);
		\draw [style=wire] (28) to (29);
		\draw [style=wire] (30) to (27);
		\draw [style=wire] (17) to (24);
	\end{pgfonlayer}
\end{tikzpicture}
}%
 \end{gathered}\end{equation}
 
As explained in \cite{schalk2004categorical}, the monoidal coalgebra modality coherences are precisely what is required so that tensor product of the base category becomes a product in the co-Eilenberg-Moore category. And in the presence of finite products, monoidal coalgebra modalities can be equivalently be described using the Seely isomorphisms -- which we discuss in Section \ref{Seelysec}. 

\section{Additive Bialgebra Modalities}

In this section we introduce the notion of an {\em additive} bialgebra modality. In the presence of additive structure, additive bialgebra modalities are in bijective correspondence to monoidal coalgebra modalities (Theorem \ref{addlinaddbialg}). In particular, in Section \ref{equivalencesec}, we will show that for additive bialgebra modalities, deriving transformations and coderiction maps are equivalent. The majority of the proofs of this section, due to their length, can be found in the appendix. 

We begin by recalling additive structure by starting with the notion of an additive category.  Here we mean ``additive'' in the sense of being commutative monoid enriched: we do not assume negatives nor do we assume biproducts (this differs from the usage in \cite{mac2013categories} for example). This allows many important examples such as the category of sets and relations or the category of modules for a commutative rig\footnote{Rigs are also known as a semirings: they are ri{\em n\/}gs without {\em n\/}egatives.}.

\begin{definition}\label{addcatdef} An \textbf{additive category} \cite{blute2006differential} is a commutative monoid enriched category, that is, a category in which each hom-set is a commutative monoid with an addition operation $+$ and a zero $0$, and such that composition preserves the additive structure, that is:
\begin{align*}
k(f\!+\!g)h\!=kfh\!+\!kgh && k0h=0
\end{align*}
An \textbf{additive symmetric monoidal category} \cite{blute2006differential} is a symmetric monoidal category which is also an additive category in which the tensor product is compatible with the additive structure in the sense that: 
\begin{align*}
k \otimes (f\!+\!g)\otimes h\!= \!k\otimes\!f\otimes h \!+ \!k\otimes\!g\otimes h && k \otimes 0 \otimes h\!=\!0 
\end{align*}
\end{definition}

In \cite{blute2006differential}, it was observed that if a coalgebra modality came equipped with a natural bialgebra structure and a codereliction then one could obtain a deriving transformation (more on this in Section \ref{codersec}). 

\begin{definition} A \textbf{bialgebra modality} \cite{blute2006differential} on an additive symmetric monoidal category is a septuple $(\oc,  \delta, \varepsilon, \Delta, e, \nabla, u)$ consisting of a coalgebra modality $(\oc,  \delta, \varepsilon, \Delta, e)$, a natural transformation $\nabla: \oc A \otimes \oc A \to \oc A$, and a natural transformation $u: K \to \oc A$ such that $(\oc A, \nabla, u)$ is a commutative monoid, that is, the dual diagrams of (\ref{cocomdiagrams}) commute, and $(\oc A, \nabla, u, \Delta, e)$ is a bialgebra, that is, the following diagrams commute: 
\begin{equation}\label{}\begin{gathered}\xymatrixcolsep{1.5pc}\xymatrix{\oc A \otimes \oc A \ar[dr]_{e \otimes e} \ar[r]^-{\nabla} & \oc A  \ar[d]^-{e}& K \ar[dr]_-{u \otimes u} \ar[r]^-{u} & \oc A \ar[d]^-{\Delta} & K \ar@{=}[dr] \ar[r]^-{u} & \oc A \ar[d]^-{e} & \oc A \otimes \oc A \ar[dd]_-{\nabla} \ar[rr]^-{\Delta \otimes \Delta}&& \oc A \otimes \oc A \otimes \oc A \otimes \oc A \ar[d]^-{1 \otimes \sigma \otimes 1}   \\ 
&K&& \oc A \otimes \oc A && K  & && \oc A \otimes \oc A \otimes \oc A \otimes \oc A \ar[d]^-{\nabla \otimes \nabla} \\
&&&&&&   \oc A \ar[rr]_-{\Delta} && \oc A \otimes \oc A
  } \end{gathered}\end{equation}
and $\varepsilon$ is compatible with $\nabla$ in the sense that the following diagram commutes: 
\begin{equation}\label{}\begin{gathered} \xymatrixcolsep{5pc}\xymatrix{ \oc A \otimes \oc A \ar[dr]_-{\varepsilon \otimes e+ e \otimes \varepsilon} \ar[r]^-{\nabla} & \oc A \ar[d]^-{\varepsilon} \\
& A 
  } \end{gathered}\end{equation}
\end{definition}

In the graphical calculus, the extra bialgebra modality identities are drawn as follows: 
\begin{equation}\label{moneq}\begin{gathered} 
\resizebox{!}{2cm}{%
\begin{tikzpicture}
	\begin{pgfonlayer}{nodelayer}
		\node [style=port] (0) at (6.75, 0.75) {};
		\node [style=duplicate] (1) at (6.75, 1.75) {$\nabla$};
		\node [style=component] (2) at (6.25, 2.75) {$u$};
		\node [style=port] (3) at (7.25, 3.25) {};
		\node [style=port] (4) at (8.25, 0.75) {};
		\node [style=port] (5) at (8.25, 3.25) {};
		\node [style=duplicate] (6) at (9.75, 1.75) {$\nabla$};
		\node [style=port] (7) at (9.25, 3.25) {};
		\node [style=component] (8) at (10.25, 2.75) {$u$};
		\node [style=port] (9) at (9.75, 0.75) {};
		\node [style=port] (10) at (7.5, 2.25) {$=$};
		\node [style=port] (11) at (8.75, 2.25) {$=$};
		\node [style=port] (12) at (1.5, 0.5) {};
		\node [style=port] (13) at (1.5, 3.5) {};
		\node [style=port] (14) at (0.5, 3.5) {};
		\node [style=duplicate] (15) at (1, 2.5) {$\nabla$};
		\node [style=duplicate] (16) at (1.5, 1.5) {$\nabla$};
		\node [style=port] (17) at (2.25, 3.5) {};
		\node [style=port] (18) at (4, 0.5) {};
		\node [style=port] (19) at (4, 3.5) {};
		\node [style=port] (20) at (5, 3.5) {};
		\node [style=duplicate] (21) at (4.5, 2.5) {$\nabla$};
		\node [style=duplicate] (22) at (4, 1.5) {$\nabla$};
		\node [style=port] (23) at (3.25, 3.5) {};
		\node [style=port] (24) at (2.75, 2) {$=$};
		\node [style=duplicate] (25) at (12.5, 1.75) {$\nabla$};
		\node [style=port] (26) at (13.25, 3.25) {};
		\node [style=port] (27) at (11.75, 3.25) {};
		\node [style=port] (28) at (12.5, 0.75) {};
		\node [style=duplicate] (29) at (15, 1.5) {$\nabla$};
		\node [style=port] (30) at (15.75, 3.25) {};
		\node [style=port] (31) at (14.25, 3.25) {};
		\node [style=port] (32) at (15, 0.75) {};
		\node [style=port] (33) at (13.75, 2.25) {$=$};
	\end{pgfonlayer}
	\begin{pgfonlayer}{edgelayer}
		\draw [style=wire, bend left=15, looseness=1.25] (1) to (2);
		\draw [style=wire, bend right=15] (1) to (3);
		\draw [style=wire] (0) to (1);
		\draw [style=wire] (4) to (5);
		\draw [style=wire, bend right=15, looseness=1.25] (6) to (8);
		\draw [style=wire, bend left=15] (6) to (7);
		\draw [style=wire] (9) to (6);
		\draw [style=wire, bend left=15, looseness=1.25] (16) to (15);
		\draw [style=wire, bend right=15] (16) to (17);
		\draw [style=wire] (12) to (16);
		\draw [style=wire, bend left=15] (15) to (14);
		\draw [style=wire, bend right=15] (15) to (13);
		\draw [style=wire, bend right=15, looseness=1.25] (22) to (21);
		\draw [style=wire, bend left=15] (22) to (23);
		\draw [style=wire] (18) to (22);
		\draw [style=wire, bend right=15] (21) to (20);
		\draw [style=wire, bend left=15] (21) to (19);
		\draw [style=wire, in=-90, out=39] (25) to (26);
		\draw [style=wire] (28) to (25);
		\draw [style=wire, in=-90, out=141] (25) to (27);
		\draw [style=wire] (32) to (29);
		\draw [style=wire, in=-90, out=45, looseness=1.50] (29) to (31);
		\draw [style=wire, in=-90, out=135, looseness=1.50] (29) to (30);
	\end{pgfonlayer}
\end{tikzpicture}
}%
 \end{gathered}\end{equation}
\begin{equation}\label{bialgeq}\begin{gathered} 
\resizebox{!}{2.5cm}{%
\begin{tikzpicture}
	\begin{pgfonlayer}{nodelayer}
		\node [style=duplicate] (3) at (2.5, -1) {$\Delta$};
		\node [style=port] (4) at (3.25, -2) {};
		\node [style=port] (5) at (1.75, -2) {};
		\node [style=port] (6) at (6.5, -2) {};
		\node [style=port] (7) at (6.5, 2) {};
		\node [style=duplicate] (8) at (6.5, 1) {$\Delta$};
		\node [style=port] (9) at (4.75, 2) {};
		\node [style=duplicate] (10) at (6.5, -1) {$\nabla$};
		\node [style=duplicate] (11) at (4.75, -1) {$\nabla$};
		\node [style=duplicate] (12) at (4.75, 1) {$\Delta$};
		\node [style=port] (13) at (4.75, -2) {};
		\node [style=duplicate] (14) at (9.5, 0) {$\Delta$};
		\node [style=port] (15) at (10.25, -1) {};
		\node [style=port] (16) at (8.75, -1) {};
		\node [style=component] (17) at (9.5, 1) {$u$};
		\node [style=port] (18) at (12.5, -1) {};
		\node [style=port] (19) at (11.5, -1) {};
		\node [style=component] (20) at (11.5, 0.5) {$u$};
		\node [style=component] (21) at (12.5, 0.5) {$u$};
		\node [style=component] (30) at (19.75, -1) {$e$};
		\node [style=component] (31) at (19.75, 0.5) {$u$};
		\node [style=port] (32) at (3.5, -0.5) {$=$};
		\node [style=port] (33) at (10.75, 0) {$=$};
		\node [style=port] (35) at (20.5, -0.25) {$=$};
		\node [style=duplicate] (36) at (2.5, 0) {$\nabla$};
		\node [style=port] (37) at (3.25, 1) {};
		\node [style=port] (38) at (1.75, 1) {};
		\node [style=duplicate] (39) at (14.75, 0) {$\nabla$};
		\node [style=port] (40) at (15.5, 1) {};
		\node [style=port] (41) at (14, 1) {};
		\node [style=component] (42) at (14.75, -1) {$e$};
		\node [style=port] (43) at (17.75, 1) {};
		\node [style=port] (44) at (16.75, 1) {};
		\node [style=component] (45) at (16.75, -0.5) {$e$};
		\node [style=component] (46) at (17.75, -0.5) {$e$};
		\node [style=port] (47) at (16, 0) {$=$};
	\end{pgfonlayer}
	\begin{pgfonlayer}{edgelayer}
		\draw [style=wire, bend right] (3) to (5);
		\draw [style=wire, bend left] (3) to (4);
		\draw [style=wire] (9) to (12);
		\draw [style=wire] (11) to (13);
		\draw [style=wire, in=150, out=-135, looseness=1.25] (12) to (11);
		\draw [style=wire] (7) to (8);
		\draw [style=wire] (10) to (6);
		\draw [style=wire, in=30, out=-45, looseness=1.25] (8) to (10);
		\draw [style=wire, in=15, out=-150] (8) to (11);
		\draw [style=wire, in=165, out=-30] (12) to (10);
		\draw [style=wire, bend right] (14) to (16);
		\draw [style=wire, bend left] (14) to (15);
		\draw [style=wire] (17) to (14);
		\draw [style=wire] (20) to (19);
		\draw [style=wire] (21) to (18);
		\draw [style=wire] (31) to (30);
		\draw [style=wire, bend left] (36) to (38);
		\draw [style=wire, bend right] (36) to (37);
		\draw [style=wire] (36) to (3);
		\draw [style=wire, bend left] (39) to (41);
		\draw [style=wire, bend right] (39) to (40);
		\draw [style=wire] (42) to (39);
		\draw [style=wire] (45) to (44);
		\draw [style=wire] (46) to (43);
	\end{pgfonlayer}
\end{tikzpicture}
}%
 \end{gathered}\end{equation}
\begin{equation}\label{nablaep}\begin{gathered} 
\resizebox{!}{2cm}{%
\begin{tikzpicture}
	\begin{pgfonlayer}{nodelayer}
		\node [style=port] (0) at (0.5, 3) {};
		\node [style=port] (1) at (1.25, 0) {};
		\node [style=duplicate] (2) at (1.25, 2) {$\nabla$};
		\node [style=component] (3) at (1.25, 1) {$\varepsilon$};
		\node [style=port] (4) at (2, 3) {};
		\node [style=component] (5) at (3.25, 1.5) {$\varepsilon$};
		\node [style=port] (6) at (4.25, 3) {};
		\node [style=component] (7) at (4.25, 1.5) {$e$};
		\node [style=port] (8) at (3.25, 3) {};
		\node [style=port] (9) at (3.25, 0) {};
		\node [style=component] (10) at (6.75, 1.5) {$\varepsilon$};
		\node [style=port] (11) at (5.75, 3) {};
		\node [style=component] (12) at (5.75, 1.5) {$e$};
		\node [style=port] (13) at (6.75, 3) {};
		\node [style=port] (14) at (6.75, 0) {};
		\node [style=port] (15) at (2.25, 1.5) {$=$};
		\node [style=port] (16) at (5, 1.5) {$+$};
	\end{pgfonlayer}
	\begin{pgfonlayer}{edgelayer}
		\draw [style=wire] (1) to (3);
		\draw [style=wire, bend right] (2) to (4);
		\draw [style=wire, bend left] (2) to (0);
		\draw [style=wire] (3) to (2);
		\draw [style=wire] (7) to (6);
		\draw [style=wire] (9) to (5);
		\draw [style=wire] (5) to (8);
		\draw [style=wire] (12) to (11);
		\draw [style=wire] (14) to (10);
		\draw [style=wire] (10) to (13);
	\end{pgfonlayer}
\end{tikzpicture}
}%
 \end{gathered}\end{equation}

By the naturality of $\nabla$, $\Delta$, $u$ and $e$ we note that for every map $f$, $\oc(f)$ is both a monoid and comonoid morphism. In the original definition of a bialgebra modality in \cite{blute2006differential} it was also required that $u \varepsilon =0$; however this is provable: 

\begin{lemma}\label{epu} For a bialgebra modality $(\oc,  \delta, \varepsilon, \Delta, e, \nabla, u)$, the following diagram commutes: 
\begin{align*}
\begin{array}[c]{c}
 \xymatrixcolsep{5pc}\xymatrix{ K \ar[r]^-{u} \ar[dr]_-{0} & \oc A \ar[d]^-{\varepsilon} \\
& A 
  }
   \end{array} && 
   \begin{array}[c]{c}
\resizebox{!}{1.5cm}{%
\begin{tikzpicture}
	\begin{pgfonlayer}{nodelayer}
		\node [style=port] (0) at (0.5, 1) {};
		\node [style=component] (1) at (0.5, 3) {$u$};
		\node [style=component] (2) at (0.5, 1.75) {$\varepsilon$};
		\node [style=port] (3) at (2, 2.25) {$0$};
		\node [style=port] (4) at (1.25, 2.25) {$=$};
	\end{pgfonlayer}
	\begin{pgfonlayer}{edgelayer}
		\draw [style=wire] (0) to (2);
		\draw [style=wire] (2) to (1);
	\end{pgfonlayer}
\end{tikzpicture}
}%
   \end{array}
\end{align*}
\end{lemma}
\begin{proof} By the naturality of $u$ and $\varepsilon$, and the additive structure we have the following: 
\[\resizebox{!}{2cm}{%
\begin{tikzpicture}
	\begin{pgfonlayer}{nodelayer}
		\node [style=port] (0) at (0.5, 0.25) {};
		\node [style=component] (1) at (0.5, 2.25) {$u$};
		\node [style=component] (2) at (0.5, 1) {$\varepsilon$};
		\node [style=port] (3) at (1.75, 1.75) {$=$};
		\node [style=port] (4) at (3, 0.25) {};
		\node [style=function] (5) at (3, 2) {$0$};
		\node [style=component] (6) at (3, 1) {$\varepsilon$};
		\node [style=component] (7) at (3, 3) {$u$};
		\node [style=port] (8) at (4.25, 1.75) {$=$};
		\node [style=port] (9) at (1.75, 1.5) {\small Nat. of $u$};
		\node [style=port] (10) at (4.25, 1.5) {\small Nat of $\varepsilon$};
		\node [style=port] (11) at (5.25, 0.25) {};
		\node [style=component] (12) at (5.25, 2) {$\varepsilon$};
		\node [style=component] (13) at (5.25, 1) {$0$};
		\node [style=component] (14) at (5.25, 3) {$u$};
		\node [style=port] (15) at (6.25, 1.75) {$=$};
		\node [style=port] (16) at (7, 1.75) {$0$};
	\end{pgfonlayer}
	\begin{pgfonlayer}{edgelayer}
		\draw [style=wire] (0) to (2);
		\draw [style=wire] (2) to (1);
		\draw [style=wire] (4) to (6);
		\draw [style=wire] (6) to (5);
		\draw [style=wire] (5) to (7);
		\draw [style=wire] (11) to (13);
		\draw [style=wire] (13) to (12);
		\draw [style=wire] (12) to (14);
	\end{pgfonlayer}
\end{tikzpicture}
}%
\]
\end{proof} 

Additive bialgebra modalities are bialgebra modalities such that the additive structure of the category and the natural bialgebra structure of the bialgebra modality are compatible via bialgebra convolution. 

\begin{definition} \normalfont An \textbf{additive bialgebra modality} on an additive symmetric monoidal category is a bialgebra modality $(\oc, \delta, \varepsilon, \Delta, e, \nabla, u)$ which is compatible with the additive structure in the sense that the following diagrams commute (for any parallel maps $f$ and $g$):
\begin{equation}\label{}\begin{gathered} \xymatrixcolsep{5pc}\xymatrix{ \oc A \ar[d]_-{\Delta} \ar[r]^-{\oc(f+g)} & \oc B & \oc A \ar[dr]_-{e} \ar[rr]^-{\oc(0)} && \oc B\\
\oc A \otimes \oc A \ar[r]_-{\oc(f) \otimes \oc(g)} & \oc B \otimes \oc B \ar[u]_-{\nabla} & & K \ar[ur]_-{u}
  } \end{gathered}\end{equation}
\end{definition}

In the graphical calculus, the additive bialgebra modality identities are drawn as follows: 
\begin{equation}\label{addbialgeq}\begin{gathered} 
\resizebox{!}{3cm}{%
\begin{tikzpicture}
	\begin{pgfonlayer}{nodelayer}
		\node [style=port] (0) at (3.5, 6.25) {};
		\node [style=duplicate] (1) at (3.5, 5.25) {$\Delta$};
		\node [style=duplicate] (2) at (3.5, 2.75) {$\nabla$};
		\node [style=port] (3) at (3.5, 1.75) {};
		\node [style=function2] (4) at (2.75, 4) {$f$};
		\node [style=function] (5) at (4.25, 4) {$g$};
		\node [style=port] (6) at (2, 4) {$=$};
		\node [style=port] (7) at (0.25, 6) {};
		\node [style=port] (8) at (0.25, 2) {};
		\node [style=function] (10) at (6.75, 4) {$0$};
		\node [style=port] (11) at (6.75, 5) {};
		\node [style=port] (12) at (6.75, 3) {};
		\node [style=port] (13) at (7.75, 4) {$=$};
		\node [style=port] (14) at (8.75, 5.5) {};
		\node [style=port] (15) at (8.75, 2.5) {};
		\node [style=component] (16) at (8.75, 3.5) {$u$};
		\node [style=component] (17) at (8.75, 4.5) {$e$};
		\node [style=port] (18) at (-1, 4.75) {};
		\node [style=port] (19) at (1.5, 4.75) {};
		\node [style=port] (20) at (-1, 3.25) {};
		\node [style=port] (21) at (1.5, 3.25) {};
		\node [style=port] (22) at (-0.5, 4.75) {};
		\node [style=port] (23) at (-0.5, 3.25) {};
		\node [style=port] (24) at (1, 4.75) {};
		\node [style=port] (25) at (1, 3.25) {};
		\node [style=port] (26) at (0.25, 4) {$+$};
		\node [style=component] (27) at (-0.5, 4) {$f$};
		\node [style=component] (28) at (1, 4) {$g$};
		\node [style=port] (29) at (0.25, 4.75) {};
		\node [style=port] (30) at (0.25, 3.25) {};
	\end{pgfonlayer}
	\begin{pgfonlayer}{edgelayer}
		\draw [style=wire] (0) to (1);
		\draw [style=wire] (2) to (3);
		\draw [style=wire, in=90, out=-150] (1) to (4);
		\draw [style=wire, in=150, out=-90] (4) to (2);
		\draw [style=wire, in=30, out=-90] (5) to (2);
		\draw [style=wire, in=90, out=-30] (1) to (5);
		\draw [style=wire] (11) to (10);
		\draw [style=wire] (10) to (12);
		\draw [style=wire] (14) to (17);
		\draw [style=wire] (16) to (15);
		\draw [style=wire] (18) to (20);
		\draw [style=wire] (20) to (21);
		\draw [style=wire] (18) to (19);
		\draw [style=wire] (19) to (21);
		\draw [style=wire] (22) to (27);
		\draw [style=wire] (23) to (27);
		\draw [style=wire] (28) to (25);
		\draw [style=wire] (28) to (24);
		\draw [style=wire] (7) to (29);
		\draw [style=wire] (30) to (8);
	\end{pgfonlayer}
\end{tikzpicture}
}%
 \end{gathered}\end{equation}

We now explain the relation between additive bialgebra modalities and monoidal coalgebra modalities. 

\begin{definition} An \textbf{additive linear category} is a linear category which is also an additive symmetric monoidal category. 
\end{definition}

The monoidal coalgebra modality of an additive linear category induces an additive bialgebra modality where $\nabla$ and $u$ are respectively:
\begin{align*}
\nabla := \xymatrixcolsep{2pc}\xymatrix{\oc A \otimes \oc A  \ar[r]^-{\delta \otimes \delta} & \oc \oc A \otimes \oc \oc A \ar[r]^-{m_\otimes} & \oc(\oc A \otimes \oc A) \ar[rr]^-{\oc(\varepsilon \otimes e + e\otimes \varepsilon)} & & \oc A  } && u :=\xymatrixcolsep{2pc}\xymatrix{K \ar[r]^-{m_K} & \oc K \ar[r]^-{\oc(0)} & \oc A}
\end{align*}

\begin{proposition}\label{monoidaltoaddbialg} The monoidal coalgebra modality of an additive linear category is an additive bialgebra modality with $\nabla$ and $u$ defined as above. 
\end{proposition} 
\begin{proof} See Appendix \ref{montoadd}.
\end{proof} 

The monoidal coalgebra modality structure and the bialgebra modality structure are compatible in the following sense: 

\begin{proposition}\label{monoidalnabla} \cite[Theorem 3.1]{fiore2007differential} In an additive linear category, $u$ and $\nabla$ are $\oc$-coalgbera morphisms, that is, the following diagram commutes: 
  \[  \xymatrixcolsep{5pc}\xymatrix{ \oc A \otimes \oc A \ar[d]_-{\nabla} \ar[r]^-{\delta \otimes \delta} & \oc \oc A \otimes \oc \oc A \ar[r]^-{m_\otimes} & \oc(\oc A \otimes \oc A) \ar[d]^-{\oc(\nabla)} & K \ar[d]_-{u} \ar[r]^-{m_K} & \oc K \ar[d]^-{\oc(u)} \\
  \oc A \ar[rr]_-{\delta} && \oc \oc A & \oc A \ar[r]_-{\delta} & \oc \oc A
  } \]
and also the following diagrams commute: 
  \[  \xymatrixcolsep{3pc}\xymatrix{ \oc A \ar[d]_-{e}  \ar[r]^-{1 \otimes u} & \oc A \otimes \oc B \ar[d]^-{m_\otimes} & \oc A \otimes \oc B \otimes \oc B \ar[dd]_-{1 \otimes \nabla} \ar[r]^-{\Delta \otimes 1 \otimes 1} & \oc A \otimes \oc A \otimes \oc B \otimes \oc B \ar[r]^-{1 \otimes \sigma \otimes 1} & \oc A \otimes \oc B \otimes \oc A \otimes \oc B \ar[d]^-{m_ \otimes \otimes m_\otimes} \\
  K \ar[r]_-{u}  & \oc(A \otimes B) & & & \oc(A \otimes B) \otimes \oc(A \otimes B) \ar[d]_-{\nabla} \\
  && \oc A \otimes \oc B \ar[rr]_-{m_\otimes} && \oc(A \otimes B)
  } \]
\end{proposition} 
\begin{proof} See Appendix \ref{Propappendix}
\end{proof} 

In the graphical calculus, the equalities of the above proposition are drawn as follows: 
\begin{equation}\label{nablam1}\begin{gathered} 
\resizebox{!}{3cm}{%
\begin{tikzpicture}
	\begin{pgfonlayer}{nodelayer}
		\node [style=duplicate] (0) at (1.75, 0) {$\nabla$};
		\node [style=port] (1) at (2.5, 1) {};
		\node [style=port] (2) at (1, 1) {};
		\node [style=port] (3) at (0.5, 1) {};
		\node [style=port] (4) at (3, 1) {};
		\node [style=port] (5) at (3, -0.5) {};
		\node [style=port] (6) at (0.5, -0.5) {};
		\node [style=port] (8) at (1.75, -1.25) {};
		\node [style=port] (10) at (1.75, 1) {};
		\node [style=function] (17) at (5.25, -0.25) {$\delta$};
		\node [style=component] (18) at (5.25, 0.75) {$u$};
		\node [style=port] (20) at (5.25, -1.25) {};
		\node [style=component] (21) at (6.75, 0.75) {$m$};
		\node [style=function] (22) at (6.75, -0.25) {$u$};
		\node [style=port] (23) at (6.75, -1.25) {};
		\node [style=port] (25) at (-0.25, 0) {$=$};
		\node [style=port] (26) at (6, 0.25) {$=$};
		\node [style=component] (27) at (2.5, 2.75) {$\delta$};
		\node [style=component] (28) at (1, 2.75) {$\delta$};
		\node [style=function2] (29) at (1.75, 1.75) {$\bigotimes$};
		\node [style=port] (30) at (1.75, 1) {};
		\node [style=port] (31) at (1, 3.5) {};
		\node [style=port] (32) at (2.5, 3.5) {};
		\node [style=duplicate] (33) at (-1.25, 1) {$\nabla$};
		\node [style=port] (34) at (-0.5, 2) {};
		\node [style=port] (35) at (-2, 2) {};
		\node [style=component] (36) at (-1.25, 0) {$\delta$};
		\node [style=port] (37) at (-1.25, -1.25) {};
	\end{pgfonlayer}
	\begin{pgfonlayer}{edgelayer}
		\draw [style=wire, bend left] (0) to (2);
		\draw [style=wire, bend right] (0) to (1);
		\draw [style=wire] (6) to (5);
		\draw [style=wire] (6) to (3);
		\draw [style=wire] (3) to (4);
		\draw [style=wire] (5) to (4);
		\draw [style=wire] (18) to (17);
		\draw [style=wire] (22) to (23);
		\draw [style=wire] (17) to (20);
		\draw [style=wire] (0) to (8);
		\draw [style=wire] (29) to (30);
		\draw [style=wire, in=0, out=-90, looseness=1.75] (27) to (29);
		\draw [style=wire, in=180, out=-90, looseness=1.75] (28) to (29);
		\draw [style=wire] (31) to (28);
		\draw [style=wire] (32) to (27);
		\draw [style=wire, bend left] (33) to (35);
		\draw [style=wire, bend right] (33) to (34);
		\draw [style=wire] (33) to (36);
		\draw [style=wire] (36) to (37);
		\draw [style=wire] (21) to (22);
	\end{pgfonlayer}
\end{tikzpicture}
}%
 \end{gathered}\end{equation}
 \begin{equation}\label{nablam2}\begin{gathered} 
 \resizebox{!}{3cm}{%
\begin{tikzpicture}
	\begin{pgfonlayer}{nodelayer}
		\node [style=function2] (0) at (7.25, 0) {$\bigotimes$};
		\node [style=port] (1) at (6.25, 2) {};
		\node [style=port] (2) at (8, 1) {};
		\node [style=port] (3) at (7.25, -1) {};
		\node [style=duplicate] (4) at (8, 1) {$\nabla$};
		\node [style=port] (5) at (8.75, 2) {};
		\node [style=port] (6) at (7.25, 2) {};
		\node [style=port] (7) at (9.25, 0.25) {$=$};
		\node [style=port] (9) at (12.25, 2.75) {};
		\node [style=port] (10) at (10.5, 2.75) {};
		\node [style=function2] (11) at (12.5, -0.25) {$\bigotimes$};
		\node [style=function2] (12) at (10.5, -0.25) {$\bigotimes$};
		\node [style=duplicate] (13) at (10.5, 1.75) {$\Delta$};
		\node [style=port] (15) at (13.5, 2.75) {};
		\node [style=port] (16) at (15.75, 1) {};
		\node [style=component] (17) at (17.25, 1) {$u$};
		\node [style=function2] (18) at (16.5, 0) {$\bigotimes$};
		\node [style=port] (19) at (16.5, -1) {};
		\node [style=port] (20) at (18, 0) {$=$};
		\node [style=port] (21) at (19, 1.5) {};
		\node [style=port] (22) at (19, -1.5) {};
		\node [style=component] (23) at (19, -0.5) {$u$};
		\node [style=component] (24) at (19, 0.5) {$e$};
		\node [style=duplicate] (25) at (11.5, -1.25) {$\nabla$};
		\node [style=port] (26) at (11.5, -2) {};
	\end{pgfonlayer}
	\begin{pgfonlayer}{edgelayer}
		\draw [style=wire, in=-90, out=180, looseness=1.25] (0) to (1);
		\draw [style=wire, in=0, out=-90, looseness=1.50] (2) to (0);
		\draw [style=wire] (0) to (3);
		\draw [style=wire, bend left] (4) to (6);
		\draw [style=wire, bend right] (4) to (5);
		\draw [style=wire] (10) to (13);
		\draw [style=wire, in=180, out=-135, looseness=1.25] (13) to (12);
		\draw [style=wire, in=180, out=-30] (13) to (11);
		\draw [style=wire, in=0, out=-90, looseness=0.75] (9) to (12);
		\draw [style=wire, in=0, out=-90, looseness=0.75] (15) to (11);
		\draw [style=wire, in=-90, out=180, looseness=1.50] (18) to (16);
		\draw [style=wire, in=0, out=-90, looseness=1.50] (17) to (18);
		\draw [style=wire] (18) to (19);
		\draw [style=wire] (21) to (24);
		\draw [style=wire] (23) to (22);
		\draw [style=wire, in=135, out=-90, looseness=1.25] (12) to (25);
		\draw [style=wire, in=45, out=-90, looseness=1.25] (11) to (25);
		\draw [style=wire] (25) to (26);
	\end{pgfonlayer}
\end{tikzpicture}
}%
 \end{gathered}\end{equation}

Conversly an additive bialgebra modality induces a monoidal coalgebra modality.  The monoidal structure $m_\otimes$ and $m_K$ are defined 
respectively as follows: 
\[\xymatrixcolsep{5pc}\xymatrix{\oc A \!\otimes\! \oc B \ar[d]_{m_\otimes}^{~~:=} \ar[r]^-{\delta \otimes \delta} & \oc A \!\otimes\! \oc B \ar[r]^-{\oc(1 \otimes u) \otimes \oc(u \otimes 1)} & \oc(\oc A \!\otimes\!\oc B) \!\otimes\! \oc(\oc A \!\otimes\! \oc B) \ar[r]^-{\nabla} & \oc(\oc A \!\otimes\! \oc B) \ar[d]^-{\delta} & ~\\
\oc(A \!\otimes\! B) & \oc(\oc A \!\otimes\! \oc B) \ar[l]^-{\oc(\varepsilon \otimes \varepsilon)}& \oc(\oc(\oc A \!\otimes\! \oc B) \!\otimes\! \oc(\oc A \!\otimes\! \oc B)) \ar[l]^-{\oc(\oc(\varepsilon \otimes e) \otimes \oc(e \otimes \varepsilon))} & \oc \oc(\oc A \!\otimes\! \oc B) \ar[l]^-{\oc(\Delta)}  
  }\]
\[\xymatrixcolsep{5pc}\xymatrix{K \ar[r]^-{u} \ar[d]_{m_K}^{~~:=}  & \oc K \ar[d]^-{\delta} \\   \oc(K) & \oc \oc K \ar[l]^-{\oc(e)}   }\]

\begin{proposition}\label{addbialgtomonoidal} Every additive bialgebra modality is a monoidal coalgebra modality. 
\end{proposition} 
\begin{proof} See Appendix \ref{addtomon}. \end{proof} 

These constructions, between additive bialgebra modalities and monoidal coalgebra modalities, are in fact inverses of each other. 

\begin{theorem}\label{addlinaddbialg} For an additive symmetric monoidal category, monoidal coalgebra modalities correspond bijectively to additive bialgebra modalities. Therefore, the following are equivalent:
\begin{enumerate}[{\em (i)}]
\item An additive linear category;
\item An additive symmetric monoidal category with an additive bialgebra modality. 
\end{enumerate}
\end{theorem} 
\begin{proof} See Appendix \ref{addmonequiv}. \end{proof} 

\section{Differential Categories}

In this section we review (tensor) differential categories, which are are structures over additive symmetric monoidal categories with a coalgebra modality. In particular we revisit the axioms of a deriving transformation for coalgebra modalities and bialgebra modalities. For a full detailed introduction to differential categories, we refer the reader to \cite{blute2006differential,blute2015cartesian}.

\begin{definition} A \textbf{differential category} is an additive symmetric monoidal category with a coalgebra modality $(\oc, \delta, \varepsilon, \Delta, e)$ which comes equipped with a \textbf{deriving transformation} \cite{blute2006differential}, that is, a natural transformation $\mathsf{d}: \oc A \otimes A \to \oc A$ such that the following diagrams commute: 
\begin{description}
\item[{\bf [d.1]}] Constant Rule:
  \[  \xymatrixcolsep{5pc}\xymatrix{\oc A \otimes A \ar[dr]_-{0} \ar[r]^-{\mathsf{d}} & \oc A \ar[d]^-{e} \\
  & K
  } \]
\item[{\bf [d.2]}] Leibniz Rule (or Product Rule):
  \[  \xymatrixcolsep{7pc}\xymatrix{\oc A \otimes A \ar[d]_-{\Delta \otimes 1} \ar[r]^-{\mathsf{d}} & \oc A \ar[d]^-{\Delta} \\
    \oc A \otimes \oc A \otimes A \ar[r]_-{(1 \otimes \mathsf{d}) + (1 \otimes \sigma)(\mathsf{d} \otimes 1)} & \oc A \otimes \oc A
  } \]
\item[{\bf [d.3]}] Linear Rule:
  \[  \xymatrixcolsep{5pc}\xymatrix{\oc A \otimes A \ar[r]^-{\mathsf{d}} \ar[dr]_-{e \otimes 1} & \oc A \ar[d]^-{\varepsilon} \\
  & A   
  } \] 
\item[{\bf [d.4]}] Chain Rule: 
  \[  \xymatrixcolsep{3pc}\xymatrix{\oc A \otimes A \ar[d]_-{\Delta \otimes 1} \ar[rr]^-{\mathsf{d}} && \oc A \ar[d]^-{\delta} \\
    \oc A \otimes \oc A \otimes A \ar[r]_-{\delta \otimes \mathsf{d}} & \oc \oc A \otimes \oc A \ar[r]_-{\mathsf{d}} & \oc \oc A
  } \]
     \item[{\bf [d.5]}] Interchange Rule:
       \[  \xymatrixcolsep{5pc}\xymatrix{\oc A \otimes A \otimes A  \ar[d]_-{\mathsf{d} \otimes 1}\ar[r]^-{1 \otimes \sigma} & \oc A \otimes A \otimes A \ar[r]^-{\mathsf{d} \otimes 1} & \oc A \otimes A \ar[d]^-{\mathsf{d}} \\
         \oc A \otimes A \ar[rr]_-{\mathsf{d}} && \oc A 
  } \]
\end{description}
\end{definition}

In the graphical calculus, the deriving transformation $\mathsf{d}$ is represented as: 
\[\mathsf{d}:= \begin{array}[c]{c} \resizebox{!}{1cm}{%
\begin{tikzpicture}
	\begin{pgfonlayer}{nodelayer}
		\node [style=port] (0) at (1.75, 2.75) {};
		\node [style=port] (1) at (1.25, 1.25) {};
		\node [style=integral] (2) at (1.25, 2) {{\bf =\!=\!=\!=}};
		\node [style=port] (3) at (0.75, 2.75) {};
	\end{pgfonlayer}
	\begin{pgfonlayer}{edgelayer}
		\draw [style=wire, bend right] (2) to (0);
		\draw [style=wire] (1) to (2);
		\draw [style=wire, bend left] (2) to (3);
	\end{pgfonlayer}
\end{tikzpicture}}
   \end{array}\]
and so the deriving transformation axioms \textbf{[d.1]}-\textbf{[d.5]} are drawn as follows: 
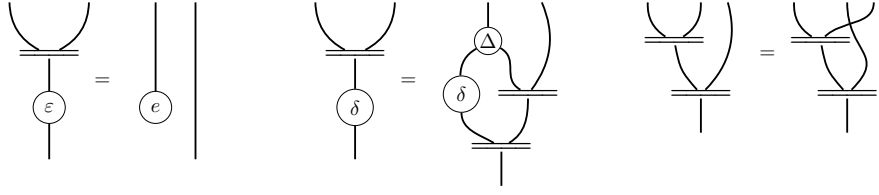
\[\resizebox{!}{2cm}{%
\begin{tikzpicture}
	\begin{pgfonlayer}{nodelayer}
		\node [style=component] (0) at (1.25, 0) {$e$};
		\node [style=port] (1) at (2, 2) {};
		\node [style=differential] (2) at (1.25, 1) {{\bf =\!=\!=\!=}};
		\node [style=port] (3) at (0.5, 2) {};
		\node [style=port] (4) at (2.25, 1) {$=$};
		\node [style=port] (5) at (3, 1) {$0$};
		\node [style=differential] (6) at (7, 1) {{\bf =\!=\!=\!=}};
		\node [style=port] (7) at (7.75, 2) {};
		\node [style=port] (8) at (6.25, 2) {};
		\node [style=port] (9) at (6.25, -0.75) {};
		\node [style=duplicate] (10) at (7, 0.25) {$\Delta$};
		\node [style=port] (11) at (7.75, -0.75) {};
		\node [style=port] (12) at (9.5, 2) {};
		\node [style=differential] (13) at (9, 0) {{\bf =\!=\!=\!=}};
		\node [style=port] (14) at (10.75, 2) {};
		\node [style=duplicate] (15) at (9.5, 1) {$\Delta$};
		\node [style=port] (16) at (9, -0.75) {};
		\node [style=port] (17) at (10.5, -0.75) {};
		\node [style=port] (18) at (13.5, 2) {};
		\node [style=differential] (19) at (13.25, 0) {{\bf =\!=\!=\!=}};
		\node [style=port] (20) at (13.25, -0.75) {};
		\node [style=port] (21) at (11.75, -0.75) {};
		\node [style=port] (22) at (12.25, 2) {};
		\node [style=duplicate] (23) at (12.25, 1) {$\Delta$};
		\node [style=port] (24) at (8, 0.5) {$=$};
		\node [style=port] (25) at (11.25, 0.5) {$+$};
	\end{pgfonlayer}
	\begin{pgfonlayer}{edgelayer}
		\draw [style=wire, bend right] (2) to (1);
		\draw [style=wire, bend left] (2) to (3);
		\draw [style=wire] (2) to (0);
		\draw [style=wire, bend right] (6) to (7);
		\draw [style=wire, bend left] (6) to (8);
		\draw [style=wire, bend right] (10) to (9);
		\draw [style=wire, bend left] (10) to (11);
		\draw [style=wire] (6) to (10);
		\draw [style=wire, in=-90, out=30] (13) to (14);
		\draw [style=wire, in=150, out=-150, looseness=1.50] (15) to (13);
		\draw [style=wire] (12) to (15);
		\draw [style=wire] (13) to (16);
		\draw [style=wire, bend left, looseness=1.25] (15) to (17);
		\draw [style=wire, bend right] (19) to (18);
		\draw [style=wire, in=91, out=-135, looseness=0.75] (23) to (21);
		\draw [style=wire, in=150, out=-30] (23) to (19);
		\draw [style=wire] (22) to (23);
		\draw [style=wire] (19) to (20);
	\end{pgfonlayer}
\end{tikzpicture}
}%
\]
\[\resizebox{!}{2.5cm}{%
\begin{tikzpicture}
	\begin{pgfonlayer}{nodelayer}
		\node [style=port] (28) at (2.5, 1.75) {};
		\node [style=port] (29) at (3.25, 1.75) {};
		\node [style=port] (30) at (3.25, -1.25) {};
		\node [style=component] (31) at (2.5, -0.25) {$e$};
		\node [style=port] (32) at (1.25, 1.75) {};
		\node [style=port] (33) at (-0.25, 1.75) {};
		\node [style=component] (34) at (0.5, -0.25) {$\varepsilon$};
		\node [style=differential] (35) at (0.5, 0.75) {{\bf =\!=\!=\!=}};
		\node [style=port] (36) at (0.5, -1.25) {};
		\node [style=port] (37) at (1.5, 0.25) {$=$};
		\node [style=port] (38) at (7, 1.75) {};
		\node [style=differential] (39) at (6.25, 0.75) {{\bf =\!=\!=\!=}};
		\node [style=port] (40) at (6.25, -1.25) {};
		\node [style=port] (41) at (5.5, 1.75) {};
		\node [style=component] (42) at (6.25, -0.25) {$\delta$};
		\node [style=component] (43) at (8.25, 0) {$\delta$};
		\node [style=duplicate] (44) at (8.75, 1) {$\Delta$};
		\node [style=port] (45) at (8.75, 1.75) {};
		\node [style=differential] (46) at (9.5, 0) {{\bf =\!=\!=\!=}};
		\node [style=port] (47) at (9.75, 1.75) {};
		\node [style=differential] (48) at (9, -1) {{\bf =\!=\!=\!=}};
		\node [style=port] (49) at (9, -1.75) {};
		\node [style=port] (50) at (7.25, 0.25) {$=$};
		\node [style=port] (51) at (13.25, 1.75) {};
		\node [style=port] (52) at (12.75, 1.75) {};
		\node [style=port] (53) at (12.75, -0.75) {};
		\node [style=port] (54) at (11.75, 1.75) {};
		\node [style=codifferential] (55) at (12.75, 0) {{\bf =\!=\!=\!=}};
		\node [style=codifferential] (56) at (12.25, 1) {{\bf =\!=\!=\!=}};
		\node [style=port] (57) at (14.5, 1.75) {};
		\node [style=codifferential] (58) at (15, 1) {{\bf =\!=\!=\!=}};
		\node [style=codifferential] (59) at (15.5, 0) {{\bf =\!=\!=\!=}};
		\node [style=port] (60) at (15.5, -0.75) {};
		\node [style=port] (61) at (15.5, 1.75) {};
		\node [style=port] (62) at (16, 1.75) {};
		\node [style=port] (63) at (14, 0.75) {$=$};
	\end{pgfonlayer}
	\begin{pgfonlayer}{edgelayer}
		\draw [style=wire] (28) to (31);
		\draw [style=wire] (29) to (30);
		\draw [style=wire, bend right] (35) to (32);
		\draw [style=wire, bend left] (35) to (33);
		\draw [style=wire] (35) to (34);
		\draw [style=wire] (34) to (36);
		\draw [style=wire, bend right] (39) to (38);
		\draw [style=wire, bend left] (39) to (41);
		\draw [style=wire] (39) to (42);
		\draw [style=wire] (42) to (40);
		\draw [style=wire, bend right] (44) to (43);
		\draw [style=wire] (45) to (44);
		\draw [style=wire, bend right] (46) to (47);
		\draw [style=wire, in=150, out=-30, looseness=1.25] (44) to (46);
		\draw [style=wire] (48) to (49);
		\draw [style=wire, in=30, out=-90] (46) to (48);
		\draw [style=wire, in=150, out=-90] (43) to (48);
		\draw [style=wire, bend left=15, looseness=1.25] (55) to (56);
		\draw [style=wire, bend right] (55) to (51);
		\draw [style=wire] (53) to (55);
		\draw [style=wire, bend left] (56) to (54);
		\draw [style=wire, bend right] (56) to (52);
		\draw [style=wire, bend left=15, looseness=1.25] (59) to (58);
		\draw [style=wire] (60) to (59);
		\draw [style=wire, bend left] (58) to (57);
		\draw [style=wire, in=45, out=-90] (62) to (58);
		\draw [style=wire, in=45, out=-90, looseness=1.50] (61) to (59);
	\end{pgfonlayer}
\end{tikzpicture}
}%
\]

The coKleisli maps of the coalgebra modality of a differential category are important: these maps are of the form $f: \oc A \to B$ are considered as \emph{smooth maps}. Amongst these are the {\bf linear maps} $\varepsilon g : \oc A \to B$ where $g: A \to B$. The differential of a smooth map $f: \oc A \to B$ is the map $\mathsf{D}[f]: \oc A \otimes A \to B$ defined by precomposing with the deriving transformation, $\mathsf{D}[f]=\mathsf{d}f$. The first axiom {\bf [d.1]} states that the derivative of a constant map is zero. The second axiom {\bf [d.2]} is the Leibniz rule or the product rule for differentiation. The third axiom {\bf [d.3]} says that the derivative of a linear map is constant. The fourth axiom {\bf [d.4]} is the chain rule. The last axiom {\bf [d.5]} is the interchange law, which naively states that differentiating with respect to $x$ then $y$ is the same as differentiation with respect to $y$ then $x$. It should be noted that {\bf [d.5]} was not originally a requirement in \cite{blute2006differential} but was later added to the definition to ensure that the coKleisli category of a differential category was a Cartesian differential category \cite{blute2009cartesian}. 

Our first revision of differential categories is that the constant rule {\bf [d.1]} is in fact derivable: 

\begin{lemma}\label{dericon} For a coalgebra modality $(\oc, \delta, \varepsilon, \Delta, e)$ on additive symmetric monoidal category, any natural transformation $\mathsf{d}: \oc A \otimes A \to \oc A$ satisfies the constant rule {\bf [d.1]}. 
\end{lemma}
\begin{proof} By naturality of $e$ and $\mathsf{d}$, and the additive structure, we have the following equalities: 
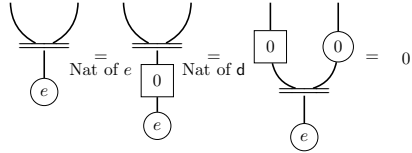
\[\resizebox{!}{2cm}{%
\begin{tikzpicture}
	\begin{pgfonlayer}{nodelayer}
		\node [style=component] (0) at (1.25, 0) {$e$};
		\node [style=port] (1) at (2, 2) {};
		\node [style=differential] (2) at (1.25, 1) {{\bf =\!=\!=\!=}};
		\node [style=port] (3) at (0.5, 2) {};
		\node [style=port] (4) at (8.5, 0.75) {$=$};
		\node [style=port] (5) at (9.25, 0.75) {$0$};
		\node [style=component] (6) at (3.75, -0.75) {$e$};
		\node [style=port] (7) at (4.5, 2) {};
		\node [style=differential] (8) at (3.75, 1) {{\bf =\!=\!=\!=}};
		\node [style=port] (9) at (3, 2) {};
		\node [style=function] (10) at (3.75, 0.25) {$0$};
		\node [style=port] (11) at (2.5, 0.75) {$=$};
		\node [style=port] (12) at (5, 0.75) {$=$};
		\node [style=component] (13) at (7, -1) {$e$};
		\node [style=component] (14) at (7.75, 1) {$0$};
		\node [style=differential] (15) at (7, 0) {{\bf =\!=\!=\!=}};
		\node [style=function] (16) at (6.25, 1) {$0$};
		\node [style=port] (17) at (6.25, 2) {};
		\node [style=port] (18) at (7.75, 2) {};
		\node [style=port] (19) at (2.5, 0.5) {Nat of $e$};
		\node [style=port] (20) at (5, 0.5) {Nat of $\mathsf{d}$};
	\end{pgfonlayer}
	\begin{pgfonlayer}{edgelayer}
		\draw [style=wire, bend right] (2) to (1);
		\draw [style=wire, bend left] (2) to (3);
		\draw [style=wire] (2) to (0);
		\draw [style=wire, bend right] (8) to (7);
		\draw [style=wire, bend left] (8) to (9);
		\draw [style=wire] (8) to (10);
		\draw [style=wire] (10) to (6);
		\draw [style=wire, bend right] (15) to (14);
		\draw [style=wire, bend left] (15) to (16);
		\draw [style=wire] (15) to (13);
		\draw [style=wire] (18) to (14);
		\draw [style=wire] (17) to (16);
	\end{pgfonlayer}
\end{tikzpicture}
}%
\]
\end{proof} 

\begin{corollary}\label{corderi} For a coalgebra modality $(\oc, \delta, \varepsilon, \Delta, e)$ on additive symmetric monoidal category, the following are equivalent for a natural transformation $\mathsf{d}: \oc A \otimes A \to \oc A$: 
\begin{enumerate}[{\em (i)}]
\item $\mathsf{d}$ is a deriving transformation; 
\item $\mathsf{d}$ satisfies the product rule {\bf [d.2]}, the linear rule {\bf [d.3]}, the chain rule {\bf [d.4]}, and the interchange rule {\bf [d.5]}. 
\end{enumerate}
\end{corollary} 

Let us now consider the relation between deriving transformations and bialgebra modalities. This is captured by the $\nabla$-rule \cite{blute2006differential}:

\begin{definition} For a bialgebra modality $(\oc,  \delta, \varepsilon, \Delta, e, \nabla, u)$, a natural transformation ${\mathsf{d}: \oc A \otimes A \to \oc A}$ is said to satisfy the \textbf{$\nabla$-rule} if the following diagram commutes: \\

\noindent {\bf [d.$\nabla$]} $\nabla$-Rule: 
\begin{align*}
\begin{array}[c]{c}
\xymatrixcolsep{5pc}\xymatrix{\oc A \otimes \oc A \otimes A \ar[d]_-{1 \otimes \mathsf{d}} \ar[r]^-{\nabla \otimes 1} & \oc A \otimes A \ar[d]^-{\mathsf{d}} \\
\oc A \otimes \oc A \ar[r]_-{\nabla} & \oc A   
  } 
   \end{array} &&
   \begin{array}[c]{c}
\resizebox{!}{2cm}{%
\begin{tikzpicture}
	\begin{pgfonlayer}{nodelayer}
		\node [style=port] (0) at (0.5, 3) {};
		\node [style=port] (1) at (2.5, 3) {};
		\node [style=codifferential] (2) at (1.5, 0.75) {{\bf =\!=\!=}};
		\node [style=port] (3) at (1.5, -0.25) {};
		\node [style=port] (4) at (1.5, 3) {};
		\node [style=duplicate] (5) at (1, 1.75) {$\nabla$};
		\node [style=port] (6) at (4.5, -0.25) {};
		\node [style=port] (7) at (4.5, 3) {};
		\node [style=duplicate] (8) at (4.5, 0.75) {$\nabla$};
		\node [style=port] (9) at (5.5, 3) {};
		\node [style=codifferential] (10) at (5, 1.75) {{\bf =\!=\!=}};
		\node [style=port] (11) at (3.5, 3) {};
		\node [style=port] (12) at (3, 1.5) {$=$};
	\end{pgfonlayer}
	\begin{pgfonlayer}{edgelayer}
		\draw [style=wire, bend left=15, looseness=1.25] (2) to (5);
		\draw [style=wire, bend right=15] (2) to (1);
		\draw [style=wire] (3) to (2);
		\draw [style=wire, bend left=15] (5) to (0);
		\draw [style=wire, bend right=15] (5) to (4);
		\draw [style=wire, bend right=15, looseness=1.25] (8) to (10);
		\draw [style=wire, bend left=15] (8) to (11);
		\draw [style=wire] (6) to (8);
		\draw [style=wire, bend right=15] (10) to (9);
		\draw [style=wire, bend left=15] (10) to (7);
	\end{pgfonlayer}
\end{tikzpicture}
}%
   \end{array}
\end{align*}
\end{definition}

We observe that the $\nabla$-rule implies the interchange rule:

\begin{lemma}\label{diden} For a bialgebra modality $(\oc,  \delta, \varepsilon, \Delta, e, \nabla, u)$, for any natural transformation ${\mathsf{d}: \oc A \otimes A \to \oc A}$ which satisfies the $\nabla$-rule, {\bf [d.$\nabla$]}, the following diagram commutes: 
\begin{align*}
\begin{array}[c]{c}
 \xymatrixcolsep{5pc}\xymatrix{\oc A \otimes A \ar[drr]_-{\mathsf{d}}\ar[r]^-{1 \otimes u \otimes 1} & \oc A \otimes \oc A \otimes  A \ar[r]^-{1 \otimes \mathsf{d}} & \oc A \otimes \oc A \ar[d]^-{\nabla} \\
 & & \oc A
  } 
   \end{array}&&
   \begin{array}[c]{c}
\resizebox{!}{2.5cm}{%
\begin{tikzpicture}
	\begin{pgfonlayer}{nodelayer}
		\node [style=port] (6) at (3.75, 0.25) {};
		\node [style=component] (7) at (3.75, 3.5) {$u$};
		\node [style=duplicate] (8) at (3.75, 1.25) {$\nabla$};
		\node [style=port] (9) at (4.75, 3.5) {};
		\node [style=codifferential] (10) at (4.25, 2.25) {{\bf =\!=\!=}};
		\node [style=port] (11) at (2.75, 3.5) {};
		\node [style=port] (13) at (1.5, 2.25) {};
		\node [style=port] (14) at (1, 0.75) {};
		\node [style=integral] (15) at (1, 1.5) {{\bf =\!=\!=\!=}};
		\node [style=port] (16) at (0.5, 2.25) {};
		\node [style=port] (17) at (2.5, 1.5) {$=$};
	\end{pgfonlayer}
	\begin{pgfonlayer}{edgelayer}
		\draw [style=wire, bend right=15, looseness=1.25] (8) to (10);
		\draw [style=wire, bend left=15] (8) to (11);
		\draw [style=wire] (6) to (8);
		\draw [style=wire, bend right=15] (10) to (9);
		\draw [style=wire, bend left=15] (10) to (7);
		\draw [style=wire, bend right] (15) to (13);
		\draw [style=wire] (14) to (15);
		\draw [style=wire, bend left] (15) to (16);
	\end{pgfonlayer}
\end{tikzpicture}
}%
   \end{array}
\end{align*}
\end{lemma}
\begin{proof} Using {\bf [d.$\nabla$]} and the monoid unit identity, we obtain the following: 
\[\resizebox{!}{2.5cm}{%
\begin{tikzpicture}
	\begin{pgfonlayer}{nodelayer}
		\node [style=port] (0) at (0.5, 3) {};
		\node [style=port] (1) at (2.5, 3) {};
		\node [style=codifferential] (2) at (1.5, 0.75) {{\bf =\!=\!=}};
		\node [style=port] (3) at (1.5, -0.25) {};
		\node [style=component] (4) at (1.5, 3) {$u$};
		\node [style=duplicate] (5) at (1, 1.75) {$\nabla$};
		\node [style=port] (6) at (4.5, -0.25) {};
		\node [style=component] (7) at (4.5, 3) {$u$};
		\node [style=duplicate] (8) at (4.5, 0.75) {$\nabla$};
		\node [style=port] (9) at (5.5, 3) {};
		\node [style=codifferential] (10) at (5, 1.75) {{\bf =\!=\!=}};
		\node [style=port] (11) at (3.5, 3) {};
		\node [style=port] (12) at (3, 1.5) {$=$};
		\node [style=port] (13) at (-1, 2.25) {};
		\node [style=port] (14) at (-1.5, 0.75) {};
		\node [style=integral] (15) at (-1.5, 1.5) {{\bf =\!=\!=\!=}};
		\node [style=port] (16) at (-2, 2.25) {};
		\node [style=port] (17) at (0, 1.5) {$=$};
		\node [style=port] (18) at (3, 1) {\textbf{[d.$\nabla$]}};
		\node [style=port] (19) at (0, 1) {(\ref{moneq})};
	\end{pgfonlayer}
	\begin{pgfonlayer}{edgelayer}
		\draw [style=wire, bend left=15, looseness=1.25] (2) to (5);
		\draw [style=wire, bend right=15] (2) to (1);
		\draw [style=wire] (3) to (2);
		\draw [style=wire, bend left=15] (5) to (0);
		\draw [style=wire, bend right=15] (5) to (4);
		\draw [style=wire, bend right=15, looseness=1.25] (8) to (10);
		\draw [style=wire, bend left=15] (8) to (11);
		\draw [style=wire] (6) to (8);
		\draw [style=wire, bend right=15] (10) to (9);
		\draw [style=wire, bend left=15] (10) to (7);
		\draw [style=wire, bend right] (15) to (13);
		\draw [style=wire] (14) to (15);
		\draw [style=wire, bend left] (15) to (16);
	\end{pgfonlayer}
\end{tikzpicture}
}%
\]
\end{proof} 

\begin{lemma} For a bialgebra modality $(\oc,  \delta, \varepsilon, \Delta, e, \nabla, u)$, any natural transformation ${\mathsf{d}: \oc A \otimes A \to \oc A}$ which satisfies the $\nabla$-rule, {\bf [d.$\nabla$]}, also satisfies the interchange rule, {\bf [d.5]}. \end{lemma}  
\begin{proof} Using Lemma \ref{diden}, and both associativity and commutativity of the multiplication, we have the following equality:
\[\resizebox{!}{3cm}{%
\begin{tikzpicture}
	\begin{pgfonlayer}{nodelayer}
		\node [style=port] (0) at (2.5, 7.25) {};
		\node [style=port] (1) at (1.5, 7.25) {};
		\node [style=port] (2) at (1.5, 4) {};
		\node [style=port] (3) at (0.5, 7.25) {};
		\node [style=codifferential] (4) at (1.5, 5) {{\bf =\!=\!=}};
		\node [style=codifferential] (5) at (1, 6) {{\bf =\!=\!=}};
		\node [style=port] (6) at (4.75, 3.5) {};
		\node [style=duplicate] (7) at (4.75, 4.5) {$\nabla$};
		\node [style=port] (8) at (3.75, 7.25) {};
		\node [style=codifferential] (9) at (5, 6.25) {{\bf =\!=\!=}};
		\node [style=port] (10) at (5.5, 7.25) {};
		\node [style=component] (11) at (4.5, 7) {$u$};
		\node [style=duplicate] (12) at (4.25, 5.25) {$\nabla$};
		\node [style=port] (13) at (7.25, 7.25) {};
		\node [style=codifferential] (14) at (6.75, 6.25) {{\bf =\!=\!=}};
		\node [style=component] (15) at (6.25, 7) {$u$};
		\node [style=port] (16) at (3, 5.75) {$=$};
		\node [style=port] (17) at (7.75, 5.75) {$=$};
		\node [style=port] (18) at (7.75, 5.25) {(\ref{moneq})};
		\node [style=port] (19) at (9.5, 3.5) {};
		\node [style=duplicate] (20) at (9.5, 4.5) {$\nabla$};
		\node [style=port] (21) at (8.5, 7.25) {};
		\node [style=codifferential] (22) at (9.75, 6.25) {{\bf =\!=\!=}};
		\node [style=port] (23) at (10.25, 7.25) {};
		\node [style=component] (24) at (9.25, 7) {$u$};
		\node [style=port] (25) at (12, 7.25) {};
		\node [style=codifferential] (26) at (11.5, 6.25) {{\bf =\!=\!=}};
		\node [style=component] (27) at (11, 7) {$u$};
		\node [style=duplicate] (28) at (10.5, 5.25) {$\nabla$};
		\node [style=port] (29) at (14.25, 2.75) {};
		\node [style=duplicate] (30) at (14.25, 3.5) {$\nabla$};
		\node [style=port] (31) at (13.25, 7.25) {};
		\node [style=codifferential] (32) at (16.5, 6.25) {{\bf =\!=\!=}};
		\node [style=port] (33) at (17, 7.25) {};
		\node [style=component] (34) at (16, 7) {$u$};
		\node [style=port] (35) at (15.25, 7.25) {};
		\node [style=codifferential] (36) at (14.75, 6.25) {{\bf =\!=\!=}};
		\node [style=component] (37) at (14.25, 7) {$u$};
		\node [style=duplicate] (38) at (15.5, 4.25) {$\nabla$};
		\node [style=port] (39) at (12.5, 5.75) {$=$};
		\node [style=port] (40) at (12.5, 5.25) {(\ref{moneq})};
		\node [style=port] (41) at (19.75, 2.75) {};
		\node [style=duplicate] (42) at (19.75, 3.5) {$\nabla$};
		\node [style=port] (43) at (18.25, 7.25) {};
		\node [style=codifferential] (44) at (21.5, 6.25) {{\bf =\!=\!=}};
		\node [style=port] (45) at (22, 7.25) {};
		\node [style=component] (46) at (21, 7) {$u$};
		\node [style=port] (47) at (20.25, 7.25) {};
		\node [style=codifferential] (48) at (19.75, 6.25) {{\bf =\!=\!=}};
		\node [style=component] (49) at (19.25, 7) {$u$};
		\node [style=duplicate] (50) at (18.75, 4.75) {$\nabla$};
		\node [style=port] (51) at (17.5, 5.75) {$=$};
		\node [style=port] (52) at (17.5, 5.25) {(\ref{moneq})};
		\node [style=port] (53) at (25.25, 7.25) {};
		\node [style=port] (54) at (24.25, 7.25) {};
		\node [style=differential] (55) at (24.25, 5) {{\bf =\!=\!=}};
		\node [style=differential] (56) at (23.75, 6) {{\bf =\!=\!=}};
		\node [style=port] (57) at (24.25, 4) {};
		\node [style=port] (58) at (23.25, 7.25) {};
		\node [style=port] (59) at (22.5, 5.75) {$=$};
		\node [style=port] (60) at (22.5, 5.25) {Lem \ref{diden}};
		\node [style=port] (61) at (3, 5.25) {Lem \ref{diden}};
		\node [style=port] (62) at (7.75, 4.75) {assoc.};
		\node [style=port] (63) at (12.5, 4.75) {comm.};
		\node [style=port] (64) at (17.5, 4.75) {assoc.};
	\end{pgfonlayer}
	\begin{pgfonlayer}{edgelayer}
		\draw [style=wire, bend left=15, looseness=1.25] (4) to (5);
		\draw [style=wire, bend right] (4) to (0);
		\draw [style=wire] (2) to (4);
		\draw [style=wire, bend left] (5) to (3);
		\draw [style=wire, bend right] (5) to (1);
		\draw [style=wire] (6) to (7);
		\draw [style=wire, bend left, looseness=1.25] (10) to (9);
		\draw [style=wire, bend right] (11) to (9);
		\draw [style=wire, in=150, out=-90, looseness=1.25] (12) to (7);
		\draw [style=wire, bend left, looseness=1.25] (13) to (14);
		\draw [style=wire, bend right] (15) to (14);
		\draw [style=wire, in=-90, out=15] (7) to (14);
		\draw [style=wire, in=-90, out=45, looseness=1.50] (12) to (9);
		\draw [style=wire, in=-105, out=135, looseness=0.75] (12) to (8);
		\draw [style=wire] (19) to (20);
		\draw [style=wire, bend left, looseness=1.25] (23) to (22);
		\draw [style=wire, bend right] (24) to (22);
		\draw [style=wire, bend left, looseness=1.25] (25) to (26);
		\draw [style=wire, bend right] (27) to (26);
		\draw [style=wire, in=-90, out=127, looseness=1.25] (28) to (22);
		\draw [style=wire, in=-90, out=150] (20) to (21);
		\draw [style=wire, in=-90, out=37] (20) to (28);
		\draw [style=wire, in=-90, out=45, looseness=1.25] (28) to (26);
		\draw [style=wire] (29) to (30);
		\draw [style=wire, bend left, looseness=1.25] (33) to (32);
		\draw [style=wire, bend right] (34) to (32);
		\draw [style=wire, bend left, looseness=1.25] (35) to (36);
		\draw [style=wire, bend right] (37) to (36);
		\draw [style=wire, in=-90, out=127, looseness=1.25] (38) to (32);
		\draw [style=wire, in=-90, out=150] (30) to (31);
		\draw [style=wire, in=-90, out=37] (30) to (38);
		\draw [style=wire, in=-90, out=45, looseness=1.25] (38) to (36);
		\draw [style=wire] (41) to (42);
		\draw [style=wire, bend left, looseness=1.25] (45) to (44);
		\draw [style=wire, bend right] (46) to (44);
		\draw [style=wire, bend left, looseness=1.25] (47) to (48);
		\draw [style=wire, bend right] (49) to (48);
		\draw [style=wire, in=150, out=-105, looseness=0.75] (43) to (50);
		\draw [style=wire, in=-90, out=30] (50) to (44);
		\draw [style=wire, in=30, out=-105, looseness=1.50] (48) to (42);
		\draw [style=wire, in=-90, out=165] (42) to (50);
		\draw [style=wire, bend left=15, looseness=1.25] (55) to (56);
		\draw [style=wire] (57) to (55);
		\draw [style=wire, bend left=15] (56) to (58);
		\draw [style=wire, in=45, out=-90, looseness=1.50] (54) to (55);
		\draw [style=wire, bend left=15] (53) to (56);
	\end{pgfonlayer}
\end{tikzpicture}
}%
\]
\end{proof} 
\begin{corollary} For a bialgebra modality $(\oc, \delta, \varepsilon, \Delta, e, \nabla, u)$, the following are equivalent for a natural transformation $\mathsf{d}: \oc A \otimes A \to \oc A$: 
\begin{enumerate}[{\em (i)}]
\item $\mathsf{d}$ is a deriving transformation which satisfies the $\nabla$-rule {\bf [d.$\nabla$]}; 
\item $\mathsf{d}$ satisfies the product rule {\bf [d.2]}, the linear rule {\bf [d.3]}, the chain rule {\bf [d.4]}, and the $\nabla$-rule {\bf [d.$\nabla$]}.  
\end{enumerate}
\end{corollary} 

\section{Coderelictions}\label{codersec}

For a bialgebra modality there is a natural alternative way to introduce differentiation:

\begin{definition}\label{cordef} \normalfont A \textbf{codereliction} \cite{blute2006differential} for a bialgebra modality $(\oc,  \delta, \varepsilon, \Delta, e, \nabla, u)$ is a natural transformation $\eta: A \to \oc A$, such that the following diagrams commute: 
\begin{enumerate}[{\bf [dC.1]}]
 \item Constant Rule: 
  \[  \xymatrixcolsep{5pc}\xymatrix{A \ar[dr]_-{0} \ar[r]^-{\eta} & \oc A \ar[d]^-{e}Ê\\
  & K
  } \]   
\item Product Rule: 
  \[  \xymatrixcolsep{5pc}\xymatrix{A \ar[dr]_-{\eta \otimes u + u \otimes \eta} \ar[r]^-{\eta} & \oc A \ar[d]^-{\Delta} \\
  & \oc A \otimes \oc A  
  } \]
\item Linear Rule:
  \[  \xymatrixcolsep{5pc}\xymatrix{A \ar@{=}[dr]^-{}  \ar[r]^-{\eta} & \oc A \ar[d]^-{\varepsilon} \\
  & A  
  } \] 
\item Chain Rule: 
  \[  \xymatrixcolsep{3pc}\xymatrix{\oc A \otimes A \ar[d]_-{\Delta \otimes \eta}  \ar[r]^-{1 \otimes \eta} & \oc A \otimes \oc A \ar[r]^-{\nabla} & \oc A \ar[dd]^-{\delta} \\
    \oc A \otimes \oc A \otimes \oc A \ar[d]_-{1 \otimes \nabla} \\
     \oc A \otimes \oc A \ar[r]_-{\delta \otimes \eta} & \oc \oc A \otimes \oc \oc A \ar[r]_-{\nabla} & \oc  \oc A 
  } \]
\end{enumerate}
\end{definition}

In the graphical calculus, the codereliction axioms are drawn as follows: 
\[\resizebox{!}{3cm}{%
\begin{tikzpicture}
	\begin{pgfonlayer}{nodelayer}
		\node [style=port] (29) at (4.5, 3.75) {};
		\node [style=port] (30) at (4.5, 6.75) {};
		\node [style=component] (34) at (3, 5.75) {$\eta$};
		\node [style=port] (36) at (3, 6.75) {};
		\node [style=port] (37) at (3.75, 5.25) {$=$};
		\node [style=port] (50) at (9, 4.75) {$=$};
		\node [style=port] (51) at (-9, 6.75) {};
		\node [style=component] (52) at (-9, 6) {$\eta$};
		\node [style=component] (53) at (-9, 5) {$e$};
		\node [style=port] (54) at (-8.25, 5.5) {$=$};
		\node [style=port] (55) at (-7.5, 5.5) {$0$};
		\node [style=port] (57) at (-5.75, 3.75) {};
		\node [style=port] (58) at (-5, 6.75) {};
		\node [style=duplicate] (59) at (-5, 4.75) {$\Delta$};
		\node [style=component] (60) at (-5, 5.75) {$\eta$};
		\node [style=port] (61) at (-4.25, 3.75) {};
		\node [style=component] (62) at (-3, 5.25) {$\eta$};
		\node [style=port] (63) at (-2, 3.75) {};
		\node [style=component] (64) at (-2, 5.25) {$u$};
		\node [style=port] (65) at (-3, 3.75) {};
		\node [style=port] (66) at (-3, 6.75) {};
		\node [style=component] (67) at (0.5, 5.25) {$\eta$};
		\node [style=port] (68) at (-0.5, 3.75) {};
		\node [style=component] (69) at (-0.5, 5.25) {$u$};
		\node [style=port] (70) at (0.5, 3.75) {};
		\node [style=port] (71) at (0.5, 6.75) {};
		\node [style=port] (72) at (-4, 5.25) {$=$};
		\node [style=port] (73) at (-1.25, 5.25) {$+$};
		\node [style=component] (75) at (3, 4.75) {$\varepsilon$};
		\node [style=port] (76) at (3, 3.75) {};
		\node [style=port] (77) at (10.5, 7) {};
		\node [style=duplicate] (78) at (10.5, 6) {$\Delta$};
		\node [style=component] (79) at (10, 5) {$\delta$};
		\node [style=port] (80) at (11.75, 7) {};
		\node [style=component] (81) at (11.75, 6) {$\eta$};
		\node [style=duplicate] (82) at (11.25, 5) {$\nabla$};
		\node [style=duplicate] (83) at (10.75, 3.25) {$\nabla$};
		\node [style=component] (84) at (11.25, 4) {$\eta$};
		\node [style=port] (85) at (10.75, 2.5) {};
		\node [style=component] (86) at (8.25, 5.75) {$\eta$};
		\node [style=port] (87) at (6.75, 7) {};
		\node [style=port] (88) at (8.25, 7) {};
		\node [style=component] (89) at (7.75, 3.25) {$\delta$};
		\node [style=duplicate] (90) at (7.75, 4.25) {$\nabla$};
		\node [style=port] (91) at (7.75, 2.5) {};
	\end{pgfonlayer}
	\begin{pgfonlayer}{edgelayer}
		\draw [style=wire] (29) to (30);
		\draw [style=wire] (34) to (36);
		\draw [style=wire] (51) to (52);
		\draw [style=wire] (52) to (53);
		\draw [style=wire] (58) to (60);
		\draw [style=wire, bend left] (59) to (61);
		\draw [style=wire, bend right] (59) to (57);
		\draw [style=wire] (60) to (59);
		\draw [style=wire] (64) to (63);
		\draw [style=wire] (66) to (62);
		\draw [style=wire] (62) to (65);
		\draw [style=wire] (69) to (68);
		\draw [style=wire] (71) to (67);
		\draw [style=wire] (67) to (70);
		\draw [style=wire] (75) to (76);
		\draw [style=wire] (34) to (75);
		\draw [style=wire, bend right] (78) to (79);
		\draw [style=wire] (77) to (78);
		\draw [style=wire, bend right=15, looseness=1.25] (82) to (81);
		\draw [style=wire] (81) to (80);
		\draw [style=wire, in=135, out=-30] (78) to (82);
		\draw [style=wire, bend right=15, looseness=1.25] (83) to (84);
		\draw [style=wire] (82) to (84);
		\draw [style=wire, bend right=15, looseness=1.25] (79) to (83);
		\draw [style=wire] (83) to (85);
		\draw [style=wire, bend right=15, looseness=1.25] (90) to (86);
		\draw [style=wire, bend left=15] (90) to (87);
		\draw [style=wire] (89) to (90);
		\draw [style=wire] (89) to (91);
		\draw [style=wire] (88) to (86);
	\end{pgfonlayer}
\end{tikzpicture}
}%
\]

As for the constant rule for the deriving transformation, the constant rule {\bf [dC.1]} for a codereliction can be derived: 

\begin{lemma} For a bialgebra modality $(\oc,  \delta, \varepsilon, \Delta, e, \nabla, u)$, any natural transformation $\eta: A \to \oc A$ satisfies the constant rule {\bf [dC.1]}.
\end{lemma}
\begin{proof} By naturality of $e$ and $\eta$, and the additive structure, we have the following equalities: 
\[\resizebox{!}{2cm}{%
\begin{tikzpicture}
	\begin{pgfonlayer}{nodelayer}
		\node [style=port] (51) at (-9.25, 6.75) {};
		\node [style=component] (52) at (-9.25, 6) {$\eta$};
		\node [style=component] (53) at (-9.25, 5) {$e$};
		\node [style=port] (54) at (-3.25, 5.5) {$=$};
		\node [style=port] (55) at (-2.5, 5.5) {$0$};
		\node [style=port] (56) at (-6.75, 6.75) {};
		\node [style=component] (57) at (-6.75, 6) {$\eta$};
		\node [style=component] (58) at (-6.75, 4) {$e$};
		\node [style=function] (59) at (-6.75, 5) {$0$};
		\node [style=port] (60) at (-8, 5.5) {$=$};
		\node [style=port] (61) at (-8, 5) {Nat of $e$};
		\node [style=port] (62) at (-4.25, 6.75) {};
		\node [style=component] (63) at (-4.25, 6) {$0$};
		\node [style=component] (64) at (-4.25, 4) {$e$};
		\node [style=component] (65) at (-4.25, 5) {$\eta$};
		\node [style=port] (66) at (-5.5, 5.5) {$=$};
		\node [style=port] (67) at (-5.5, 5) {Nat of $\eta$};
	\end{pgfonlayer}
	\begin{pgfonlayer}{edgelayer}
		\draw [style=wire] (51) to (52);
		\draw [style=wire] (52) to (53);
		\draw [style=wire] (56) to (57);
		\draw [style=wire] (57) to (59);
		\draw [style=wire] (59) to (58);
		\draw [style=wire] (62) to (63);
		\draw [style=wire] (63) to (65);
		\draw [style=wire] (65) to (64);
	\end{pgfonlayer}
\end{tikzpicture}
}%
\]
\end{proof} 

\begin{corollary} For a bialgebra modality $(\oc,  \delta, \varepsilon, \Delta, e, \nabla, u)$, the following are equivalent for a natural transformation $\eta: A \to \oc A$: 
\begin{enumerate}[{\em (i)}]
\item $\eta$ is a codereliction; 
\item $\eta$ satisfies the product rule {\bf [dC.2]}, the linear rule {\bf [dC.3]}, and the chain rule {\bf [dC.4]}. 
\end{enumerate}
\end{corollary} 

In \cite{fiore2007differential} an alternative axiom for the chain rule {\bf [$\text{dC.4}$]} is used:
\begin{enumerate}[{\bf [$\text{dC.4}^\prime$]}] 
\item Alternative Chain Rule: 
\begin{align*}
\begin{array}[c]{c}
 \xymatrixcolsep{5pc}\xymatrix{A \ar[d]_-{u \otimes \eta} \ar[rr]^-{\eta} && \oc A \ar[d]^-{\delta}   \\
 \oc A \otimes \oc A \ar[r]_-{\delta \otimes \eta} & \oc \oc A \otimes \oc \oc A \ar[r]_-{\nabla} & \oc \oc A 
  } 
   \end{array} && 
   \begin{array}[c]{c}
\resizebox{!}{2cm}{%
\begin{tikzpicture}
	\begin{pgfonlayer}{nodelayer}
		\node [style=port] (0) at (0.5, -0.25) {};
		\node [style=component] (1) at (0.5, 1) {$\delta$};
		\node [style=port] (2) at (0.5, 3) {};
		\node [style=component] (3) at (0.5, 2) {$\eta$};
		\node [style=port] (4) at (1.25, 1.5) {$=$};
		\node [style=duplicate] (5) at (2.75, 0.5) {$\nabla$};
		\node [style=component] (6) at (3.25, 1.25) {$\eta$};
		\node [style=port] (7) at (2.75, -0.25) {};
		\node [style=duplicate] (8) at (2.75, 0.5) {$\nabla$};
		\node [style=component] (9) at (2.25, 1.25) {$\delta$};
		\node [style=component] (10) at (2.25, 2.25) {$u$};
		\node [style=component] (11) at (3.25, 2.25) {$\eta$};
		\node [style=port] (12) at (3.25, 3) {};
	\end{pgfonlayer}
	\begin{pgfonlayer}{edgelayer}
		\draw [style=wire] (2) to (3);
		\draw [style=wire] (3) to (1);
		\draw [style=wire] (1) to (0);
		\draw [style=wire, bend right=15, looseness=1.25] (5) to (6);
		\draw [style=wire] (5) to (7);
		\draw [style=wire, bend left=15, looseness=1.25] (8) to (9);
		\draw [style=wire] (12) to (11);
		\draw [style=wire] (11) to (6);
		\draw [style=wire] (10) to (9);
	\end{pgfonlayer}
\end{tikzpicture}
}%
   \end{array}
\end{align*}
\end{enumerate}
In a monoidal storage category -- the setting assumed in \cite{fiore2007differential} -- {\bf [dC.4]} and {\bf [$\text{dC.4}^\prime$]} are equivalent. However in the setting of a mere bialgebra modality, it is clear that {\bf [dC.4]} implies {\bf [$\text{dC.4}^\prime$]}: the reverse implication, however, does not appear to hold.  Thus, at this stage we prove the implication in one direction: 

\begin{lemma}\label{dc4'} For a bialgebra modality $(\oc,  \delta, \varepsilon, \Delta, e, \nabla, u)$, any natural transformation $\eta: A \to \oc A$ which satisfies the chain rule {\bf [dC.4]} also satisfies the alternative chain rule {\bf [$\text{dC.4}^\prime$]}. 
\end{lemma} 

\begin{proof} The bialgebra structure gives the following chain of equalities:  
\[\resizebox{!}{3cm}{%
\begin{tikzpicture}
	\begin{pgfonlayer}{nodelayer}
		\node [style=port] (0) at (-5.25, -0.25) {};
		\node [style=component] (1) at (-5.25, 1) {$\delta$};
		\node [style=port] (2) at (-5.25, 3) {};
		\node [style=component] (3) at (-5.25, 2) {$\eta$};
		\node [style=port] (4) at (-4.25, 1.5) {$=$};
		\node [style=duplicate] (5) at (6.25, -0.25) {$\nabla$};
		\node [style=component] (6) at (6.75, 0.5) {$\eta$};
		\node [style=port] (7) at (6.25, -1) {};
		\node [style=duplicate] (8) at (6.25, -0.25) {$\nabla$};
		\node [style=component] (9) at (5.75, 0.5) {$\delta$};
		\node [style=component] (10) at (5.75, 1.5) {$u$};
		\node [style=component] (11) at (6.75, 1.5) {$\eta$};
		\node [style=port] (12) at (6.75, 2.25) {};
		\node [style=component] (13) at (-3, 0.5) {$\delta$};
		\node [style=port] (14) at (-2.5, 3.5) {};
		\node [style=component] (15) at (-2.5, 2.5) {$\eta$};
		\node [style=port] (16) at (-3, -0.25) {};
		\node [style=component] (17) at (-3.5, 2.75) {$u$};
		\node [style=duplicate] (18) at (-3, 1.5) {$\nabla$};
		\node [style=port] (19) at (-4.25, 1) {(\ref{moneq})};
		\node [style=component] (20) at (-0.5, 3.5) {$u$};
		\node [style=duplicate] (21) at (-0.5, 2.5) {$\Delta$};
		\node [style=component] (22) at (-1, 1.5) {$\delta$};
		\node [style=port] (23) at (0.75, 3.5) {};
		\node [style=component] (24) at (0.75, 2.5) {$\eta$};
		\node [style=duplicate] (25) at (0.25, 1.5) {$\nabla$};
		\node [style=duplicate] (26) at (-0.25, -0.25) {$\nabla$};
		\node [style=component] (27) at (0.25, 0.5) {$\eta$};
		\node [style=port] (28) at (-0.25, -1) {};
		\node [style=port] (29) at (-2, 1.5) {$=$};
		\node [style=port] (30) at (-2, 1) {\textbf{[dC.4]}};
		\node [style=duplicate] (31) at (3.5, 1.25) {$\nabla$};
		\node [style=port] (32) at (4, 3.5) {};
		\node [style=component] (33) at (2.25, 0.75) {$\delta$};
		\node [style=component] (34) at (2.25, 2) {$u$};
		\node [style=component] (35) at (4, 2.25) {$\eta$};
		\node [style=port] (36) at (3, -1) {};
		\node [style=duplicate] (37) at (3, -0.25) {$\nabla$};
		\node [style=component] (38) at (3.5, 0.5) {$\eta$};
		\node [style=component] (39) at (3, 2) {$u$};
		\node [style=port] (40) at (1.25, 1.5) {$=$};
		\node [style=port] (41) at (1.25, 1) {(\ref{bialgeq})};
		\node [style=port] (42) at (4.75, 1.5) {$=$};
		\node [style=port] (43) at (4.75, 1) {(\ref{moneq})};
	\end{pgfonlayer}
	\begin{pgfonlayer}{edgelayer}
		\draw [style=wire] (2) to (3);
		\draw [style=wire] (3) to (1);
		\draw [style=wire] (1) to (0);
		\draw [style=wire, bend right=15, looseness=1.25] (5) to (6);
		\draw [style=wire] (5) to (7);
		\draw [style=wire, bend left=15, looseness=1.25] (8) to (9);
		\draw [style=wire] (12) to (11);
		\draw [style=wire] (11) to (6);
		\draw [style=wire] (10) to (9);
		\draw [style=wire, bend right=15, looseness=1.25] (18) to (15);
		\draw [style=wire, bend left=15] (18) to (17);
		\draw [style=wire] (13) to (18);
		\draw [style=wire] (13) to (16);
		\draw [style=wire] (14) to (15);
		\draw [style=wire, bend right] (21) to (22);
		\draw [style=wire] (20) to (21);
		\draw [style=wire, bend right=15, looseness=1.25] (25) to (24);
		\draw [style=wire] (24) to (23);
		\draw [style=wire, in=135, out=-30] (21) to (25);
		\draw [style=wire, bend right=15, looseness=1.25] (26) to (27);
		\draw [style=wire] (25) to (27);
		\draw [style=wire, bend right=15, looseness=1.25] (22) to (26);
		\draw [style=wire] (26) to (28);
		\draw [style=wire, bend right=15, looseness=1.25] (31) to (35);
		\draw [style=wire] (35) to (32);
		\draw [style=wire, bend right=15, looseness=1.25] (37) to (38);
		\draw [style=wire] (31) to (38);
		\draw [style=wire, bend right=15, looseness=1.25] (33) to (37);
		\draw [style=wire] (37) to (36);
		\draw [style=wire, bend right=15] (39) to (31);
		\draw [style=wire] (34) to (33);
	\end{pgfonlayer}
\end{tikzpicture}
}%
\]
\end{proof}

All deriving transformations which satisfy the $\nabla$-rule {\bf [d.$\nabla$]} induce a codereliction defined as: 
$$\eta := \xymatrixcolsep{2pc}\xymatrix{A \ar[r]^-{u \otimes 1} & \oc A \otimes A \ar[r]^-{\mathsf{d}} & \oc A 
  }~~~~~~~~~~\eta :=    \begin{array}[c]{c}\resizebox{!}{1.5cm}{%
  \begin{tikzpicture}
	\begin{pgfonlayer}{nodelayer}
		\node [style={circle, draw}] (0) at (0.5, 2) {$u$};
		\node [style=differential] (1) at (1, 1) {{\bf =\!=\!=}};
		\node [style=port] (2) at (1.75, 2.5) {};
		\node [style=port] (3) at (1, 0.25) {};
	\end{pgfonlayer}
	\begin{pgfonlayer}{edgelayer}
		\draw [style=wire, bend left, looseness=1.00] (1) to (0);
		\draw [style=wire, bend right, looseness=1.00] (1) to (2);
		\draw [style=wire] (1) to (3);
	\end{pgfonlayer}
\end{tikzpicture}
}
   \end{array}$$
Conversely, every codereliction induces a deriving transformation which satisfies the $\nabla$-rule:
  $$\mathsf{d}:=  \xymatrixcolsep{2pc}\xymatrix{\oc A \otimes A \ar[r]^-{1 \otimes \eta} & \oc A \otimes \oc A \ar[r]^-{\nabla} & \oc A 
  } ~~~~~~~~~~ \mathsf{d}:=    \begin{array}[c]{c}\resizebox{!}{1.5cm}{%
\begin{tikzpicture}
	\begin{pgfonlayer}{nodelayer}
		\node [style={circle, draw}] (0) at (2, 2.5) {$\eta$};
		\node [style=duplicate] (1) at (1.25, 1.5) {$\nabla$};
		\node [style=port] (2) at (0.5, 3.25) {};
		\node [style=port] (3) at (1.25, 0.75) {};
		\node [style=port] (4) at (2, 3.25) {};
	\end{pgfonlayer}
	\begin{pgfonlayer}{edgelayer}
		\draw [style=wire] (3) to (1);
		\draw [style=wire, in=-90, out=0, looseness=1.25] (1) to (0);
		\draw [style=wire] (0) to (4);
		\draw [style=wire, in=-90, out=180, looseness=1.00] (1) to (2);
	\end{pgfonlayer}
\end{tikzpicture}}
   \end{array}$$
Using the same proof in \cite{blute2006differential} -- which was for monoidal storage categories -- it  is easily seen that:

\begin{theorem} \label{bialgdiff} \cite[Theorem 4.12]{blute2006differential}  For an additive symmetric monoidal category with a bialgebra modality, deriving transformations which satisfy the $\nabla$-rule {\bf [d.$\nabla$]} are in bijective correspondence to coderelictions by $\mathsf{d} \longmapsto \eta := (u \otimes 1)\mathsf{d}$ and $\eta  \longmapsto  \mathsf{d} := (1 \otimes \eta)\nabla$. 
\end{theorem}

%%%%%%%%%%%%%%%%%%%%%%%%%%%%%%%%%%%%%%%%%%%%%%%%%%%%%%%%%%%%%%%%%%%%%%%

\section{Differentiation for Additive Bialgebra Modalities}\label{equivalencesec}

%%%%%%%%%%%%%%%%%%%%%%%%%%%%%%%%%%%%%%%%%%%%%%%%%%%%%%%%%%%%%%%%%%%%%%%

In this section we prove that for additive bialgebra modalities, there is only one notion of differentiation. 

First we examine coderelictions for additive bialgebra modalities. When a natural transformation $\eta$ is a section of $\varepsilon$, that is, $\eta$ satisfies {\bf [dC.3]}, we can define four natural transformations: 
\begin{align*}
\mathsf{p}_0= \varepsilon \otimes e: \oc A \otimes \oc B \to A && \mathsf{p}_1= e \otimes \varepsilon: \oc A \otimes \oc B \to B \\
\mathsf{i}_0= \eta \otimes u: A \to \oc A \otimes \oc B &&\mathsf{i}_1= u \otimes \eta: B \to \oc A \otimes \oc B
\end{align*}
Notice that since $\eta$ satisfies the constant rule {\bf [dC.1]} and the linear rule {\bf [dC.3]}, then from the properties of a bialgebra modality, we have:
\begin{equation}\label{pi1}\begin{gathered} 
\mathsf{i}_j\mathsf{p}_k=\begin{cases} 0 & \text{ if } j\neq k \\
1 & \text{ if } j=k
\end{cases}
 \end{gathered}\end{equation}
which is reminiscent of the identities satisfied by the projection and injection maps of a biproduct. These maps will be key to the proof of Lemma \ref{dc3dc4} below. For an additive bialgebra modality, since $\oc(0)=eu$, this means that we have the following:
\begin{equation}\label{pi2}\begin{gathered} 
\oc(\mathsf{i}_j)\oc(\mathsf{p}_k)=\begin{cases} eu & \text{ if } j\neq k \\
1 & \text{ if } j=k
\end{cases}
 \end{gathered}\end{equation}
This allows the derivation of the following useful identity: 

\begin{lemma}\label{seelyish} For an additive bialgebra modality $(\oc, \delta, \varepsilon, \Delta, e, \nabla, u)$ and a natural transformation $\eta: A \to \oc A$ which satisfies the linear rule {\bf [dC.3]}, the following diagram commutes: 
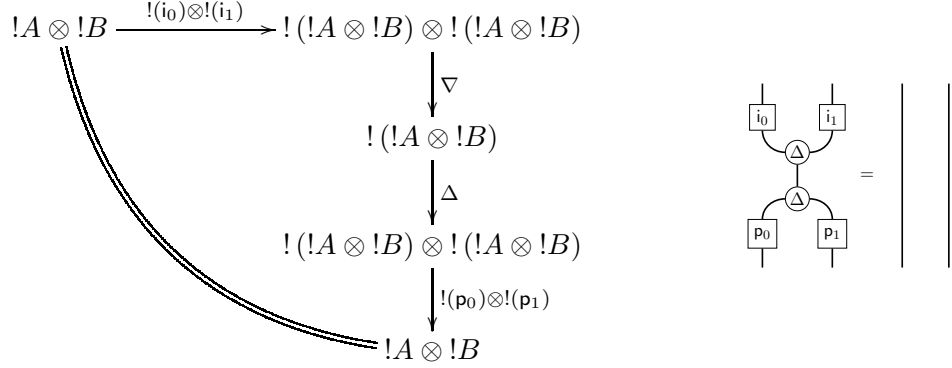
\begin{align*}
\begin{array}[c]{c}
\xymatrixcolsep{5pc}\xymatrix{ \oc A \otimes \oc B \ar@{=}@/_3pc/[dddr]^-{}  \ar[r]^-{\oc(\mathsf{i}_0) \otimes \oc(\mathsf{i}_1)} & \oc\left( \oc A \otimes \oc B \right)  \otimes \oc\left( \oc A \otimes \oc B \right)   \ar[d]^-{\nabla} \\
& \oc\left( \oc A \otimes \oc B \right)   \ar[d]^-{\Delta} \\
& \oc\left( \oc A \otimes \oc B \right)  \otimes \oc\left( \oc A \otimes \oc B \right)   \ar[d]^-{\oc(\mathsf{p}_0) \otimes \oc(\mathsf{p}_1)} \\
& \oc A \otimes \oc B
  } 
   \end{array} &&
   \begin{array}[c]{c}
   \resizebox{!}{2.5cm}{%
   \begin{tikzpicture}
	\begin{pgfonlayer}{nodelayer}
		\node [style=duplicate] (0) at (1.25, 0.75) {$\Delta$};
		\node [style=function2] (1) at (2, 0) {$\mathsf{p}_1$};
		\node [style=port] (2) at (2, -0.75) {};
		\node [style=port] (3) at (0.5, -0.75) {};
		\node [style=function2] (4) at (0.5, 0) {$\mathsf{p}_0$};
		\node [style=duplicate] (5) at (1.25, 1.75) {$\Delta$};
		\node [style=function2] (6) at (2, 2.5) {$\mathsf{i}_1$};
		\node [style=port] (7) at (2, 3.25) {};
		\node [style=port] (8) at (0.5, 3.25) {};
		\node [style=function2] (9) at (0.5, 2.5) {$\mathsf{i}_0$};
		\node [style=port] (11) at (3.5, 3.25) {};
		\node [style=port] (12) at (4.5, 3.25) {};
		\node [style=port] (13) at (3.5, -0.75) {};
		\node [style=port] (14) at (4.5, -0.75) {};
		\node [style=port] (15) at (2.75, 1.25) {$=$};
	\end{pgfonlayer}
	\begin{pgfonlayer}{edgelayer}
		\draw [style=wire] (1) to (2);
		\draw [style=wire] (4) to (3);
		\draw [style=wire, in=90, out=180] (0) to (4);
		\draw [style=wire, in=90, out=0] (0) to (1);
		\draw [style=wire] (6) to (7);
		\draw [style=wire] (9) to (8);
		\draw [style=wire, in=-90, out=180] (5) to (9);
		\draw [style=wire, in=-90, out=0] (5) to (6);
		\draw [style=wire] (5) to (0);
		\draw [style=wire] (11) to (13);
		\draw [style=wire] (12) to (14);
	\end{pgfonlayer}
\end{tikzpicture}
}%
   \end{array}
\end{align*}
\end{lemma}
\begin{proof} By the bialgebra identities, naturality of $\Delta$ and $\nabla$, and the $\mathsf{p}_i$ and $\mathsf{i}_j$ identities, we have the following equality: 
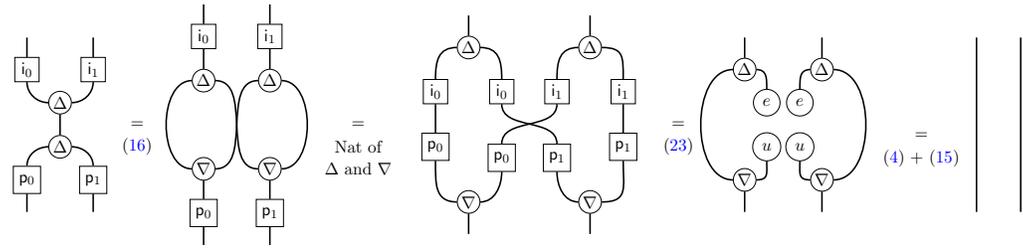
\[\resizebox{!}{3.25cm}{%
\begin{tikzpicture}
	\begin{pgfonlayer}{nodelayer}
		\node [style=duplicate] (0) at (-2.5, -5.25) {$\Delta$};
		\node [style=function2] (1) at (-1.75, -6) {$\mathsf{p}_1$};
		\node [style=port] (2) at (-1.75, -6.75) {};
		\node [style=port] (3) at (-3.25, -6.75) {};
		\node [style=function2] (4) at (-3.25, -6) {$\mathsf{p}_0$};
		\node [style=duplicate] (5) at (-2.5, -4.25) {$\Delta$};
		\node [style=function2] (6) at (-1.75, -3.5) {$\mathsf{i}_1$};
		\node [style=port] (7) at (-1.75, -2.75) {};
		\node [style=port] (8) at (-3.25, -2.75) {};
		\node [style=function2] (9) at (-3.25, -3.5) {$\mathsf{i}_0$};
		\node [style=port] (15) at (-0.75, -4.75) {$=$};
		\node [style=function2] (16) at (2.25, -6.75) {$\mathsf{p}_1$};
		\node [style=port] (17) at (2.25, -7.5) {};
		\node [style=port] (18) at (0.75, -7.5) {};
		\node [style=function2] (19) at (0.75, -6.75) {$\mathsf{p}_0$};
		\node [style=function2] (20) at (2.25, -2.75) {$\mathsf{i}_1$};
		\node [style=port] (21) at (2.25, -2) {};
		\node [style=port] (22) at (0.75, -2) {};
		\node [style=function2] (23) at (0.75, -2.75) {$\mathsf{i}_0$};
		\node [style=duplicate] (24) at (0.75, -5.75) {$\nabla$};
		\node [style=duplicate] (25) at (2.25, -3.75) {$\Delta$};
		\node [style=duplicate] (26) at (0.75, -3.75) {$\Delta$};
		\node [style=duplicate] (27) at (2.25, -5.75) {$\nabla$};
		\node [style=port] (28) at (-0.75, -5.25) {(\ref{bialgeq})};
		\node [style=port] (29) at (9.5, -7.25) {};
		\node [style=port] (30) at (6.75, -7.25) {};
		\node [style=port] (31) at (9.5, -2.25) {};
		\node [style=port] (32) at (6.75, -2.25) {};
		\node [style=duplicate] (33) at (6.75, -6.5) {$\nabla$};
		\node [style=duplicate] (34) at (9.5, -3) {$\Delta$};
		\node [style=duplicate] (35) at (6.75, -3) {$\Delta$};
		\node [style=duplicate] (36) at (9.5, -6.5) {$\nabla$};
		\node [style=function2] (37) at (7.5, -4) {$\mathsf{i}_0$};
		\node [style=function2] (38) at (8.75, -4) {$\mathsf{i}_1$};
		\node [style=function2] (39) at (10.25, -4) {$\mathsf{i}_1$};
		\node [style=function2] (40) at (10.25, -5.25) {$\mathsf{p}_1$};
		\node [style=function2] (41) at (6, -5.25) {$\mathsf{p}_0$};
		\node [style=function2] (42) at (6, -4) {$\mathsf{i}_0$};
		\node [style=function2] (43) at (8.75, -5.5) {$\mathsf{p}_1$};
		\node [style=function2] (44) at (7.5, -5.5) {$\mathsf{p}_0$};
		\node [style=port] (45) at (4.25, -4.75) {$=$};
		\node [style=port] (46) at (4.25, -5.25) {Nat of };
		\node [style=port] (47) at (13, -6.75) {};
		\node [style=port] (48) at (13, -2.75) {};
		\node [style=duplicate] (49) at (13, -6) {$\nabla$};
		\node [style=duplicate] (50) at (13, -3.5) {$\Delta$};
		\node [style=component] (51) at (13.5, -4.25) {$e$};
		\node [style=component] (52) at (13.5, -5.25) {$u$};
		\node [style=component] (53) at (14.25, -4.25) {$e$};
		\node [style=port] (54) at (14.75, -6.75) {};
		\node [style=component] (55) at (14.25, -5.25) {$u$};
		\node [style=port] (56) at (14.75, -2.75) {};
		\node [style=duplicate] (57) at (14.75, -3.5) {$\Delta$};
		\node [style=duplicate] (58) at (14.75, -6) {$\nabla$};
		\node [style=port] (59) at (11.5, -4.75) {$=$};
		\node [style=port] (60) at (11.5, -5.25) {(\ref{pi2})};
		\node [style=port] (61) at (17, -5) {$=$};
		\node [style=port] (62) at (17, -5.5) {(\ref{coalgeq}) + (\ref{moneq})};
		\node [style=port] (63) at (18.25, -2.75) {};
		\node [style=port] (64) at (19.25, -2.75) {};
		\node [style=port] (65) at (18.25, -6.75) {};
		\node [style=port] (66) at (19.25, -6.75) {};
		\node [style=port] (67) at (4.25, -5.75) {$\Delta$ and $\nabla$};
	\end{pgfonlayer}
	\begin{pgfonlayer}{edgelayer}
		\draw [style=wire] (1) to (2);
		\draw [style=wire] (4) to (3);
		\draw [style=wire, in=90, out=180] (0) to (4);
		\draw [style=wire, in=90, out=0] (0) to (1);
		\draw [style=wire] (6) to (7);
		\draw [style=wire] (9) to (8);
		\draw [style=wire, in=-90, out=180] (5) to (9);
		\draw [style=wire, in=-90, out=0] (5) to (6);
		\draw [style=wire] (5) to (0);
		\draw [style=wire] (16) to (17);
		\draw [style=wire] (19) to (18);
		\draw [style=wire] (20) to (21);
		\draw [style=wire] (23) to (22);
		\draw [style=wire, in=180, out=180] (26) to (24);
		\draw [style=wire, in=0, out=0] (25) to (27);
		\draw [style=wire, in=0, out=180] (25) to (24);
		\draw [style=wire, in=180, out=0] (26) to (27);
		\draw [style=wire] (23) to (26);
		\draw [style=wire] (20) to (25);
		\draw [style=wire] (24) to (19);
		\draw [style=wire] (27) to (16);
		\draw [style=wire] (32) to (35);
		\draw [style=wire] (31) to (34);
		\draw [style=wire] (33) to (30);
		\draw [style=wire] (36) to (29);
		\draw [style=wire, in=90, out=0, looseness=1.25] (35) to (37);
		\draw [style=wire, in=90, out=180, looseness=1.25] (34) to (38);
		\draw [style=wire, in=90, out=180, looseness=1.25] (35) to (42);
		\draw [style=wire] (42) to (41);
		\draw [style=wire, in=180, out=-90, looseness=1.25] (41) to (33);
		\draw [style=wire] (39) to (40);
		\draw [style=wire, in=0, out=-90, looseness=1.50] (40) to (36);
		\draw [style=wire, in=90, out=0, looseness=1.50] (34) to (39);
		\draw [style=wire, in=90, out=-90] (37) to (43);
		\draw [style=wire, in=90, out=-90] (38) to (44);
		\draw [style=wire, in=180, out=-90] (43) to (36);
		\draw [style=wire, in=0, out=-90] (44) to (33);
		\draw [style=wire, in=180, out=180] (50) to (49);
		\draw [style=wire] (48) to (50);
		\draw [style=wire] (49) to (47);
		\draw [style=wire, in=0, out=-90, looseness=1.50] (52) to (49);
		\draw [style=wire, in=0, out=90, looseness=1.25] (51) to (50);
		\draw [style=wire, in=0, out=0] (57) to (58);
		\draw [style=wire] (56) to (57);
		\draw [style=wire] (58) to (54);
		\draw [style=wire, in=180, out=-90, looseness=1.50] (55) to (58);
		\draw [style=wire, in=180, out=90, looseness=1.25] (53) to (57);
		\draw [style=wire] (63) to (65);
		\draw [style=wire] (64) to (66);
	\end{pgfonlayer}
\end{tikzpicture}
}%
\]
\end{proof} 

For an additive bialgebra modality, the linear rule {\bf [dC.3]} implies the product rule {\bf [dC.2]}: 

\begin{proposition}\label{dc3dc4} For an additive bialgebra modality $(\oc, \delta, \varepsilon, \Delta, e, \nabla, u)$, any natural transformation $\eta: A \to \oc A$ which satisfies the linear rule {\bf [dC.3]}, also satisfies the product rule {\bf [dC.2]}. 
\end{proposition}
\begin{proof} Notice that by naturality of $\eta$, we have that:
\begin{equation}\label{eta+}\begin{gathered} 
\eta\oc(f+g)=(f+g)\eta=f \eta + g \eta =\eta\oc(f)+\eta\oc(g)
 \end{gathered}\end{equation}
Then using the $\mathsf{i}_j$ and $\mathsf{p}_k$ identities, we obtain the following: 
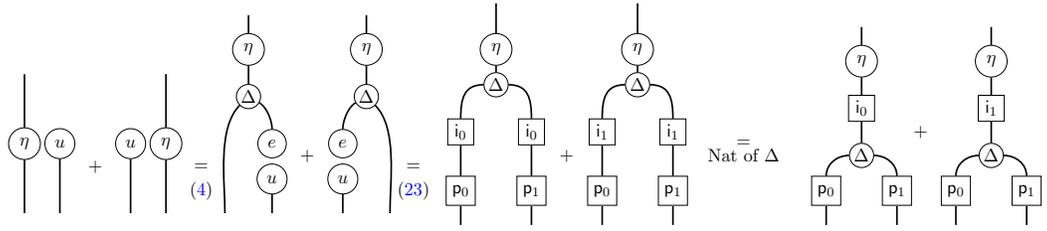
\[\resizebox{!}{3cm}{%
\begin{tikzpicture}
	\begin{pgfonlayer}{nodelayer}
		\node [style=component] (0) at (0, 0.5) {$\eta$};
		\node [style=port] (1) at (0.75, -1) {};
		\node [style=component] (2) at (0.75, 0.5) {$u$};
		\node [style=port] (3) at (0, -1) {};
		\node [style=port] (4) at (0, 2) {};
		\node [style=component] (6) at (2.25, 0.5) {$u$};
		\node [style=port] (7) at (2.25, -1) {};
		\node [style=port] (8) at (3, 2) {};
		\node [style=component] (9) at (3, 0.5) {$\eta$};
		\node [style=port] (10) at (3, -1) {};
		\node [style=duplicate] (11) at (4.75, 1.5) {$\Delta$};
		\node [style=component] (12) at (4.75, 2.5) {$\eta$};
		\node [style=port] (13) at (4.75, 3.25) {};
		\node [style=port] (14) at (4.25, -1) {};
		\node [style=component] (15) at (5.25, 0.5) {$e$};
		\node [style=component] (16) at (5.25, -0.25) {$u$};
		\node [style=port] (17) at (5.25, -1) {};
		\node [style=duplicate] (25) at (13, 1.75) {$\Delta$};
		\node [style=component] (26) at (13, 2.5) {$\eta$};
		\node [style=port] (27) at (13, 3.5) {};
		\node [style=function2] (28) at (13.75, 0.75) {$\mathsf{i}_1$};
		\node [style=function2] (29) at (13.75, -0.5) {$\mathsf{p}_1$};
		\node [style=port] (30) at (13.75, -1.25) {};
		\node [style=port] (31) at (12.25, -1.25) {};
		\node [style=function2] (32) at (12.25, -0.5) {$\mathsf{p}_0$};
		\node [style=function2] (33) at (12.25, 0.75) {$\mathsf{i}_1$};
		\node [style=duplicate] (34) at (10, 1.75) {$\Delta$};
		\node [style=component] (35) at (10, 2.5) {$\eta$};
		\node [style=port] (36) at (10, 3.5) {};
		\node [style=function2] (37) at (10.75, 0.75) {$\mathsf{i}_0$};
		\node [style=function2] (38) at (10.75, -0.5) {$\mathsf{p}_1$};
		\node [style=port] (39) at (10.75, -1.25) {};
		\node [style=port] (40) at (9.25, -1.25) {};
		\node [style=function2] (41) at (9.25, -0.5) {$\mathsf{p}_0$};
		\node [style=function2] (42) at (9.25, 0.75) {$\mathsf{i}_0$};
		\node [style=duplicate] (43) at (17.75, 0.25) {$\Delta$};
		\node [style=component] (44) at (17.75, 2.25) {$\eta$};
		\node [style=port] (45) at (17.75, 3) {};
		\node [style=function2] (46) at (18.5, -0.5) {$\mathsf{p}_1$};
		\node [style=port] (47) at (18.5, -1.25) {};
		\node [style=port] (48) at (17, -1.25) {};
		\node [style=function2] (49) at (17, -0.5) {$\mathsf{p}_0$};
		\node [style=function2] (50) at (17.75, 1.25) {$\mathsf{i}_0$};
		\node [style=duplicate] (51) at (20.5, 0.25) {$\Delta$};
		\node [style=component] (52) at (20.5, 2.25) {$\eta$};
		\node [style=port] (53) at (20.5, 3) {};
		\node [style=function2] (54) at (21.25, -0.5) {$\mathsf{p}_1$};
		\node [style=port] (55) at (21.25, -1.25) {};
		\node [style=port] (56) at (19.75, -1.25) {};
		\node [style=function2] (57) at (19.75, -0.5) {$\mathsf{p}_0$};
		\node [style=function2] (58) at (20.5, 1.25) {$\mathsf{i}_1$};
		\node [style=port] (81) at (1.5, 0) {$+$};
		\node [style=port] (82) at (3.75, 0) {$=$};
		\node [style=duplicate] (83) at (7.25, 1.5) {$\Delta$};
		\node [style=component] (84) at (7.25, 2.5) {$\eta$};
		\node [style=port] (85) at (7.25, 3.25) {};
		\node [style=port] (86) at (7.75, -1) {};
		\node [style=component] (87) at (6.75, 0.5) {$e$};
		\node [style=component] (88) at (6.75, -0.25) {$u$};
		\node [style=port] (89) at (6.75, -1) {};
		\node [style=port] (90) at (6, 0.25) {$+$};
		\node [style=port] (91) at (8.25, 0) {$=$};
		\node [style=port] (92) at (11.5, 0.25) {$+$};
		\node [style=port] (93) at (15.25, 0.5) {$=$};
		\node [style=port] (94) at (19, 0.75) {$+$};
		\node [style=port] (95) at (3.75, -0.5) {(\ref{coalgeq})};
		\node [style=port] (96) at (8.25, -0.5) {(\ref{pi2})};
		\node [style=port] (97) at (15.25, 0.25) {Nat of $\Delta$};
	\end{pgfonlayer}
	\begin{pgfonlayer}{edgelayer}
		\draw [style=wire] (2) to (1);
		\draw [style=wire] (4) to (0);
		\draw [style=wire] (0) to (3);
		\draw [style=wire] (6) to (7);
		\draw [style=wire] (8) to (9);
		\draw [style=wire] (9) to (10);
		\draw [style=wire] (13) to (12);
		\draw [style=wire] (12) to (11);
		\draw [style=wire, in=90, out=-135, looseness=0.75] (11) to (14);
		\draw [style=wire, in=90, out=-30, looseness=0.75] (11) to (15);
		\draw [style=wire] (16) to (17);
		\draw [style=wire] (27) to (26);
		\draw [style=wire] (26) to (25);
		\draw [style=wire, in=90, out=0, looseness=1.25] (25) to (28);
		\draw [style=wire] (28) to (29);
		\draw [style=wire] (29) to (30);
		\draw [style=wire] (33) to (32);
		\draw [style=wire] (32) to (31);
		\draw [style=wire, in=90, out=180, looseness=1.25] (25) to (33);
		\draw [style=wire] (36) to (35);
		\draw [style=wire] (35) to (34);
		\draw [style=wire, in=90, out=0, looseness=1.25] (34) to (37);
		\draw [style=wire] (37) to (38);
		\draw [style=wire] (38) to (39);
		\draw [style=wire] (42) to (41);
		\draw [style=wire] (41) to (40);
		\draw [style=wire, in=90, out=180, looseness=1.25] (34) to (42);
		\draw [style=wire] (45) to (44);
		\draw [style=wire] (46) to (47);
		\draw [style=wire] (49) to (48);
		\draw [style=wire, in=90, out=180] (43) to (49);
		\draw [style=wire, in=90, out=0] (43) to (46);
		\draw [style=wire] (50) to (43);
		\draw [style=wire] (44) to (50);
		\draw [style=wire] (53) to (52);
		\draw [style=wire] (54) to (55);
		\draw [style=wire] (57) to (56);
		\draw [style=wire, in=90, out=180] (51) to (57);
		\draw [style=wire, in=90, out=0] (51) to (54);
		\draw [style=wire] (58) to (51);
		\draw [style=wire] (52) to (58);
		\draw [style=wire] (85) to (84);
		\draw [style=wire] (84) to (83);
		\draw [style=wire, in=90, out=-45, looseness=0.75] (83) to (86);
		\draw [style=wire, in=90, out=-150, looseness=0.75] (83) to (87);
		\draw [style=wire] (88) to (89);
	\end{pgfonlayer}
\end{tikzpicture}
}%
\]
\[\resizebox{!}{3.5cm}{%
\begin{tikzpicture}
	\begin{pgfonlayer}{nodelayer}
		\node [style=duplicate] (0) at (18.5, 0.5) {$\Delta$};
		\node [style=component] (1) at (18.5, 3) {$\eta$};
		\node [style=port] (2) at (18.5, 3.75) {};
		\node [style=function2] (3) at (19.25, -0.25) {$\mathsf{p}_1$};
		\node [style=port] (4) at (19.25, -1) {};
		\node [style=port] (5) at (17.75, -1) {};
		\node [style=function2] (6) at (17.75, -0.25) {$\mathsf{p}_0$};
		\node [style=function2] (7) at (18.5, 1.75) {$\mathsf{i}_0 + \mathsf{i}_1$};
		\node [style=duplicate] (8) at (24.25, 2.25) {$\Delta$};
		\node [style=component] (9) at (24.25, 3.25) {$\eta$};
		\node [style=port] (10) at (24.25, 4) {};
		\node [style=function2] (11) at (25.5, 1) {$\mathsf{i}_0 + \mathsf{i}_1$};
		\node [style=function2] (12) at (25.5, -0.5) {$\mathsf{p}_1$};
		\node [style=port] (13) at (25.5, -1.25) {};
		\node [style=port] (14) at (23, -1.25) {};
		\node [style=function2] (15) at (23, -0.5) {$\mathsf{p}_0$};
		\node [style=function2] (16) at (23, 1) {$\mathsf{i}_0 + \mathsf{i}_1$};
		\node [style=port] (17) at (28.25, -0.75) {};
		\node [style=port] (18) at (29, 2.25) {};
		\node [style=duplicate] (19) at (29, 0.25) {$\Delta$};
		\node [style=component] (20) at (29, 1.25) {$\eta$};
		\node [style=port] (21) at (29.75, -0.75) {};
		\node [style=port] (22) at (16.5, 0.75) {$=$};
		\node [style=port] (23) at (21, 0.75) {$=$};
		\node [style=port] (24) at (27.25, 0.5) {$=$};
		\node [style=port] (25) at (16.5, 0.25) {(\ref{eta+})};
		\node [style=port] (26) at (21, 0.5) {Nat of $\Delta$};
		\node [style=port] (27) at (27.25, 0) {(\ref{pi1})};
	\end{pgfonlayer}
	\begin{pgfonlayer}{edgelayer}
		\draw [style=wire] (2) to (1);
		\draw [style=wire] (3) to (4);
		\draw [style=wire] (6) to (5);
		\draw [style=wire, in=90, out=180] (0) to (6);
		\draw [style=wire, in=90, out=0] (0) to (3);
		\draw [style=wire] (7) to (0);
		\draw [style=wire] (1) to (7);
		\draw [style=wire] (10) to (9);
		\draw [style=wire] (9) to (8);
		\draw [style=wire, in=90, out=0, looseness=1.25] (8) to (11);
		\draw [style=wire] (11) to (12);
		\draw [style=wire] (12) to (13);
		\draw [style=wire] (16) to (15);
		\draw [style=wire] (15) to (14);
		\draw [style=wire, in=90, out=180, looseness=1.25] (8) to (16);
		\draw [style=wire] (18) to (20);
		\draw [style=wire, bend left] (19) to (21);
		\draw [style=wire, bend right] (19) to (17);
		\draw [style=wire] (20) to (19);
	\end{pgfonlayer}
\end{tikzpicture}
}%
\]
\end{proof}  

\begin{corollary} For an additive bialgebra modality $(\oc,  \delta, \varepsilon, \Delta, e, \nabla, u)$, the following are equivalent for a natural transformation $\eta: A \to \oc A$: 
\begin{enumerate}[{\em (i)}]
\item $\eta$ is a codereliction; 
\item $\eta$ satisfies the linear rule {\bf [dC.3]} and the chain rule {\bf [dC.4]}. 
\end{enumerate}
\end{corollary} 

In \cite{fiore2007differential}, for a monoidal coalgebra modality, Fiore introduced another axiom relating $\eta$ to the monoidal structure: 
\begin{description}
\item[{\bf [dC.m]}] Monoidal Rule: 
\begin{align*}
\begin{array}[c]{c}
\xymatrixcolsep{5pc}\xymatrix{\oc A \otimes B \ar[d]_-{\varepsilon \otimes 1} \ar[r]^-{1 \otimes \eta} & \oc A \otimes \oc B \ar[d]^-{m_\otimes} \\
  A \otimes B \ar[r]_-{\eta} & \oc(A \otimes B)
  } 
   \end{array} &&
   \begin{array}[c]{c}
\resizebox{!}{2.25cm}{%
\begin{tikzpicture}
	\begin{pgfonlayer}{nodelayer}
		\node [style=component] (0) at (2, 1.25) {$\eta$};
		\node [style=port] (1) at (1.25, -1.75) {};
		\node [style=port] (2) at (0.5, 2) {};
		\node [style=function2] (3) at (1.25, 0) {$\bigotimes$};
		\node [style=port] (4) at (2, 2) {};
		\node [style=port] (5) at (4, 0.25) {$\bigotimes$};
		\node [style=port] (6) at (3.25, 2) {};
		\node [style=port] (7) at (4.75, 2) {};
		\node [style=component] (8) at (4, -0.75) {$\eta$};
		\node [style=port] (9) at (4, -1.75) {};
		\node [style=component] (10) at (3.25, 1.25) {$\varepsilon$};
		\node [style=port] (11) at (2.5, -0.5) {$=$};
	\end{pgfonlayer}
	\begin{pgfonlayer}{edgelayer}
		\draw [style=wire, in=-90, out=165] (3) to (2);
		\draw [style=wire] (4) to (0);
		\draw [style=wire, in=15, out=-90] (0) to (3);
		\draw [style=wire] (3) to (1);
		\draw [style=wire, in=-90, out=15] (5) to (7);
		\draw [style=wire] (5) to (8);
		\draw [style=wire] (8) to (9);
		\draw [style=wire] (6) to (10);
		\draw [style=wire, in=165, out=-90] (10) to (5);
	\end{pgfonlayer}
\end{tikzpicture}
}%
   \end{array}
\end{align*}
\end{description}

However, it turns out that coderelictions for the additive bialgebra modalities always satisfy the monoidal rule {\bf [dC.m]}: 

\begin{proposition}\label{etamonoidal} For the induced additive bialgabra modality of an additive linear category: all coderelictions satisfy the monoidal rule {\bf [dC.m]}. 
\end{proposition} 
\begin{proof} By Lemma \ref{seelyish} and the fact that $\Delta$ is a $\oc$-coalgebra morphism, we first obtain the following:
\begin{equation}\label{mcoder}\begin{gathered} 
\resizebox{!}{6cm}{%
\begin{tikzpicture}
	\begin{pgfonlayer}{nodelayer}
		\node [style=port] (0) at (-1.25, 0.5) {};
		\node [style=port] (1) at (-2.75, 0.5) {};
		\node [style=function2] (2) at (-2, -0.5) {$\bigotimes$};
		\node [style=port] (3) at (-2, -1.25) {};
		\node [style=function2] (4) at (4.25, 3.25) {$\mathsf{i}_0$};
		\node [style=port] (5) at (5.75, 4) {};
		\node [style=function2] (6) at (5.75, 0.75) {$\mathsf{p}_1$};
		\node [style=port] (7) at (4.25, 4) {};
		\node [style=duplicate] (8) at (5, 1.5) {$\Delta$};
		\node [style=function2] (9) at (5.75, 3.25) {$\mathsf{i}_1$};
		\node [style=duplicate] (10) at (5, 2.5) {$\Delta$};
		\node [style=component] (11) at (4.25, -0.25) {$\delta$};
		\node [style=function2] (12) at (4.25, 0.75) {$\mathsf{p}_0$};
		\node [style=component] (13) at (5.75, -0.25) {$\delta$};
		\node [style=function] (14) at (4.25, -1.25) {$\varepsilon$};
		\node [style=function2] (15) at (5, -2.25) {$\bigotimes$};
		\node [style=function] (16) at (5.75, -1.25) {$\varepsilon$};
		\node [style=port] (17) at (5, -3) {};
		\node [style=port] (18) at (10.5, -3.5) {};
		\node [style=component] (19) at (9, -2.5) {$\varepsilon$};
		\node [style=function2] (20) at (9, -1.25) {$\mathsf{p}_0$};
		\node [style=port] (21) at (8.5, -0.25) {};
		\node [style=port] (22) at (9, -3.5) {};
		\node [style=port] (23) at (8.5, -3.5) {};
		\node [style=port] (25) at (11, -3.5) {};
		\node [style=port] (26) at (11, -0.25) {};
		\node [style=function2] (27) at (10.5, -1.25) {$\mathsf{p}_1$};
		\node [style=function2] (28) at (9, 4.25) {$\mathsf{i}_0$};
		\node [style=function2] (29) at (10.5, 4.25) {$\mathsf{i}_1$};
		\node [style=port] (30) at (9.75, -4.5) {};
		\node [style=component] (31) at (10.5, -2.5) {$\varepsilon$};
		\node [style=port] (32) at (10.5, 5) {};
		\node [style=port] (33) at (9, 5) {};
		\node [style=duplicate] (34) at (9.75, 3.25) {$\nabla$};
		\node [style=port] (35) at (9.75, -3.5) {};
		\node [style=component] (36) at (10.5, 1.5) {$\delta$};
		\node [style=component] (37) at (9, 1.5) {$\delta$};
		\node [style=function2] (38) at (9.75, 0.5) {$\bigotimes$};
		\node [style=duplicate] (39) at (9.75, 2.5) {$\Delta$};
		\node [style=function2] (40) at (15.25, -1.25) {$\mathsf{p}_1$};
		\node [style=port] (41) at (14.5, -4.25) {};
		\node [style=function2] (42) at (15.25, 3) {$\mathsf{i}_1$};
		\node [style=port] (43) at (15.75, -3.25) {};
		\node [style=port] (44) at (14.5, -3.25) {};
		\node [style=component] (45) at (14.5, 1) {$\delta$};
		\node [style=port] (46) at (15.25, -3.25) {};
		\node [style=component] (47) at (13.75, -2.25) {$\varepsilon$};
		\node [style=port] (48) at (13.25, 0.25) {};
		\node [style=port] (49) at (15.75, 0.25) {};
		\node [style=duplicate] (50) at (14.5, 2) {$\nabla$};
		\node [style=port] (51) at (13.25, -3.25) {};
		\node [style=duplicate] (52) at (14.5, -0.25) {$\Delta$};
		\node [style=function2] (53) at (13.75, 3) {$\mathsf{i}_0$};
		\node [style=port] (54) at (13.75, -3.25) {};
		\node [style=function2] (55) at (13.75, -1.25) {$\mathsf{p}_0$};
		\node [style=component] (56) at (15.25, -2.25) {$\varepsilon$};
		\node [style=port] (57) at (13.75, 3.75) {};
		\node [style=port] (58) at (15.25, 3.75) {};
		\node [style=port] (60) at (1.75, 1.25) {};
		\node [style=port] (62) at (0.25, 1.25) {};
		\node [style=component] (66) at (0.25, 0.25) {$\delta$};
		\node [style=component] (68) at (1.75, 0.25) {$\delta$};
		\node [style=function] (69) at (0.25, -0.75) {$\varepsilon$};
		\node [style=function2] (70) at (1, -1.75) {$\bigotimes$};
		\node [style=function] (71) at (1.75, -0.75) {$\varepsilon$};
		\node [style=port] (72) at (1, -2.5) {};
		\node [style=port] (73) at (-0.75, -0.25) {$=$};
		\node [style=port] (74) at (-0.75, -0.75) {(\ref{comonad})};
		\node [style=port] (75) at (3, -0.25) {$=$};
		\node [style=port] (76) at (3, -0.75) {Lem \ref{seelyish}};
		\node [style=port] (77) at (9, -0.25) {};
		\node [style=port] (78) at (10.5, -0.25) {};
		\node [style=port] (79) at (9.75, -0.25) {};
		\node [style=port] (80) at (7.25, -0.25) {$=$};
		\node [style=port] (81) at (7.25, -0.5) {Nat of};
		\node [style=port] (82) at (7.25, -1) {$\delta$ and $m_\otimes$};
		\node [style=port] (83) at (12, -0.25) {$=$};
		\node [style=port] (84) at (12, -0.75) {(\ref{mcm2})};
	\end{pgfonlayer}
	\begin{pgfonlayer}{edgelayer}
		\draw [style=wire, in=-90, out=180, looseness=1.50] (2) to (1);
		\draw [style=wire, in=0, out=-90, looseness=1.50] (0) to (2);
		\draw [style=wire] (2) to (3);
		\draw [style=wire] (6) to (13);
		\draw [style=wire] (12) to (11);
		\draw [style=wire, in=90, out=180] (8) to (12);
		\draw [style=wire, in=90, out=0] (8) to (6);
		\draw [style=wire] (9) to (5);
		\draw [style=wire] (4) to (7);
		\draw [style=wire, in=-90, out=180] (10) to (4);
		\draw [style=wire, in=-90, out=0] (10) to (9);
		\draw [style=wire] (10) to (8);
		\draw [style=wire] (15) to (17);
		\draw [style=wire, in=0, out=-90, looseness=1.75] (16) to (15);
		\draw [style=wire, in=180, out=-90, looseness=1.75] (14) to (15);
		\draw [style=wire] (11) to (14);
		\draw [style=wire] (13) to (16);
		\draw [style=wire] (22) to (19);
		\draw [style=wire] (18) to (31);
		\draw [style=wire] (21) to (23);
		\draw [style=wire] (21) to (26);
		\draw [style=wire] (26) to (25);
		\draw [style=wire] (23) to (25);
		\draw [style=wire] (35) to (30);
		\draw [style=wire, bend left] (34) to (28);
		\draw [style=wire, bend right] (34) to (29);
		\draw [style=wire] (20) to (19);
		\draw [style=wire] (27) to (31);
		\draw [style=wire] (33) to (28);
		\draw [style=wire] (32) to (29);
		\draw [style=wire, bend right] (39) to (37);
		\draw [style=wire, bend left] (39) to (36);
		\draw [style=wire, in=0, out=-90, looseness=1.75] (36) to (38);
		\draw [style=wire, in=180, out=-90, looseness=1.75] (37) to (38);
		\draw [style=wire] (34) to (39);
		\draw [style=wire] (54) to (47);
		\draw [style=wire] (46) to (56);
		\draw [style=wire] (48) to (51);
		\draw [style=wire] (48) to (49);
		\draw [style=wire] (49) to (43);
		\draw [style=wire] (51) to (43);
		\draw [style=wire] (44) to (41);
		\draw [style=wire, bend left] (50) to (53);
		\draw [style=wire, bend right] (50) to (42);
		\draw [style=wire] (50) to (45);
		\draw [style=wire, bend right] (52) to (55);
		\draw [style=wire, bend left] (52) to (40);
		\draw [style=wire] (55) to (47);
		\draw [style=wire] (40) to (56);
		\draw [style=wire] (45) to (52);
		\draw [style=wire] (57) to (53);
		\draw [style=wire] (58) to (42);
		\draw [style=wire] (70) to (72);
		\draw [style=wire, in=0, out=-90, looseness=1.75] (71) to (70);
		\draw [style=wire, in=180, out=-90, looseness=1.75] (69) to (70);
		\draw [style=wire] (66) to (69);
		\draw [style=wire] (68) to (71);
		\draw [style=wire] (62) to (66);
		\draw [style=wire] (60) to (68);
		\draw [style=wire] (77) to (20);
		\draw [style=wire] (78) to (27);
		\draw [style=wire] (38) to (79);
	\end{pgfonlayer}
\end{tikzpicture}}%
 \end{gathered}\end{equation}
Expressing $m_\otimes$ as above, then by the linear rule {\bf [dC.3]}, chain rule {\bf [dC.4]}, the naturality of $u$ and $\nabla$, and Proposition \ref{monoidalnabla}, we have the following equality: 
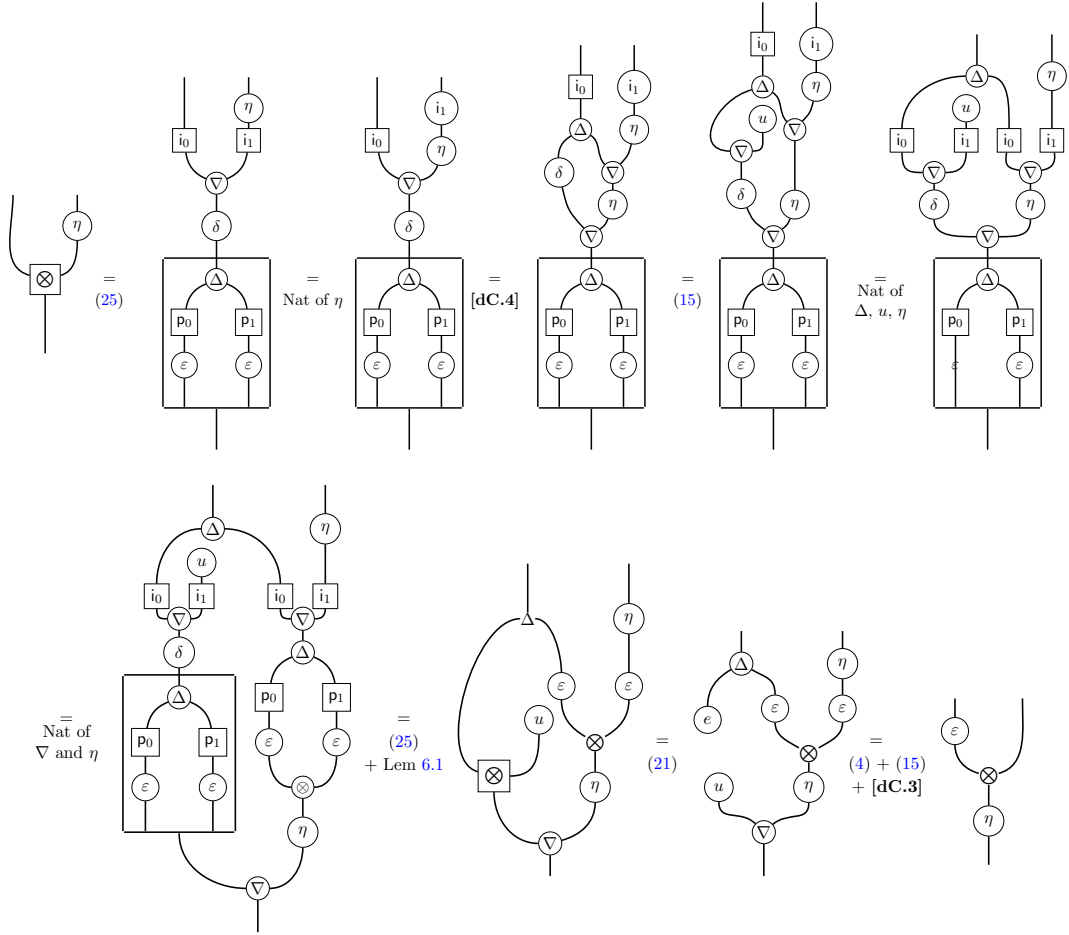
\[\resizebox{!}{6cm}{%
\begin{tikzpicture}
	\begin{pgfonlayer}{nodelayer}
		\node [style=function2] (40) at (15.25, -1.25) {$\mathsf{p}_1$};
		\node [style=port] (41) at (14.5, -4.25) {};
		\node [style=function2] (42) at (15.25, 3) {$\mathsf{i}_1$};
		\node [style=port] (43) at (15.75, -3.25) {};
		\node [style=port] (44) at (14.5, -3.25) {};
		\node [style=component] (45) at (14.5, 1) {$\delta$};
		\node [style=port] (46) at (15.25, -3.25) {};
		\node [style=component] (47) at (13.75, -2.25) {$\varepsilon$};
		\node [style=port] (48) at (13.25, 0.25) {};
		\node [style=port] (49) at (15.75, 0.25) {};
		\node [style=duplicate] (50) at (14.5, 2) {$\nabla$};
		\node [style=port] (51) at (13.25, -3.25) {};
		\node [style=duplicate] (52) at (14.5, -0.25) {$\Delta$};
		\node [style=function2] (53) at (13.75, 3) {$\mathsf{i}_0$};
		\node [style=port] (54) at (13.75, -3.25) {};
		\node [style=function2] (55) at (13.75, -1.25) {$\mathsf{p}_0$};
		\node [style=component] (56) at (15.25, -2.25) {$\varepsilon$};
		\node [style=port] (57) at (13.75, 4.5) {};
		\node [style=port] (83) at (12, -0.25) {$=$};
		\node [style=port] (84) at (12, -0.75) {(\ref{mcoder})};
		\node [style=component] (85) at (11.25, 1) {$\eta$};
		\node [style=port] (86) at (10.5, -2) {};
		\node [style=port] (87) at (9.75, 1.75) {};
		\node [style=function2] (88) at (10.5, -0.25) {$\bigotimes$};
		\node [style=port] (89) at (11.25, 1.75) {};
		\node [style=component] (90) at (15.25, 3.75) {$\eta$};
		\node [style=port] (91) at (15.25, 4.5) {};
		\node [style=port] (92) at (16.75, -0.25) {$=$};
		\node [style=port] (93) at (16.75, -0.75) {Nat of $\eta$};
		\node [style=function2] (94) at (19.75, -1.25) {$\mathsf{p}_1$};
		\node [style=port] (95) at (19, -4.25) {};
		\node [style=component] (96) at (19.75, 2.75) {$\eta$};
		\node [style=port] (97) at (20.25, -3.25) {};
		\node [style=port] (98) at (19, -3.25) {};
		\node [style=component] (99) at (19, 1) {$\delta$};
		\node [style=port] (100) at (19.75, -3.25) {};
		\node [style=component] (101) at (18.25, -2.25) {$\varepsilon$};
		\node [style=port] (102) at (17.75, 0.25) {};
		\node [style=port] (103) at (20.25, 0.25) {};
		\node [style=duplicate] (104) at (19, 2) {$\nabla$};
		\node [style=port] (105) at (17.75, -3.25) {};
		\node [style=duplicate] (106) at (19, -0.25) {$\Delta$};
		\node [style=function2] (107) at (18.25, 3) {$\mathsf{i}_0$};
		\node [style=port] (108) at (18.25, -3.25) {};
		\node [style=function2] (109) at (18.25, -1.25) {$\mathsf{p}_0$};
		\node [style=component] (110) at (19.75, -2.25) {$\varepsilon$};
		\node [style=port] (111) at (18.25, 4.5) {};
		\node [style=component] (112) at (19.75, 3.75) {$\mathsf{i}_1$};
		\node [style=port] (113) at (19.75, 4.5) {};
		\node [style=duplicate] (114) at (23.75, 2.25) {$\nabla$};
		\node [style=component] (115) at (23.75, 1.5) {$\eta$};
		\node [style=duplicate] (116) at (23, 3.25) {$\Delta$};
		\node [style=duplicate] (117) at (23.25, 0.75) {$\nabla$};
		\node [style=component] (118) at (24.25, 3.25) {$\eta$};
		\node [style=component] (119) at (22.5, 2.25) {$\delta$};
		\node [style=function2] (120) at (23, 4.25) {$\mathsf{i}_0$};
		\node [style=port] (121) at (23, 5.25) {};
		\node [style=component] (122) at (24.25, 4.25) {$\mathsf{i}_1$};
		\node [style=port] (123) at (24.25, 5.25) {};
		\node [style=port] (124) at (22, 0.25) {};
		\node [style=port] (125) at (24.5, 0.25) {};
		\node [style=function2] (126) at (24, -1.25) {$\mathsf{p}_1$};
		\node [style=port] (127) at (23.25, -4.25) {};
		\node [style=duplicate] (128) at (23.25, -0.25) {$\Delta$};
		\node [style=port] (129) at (24, -3.25) {};
		\node [style=component] (130) at (24, -2.25) {$\varepsilon$};
		\node [style=component] (131) at (22.5, -2.25) {$\varepsilon$};
		\node [style=port] (132) at (24.5, -3.25) {};
		\node [style=function2] (133) at (22.5, -1.25) {$\mathsf{p}_0$};
		\node [style=port] (134) at (22, -3.25) {};
		\node [style=port] (135) at (22.5, -3.25) {};
		\node [style=port] (136) at (23.25, -3.25) {};
		\node [style=port] (137) at (21, -0.25) {$=$};
		\node [style=port] (138) at (21, -0.75) {\textbf{[dC.4]}};
		\node [style=duplicate] (166) at (28, 3.25) {$\nabla$};
		\node [style=component] (167) at (28, 1.5) {$\eta$};
		\node [style=duplicate] (168) at (27.25, 4.25) {$\Delta$};
		\node [style=duplicate] (169) at (27.5, 0.75) {$\nabla$};
		\node [style=component] (170) at (28.5, 4.25) {$\eta$};
		\node [style=component] (171) at (26.75, 1.75) {$\delta$};
		\node [style=function2] (172) at (27.25, 5.25) {$\mathsf{i}_0$};
		\node [style=port] (173) at (27.25, 6.25) {};
		\node [style=component] (174) at (28.5, 5.25) {$\mathsf{i}_1$};
		\node [style=port] (175) at (28.5, 6.25) {};
		\node [style=port] (176) at (26.25, 0.25) {};
		\node [style=port] (177) at (28.75, 0.25) {};
		\node [style=function2] (178) at (28.25, -1.25) {$\mathsf{p}_1$};
		\node [style=port] (179) at (27.5, -4.25) {};
		\node [style=duplicate] (180) at (27.5, -0.25) {$\Delta$};
		\node [style=port] (181) at (28.25, -3.25) {};
		\node [style=component] (182) at (28.25, -2.25) {$\varepsilon$};
		\node [style=component] (183) at (26.75, -2.25) {$\varepsilon$};
		\node [style=port] (184) at (28.75, -3.25) {};
		\node [style=function2] (185) at (26.75, -1.25) {$\mathsf{p}_0$};
		\node [style=port] (186) at (26.25, -3.25) {};
		\node [style=port] (187) at (26.75, -3.25) {};
		\node [style=port] (188) at (27.5, -3.25) {};
		\node [style=duplicate] (189) at (26.75, 2.75) {$\nabla$};
		\node [style=component] (190) at (27.25, 3.5) {$u$};
		\node [style=port] (191) at (25.5, -0.25) {$=$};
		\node [style=port] (192) at (25.5, -0.75) {(\ref{moneq})};
		\node [style=component] (193) at (33.5, 1.5) {$\eta$};
		\node [style=duplicate] (194) at (32.25, 4.5) {$\Delta$};
		\node [style=duplicate] (195) at (32.5, 0.75) {$\nabla$};
		\node [style=component] (196) at (31.25, 1.5) {$\delta$};
		\node [style=port] (197) at (32.25, 5.5) {};
		\node [style=component] (198) at (34, 4.5) {$\eta$};
		\node [style=port] (199) at (34, 5.5) {};
		\node [style=function2] (200) at (30.5, 3) {$\mathsf{i}_0$};
		\node [style=function2] (201) at (32, 3) {$\mathsf{i}_1$};
		\node [style=component] (202) at (32, 3.75) {$u$};
		\node [style=duplicate] (203) at (31.25, 2.25) {$\nabla$};
		\node [style=function2] (204) at (33, 3) {$\mathsf{i}_0$};
		\node [style=function2] (205) at (34, 3) {$\mathsf{i}_1$};
		\node [style=duplicate] (206) at (33.5, 2.25) {$\nabla$};
		\node [style=port] (207) at (31.25, 0.25) {};
		\node [style=port] (208) at (33.75, 0.25) {};
		\node [style=function2] (209) at (33.25, -1.25) {$\mathsf{p}_1$};
		\node [style=port] (210) at (32.5, -4.25) {};
		\node [style=duplicate] (211) at (32.5, -0.25) {$\Delta$};
		\node [style=port] (212) at (33.25, -3.25) {};
		\node [style=component] (213) at (33.25, -2.25) {$\varepsilon$};
		\node [style=none] (214) at (31.75, -2.25) {$\varepsilon$};
		\node [style=port] (215) at (33.75, -3.25) {};
		\node [style=function2] (216) at (31.75, -1.25) {$\mathsf{p}_0$};
		\node [style=port] (217) at (31.25, -3.25) {};
		\node [style=port] (218) at (31.75, -3.25) {};
		\node [style=port] (219) at (32.5, -3.25) {};
		\node [style=port] (220) at (30, -0.25) {$=$};
		\node [style=port] (221) at (30, -0.5) {Nat of};
		\node [style=port] (222) at (30, -1) {$\Delta$, $u$, $\eta$};
	\end{pgfonlayer}
	\begin{pgfonlayer}{edgelayer}
		\draw [style=wire] (54) to (47);
		\draw [style=wire] (46) to (56);
		\draw [style=wire] (48) to (51);
		\draw [style=wire] (48) to (49);
		\draw [style=wire] (49) to (43);
		\draw [style=wire] (51) to (43);
		\draw [style=wire] (44) to (41);
		\draw [style=wire, bend left] (50) to (53);
		\draw [style=wire, bend right] (50) to (42);
		\draw [style=wire] (50) to (45);
		\draw [style=wire, bend right] (52) to (55);
		\draw [style=wire, bend left] (52) to (40);
		\draw [style=wire] (55) to (47);
		\draw [style=wire] (40) to (56);
		\draw [style=wire] (45) to (52);
		\draw [style=wire] (57) to (53);
		\draw [style=wire, in=-90, out=165] (88) to (87);
		\draw [style=wire] (89) to (85);
		\draw [style=wire, in=15, out=-90] (85) to (88);
		\draw [style=wire] (88) to (86);
		\draw [style=wire] (91) to (90);
		\draw [style=wire] (90) to (42);
		\draw [style=wire] (108) to (101);
		\draw [style=wire] (100) to (110);
		\draw [style=wire] (102) to (105);
		\draw [style=wire] (102) to (103);
		\draw [style=wire] (103) to (97);
		\draw [style=wire] (105) to (97);
		\draw [style=wire] (98) to (95);
		\draw [style=wire, bend left] (104) to (107);
		\draw [style=wire, bend right] (104) to (96);
		\draw [style=wire] (104) to (99);
		\draw [style=wire, bend right] (106) to (109);
		\draw [style=wire, bend left] (106) to (94);
		\draw [style=wire] (109) to (101);
		\draw [style=wire] (94) to (110);
		\draw [style=wire] (99) to (106);
		\draw [style=wire] (111) to (107);
		\draw [style=wire] (113) to (112);
		\draw [style=wire] (112) to (96);
		\draw [style=wire, bend right] (116) to (119);
		\draw [style=wire, bend right=15, looseness=1.25] (114) to (118);
		\draw [style=wire, in=135, out=-30] (116) to (114);
		\draw [style=wire, bend right=15, looseness=1.25] (117) to (115);
		\draw [style=wire] (114) to (115);
		\draw [style=wire, bend right=15, looseness=1.25] (119) to (117);
		\draw [style=wire] (120) to (121);
		\draw [style=wire] (120) to (116);
		\draw [style=wire] (122) to (123);
		\draw [style=wire] (122) to (118);
		\draw [style=wire] (135) to (131);
		\draw [style=wire] (129) to (130);
		\draw [style=wire] (124) to (134);
		\draw [style=wire] (124) to (125);
		\draw [style=wire] (125) to (132);
		\draw [style=wire] (134) to (132);
		\draw [style=wire] (136) to (127);
		\draw [style=wire, bend right] (128) to (133);
		\draw [style=wire, bend left] (128) to (126);
		\draw [style=wire] (133) to (131);
		\draw [style=wire] (126) to (130);
		\draw [style=wire] (117) to (128);
		\draw [style=wire, bend right=15, looseness=1.25] (166) to (170);
		\draw [style=wire, in=135, out=-30] (168) to (166);
		\draw [style=wire, bend right=15, looseness=1.25] (169) to (167);
		\draw [style=wire] (166) to (167);
		\draw [style=wire, bend right=15, looseness=1.25] (171) to (169);
		\draw [style=wire] (172) to (173);
		\draw [style=wire] (172) to (168);
		\draw [style=wire] (174) to (175);
		\draw [style=wire] (174) to (170);
		\draw [style=wire] (187) to (183);
		\draw [style=wire] (181) to (182);
		\draw [style=wire] (176) to (186);
		\draw [style=wire] (176) to (177);
		\draw [style=wire] (177) to (184);
		\draw [style=wire] (186) to (184);
		\draw [style=wire] (188) to (179);
		\draw [style=wire, bend right] (180) to (185);
		\draw [style=wire, bend left] (180) to (178);
		\draw [style=wire] (185) to (183);
		\draw [style=wire] (178) to (182);
		\draw [style=wire] (169) to (180);
		\draw [style=wire] (189) to (171);
		\draw [style=wire, in=180, out=-165, looseness=1.50] (168) to (189);
		\draw [style=wire, in=0, out=-90, looseness=1.50] (190) to (189);
		\draw [style=wire, in=-90, out=0, looseness=1.25] (195) to (193);
		\draw [style=wire, in=180, out=-90] (196) to (195);
		\draw [style=wire] (198) to (199);
		\draw [style=wire, in=90, out=180] (194) to (200);
		\draw [style=wire] (197) to (194);
		\draw [style=wire] (203) to (196);
		\draw [style=wire, bend right=45, looseness=1.50] (200) to (203);
		\draw [style=wire, bend left=45, looseness=1.25] (201) to (203);
		\draw [style=wire] (202) to (201);
		\draw [style=wire, bend right=45, looseness=1.50] (204) to (206);
		\draw [style=wire, bend left=45, looseness=1.25] (205) to (206);
		\draw [style=wire] (206) to (193);
		\draw [style=wire, in=90, out=0] (194) to (204);
		\draw [style=wire, in=90, out=-90, looseness=1.25] (198) to (205);
		\draw [style=wire] (218) to (214.center);
		\draw [style=wire] (212) to (213);
		\draw [style=wire] (207) to (217);
		\draw [style=wire] (207) to (208);
		\draw [style=wire] (208) to (215);
		\draw [style=wire] (217) to (215);
		\draw [style=wire] (219) to (210);
		\draw [style=wire, bend right] (211) to (216);
		\draw [style=wire, bend left] (211) to (209);
		\draw [style=wire] (216) to (214.center);
		\draw [style=wire] (209) to (213);
		\draw [style=wire] (195) to (211);
	\end{pgfonlayer}
\end{tikzpicture}
}%
\] 
\[\resizebox{!}{6cm}{%
\begin{tikzpicture}
	\begin{pgfonlayer}{nodelayer}
		\node [style=differential] (1) at (12.5, 0.75) {$\bigotimes$};
		\node [style=port] (2) at (11.75, 2.5) {};
		\node [style=port] (3) at (13.25, 2.5) {};
		\node [style=component] (4) at (12.5, -0.25) {$\eta$};
		\node [style=port] (5) at (12.5, -1.25) {};
		\node [style=component] (6) at (11.75, 1.75) {$\varepsilon$};
		\node [style=duplicate] (7) at (7, 3.25) {$\Delta$};
		\node [style=port] (8) at (7, 4) {};
		\node [style=component] (9) at (9.25, 3.25) {$\eta$};
		\node [style=port] (10) at (9.25, 4) {};
		\node [style=component] (11) at (6.5, 0.5) {$u$};
		\node [style=duplicate] (12) at (7.5, -0.5) {$\nabla$};
		\node [style=component] (13) at (7.75, 2.25) {$\varepsilon$};
		\node [style=component] (14) at (9.25, 2.25) {$\varepsilon$};
		\node [style=component] (15) at (8.5, 0.5) {$\eta$};
		\node [style=port] (16) at (7.5, -1.5) {};
		\node [style=port] (17) at (8.5, 1.25) {$\bigotimes$};
		\node [style=component] (18) at (6.25, 2) {$e$};
		\node [style=duplicate] (19) at (-4.75, 6.25) {$\Delta$};
		\node [style=duplicate] (20) at (-3.75, -1.75) {$\nabla$};
		\node [style=component] (21) at (-5.5, 3.5) {$\delta$};
		\node [style=port] (22) at (-4.75, 7.25) {};
		\node [style=component] (23) at (-2.25, 6.25) {$\eta$};
		\node [style=port] (24) at (-2.25, 7.25) {};
		\node [style=function2] (25) at (-6, 4.75) {$\mathsf{i}_0$};
		\node [style=function2] (26) at (-5, 4.75) {$\mathsf{i}_1$};
		\node [style=component] (27) at (-5, 5.5) {$u$};
		\node [style=duplicate] (28) at (-5.5, 4.25) {$\nabla$};
		\node [style=function2] (29) at (-3.25, 4.75) {$\mathsf{i}_0$};
		\node [style=function2] (30) at (-2.25, 4.75) {$\mathsf{i}_1$};
		\node [style=duplicate] (31) at (-2.75, 4.25) {$\nabla$};
		\node [style=port] (32) at (-3.75, -2.75) {};
		\node [style=port] (33) at (-6.25, -0.5) {};
		\node [style=port] (34) at (-4.25, 3) {};
		\node [style=port] (35) at (-4.75, -0.5) {};
		\node [style=port] (36) at (-6.75, -0.5) {};
		\node [style=port] (37) at (-6.75, 3) {};
		\node [style=port] (38) at (-5.5, -0.5) {};
		\node [style=function2] (39) at (-6.25, 1.5) {$\mathsf{p}_0$};
		\node [style=component] (40) at (-4.75, 0.5) {$\varepsilon$};
		\node [style=port] (41) at (-4.25, -0.5) {};
		\node [style=function2] (42) at (-4.75, 1.5) {$\mathsf{p}_1$};
		\node [style=component] (43) at (-6.25, 0.5) {$\varepsilon$};
		\node [style=duplicate] (44) at (-5.5, 2.5) {$\Delta$};
		\node [style=component] (45) at (-2.75, -0.5) {$\eta$};
		\node [style=function2] (46) at (-3.5, 2.5) {$\mathsf{p}_0$};
		\node [style=component] (47) at (-2, 1.5) {$\varepsilon$};
		\node [style=function2] (48) at (-2, 2.5) {$\mathsf{p}_1$};
		\node [style=component] (49) at (-3.5, 1.5) {$\varepsilon$};
		\node [style=duplicate] (50) at (-2.75, 3.5) {$\Delta$};
		\node [style=duplicate] (51) at (-2.75, 0.5) {$\otimes$};
		\node [style=port] (52) at (2.25, 4.25) {$\Delta$};
		\node [style=port] (53) at (2.25, 5.5) {};
		\node [style=component] (54) at (4.5, 4.25) {$\eta$};
		\node [style=port] (55) at (4.5, 5.5) {};
		\node [style=component] (56) at (2.5, 2) {$u$};
		\node [style=duplicate] (57) at (2.75, -0.75) {$\nabla$};
		\node [style=component] (58) at (3, 2.75) {$\varepsilon$};
		\node [style=component] (59) at (4.5, 2.75) {$\varepsilon$};
		\node [style=component] (60) at (3.75, 0.5) {$\eta$};
		\node [style=port] (61) at (2.75, -1.5) {};
		\node [style=port] (62) at (3.75, 1.5) {$\bigotimes$};
		\node [style=function2] (63) at (1.5, 0.75) {$\bigotimes$};
		\node [style=port] (64) at (-8, 2) {$=$};
		\node [style=port] (65) at (-8, 1.75) {Nat of};
		\node [style=port] (66) at (-8, 1.25) {$\nabla$ and $\eta$};
		\node [style=port] (67) at (-0.5, 2) {$=$};
		\node [style=port] (68) at (-0.5, 1.5) {(\ref{mcoder}) };
		\node [style=port] (69) at (5.25, 1.5) {$=$};
		\node [style=port] (70) at (5.25, 1) {(\ref{nablam2})};
		\node [style=port] (71) at (10, 1.5) {$=$};
		\node [style=port] (72) at (10.25, 1) {(\ref{coalgeq}) + (\ref{moneq})};
		\node [style=port] (73) at (-0.5, 1) {+ Lem \ref{seelyish}};
		\node [style=port] (74) at (10.25, 0.5) {+ \textbf{[dC.3]}};
	\end{pgfonlayer}
	\begin{pgfonlayer}{edgelayer}
		\draw [style=wire, in=-90, out=15] (1) to (3);
		\draw [style=wire] (1) to (4);
		\draw [style=wire] (4) to (5);
		\draw [style=wire] (2) to (6);
		\draw [style=wire, in=165, out=-90] (6) to (1);
		\draw [style=wire] (9) to (10);
		\draw [style=wire] (8) to (7);
		\draw [style=wire, in=41, out=-90, looseness=1.25] (15) to (12);
		\draw [style=wire] (12) to (16);
		\draw [style=wire] (17) to (15);
		\draw [style=wire, in=150, out=-90] (13) to (17);
		\draw [style=wire, in=30, out=-90] (14) to (17);
		\draw [style=wire, in=135, out=-90, looseness=1.25] (11) to (12);
		\draw [style=wire, in=90, out=-150] (7) to (18);
		\draw [style=wire] (9) to (14);
		\draw [style=wire, in=90, out=-30, looseness=1.25] (7) to (13);
		\draw [style=wire] (23) to (24);
		\draw [style=wire, in=90, out=180] (19) to (25);
		\draw [style=wire] (22) to (19);
		\draw [style=wire] (28) to (21);
		\draw [style=wire, bend right=45, looseness=1.50] (25) to (28);
		\draw [style=wire, bend left=45, looseness=1.25] (26) to (28);
		\draw [style=wire] (27) to (26);
		\draw [style=wire, bend right=45, looseness=1.50] (29) to (31);
		\draw [style=wire, bend left=45, looseness=1.25] (30) to (31);
		\draw [style=wire, in=90, out=0] (19) to (29);
		\draw [style=wire, in=90, out=-90, looseness=1.25] (23) to (30);
		\draw [style=wire] (33) to (43);
		\draw [style=wire] (35) to (40);
		\draw [style=wire] (37) to (36);
		\draw [style=wire] (37) to (34);
		\draw [style=wire] (34) to (41);
		\draw [style=wire] (36) to (41);
		\draw [style=wire, bend right] (44) to (39);
		\draw [style=wire, bend left] (44) to (42);
		\draw [style=wire] (39) to (43);
		\draw [style=wire] (42) to (40);
		\draw [style=wire, bend right] (50) to (46);
		\draw [style=wire, bend left] (50) to (48);
		\draw [style=wire] (46) to (49);
		\draw [style=wire] (48) to (47);
		\draw [style=wire] (21) to (44);
		\draw [style=wire, in=180, out=-90] (38) to (20);
		\draw [style=wire, in=0, out=-90] (45) to (20);
		\draw [style=wire] (20) to (32);
		\draw [style=wire] (31) to (50);
		\draw [style=wire] (51) to (45);
		\draw [style=wire, in=180, out=-90] (49) to (51);
		\draw [style=wire, in=0, out=-90] (47) to (51);
		\draw [style=wire] (54) to (55);
		\draw [style=wire] (53) to (52);
		\draw [style=wire, in=15, out=-90] (60) to (57);
		\draw [style=wire] (57) to (61);
		\draw [style=wire] (62) to (60);
		\draw [style=wire, in=150, out=-90] (58) to (62);
		\draw [style=wire, in=30, out=-90] (59) to (62);
		\draw [style=wire, in=180, out=-90] (63) to (57);
		\draw [style=wire, in=0, out=-90, looseness=1.25] (56) to (63);
		\draw [style=wire, in=165, out=180, looseness=0.75] (52) to (63);
		\draw [style=wire, in=90, out=0, looseness=0.75] (52) to (58);
		\draw [style=wire] (54) to (59);
	\end{pgfonlayer}
\end{tikzpicture}
}%
\]
\end{proof} 

Conversly, the alternative chain rule {\bf [$\text{dC.4}^\prime$]} and the monoidal rule {\bf [dC.m]} imply the chain rule {\bf [dC.4]}.

\begin{lemma}\label{dCm2} For the induced additive bialgebra modality of an additive linear category, any natural transformation $\eta: A \to \oc A$ which satisfies the alternative chain rule {\bf [$\text{dC.4}^\prime$]} and the monoidal rule {\bf [dC.m]} also satisfies the chain rule {\bf [dC.4]}. 
\end{lemma}
\begin{proof} Using Proposition \ref{monoidalnabla}, the alternative chain rule {\bf [$\text{dC.4}^\prime$]}, and the bialgebra modality identities we have: 
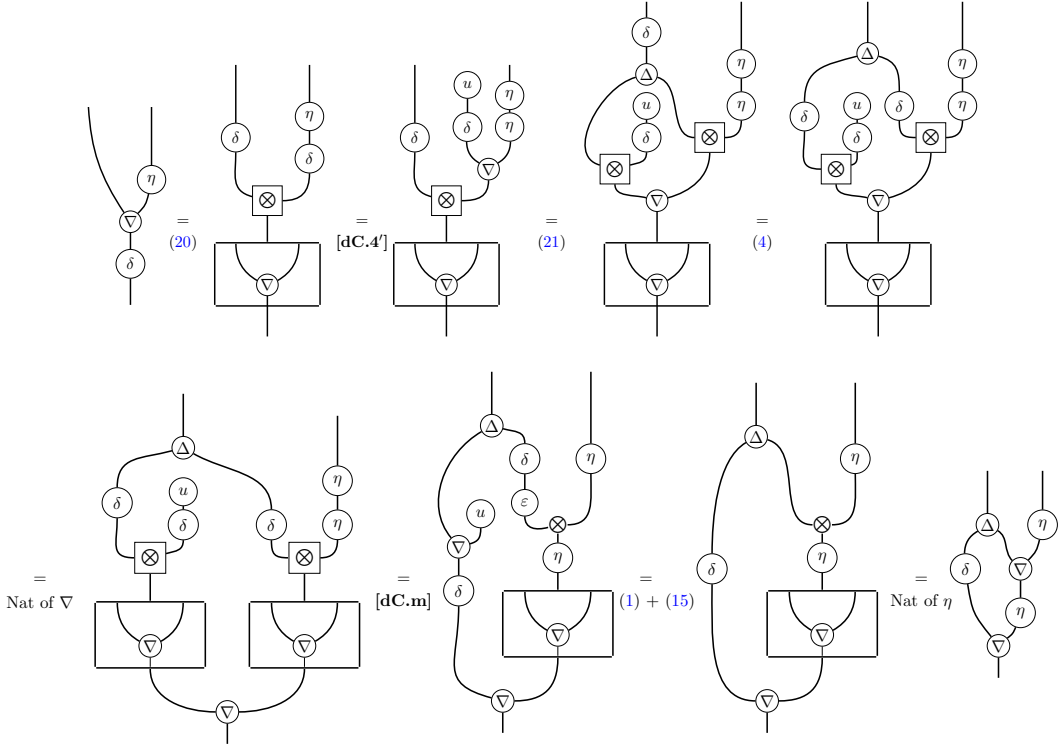
\[\resizebox{!}{4.5cm}{%
\begin{tikzpicture}
	\begin{pgfonlayer}{nodelayer}
		\node [style=component] (0) at (2, 1.75) {$\eta$};
		\node [style=port] (1) at (0.5, 3.5) {};
		\node [style=port] (2) at (2, 3.5) {};
		\node [style=component] (3) at (1.5, -0.25) {$\delta$};
		\node [style=duplicate] (4) at (1.5, 0.75) {$\nabla$};
		\node [style=port] (5) at (1.5, -1.25) {};
		\node [style=port] (6) at (2.75, 0.75) {$=$};
		\node [style=function2] (7) at (4.75, 1.25) {$\bigotimes$};
		\node [style=component] (8) at (5.75, 3.25) {$\eta$};
		\node [style=port] (9) at (5.75, 4.5) {};
		\node [style=component] (10) at (5.75, 2.25) {$\delta$};
		\node [style=component] (11) at (4, 2.75) {$\delta$};
		\node [style=port] (12) at (4, 4.5) {};
		\node [style=duplicate] (38) at (4.75, -0.75) {$\nabla$};
		\node [style=port] (39) at (5.5, 0.25) {};
		\node [style=port] (40) at (4, 0.25) {};
		\node [style=port] (41) at (3.5, 0.25) {};
		\node [style=port] (42) at (6, 0.25) {};
		\node [style=port] (43) at (6, -1.25) {};
		\node [style=port] (44) at (3.5, -1.25) {};
		\node [style=port] (45) at (4.75, -2) {};
		\node [style=port] (46) at (4.75, 0.25) {};
		\node [style=port] (47) at (4.75, 0.25) {};
		\node [style=function2] (48) at (9, 1.25) {$\bigotimes$};
		\node [style=duplicate] (49) at (10, 2) {$\nabla$};
		\node [style=component] (50) at (8.25, 2.75) {$\delta$};
		\node [style=port] (51) at (8.25, 4.5) {};
		\node [style=component] (54) at (9.5, 3) {$\delta$};
		\node [style=component] (55) at (10.5, 3) {$\eta$};
		\node [style=component] (56) at (9.5, 4) {$u$};
		\node [style=component] (57) at (10.5, 3.75) {$\eta$};
		\node [style=port] (58) at (10.5, 4.5) {};
		\node [style=duplicate] (59) at (9, -0.75) {$\nabla$};
		\node [style=port] (60) at (9.75, 0.25) {};
		\node [style=port] (61) at (8.25, 0.25) {};
		\node [style=port] (62) at (7.75, 0.25) {};
		\node [style=port] (63) at (10.25, 0.25) {};
		\node [style=port] (64) at (10.25, -1.25) {};
		\node [style=port] (65) at (7.75, -1.25) {};
		\node [style=port] (66) at (9, -2) {};
		\node [style=port] (67) at (9, 0.25) {};
		\node [style=port] (68) at (9, 0.25) {};
		\node [style=port] (69) at (2.75, 0.25) {(\ref{nablam1})};
		\node [style=port] (70) at (7, 0.75) {$=$};
		\node [style=port] (71) at (7, 0.25) {{\bf [$\text{dC.4}^\prime$]}};
		\node [style=function2] (72) at (15.25, 2.75) {$\bigotimes$};
		\node [style=duplicate] (73) at (13.75, 4.25) {$\Delta$};
		\node [style=component] (74) at (13.75, 5.25) {$\delta$};
		\node [style=duplicate] (75) at (14, 1.25) {$\nabla$};
		\node [style=component] (76) at (13.75, 2.75) {$\delta$};
		\node [style=component] (77) at (16, 3.5) {$\eta$};
		\node [style=component] (78) at (13.75, 3.5) {$u$};
		\node [style=component] (79) at (16, 4.5) {$\eta$};
		\node [style=port] (80) at (16, 6) {};
		\node [style=function2] (81) at (13, 2) {$\bigotimes$};
		\node [style=port] (84) at (13.75, 6) {};
		\node [style=duplicate] (85) at (14, -0.75) {$\nabla$};
		\node [style=port] (86) at (14.75, 0.25) {};
		\node [style=port] (87) at (13.25, 0.25) {};
		\node [style=port] (88) at (12.75, 0.25) {};
		\node [style=port] (89) at (15.25, 0.25) {};
		\node [style=port] (90) at (15.25, -1.25) {};
		\node [style=port] (91) at (12.75, -1.25) {};
		\node [style=port] (92) at (14, -2) {};
		\node [style=port] (93) at (14, 0.25) {};
		\node [style=port] (94) at (14, 0.25) {};
		\node [style=port] (95) at (11.5, 0.75) {$=$};
		\node [style=port] (96) at (11.5, 0.25) {(\ref{nablam2})};
		\node [style=function2] (97) at (20.5, 2.75) {$\bigotimes$};
		\node [style=duplicate] (98) at (19, 4.75) {$\Delta$};
		\node [style=duplicate] (99) at (19.25, 1.25) {$\nabla$};
		\node [style=component] (100) at (18.75, 2.75) {$\delta$};
		\node [style=component] (101) at (21.25, 3.5) {$\eta$};
		\node [style=component] (102) at (18.75, 3.5) {$u$};
		\node [style=component] (103) at (21.25, 4.5) {$\eta$};
		\node [style=port] (104) at (21.25, 6) {};
		\node [style=function2] (105) at (18.25, 2) {$\bigotimes$};
		\node [style=port] (108) at (19, 6) {};
		\node [style=component] (109) at (17.5, 3.25) {$\delta$};
		\node [style=component] (110) at (19.75, 3.5) {$\delta$};
		\node [style=port] (111) at (16.5, 0.75) {$=$};
		\node [style=port] (112) at (16.5, 0.25) {(\ref{coalgeq})};
		\node [style=duplicate] (113) at (19.25, -0.75) {$\nabla$};
		\node [style=port] (114) at (20, 0.25) {};
		\node [style=port] (115) at (18.5, 0.25) {};
		\node [style=port] (116) at (18, 0.25) {};
		\node [style=port] (117) at (20.5, 0.25) {};
		\node [style=port] (118) at (20.5, -1.25) {};
		\node [style=port] (119) at (18, -1.25) {};
		\node [style=port] (120) at (19.25, -2) {};
		\node [style=port] (121) at (19.25, 0.25) {};
		\node [style=port] (122) at (19.25, 0.25) {};
	\end{pgfonlayer}
	\begin{pgfonlayer}{edgelayer}
		\draw [style=wire, bend right=15, looseness=1.25] (4) to (0);
		\draw [style=wire, bend left=15] (4) to (1);
		\draw [style=wire] (3) to (4);
		\draw [style=wire] (3) to (5);
		\draw [style=wire] (2) to (0);
		\draw [style=wire] (8) to (9);
		\draw [style=wire, in=0, out=-90, looseness=1.25] (10) to (7);
		\draw [style=wire] (8) to (10);
		\draw [style=wire] (12) to (11);
		\draw [style=wire, in=165, out=-90, looseness=1.25] (11) to (7);
		\draw [style=wire, bend left] (38) to (40);
		\draw [style=wire, bend right] (38) to (39);
		\draw [style=wire] (44) to (43);
		\draw [style=wire] (44) to (41);
		\draw [style=wire] (41) to (42);
		\draw [style=wire] (43) to (42);
		\draw [style=wire] (38) to (45);
		\draw [style=wire] (7) to (47);
		\draw [style=wire, in=0, out=-90, looseness=1.25] (49) to (48);
		\draw [style=wire] (51) to (50);
		\draw [style=wire, in=165, out=-90, looseness=1.25] (50) to (48);
		\draw [style=wire] (56) to (54);
		\draw [style=wire] (58) to (57);
		\draw [style=wire] (57) to (55);
		\draw [style=wire, bend right] (54) to (49);
		\draw [style=wire, bend left] (55) to (49);
		\draw [style=wire, bend left] (59) to (61);
		\draw [style=wire, bend right] (59) to (60);
		\draw [style=wire] (65) to (64);
		\draw [style=wire] (65) to (62);
		\draw [style=wire] (62) to (63);
		\draw [style=wire] (64) to (63);
		\draw [style=wire] (59) to (66);
		\draw [style=wire] (48) to (68);
		\draw [style=wire] (74) to (73);
		\draw [style=wire, in=180, out=-15] (73) to (72);
		\draw [style=wire, in=15, out=-90] (72) to (75);
		\draw [style=wire] (78) to (76);
		\draw [style=wire] (80) to (79);
		\draw [style=wire] (79) to (77);
		\draw [style=wire, in=165, out=-165] (73) to (81);
		\draw [style=wire, in=165, out=-90] (81) to (75);
		\draw [style=wire, in=-90, out=15, looseness=1.25] (81) to (76);
		\draw [style=wire] (84) to (74);
		\draw [style=wire, in=0, out=-90, looseness=1.25] (77) to (72);
		\draw [style=wire, bend left] (85) to (87);
		\draw [style=wire, bend right] (85) to (86);
		\draw [style=wire] (91) to (90);
		\draw [style=wire] (91) to (88);
		\draw [style=wire] (88) to (89);
		\draw [style=wire] (90) to (89);
		\draw [style=wire] (85) to (92);
		\draw [style=wire] (75) to (94);
		\draw [style=wire, in=15, out=-90] (97) to (99);
		\draw [style=wire] (102) to (100);
		\draw [style=wire] (104) to (103);
		\draw [style=wire] (103) to (101);
		\draw [style=wire, in=165, out=-90] (105) to (99);
		\draw [style=wire, in=-90, out=15, looseness=1.25] (105) to (100);
		\draw [style=wire, in=0, out=-90, looseness=1.25] (101) to (97);
		\draw [style=wire] (108) to (98);
		\draw [style=wire, in=90, out=-15, looseness=1.25] (98) to (110);
		\draw [style=wire, in=90, out=-165, looseness=1.25] (98) to (109);
		\draw [style=wire, in=180, out=-90, looseness=1.50] (109) to (105);
		\draw [style=wire, in=180, out=-90, looseness=1.50] (110) to (97);
		\draw [style=wire, bend left] (113) to (115);
		\draw [style=wire, bend right] (113) to (114);
		\draw [style=wire] (119) to (118);
		\draw [style=wire] (119) to (116);
		\draw [style=wire] (116) to (117);
		\draw [style=wire] (118) to (117);
		\draw [style=wire] (113) to (120);
		\draw [style=wire] (99) to (122);
	\end{pgfonlayer}
\end{tikzpicture}
}%
\]
\[\resizebox{!}{5cm}{%
\begin{tikzpicture}
	\begin{pgfonlayer}{nodelayer}
		\node [style=function2] (0) at (5, 0.5) {$\bigotimes$};
		\node [style=duplicate] (1) at (2.25, 3) {$\Delta$};
		\node [style=component] (2) at (2.25, 1.25) {$\delta$};
		\node [style=component] (3) at (5.75, 1.25) {$\eta$};
		\node [style=component] (4) at (2.25, 2) {$u$};
		\node [style=component] (5) at (5.75, 2.25) {$\eta$};
		\node [style=port] (6) at (5.75, 3.75) {};
		\node [style=function2] (7) at (1.5, 0.5) {$\bigotimes$};
		\node [style=duplicate] (8) at (3.25, -3) {$\nabla$};
		\node [style=port] (9) at (3.25, -3.75) {};
		\node [style=port] (10) at (2.25, 4.25) {};
		\node [style=component] (11) at (0.75, 1.75) {$\delta$};
		\node [style=component] (12) at (4.25, 1.25) {$\delta$};
		\node [style=duplicate] (13) at (1.5, -1.5) {$\nabla$};
		\node [style=port] (14) at (2.25, -0.5) {};
		\node [style=port] (15) at (0.75, -0.5) {};
		\node [style=port] (16) at (0.25, -0.5) {};
		\node [style=port] (17) at (2.75, -0.5) {};
		\node [style=port] (18) at (2.75, -2) {};
		\node [style=port] (19) at (0.25, -2) {};
		\node [style=port] (20) at (1.5, -0.5) {};
		\node [style=port] (21) at (1.5, -0.5) {};
		\node [style=duplicate] (22) at (5, -1.5) {$\nabla$};
		\node [style=port] (23) at (5.75, -0.5) {};
		\node [style=port] (24) at (4.25, -0.5) {};
		\node [style=port] (25) at (3.75, -0.5) {};
		\node [style=port] (26) at (6.25, -0.5) {};
		\node [style=port] (27) at (6.25, -2) {};
		\node [style=port] (28) at (3.75, -2) {};
		\node [style=port] (29) at (5, -0.5) {};
		\node [style=port] (30) at (5, -0.5) {};
		\node [style=port] (31) at (1.5, -2) {};
		\node [style=port] (32) at (5, -2) {};
		\node [style=port] (33) at (-1, 0) {};
		\node [style=port] (34) at (-1, -0.5) {Nat of $\nabla$};
		\node [style=port] (35) at (-1, 0) {$=$};
		\node [style=port] (36) at (10.75, 1.25) {$\bigotimes$};
		\node [style=duplicate] (37) at (9.25, 3.5) {$\Delta$};
		\node [style=component] (38) at (9, 1.5) {$u$};
		\node [style=component] (39) at (11.5, 2.75) {$\eta$};
		\node [style=port] (40) at (11.5, 4.75) {};
		\node [style=duplicate] (41) at (8.5, 0.75) {$\nabla$};
		\node [style=duplicate] (42) at (9.5, -2.75) {$\nabla$};
		\node [style=port] (43) at (9.5, -3.5) {};
		\node [style=port] (44) at (9.25, 4.75) {};
		\node [style=component] (45) at (10, 2.75) {$\delta$};
		\node [style=component] (46) at (8.5, -0.25) {$\delta$};
		\node [style=component] (48) at (10.75, 0.5) {$\eta$};
		\node [style=component] (49) at (10, 1.75) {$\varepsilon$};
		\node [style=duplicate] (50) at (10.75, -1.25) {$\nabla$};
		\node [style=port] (51) at (11.5, -0.25) {};
		\node [style=port] (52) at (10, -0.25) {};
		\node [style=port] (53) at (9.5, -0.25) {};
		\node [style=port] (54) at (12, -0.25) {};
		\node [style=port] (55) at (12, -1.75) {};
		\node [style=port] (56) at (9.5, -1.75) {};
		\node [style=port] (57) at (10.75, -0.25) {};
		\node [style=port] (58) at (10.75, -0.25) {};
		\node [style=port] (59) at (10.75, -1.75) {};
		\node [style=port] (60) at (7.25, 0) {$=$};
		\node [style=port] (61) at (7.25, -0.5) {\textbf{[dC.m]}};
		\node [style=port] (62) at (20.5, 2.5) {};
		\node [style=duplicate] (63) at (20.5, 1.25) {$\Delta$};
		\node [style=component] (64) at (20, 0.25) {$\delta$};
		\node [style=port] (65) at (21.75, 2.5) {};
		\node [style=component] (66) at (21.75, 1.25) {$\eta$};
		\node [style=duplicate] (67) at (21.25, 0.25) {$\nabla$};
		\node [style=duplicate] (68) at (20.75, -1.5) {$\nabla$};
		\node [style=component] (69) at (21.25, -0.75) {$\eta$};
		\node [style=port] (70) at (20.75, -2.25) {};
		\node [style=port] (71) at (16.75, 1.25) {$\bigotimes$};
		\node [style=duplicate] (72) at (15.25, 3.25) {$\Delta$};
		\node [style=component] (74) at (17.5, 2.75) {$\eta$};
		\node [style=port] (75) at (17.5, 4.5) {};
		\node [style=duplicate] (77) at (15.5, -2.75) {$\nabla$};
		\node [style=port] (78) at (15.5, -3.5) {};
		\node [style=port] (79) at (15.25, 4.5) {};
		\node [style=component] (81) at (14.25, 0.25) {$\delta$};
		\node [style=component] (82) at (16.75, 0.5) {$\eta$};
		\node [style=duplicate] (84) at (16.75, -1.25) {$\nabla$};
		\node [style=port] (85) at (17.5, -0.25) {};
		\node [style=port] (86) at (16, -0.25) {};
		\node [style=port] (87) at (15.5, -0.25) {};
		\node [style=port] (88) at (18, -0.25) {};
		\node [style=port] (89) at (18, -1.75) {};
		\node [style=port] (90) at (15.5, -1.75) {};
		\node [style=port] (91) at (16.75, -0.25) {};
		\node [style=port] (92) at (16.75, -0.25) {};
		\node [style=port] (93) at (16.75, -1.75) {};
		\node [style=port] (94) at (12.75, 0) {$=$};
		\node [style=port] (95) at (13, -0.5) {(\ref{comonad}) + (\ref{moneq})};
		\node [style=port] (96) at (19, 0) {$=$};
		\node [style=port] (97) at (19, -0.5) {Nat of $\eta$};
	\end{pgfonlayer}
	\begin{pgfonlayer}{edgelayer}
		\draw [style=wire] (4) to (2);
		\draw [style=wire] (6) to (5);
		\draw [style=wire] (5) to (3);
		\draw [style=wire] (8) to (9);
		\draw [style=wire, in=-90, out=15, looseness=1.25] (7) to (2);
		\draw [style=wire, in=0, out=-90, looseness=1.25] (3) to (0);
		\draw [style=wire] (10) to (1);
		\draw [style=wire, in=90, out=-30] (1) to (12);
		\draw [style=wire, in=90, out=-165, looseness=1.25] (1) to (11);
		\draw [style=wire, in=180, out=-90, looseness=1.50] (11) to (7);
		\draw [style=wire, in=180, out=-90, looseness=1.50] (12) to (0);
		\draw [style=wire, bend left] (13) to (15);
		\draw [style=wire, bend right] (13) to (14);
		\draw [style=wire] (19) to (18);
		\draw [style=wire] (19) to (16);
		\draw [style=wire] (16) to (17);
		\draw [style=wire] (18) to (17);
		\draw [style=wire] (7) to (21);
		\draw [style=wire, bend left] (22) to (24);
		\draw [style=wire, bend right] (22) to (23);
		\draw [style=wire] (28) to (27);
		\draw [style=wire] (28) to (25);
		\draw [style=wire] (25) to (26);
		\draw [style=wire] (27) to (26);
		\draw [style=wire] (0) to (30);
		\draw (13) to (31);
		\draw (22) to (32);
		\draw [style=wire, in=180, out=-90] (31) to (8);
		\draw [style=wire, in=0, out=-90] (32) to (8);
		\draw [style=wire] (42) to (43);
		\draw [style=wire, in=-90, out=15, looseness=1.25] (41) to (38);
		\draw [style=wire, in=0, out=-90, looseness=1.25] (39) to (36);
		\draw [style=wire] (44) to (37);
		\draw [style=wire, in=90, out=-15, looseness=1.25] (37) to (45);
		\draw [style=wire] (41) to (46);
		\draw [style=wire, in=165, out=-90, looseness=1.50] (46) to (42);
		\draw [style=wire, in=130, out=-150] (37) to (41);
		\draw [style=wire] (36) to (48);
		\draw [style=wire] (45) to (49);
		\draw [style=wire, in=-165, out=-90, looseness=1.25] (49) to (36);
		\draw [style=wire] (40) to (39);
		\draw [style=wire, bend left] (50) to (52);
		\draw [style=wire, bend right] (50) to (51);
		\draw [style=wire] (56) to (55);
		\draw [style=wire] (56) to (53);
		\draw [style=wire] (53) to (54);
		\draw [style=wire] (55) to (54);
		\draw (50) to (59);
		\draw [style=wire] (48) to (58);
		\draw [style=wire, in=0, out=-90, looseness=1.25] (59) to (42);
		\draw [style=wire, bend right] (63) to (64);
		\draw [style=wire] (62) to (63);
		\draw [style=wire, bend right=15, looseness=1.25] (67) to (66);
		\draw [style=wire] (66) to (65);
		\draw [style=wire, in=135, out=-30] (63) to (67);
		\draw [style=wire, bend right=15, looseness=1.25] (68) to (69);
		\draw [style=wire] (67) to (69);
		\draw [style=wire, bend right=15, looseness=1.25] (64) to (68);
		\draw [style=wire] (68) to (70);
		\draw [style=wire] (77) to (78);
		\draw [style=wire, in=0, out=-90, looseness=1.25] (74) to (71);
		\draw [style=wire] (79) to (72);
		\draw [style=wire, in=180, out=-90] (81) to (77);
		\draw [style=wire] (71) to (82);
		\draw [style=wire] (75) to (74);
		\draw [style=wire, bend left] (84) to (86);
		\draw [style=wire, bend right] (84) to (85);
		\draw [style=wire] (90) to (89);
		\draw [style=wire] (90) to (87);
		\draw [style=wire] (87) to (88);
		\draw [style=wire] (89) to (88);
		\draw (84) to (93);
		\draw [style=wire] (82) to (92);
		\draw [style=wire, in=0, out=-90, looseness=1.25] (93) to (77);
		\draw [style=wire, in=90, out=-150, looseness=0.75] (72) to (81);
		\draw [style=wire, in=180, out=-15] (72) to (71);
	\end{pgfonlayer}
\end{tikzpicture}
}%
\]
\end{proof} 

\begin{corollary}\label{cor1} For the induced additive bialgebra modality of an additive linear category, the following are equivalent for a natural transformation $\eta: A \to \oc A$: 
\begin{enumerate}[{\em (i)}]
\item $\eta$ is a codereliction; 
\item $\eta$ satisfies the linear rule {\bf [dC.3]} and the chain rule {\bf [dC.4]};
\item $\eta$ satisfies the linear rule {\bf [dC.3]}, the alternative chain rule {\bf [$\text{dC.4}^\prime$]} and the monoidal rule {\bf [dC.m]}.
\end{enumerate}
\end{corollary} 

Part {\em (iii)} of the above corollary is the definition of Fiore's creation map \cite{fiore2007differential}. This shows that, for additive bialgebra modalities or equivalently monoidal coalgebra modalities, the original definition of a codereliction is equivalent to Fiore's creation map.

Turning our attention to deriving transformations for additive bialgebra modalities, we begin by noticing that satisfying the Leibniz rule is equivalent to satisfying the $\nabla$-rule: 

\begin{proposition}\label{nablaleibniz} 
For an additive bialgebra modality $(\oc,  \delta, \varepsilon, \Delta, e, \nabla, u)$, the following are equivalent for a natural transformation $\mathsf{d}: \oc A \otimes A \to \oc A$ which satisfies the linear rule {\bf [d.3]}:
\begin{enumerate}[{\em (i)}]
\item $\mathsf{d}$ satisfies the product rule {\bf [d.2]};
\item $\mathsf{d}$ satisfies the $\nabla$-rule {\bf [d.$\nabla$]}. 
\end{enumerate}
\end{proposition} 

\begin{proof} ~

\noindent
{\bf [d.$\nabla$]} $\Rightarrow$ {\bf [d.2]}: It is easy to see that since $\mathsf{d}$ satisfies {\bf [d.3]} that  $(u \otimes 1)\mathsf{d}: A \to \oc A$ satisfies {\bf [dC.3]}, the linear rule for coderelictions. 
However, by Lemma \ref{dc3dc4}, this implies that $(u \otimes 1)\mathsf{d}$ satisfies {\bf [dC.2]}, the product rule for coderelictions. Since $\mathsf{d}$ satisfies \textbf{[d.$\nabla$]}, then Lemma \ref{diden} holds. And so we have:
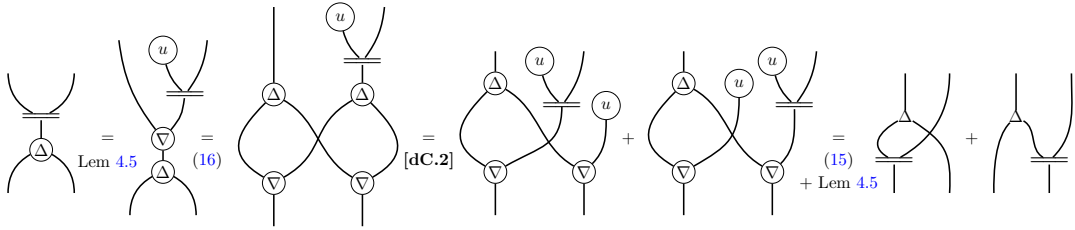
\[\resizebox{!}{3cm}{%
\begin{tikzpicture}
	\begin{pgfonlayer}{nodelayer}
		\node [style=differential] (0) at (1, 1) {{\bf =\!=\!=}};
		\node [style=port] (1) at (1.75, 2) {};
		\node [style=port] (2) at (0.25, 2) {};
		\node [style=port] (3) at (0.25, -0.75) {};
		\node [style=duplicate] (4) at (1, 0.25) {$\Delta$};
		\node [style=port] (5) at (1.75, -0.75) {};
		\node [style=port] (6) at (2.5, 0.5) {$=$};
		\node [style=port] (15) at (2.5, 0) {Lem \ref{diden}};
		\node [style=port] (16) at (4.5, -1.25) {};
		\node [style=port] (17) at (3, -1.25) {};
		\node [style=duplicate] (18) at (3.75, -0.25) {$\Delta$};
		\node [style=port] (19) at (2.75, 2.75) {};
		\node [style=duplicate] (20) at (3.75, 0.5) {$\nabla$};
		\node [style=component] (21) at (3.75, 2.5) {$u$};
		\node [style=port] (22) at (4.75, 2.75) {};
		\node [style=codifferential] (23) at (4.25, 1.5) {{\bf =\!=\!=}};
		\node [style=port] (26) at (6.25, 3.5) {};
		\node [style=component] (27) at (7.75, 3.25) {$u$};
		\node [style=port] (28) at (8.75, 3.5) {};
		\node [style=codifferential] (29) at (8.25, 2.25) {{\bf =\!=\!=}};
		\node [style=port] (30) at (8.25, -1.5) {};
		\node [style=duplicate] (31) at (6.25, -0.5) {$\nabla$};
		\node [style=duplicate] (32) at (6.25, 1.5) {$\Delta$};
		\node [style=duplicate] (33) at (8.25, -0.5) {$\nabla$};
		\node [style=duplicate] (34) at (8.25, 1.5) {$\Delta$};
		\node [style=port] (35) at (6.25, -1.5) {};
		\node [style=port] (36) at (4.75, 0.5) {$=$};
		\node [style=port] (37) at (4.75, 0) {(\ref{bialgeq})};
		\node [style=port] (38) at (11.25, 2.5) {};
		\node [style=component] (39) at (13.75, 1.25) {$u$};
		\node [style=port] (40) at (13.25, -1.25) {};
		\node [style=duplicate] (41) at (11.25, -0.25) {$\nabla$};
		\node [style=duplicate] (42) at (11.25, 1.75) {$\Delta$};
		\node [style=duplicate] (43) at (13.25, -0.25) {$\nabla$};
		\node [style=port] (44) at (12.75, 1.25) {};
		\node [style=port] (45) at (11.25, -1.25) {};
		\node [style=port] (46) at (13.25, 2.5) {};
		\node [style=codifferential] (47) at (12.75, 1.25) {{\bf =\!=\!=}};
		\node [style=component] (48) at (12.25, 2.25) {$u$};
		\node [style=port] (49) at (9.75, 0.5) {$=$};
		\node [style=port] (50) at (9.75, 0) {\textbf{[dC.2]}};
		\node [style=port] (51) at (14.25, 0.5) {$+$};
		\node [style=port] (52) at (15.5, 2.5) {};
		\node [style=component] (53) at (17.5, 2.25) {$u$};
		\node [style=port] (54) at (18.5, 2.5) {};
		\node [style=codifferential] (55) at (18, 1.25) {{\bf =\!=\!=}};
		\node [style=port] (56) at (17.5, -1.25) {};
		\node [style=duplicate] (57) at (15.5, -0.25) {$\nabla$};
		\node [style=duplicate] (58) at (15.5, 1.75) {$\Delta$};
		\node [style=duplicate] (59) at (17.5, -0.25) {$\nabla$};
		\node [style=component] (60) at (16.75, 1.75) {$u$};
		\node [style=port] (61) at (15.5, -1.25) {};
		\node [style=port] (91) at (19, 0.5) {$=$};
		\node [style=port] (92) at (19, 0) {(\ref{moneq})};
		\node [style=port] (103) at (20.5, 2) {};
		\node [style=differential] (104) at (20.25, 0) {{\bf =\!=\!=}};
		\node [style=port] (105) at (21.5, 2) {};
		\node [style=port] (106) at (20.5, 1) {$\Delta$};
		\node [style=port] (107) at (20.25, -0.75) {};
		\node [style=port] (108) at (21.5, -0.75) {};
		\node [style=port] (109) at (24.25, 2) {};
		\node [style=differential] (110) at (23.75, 0) {{\bf =\!=\!=}};
		\node [style=port] (111) at (23.75, -0.75) {};
		\node [style=port] (112) at (22.5, -0.75) {};
		\node [style=port] (113) at (23, 2) {};
		\node [style=port] (114) at (23, 1) {$\Delta$};
		\node [style=port] (115) at (22, 0.5) {$+$};
		\node [style=port] (116) at (19, -0.5) {+ Lem \ref{diden}};
	\end{pgfonlayer}
	\begin{pgfonlayer}{edgelayer}
		\draw [style=wire, bend right] (0) to (1);
		\draw [style=wire, bend left] (0) to (2);
		\draw [style=wire, bend right] (4) to (3);
		\draw [style=wire, bend left] (4) to (5);
		\draw [style=wire] (0) to (4);
		\draw [style=wire, bend left] (18) to (16);
		\draw [style=wire, bend right] (18) to (17);
		\draw [style=wire, bend right=15, looseness=1.25] (20) to (23);
		\draw [style=wire, bend left=15] (20) to (19);
		\draw [style=wire, bend left=15] (23) to (21);
		\draw [style=wire, bend right=15] (23) to (22);
		\draw [style=wire] (20) to (18);
		\draw [style=wire, bend left=15] (29) to (27);
		\draw [style=wire, bend right=15] (29) to (28);
		\draw [style=wire] (31) to (35);
		\draw [style=wire, in=150, out=-135, looseness=1.50] (32) to (31);
		\draw [style=wire] (33) to (30);
		\draw [style=wire, in=30, out=-45, looseness=1.50] (34) to (33);
		\draw [style=wire, in=30, out=-150, looseness=1.25] (34) to (31);
		\draw [style=wire, in=150, out=-30] (32) to (33);
		\draw [style=wire] (29) to (34);
		\draw [style=wire] (26) to (32);
		\draw [style=wire] (41) to (45);
		\draw [style=wire, in=150, out=-135, looseness=1.50] (42) to (41);
		\draw [style=wire] (43) to (40);
		\draw [style=wire, in=30, out=-90, looseness=1.25] (44) to (41);
		\draw [style=wire, in=150, out=-30] (42) to (43);
		\draw [style=wire] (38) to (42);
		\draw [style=wire, in=45, out=-90] (39) to (43);
		\draw [style=wire, bend left=15] (47) to (48);
		\draw [style=wire, bend right=15] (47) to (46);
		\draw [style=wire, bend left=15] (55) to (53);
		\draw [style=wire, bend right=15] (55) to (54);
		\draw [style=wire] (57) to (61);
		\draw [style=wire, in=150, out=-135, looseness=1.50] (58) to (57);
		\draw [style=wire] (59) to (56);
		\draw [style=wire, in=30, out=-90, looseness=1.25] (60) to (57);
		\draw [style=wire, in=150, out=-30] (58) to (59);
		\draw [style=wire] (52) to (58);
		\draw [style=wire, in=45, out=-90] (55) to (59);
		\draw [style=wire, in=-90, out=30] (104) to (105);
		\draw [style=wire, in=150, out=-150, looseness=1.50] (106) to (104);
		\draw [style=wire] (103) to (106);
		\draw [style=wire] (104) to (107);
		\draw [style=wire, bend left, looseness=1.25] (106) to (108);
		\draw [style=wire, in=-90, out=45] (110) to (109);
		\draw [style=wire, in=91, out=-135, looseness=0.75] (114) to (112);
		\draw [style=wire, in=150, out=-30] (114) to (110);
		\draw [style=wire] (113) to (114);
		\draw [style=wire] (110) to (111);
	\end{pgfonlayer}
\end{tikzpicture}
}%
\]
\noindent
{\bf [d.2]}  $\Rightarrow$ {\bf [d.$\nabla$]}:  By the properties of $\mathsf{i}_j$ and $\mathsf{p}_k$, Lemma \ref{seelyish}, the additive bialgebra modality identities, and the additive structure, we have that: 
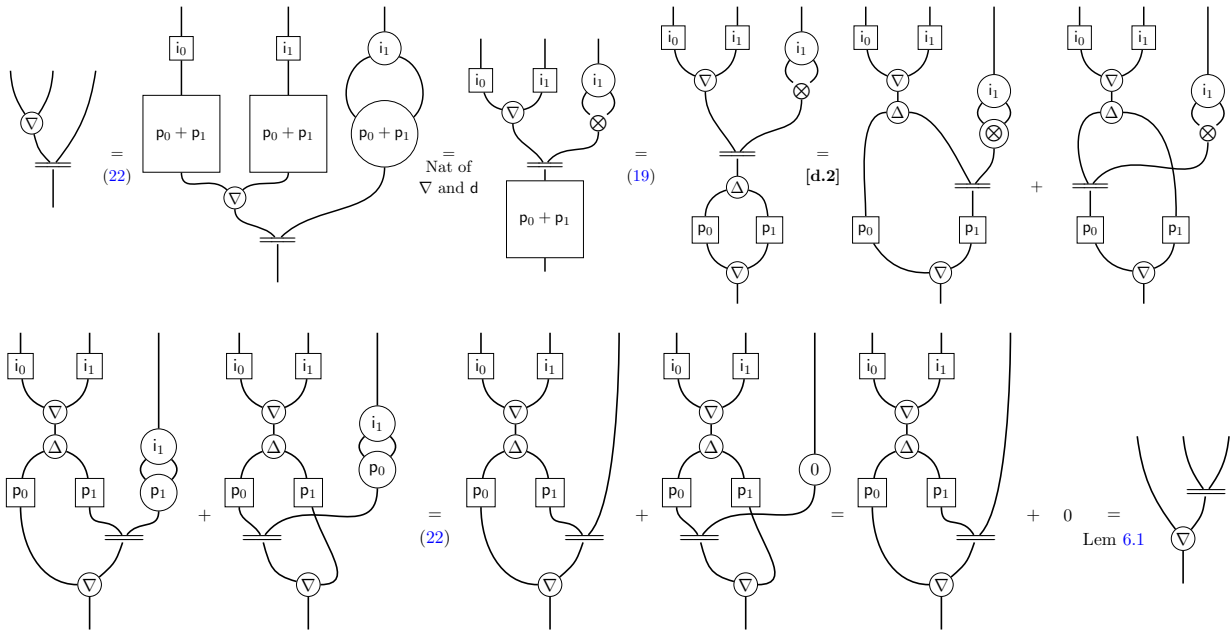
\[\resizebox{!}{4cm}{%
\begin{tikzpicture}
	\begin{pgfonlayer}{nodelayer}
		\node [style=port] (0) at (0.5, 3) {};
		\node [style=port] (1) at (2.5, 3) {};
		\node [style=codifferential] (2) at (1.5, 0.75) {{\bf =\!=\!=}};
		\node [style=port] (3) at (1.5, -0.25) {};
		\node [style=port] (4) at (1.5, 3) {};
		\node [style=duplicate] (5) at (1, 1.75) {$\nabla$};
		\node [style=port] (6) at (6.75, -2) {};
		\node [style=codifferential] (7) at (6.75, -1) {{\bf =\!=\!=}};
		\node [style=function2] (8) at (4.5, 3.5) {$\mathsf{i}_0$};
		\node [style=port] (9) at (7, 4.5) {};
		\node [style=port] (10) at (4.5, 4.5) {};
		\node [style=duplicate] (11) at (5.75, 0) {$\nabla$};
		\node [style=function2] (12) at (7, 3.5) {$\mathsf{i}_1$};
		\node [style=component] (13) at (9.25, 3.5) {$\mathsf{i}_1$};
		\node [style=port] (14) at (9.25, 4.5) {};
		\node [style=function2] (15) at (4.5, 1.5) {$\mathsf{p}_0+\mathsf{p}_1$};
		\node [style=function2] (16) at (7, 1.5) {$\mathsf{p}_0+\mathsf{p}_1$};
		\node [style=component] (17) at (9.25, 1.5) {$\mathsf{p}_0+\mathsf{p}_1$};
		\node [style=port] (18) at (3, 1) {$=$};
		\node [style=port] (19) at (3, 0.5) {(\ref{pi1})};
		\node [style=port] (20) at (13, -1.75) {};
		\node [style=codifferential] (21) at (13, 0.75) {{\bf =\!=\!=}};
		\node [style=function2] (22) at (11.5, 2.75) {$\mathsf{i}_0$};
		\node [style=port] (23) at (13, 3.75) {};
		\node [style=port] (24) at (11.5, 3.75) {};
		\node [style=duplicate] (25) at (12.25, 2) {$\nabla$};
		\node [style=function2] (26) at (13, 2.75) {$\mathsf{i}_1$};
		\node [style=port] (27) at (14.25, 1.75) {$\bigotimes$};
		\node [style=component] (28) at (14.25, 2.75) {$\mathsf{i}_1$};
		\node [style=port] (29) at (14.25, 3.75) {};
		\node [style=function2] (30) at (13, -0.5) {$\mathsf{p}_0+\mathsf{p}_1$};
		\node [style=port] (31) at (10.75, 1) {$=$};
		\node [style=port] (32) at (10.75, 0.75) {Nat of };
		\node [style=port] (33) at (10.75, 0.25) {$\nabla$ and $\mathsf{d}$};
		\node [style=port] (34) at (17.5, -2.5) {};
		\node [style=codifferential] (35) at (17.5, 1) {{\bf =\!=\!=}};
		\node [style=function2] (36) at (16, 3.75) {$\mathsf{i}_0$};
		\node [style=port] (37) at (17.5, 4.5) {};
		\node [style=port] (38) at (16, 4.5) {};
		\node [style=duplicate] (39) at (16.75, 2.75) {$\nabla$};
		\node [style=function2] (40) at (17.5, 3.75) {$\mathsf{i}_1$};
		\node [style=port] (41) at (19, 2.5) {$\bigotimes$};
		\node [style=component] (42) at (19, 3.5) {$\mathsf{i}_1$};
		\node [style=port] (43) at (19, 4.5) {};
		\node [style=function2] (44) at (18.25, -0.75) {$\mathsf{p}_1$};
		\node [style=function2] (45) at (16.75, -0.75) {$\mathsf{p}_0$};
		\node [style=duplicate] (46) at (17.5, -1.75) {$\nabla$};
		\node [style=duplicate] (47) at (17.5, 0.25) {$\Delta$};
		\node [style=port] (48) at (15.25, 1) {$=$};
		\node [style=port] (49) at (15.25, 0.5) {(\ref{addbialgeq})};
		\node [style=port] (50) at (22.25, -2.5) {};
		\node [style=codifferential] (51) at (23, 0.25) {{\bf =\!=\!=}};
		\node [style=function2] (52) at (20.5, 3.75) {$\mathsf{i}_0$};
		\node [style=port] (53) at (22, 4.5) {};
		\node [style=port] (54) at (20.5, 4.5) {};
		\node [style=duplicate] (55) at (21.25, 2.75) {$\nabla$};
		\node [style=duplicate] (56) at (21.25, 2) {$\Delta$};
		\node [style=function2] (57) at (22, 3.75) {$\mathsf{i}_1$};
		\node [style=duplicate] (58) at (23.5, 1.5) {$\bigotimes$};
		\node [style=component] (59) at (23.5, 2.5) {$\mathsf{i}_1$};
		\node [style=port] (60) at (23.5, 4.5) {};
		\node [style=function2] (61) at (23, -0.75) {$\mathsf{p}_1$};
		\node [style=function2] (62) at (20.5, -0.75) {$\mathsf{p}_0$};
		\node [style=duplicate] (63) at (22.25, -1.75) {$\nabla$};
		\node [style=port] (64) at (24.5, 0.25) {$+$};
		\node [style=port] (65) at (28.5, 4.5) {};
		\node [style=codifferential] (66) at (25.75, 0.25) {{\bf =\!=\!=}};
		\node [style=port] (67) at (25.5, 4.5) {};
		\node [style=port] (68) at (27, -2.5) {};
		\node [style=port] (69) at (28.5, 1.5) {$\bigotimes$};
		\node [style=function2] (70) at (25.75, -0.75) {$\mathsf{p}_0$};
		\node [style=port] (71) at (27, 4.5) {};
		\node [style=function2] (72) at (27, 3.75) {$\mathsf{i}_1$};
		\node [style=function2] (73) at (27.75, -0.75) {$\mathsf{p}_1$};
		\node [style=duplicate] (74) at (26.25, 2.75) {$\nabla$};
		\node [style=duplicate] (75) at (27, -1.75) {$\nabla$};
		\node [style=function2] (76) at (25.5, 3.75) {$\mathsf{i}_0$};
		\node [style=duplicate] (77) at (26.25, 2) {$\Delta$};
		\node [style=component] (78) at (28.5, 2.5) {$\mathsf{i}_1$};
		\node [style=port] (79) at (19.5, 1) {$=$};
		\node [style=port] (80) at (19.5, 0.5) {\textbf{[d.2]}};
	\end{pgfonlayer}
	\begin{pgfonlayer}{edgelayer}
		\draw [style=wire, bend left=15, looseness=1.25] (2) to (5);
		\draw [style=wire, bend right=15] (2) to (1);
		\draw [style=wire] (3) to (2);
		\draw [style=wire, bend left=15] (5) to (0);
		\draw [style=wire, bend right=15] (5) to (4);
		\draw [style=wire] (10) to (8);
		\draw [style=wire] (9) to (12);
		\draw [style=wire] (14) to (13);
		\draw [style=wire, in=135, out=-90, looseness=1.25] (15) to (11);
		\draw [style=wire, in=45, out=-90] (16) to (11);
		\draw [style=wire, bend right=60] (13) to (17);
		\draw [style=wire, in=45, out=-90] (17) to (7);
		\draw [style=wire] (12) to (16);
		\draw [style=wire] (8) to (15);
		\draw [style=wire, in=129, out=-90] (11) to (7);
		\draw [style=wire] (7) to (6);
		\draw [style=wire, bend left=60] (13) to (17);
		\draw [style=wire, bend left] (25) to (22);
		\draw [style=wire, bend right] (25) to (26);
		\draw [style=wire] (24) to (22);
		\draw [style=wire] (23) to (26);
		\draw [style=wire, in=53, out=-90, looseness=1.25] (27) to (21);
		\draw [style=wire, bend right=45, looseness=1.25] (28) to (27);
		\draw [style=wire, bend left=45, looseness=1.25] (28) to (27);
		\draw [style=wire] (29) to (28);
		\draw [style=wire, in=120, out=-90] (25) to (21);
		\draw [style=wire] (21) to (30);
		\draw [style=wire] (30) to (20);
		\draw [style=wire, bend left] (39) to (36);
		\draw [style=wire, bend right] (39) to (40);
		\draw [style=wire] (38) to (36);
		\draw [style=wire] (37) to (40);
		\draw [style=wire, in=53, out=-90, looseness=1.25] (41) to (35);
		\draw [style=wire, bend right=45, looseness=1.25] (42) to (41);
		\draw [style=wire, bend left=45, looseness=1.25] (42) to (41);
		\draw [style=wire] (43) to (42);
		\draw [style=wire, bend left] (46) to (45);
		\draw [style=wire, bend right] (46) to (44);
		\draw [style=wire] (46) to (34);
		\draw [style=wire] (35) to (47);
		\draw [style=wire, in=90, out=-165] (47) to (45);
		\draw [style=wire, in=97, out=-15, looseness=1.25] (47) to (44);
		\draw [style=wire, in=120, out=-90] (39) to (35);
		\draw [style=wire, bend left] (55) to (52);
		\draw [style=wire, bend right] (55) to (57);
		\draw [style=wire] (54) to (52);
		\draw [style=wire] (53) to (57);
		\draw [style=wire] (55) to (56);
		\draw [style=wire, in=53, out=-90, looseness=1.25] (58) to (51);
		\draw [style=wire, bend right=45, looseness=1.25] (59) to (58);
		\draw [style=wire, bend left=45, looseness=1.25] (59) to (58);
		\draw [style=wire] (60) to (59);
		\draw [style=wire, bend left] (63) to (62);
		\draw [style=wire, bend right] (63) to (61);
		\draw [style=wire, in=120, out=0, looseness=0.75] (56) to (51);
		\draw [style=wire] (63) to (50);
		\draw [style=wire, in=90, out=180, looseness=0.75] (56) to (62);
		\draw [style=wire] (51) to (61);
		\draw [style=wire, bend left] (74) to (76);
		\draw [style=wire, bend right] (74) to (72);
		\draw [style=wire] (67) to (76);
		\draw [style=wire] (71) to (72);
		\draw [style=wire] (74) to (77);
		\draw [style=wire, in=53, out=-90, looseness=1.25] (69) to (66);
		\draw [style=wire, bend right=45, looseness=1.25] (78) to (69);
		\draw [style=wire, bend left=45, looseness=1.25] (78) to (69);
		\draw [style=wire] (65) to (78);
		\draw [style=wire, bend left] (75) to (70);
		\draw [style=wire, bend right] (75) to (73);
		\draw [style=wire] (66) to (70);
		\draw [style=wire, in=89, out=0] (77) to (73);
		\draw [style=wire, in=120, out=180] (77) to (66);
		\draw [style=wire] (75) to (68);
	\end{pgfonlayer}
\end{tikzpicture}
}%
\]
\[\resizebox{!}{4cm}{%
\begin{tikzpicture}
	\begin{pgfonlayer}{nodelayer}
		\node [style=port] (0) at (2, -3) {};
		\node [style=duplicate] (1) at (2, -2) {$\nabla$};
		\node [style=codifferential] (2) at (2.75, -1) {{\bf =\!=\!=}};
		\node [style=function2] (3) at (0.5, 2.75) {$\mathsf{i}_0$};
		\node [style=function2] (4) at (0.5, 0) {$\mathsf{p}_0$};
		\node [style=function2] (5) at (2, 0) {$\mathsf{p}_1$};
		\node [style=port] (6) at (2, 3.5) {};
		\node [style=port] (7) at (0.5, 3.5) {};
		\node [style=duplicate] (8) at (1.25, 1.75) {$\nabla$};
		\node [style=duplicate] (9) at (1.25, 1) {$\Delta$};
		\node [style=function2] (10) at (2, 2.75) {$\mathsf{i}_1$};
		\node [style=component] (11) at (3.5, 0) {$\mathsf{p}_1$};
		\node [style=component] (12) at (3.5, 1) {$\mathsf{i}_1$};
		\node [style=port] (13) at (3.5, 3.5) {};
		\node [style=port] (14) at (6.75, -3) {};
		\node [style=duplicate] (15) at (6.75, -2) {$\nabla$};
		\node [style=codifferential] (16) at (5.75, -1) {{\bf =\!=\!=}};
		\node [style=function2] (17) at (5.25, 2.75) {$\mathsf{i}_0$};
		\node [style=function2] (18) at (5.25, 0) {$\mathsf{p}_0$};
		\node [style=function2] (19) at (6.75, 0) {$\mathsf{p}_1$};
		\node [style=port] (20) at (6.75, 3.5) {};
		\node [style=port] (21) at (5.25, 3.5) {};
		\node [style=duplicate] (22) at (6, 1.75) {$\nabla$};
		\node [style=duplicate] (23) at (6, 1) {$\Delta$};
		\node [style=function2] (24) at (6.75, 2.75) {$\mathsf{i}_1$};
		\node [style=component] (25) at (8.25, 0.5) {$\mathsf{p}_0$};
		\node [style=component] (26) at (8.25, 1.5) {$\mathsf{i}_1$};
		\node [style=port] (27) at (8.25, 3.5) {};
		\node [style=port] (28) at (4.5, -0.5) {$+$};
		\node [style=port] (29) at (9.5, -0.5) {};
		\node [style=port] (30) at (9.5, -0.5) {$=$};
		\node [style=port] (31) at (9.5, -1) {(\ref{pi1})};
		\node [style=port] (32) at (12, -3) {};
		\node [style=duplicate] (33) at (12, -2) {$\nabla$};
		\node [style=codifferential] (34) at (12.75, -1) {{\bf =\!=\!=}};
		\node [style=function2] (35) at (10.5, 2.75) {$\mathsf{i}_0$};
		\node [style=function2] (36) at (10.5, 0) {$\mathsf{p}_0$};
		\node [style=function2] (37) at (12, 0) {$\mathsf{p}_1$};
		\node [style=port] (38) at (12, 3.5) {};
		\node [style=port] (39) at (10.5, 3.5) {};
		\node [style=duplicate] (40) at (11.25, 1.75) {$\nabla$};
		\node [style=duplicate] (41) at (11.25, 1) {$\Delta$};
		\node [style=function2] (42) at (12, 2.75) {$\mathsf{i}_1$};
		\node [style=port] (45) at (13.5, 3.5) {};
		\node [style=port] (46) at (16.25, -3) {};
		\node [style=duplicate] (47) at (16.25, -2) {$\nabla$};
		\node [style=codifferential] (48) at (15.25, -1) {{\bf =\!=\!=}};
		\node [style=function2] (49) at (14.75, 2.75) {$\mathsf{i}_0$};
		\node [style=function2] (50) at (14.75, 0) {$\mathsf{p}_0$};
		\node [style=function2] (51) at (16.25, 0) {$\mathsf{p}_1$};
		\node [style=port] (52) at (16.25, 3.5) {};
		\node [style=port] (53) at (14.75, 3.5) {};
		\node [style=duplicate] (54) at (15.5, 1.75) {$\nabla$};
		\node [style=duplicate] (55) at (15.5, 1) {$\Delta$};
		\node [style=function2] (56) at (16.25, 2.75) {$\mathsf{i}_1$};
		\node [style=component] (57) at (17.75, 0.5) {$0$};
		\node [style=port] (59) at (17.75, 3.5) {};
		\node [style=port] (60) at (14, -0.5) {$+$};
		\node [style=port] (61) at (20.5, -3) {};
		\node [style=duplicate] (62) at (20.5, -2) {$\nabla$};
		\node [style=codifferential] (63) at (21.25, -1) {{\bf =\!=\!=}};
		\node [style=function2] (64) at (19, 2.75) {$\mathsf{i}_0$};
		\node [style=function2] (65) at (19, 0) {$\mathsf{p}_0$};
		\node [style=function2] (66) at (20.5, 0) {$\mathsf{p}_1$};
		\node [style=port] (67) at (20.5, 3.5) {};
		\node [style=port] (68) at (19, 3.5) {};
		\node [style=duplicate] (69) at (19.75, 1.75) {$\nabla$};
		\node [style=duplicate] (70) at (19.75, 1) {$\Delta$};
		\node [style=function2] (71) at (20.5, 2.75) {$\mathsf{i}_1$};
		\node [style=port] (72) at (22, 3.5) {};
		\node [style=port] (73) at (18.25, -0.5) {};
		\node [style=port] (74) at (18.25, -0.5) {$=$};
		\node [style=port] (75) at (22.5, -0.5) {$+$};
		\node [style=port] (76) at (23.25, -0.5) {$0$};
		\node [style=port] (77) at (25.75, -2) {};
		\node [style=port] (78) at (25.75, 1.25) {};
		\node [style=duplicate] (79) at (25.75, -1) {$\nabla$};
		\node [style=port] (80) at (26.75, 1.25) {};
		\node [style=codifferential] (81) at (26.25, 0) {{\bf =\!=\!=}};
		\node [style=port] (82) at (24.75, 1.25) {};
		\node [style=port] (83) at (24.25, -0.5) {};
		\node [style=port] (84) at (24.25, -0.5) {$=$};
		\node [style=port] (85) at (24.25, -1) {Lem \ref{seelyish}};
	\end{pgfonlayer}
	\begin{pgfonlayer}{edgelayer}
		\draw [style=wire, bend right=15, looseness=1.25] (1) to (2);
		\draw [style=wire] (0) to (1);
		\draw [style=wire, bend left] (8) to (3);
		\draw [style=wire, bend right] (8) to (10);
		\draw [style=wire, bend right] (9) to (4);
		\draw [style=wire, bend left] (9) to (5);
		\draw [style=wire] (7) to (3);
		\draw [style=wire] (6) to (10);
		\draw [style=wire] (8) to (9);
		\draw [style=wire, in=180, out=-90] (4) to (1);
		\draw [style=wire, in=120, out=-90, looseness=1.50] (5) to (2);
		\draw [style=wire, in=53, out=-90, looseness=1.25] (11) to (2);
		\draw [style=wire, bend right=45, looseness=1.25] (12) to (11);
		\draw [style=wire, bend left=45, looseness=1.25] (12) to (11);
		\draw [style=wire] (13) to (12);
		\draw [style=wire, in=-75, out=165, looseness=1.25] (15) to (16);
		\draw [style=wire] (14) to (15);
		\draw [style=wire, bend left] (22) to (17);
		\draw [style=wire, bend right] (22) to (24);
		\draw [style=wire, bend right] (23) to (18);
		\draw [style=wire, bend left] (23) to (19);
		\draw [style=wire] (21) to (17);
		\draw [style=wire] (20) to (24);
		\draw [style=wire] (22) to (23);
		\draw [style=wire, in=53, out=-90, looseness=1.25] (25) to (16);
		\draw [style=wire, bend right=45, looseness=1.25] (26) to (25);
		\draw [style=wire, bend left=45, looseness=1.25] (26) to (25);
		\draw [style=wire] (27) to (26);
		\draw [style=wire, in=117, out=-90] (18) to (16);
		\draw [style=wire, in=0, out=-75, looseness=1.25] (19) to (15);
		\draw [style=wire, bend right=15, looseness=1.25] (33) to (34);
		\draw [style=wire] (32) to (33);
		\draw [style=wire, bend left] (40) to (35);
		\draw [style=wire, bend right] (40) to (42);
		\draw [style=wire, bend right] (41) to (36);
		\draw [style=wire, bend left] (41) to (37);
		\draw [style=wire] (39) to (35);
		\draw [style=wire] (38) to (42);
		\draw [style=wire] (40) to (41);
		\draw [style=wire, in=180, out=-90] (36) to (33);
		\draw [style=wire, in=120, out=-90, looseness=1.50] (37) to (34);
		\draw [style=wire, in=-75, out=165, looseness=1.25] (47) to (48);
		\draw [style=wire] (46) to (47);
		\draw [style=wire, bend left] (54) to (49);
		\draw [style=wire, bend right] (54) to (56);
		\draw [style=wire, bend right] (55) to (50);
		\draw [style=wire, bend left] (55) to (51);
		\draw [style=wire] (53) to (49);
		\draw [style=wire] (52) to (56);
		\draw [style=wire] (54) to (55);
		\draw [style=wire, in=53, out=-90, looseness=1.25] (57) to (48);
		\draw [style=wire, in=117, out=-90] (50) to (48);
		\draw [style=wire, in=0, out=-75, looseness=1.25] (51) to (47);
		\draw [style=wire] (59) to (57);
		\draw [style=wire, in=45, out=-90, looseness=0.50] (45) to (34);
		\draw [style=wire, bend right=15, looseness=1.25] (62) to (63);
		\draw [style=wire] (61) to (62);
		\draw [style=wire, bend left] (69) to (64);
		\draw [style=wire, bend right] (69) to (71);
		\draw [style=wire, bend right] (70) to (65);
		\draw [style=wire, bend left] (70) to (66);
		\draw [style=wire] (68) to (64);
		\draw [style=wire] (67) to (71);
		\draw [style=wire] (69) to (70);
		\draw [style=wire, in=180, out=-90] (65) to (62);
		\draw [style=wire, in=120, out=-90, looseness=1.50] (66) to (63);
		\draw [style=wire, in=45, out=-90, looseness=0.50] (72) to (63);
		\draw [style=wire, bend right=15, looseness=1.25] (79) to (81);
		\draw [style=wire, bend left=15] (79) to (82);
		\draw [style=wire] (77) to (79);
		\draw [style=wire, bend right=15] (81) to (80);
		\draw [style=wire, bend left=15] (81) to (78);
	\end{pgfonlayer}
\end{tikzpicture}
}%
\]
\end{proof} 

\begin{corollary}For an additive bialgebra modality $(\oc, \delta, \varepsilon, \Delta, e, \nabla, u)$, the following are equivalent for a natural transformation $\mathsf{d}: \oc A \otimes A \to \oc A$: 
\begin{enumerate}[{\em (i)}]
\item $\mathsf{d}$ is a deriving transformation; 
\item $\mathsf{d}$ satisfies the product rule {\bf [d.2]}, the linear rule {\bf [d.3]}, and the chain rule {\bf [d.4]};
\item $\mathsf{d}$ satisfies the linear rule {\bf [d.3]}, the chain rule {\bf [d.4]}, and the $\nabla$-rule {\bf [d.$\nabla$]}. 
\end{enumerate}
\end{corollary} 

Therefore, we obtain the following theorem: 

\begin{theorem}\label{mainthm1} For an additive bialgebra modality, every deriving transformation satisfies the $\nabla$-rule {\bf [d.$\nabla$]} and thus is induced equivalently by a codereliction. 
\end{theorem} 

We turn our attention to the relation between the monoidal structure and the differential structure, that is, we explore deriving transformations of additive linear categories. The compatibility between a deriving transformation and the monoidal structure is described by the monoidal rule  \cite{fiore2007differential}-- this is the strength rule 
which was the subject of Fiore's addendum: 
\begin{description}
\item[{\bf [d.m]}] Monoidal Rule: 
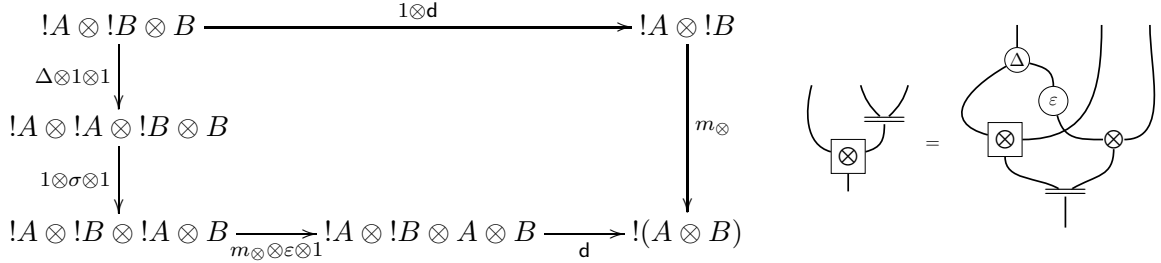
\begin{align*}
\begin{array}[c]{c}
 \xymatrixcolsep{2.5pc}\xymatrix{\oc A \otimes \oc B \otimes B \ar[d]_-{\Delta \otimes 1 \otimes 1} \ar[rr]^-{1 \otimes \mathsf{d}} && \oc A \otimes \oc B \ar[dd]^-{m_\otimes} \\
   \oc A \otimes \oc A \otimes \oc B \otimes B \ar[d]_-{1 \otimes \sigma \otimes 1} \\
   \oc A \otimes \oc B \otimes \oc A \otimes B \ar[r]_-{m_\otimes \otimes \varepsilon  \otimes 1} & \oc A \otimes \oc B \otimes A \otimes B \ar[r]_-{\mathsf{d}} & \oc(A \otimes B) 
  } 
   \end{array} &&
   \begin{array}[c]{c}
\resizebox{!}{2.75cm}{%
\begin{tikzpicture}
	\begin{pgfonlayer}{nodelayer}
		\node [style=port] (0) at (1.25, -0.75) {};
		\node [style=port] (1) at (0.5, 1.5) {};
		\node [style=function2] (2) at (1.25, 0) {$\bigotimes$};
		\node [style=codifferential] (3) at (2, 0.75) {{\bf =\!=\!=}};
		\node [style=port] (4) at (2.5, 1.5) {};
		\node [style=port] (5) at (1.5, 1.5) {};
		\node [style=port] (6) at (3, 0.25) {$=$};
		\node [style=function2] (7) at (4.5, 0.35) {$\bigotimes$};
		\node [style=duplicate] (8) at (4.75, 2) {$\Delta$};
		\node [style=port] (9) at (4.75, 2.75) {};
		\node [style=port] (10) at (6.5, 2.75) {};
		\node [style=port] (11) at (7.5, 2.75) {};
		\node [style=differential] (12) at (6.75, 0.35) {$\bigotimes$};
		\node [style=codifferential] (13) at (5.75, -0.75) {{\bf =\!=\!=}};
		\node [style=port] (14) at (5.75, -1.5) {};
		\node [style=component] (15) at (5.5, 1.15) {$\varepsilon$};
	\end{pgfonlayer}
	\begin{pgfonlayer}{edgelayer}
		\draw [style=wire, in=-105, out=165] (2) to (1);
		\draw [style=wire] (2) to (0);
		\draw [style=wire, bend right=15] (3) to (4);
		\draw [style=wire, bend left=15] (3) to (5);
		\draw [style=wire, in=15, out=-90, looseness=1.25] (3) to (2);
		\draw [style=wire] (13) to (14);
		\draw [style=wire, in=165, out=-150, looseness=1.75] (8) to (7);
		\draw [style=wire, in=150, out=-90] (7) to (13);
		\draw [style=wire, in=30, out=-90, looseness=1.25] (12) to (13);
		\draw [style=wire] (9) to (8);
		\draw [style=wire, in=0, out=-90, looseness=1.25] (10) to (7);
		\draw [style=wire, in=0, out=-90] (11) to (12);
		\draw [style=wire, in=90, out=-15, looseness=1.25] (8) to (15);
		\draw [style=wire, in=180, out=-75, looseness=1.25] (15) to (12);
	\end{pgfonlayer}
\end{tikzpicture}
}%
   \end{array}
\end{align*}
\end{description}

Fiore's creation operator \cite{fiore2007differential} was defined to satisfy the linear rule {\bf [d.3]}, the chain rule {\bf [d.4]}, the $\nabla$-rule {\bf [d.$\nabla$]}, and the monoidal rule {\bf [d.m]} -- in his addendum, he pointed out the latter was redundant. It turns out that when a natural transformation satisfies both the linear rule {\bf [d.3]} and the chain rule {\bf [d.4]}, then the monoidal rule is equivalent to both the $\nabla$-rule and the Leibniz rule:  

\begin{proposition} For the induced additive bialgebra modality of an additive linear category, the following are equivalent for a natural transformation $\mathsf{d}: \oc A \otimes A \to A$ which satisfies the linear rule {\bf [d.3]} and the chain rule {\bf [d.4]}:
\begin{enumerate}[{\em (i)}]
\item $\mathsf{d}$ satisfies the Leibniz rule {\bf [d.2]};
\item $\mathsf{d}$ satisfies the $\nabla$-rule {\bf [d.$\nabla$]};
\item $\mathsf{d}$ satisfies the monoidal rule {\bf [d.m]}
\end{enumerate}
\end{proposition}
\begin{proof} Since this is an extension of Proposition \ref{nablaleibniz}, it suffices to show that the $\nabla$-rule {\bf [d.$\nabla$]} and the monoidal rule {\bf [d.m]} are equivalent. \\
{\bf [d.$\nabla$]} $\Rightarrow$ {\bf [d.m]}: It is easy to see that since $\mathsf{d}$ satisfies the linear rule {\bf [d.3]} and the chain rule {\bf [d.4]}, $(u \otimes 1)\mathsf{d}: A \to \oc A$ satisfies the codereliction linear rule {\bf [dC.3]} and chain rule {\bf [dC.4]}, and therefore by Corollary \ref{cor1} is a codereliction and which by Proposition \ref{etamonoidal} satisfies the codereliction monoidal rule {\bf [dC.m]}. Therefore, by Lemma \ref{diden} (since $\mathsf{d}$ satisfies \textbf{[d.$\nabla$]}) and one of the identities of Proposition \ref{monoidalnabla}, we have: 
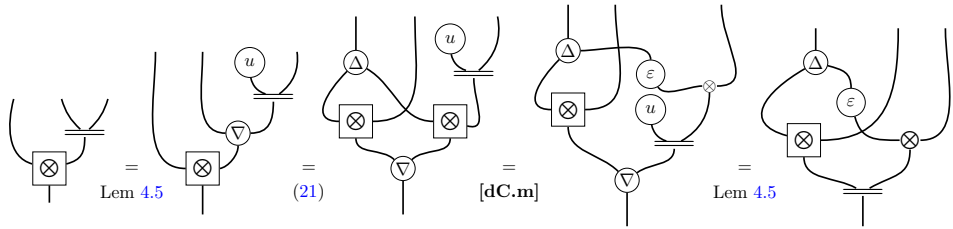
\[\resizebox{!}{3cm}{%
\begin{tikzpicture}
	\begin{pgfonlayer}{nodelayer}
		\node [style=port] (0) at (-1.75, -0.75) {};
		\node [style=port] (1) at (-2.5, 1.5) {};
		\node [style=function2] (2) at (-1.75, 0) {$\bigotimes$};
		\node [style=codifferential] (3) at (-1, 0.75) {{\bf =\!=\!=}};
		\node [style=port] (4) at (-0.5, 1.5) {};
		\node [style=port] (5) at (-1.5, 1.5) {};
		\node [style=port] (6) at (0.5, 2.5) {};
		\node [style=differential] (7) at (3, 1.5) {{\bf =\!=\!=}};
		\node [style=duplicate] (8) at (2.25, 0.75) {$\nabla$};
		\node [style=port] (9) at (1.5, 2.5) {};
		\node [style=function2] (10) at (1.5, 0) {$\bigotimes$};
		\node [style=port] (11) at (1.5, -1) {};
		\node [style=component] (12) at (2.5, 2.25) {$u$};
		\node [style=port] (13) at (3.5, 2.5) {};
		\node [style=port] (14) at (0, 0) {$=$};
		\node [style=port] (15) at (0, -0.5) {Lem \ref{diden}};
		\node [style=port] (16) at (4.75, 3.25) {};
		\node [style=function2] (17) at (4.75, 1) {$\bigotimes$};
		\node [style=duplicate] (18) at (4.75, 2.25) {$\Delta$};
		\node [style=function2] (19) at (6.75, 1) {$\bigotimes$};
		\node [style=port] (20) at (6, 3.25) {};
		\node [style=duplicate] (21) at (5.75, 0) {$\nabla$};
		\node [style=port] (22) at (5.75, -1) {};
		\node [style=component] (23) at (6.75, 2.75) {$u$};
		\node [style=differential] (24) at (7.25, 2) {{\bf =\!=\!=}};
		\node [style=port] (25) at (7.75, 3.25) {};
		\node [style=port] (26) at (3.75, 0) {$=$};
		\node [style=port] (27) at (3.75, -0.5) {(\ref{nablam2})};
		\node [style=port] (28) at (9.25, 3.5) {};
		\node [style=function2] (29) at (9.25, 1.25) {$\bigotimes$};
		\node [style=duplicate] (30) at (9.25, 2.5) {$\Delta$};
		\node [style=port] (31) at (10.25, 3.5) {};
		\node [style=duplicate] (32) at (10.5, -0.25) {$\nabla$};
		\node [style=port] (33) at (10.5, -1.25) {};
		\node [style=component] (34) at (11, 1.25) {$u$};
		\node [style=differential] (35) at (11.5, 0.5) {{\bf =\!=\!=}};
		\node [style=port] (36) at (12.25, 1.75) {$\otimes$};
		\node [style=component] (37) at (11, 2) {$\varepsilon$};
		\node [style=port] (38) at (12.5, 3.5) {};
		\node [style=port] (39) at (8, 0) {$=$};
		\node [style=port] (40) at (8, -0.5) {\textbf{[dC.m]}};
		\node [style=function2] (41) at (14.25, 0.6) {$\bigotimes$};
		\node [style=duplicate] (42) at (14.5, 2.25) {$\Delta$};
		\node [style=port] (43) at (14.5, 3) {};
		\node [style=port] (44) at (16.25, 3) {};
		\node [style=port] (45) at (17.25, 3) {};
		\node [style=differential] (46) at (16.5, 0.6) {$\bigotimes$};
		\node [style=codifferential] (47) at (15.5, -0.5) {{\bf =\!=\!=}};
		\node [style=port] (48) at (15.5, -1.25) {};
		\node [style=component] (49) at (15.25, 1.4) {$\varepsilon$};
		\node [style=port] (50) at (13, 0) {$=$};
		\node [style=port] (51) at (13, -0.5) {Lem \ref{diden}};
	\end{pgfonlayer}
	\begin{pgfonlayer}{edgelayer}
		\draw [style=wire, in=-105, out=165] (2) to (1);
		\draw [style=wire] (2) to (0);
		\draw [style=wire, bend right=15] (3) to (4);
		\draw [style=wire, bend left=15] (3) to (5);
		\draw [style=wire, in=15, out=-90, looseness=1.25] (3) to (2);
		\draw [style=wire, in=180, out=-90] (6) to (10);
		\draw [style=wire, in=180, out=-90] (9) to (8);
		\draw [style=wire] (10) to (11);
		\draw [style=wire, bend right=45, looseness=0.75] (10) to (8);
		\draw [style=wire, bend right] (7) to (13);
		\draw [style=wire, bend left] (7) to (12);
		\draw [style=wire, in=0, out=-90, looseness=1.25] (7) to (8);
		\draw [style=wire, in=180, out=-135, looseness=1.50] (18) to (17);
		\draw [style=wire, in=0, out=-90, looseness=1.25] (20) to (17);
		\draw [style=wire, in=180, out=-30] (18) to (19);
		\draw [style=wire] (16) to (18);
		\draw [style=wire, in=135, out=-90] (17) to (21);
		\draw [style=wire, in=45, out=-90] (19) to (21);
		\draw [style=wire] (21) to (22);
		\draw [style=wire, bend right] (24) to (25);
		\draw [style=wire, bend left] (24) to (23);
		\draw [style=wire, in=-15, out=-90, looseness=1.25] (24) to (19);
		\draw [style=wire, in=180, out=-135, looseness=1.50] (30) to (29);
		\draw [style=wire, in=0, out=-90, looseness=1.25] (31) to (29);
		\draw [style=wire] (28) to (30);
		\draw [style=wire, in=135, out=-90] (29) to (32);
		\draw [style=wire] (32) to (33);
		\draw [style=wire, bend right] (35) to (36);
		\draw [style=wire, bend left] (35) to (34);
		\draw [style=wire, in=37, out=-90] (35) to (32);
		\draw [style=wire, in=90, out=0, looseness=0.75] (30) to (37);
		\draw [style=wire, in=180, out=-60, looseness=0.75] (37) to (36);
		\draw [style=wire, in=0, out=-90, looseness=1.25] (38) to (36);
		\draw [style=wire] (47) to (48);
		\draw [style=wire, in=165, out=-150, looseness=1.75] (42) to (41);
		\draw [style=wire, in=150, out=-90] (41) to (47);
		\draw [style=wire, in=30, out=-90, looseness=1.25] (46) to (47);
		\draw [style=wire] (43) to (42);
		\draw [style=wire, in=0, out=-90, looseness=1.25] (44) to (41);
		\draw [style=wire, in=0, out=-90] (45) to (46);
		\draw [style=wire, in=90, out=-15, looseness=1.25] (42) to (49);
		\draw [style=wire, in=180, out=-75, looseness=1.25] (49) to (46);
	\end{pgfonlayer}
\end{tikzpicture}
}%
\]
{\bf [d.m]} $\Rightarrow$ {\bf [d.$\nabla$]}: Using the coalgebra modality identities, that $\nabla$ is a $\oc$-coalgebra morphism, the monoidal rule \textbf{[d.m]}, the chain rule \textbf{[d.4]}, the bialgebra modality compatibility between $\nabla$ and $\varepsilon$, the linear rule \textbf{[d.3]}, and the constant rule \textbf{[d.1]}, we have that: 
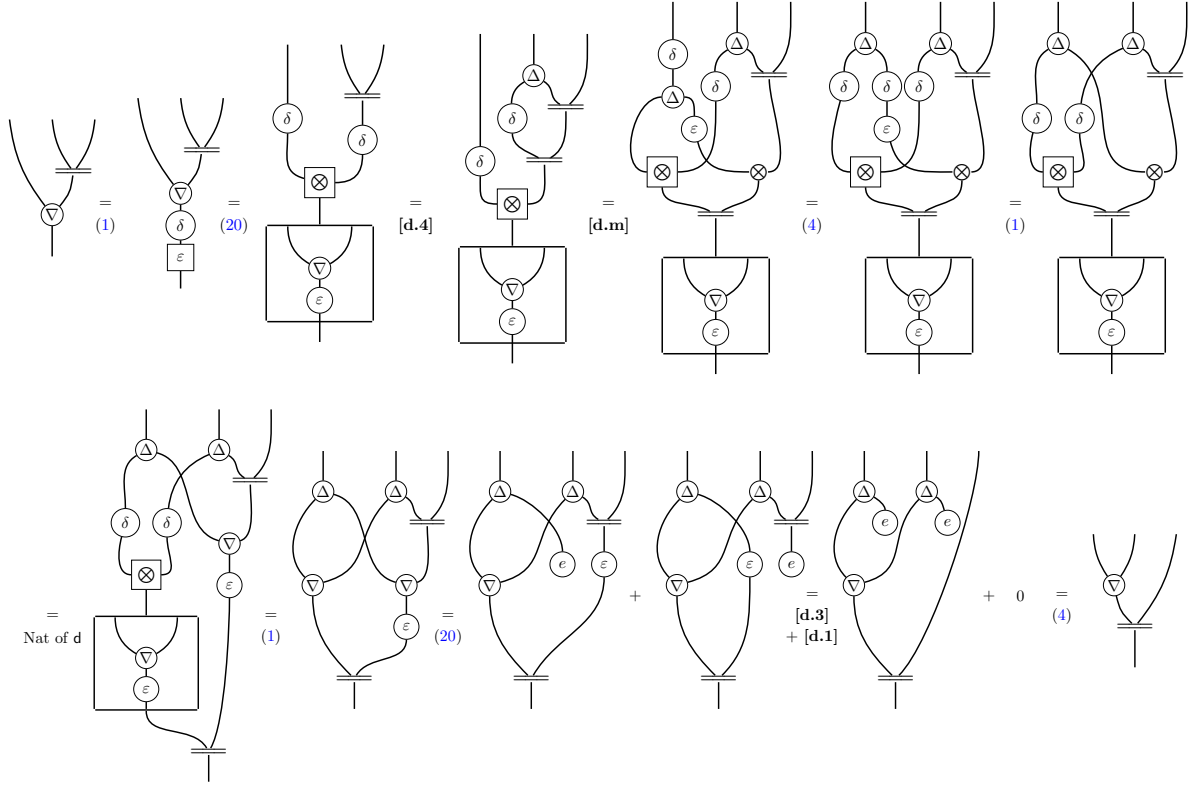
\[\resizebox{!}{5cm}{%
\begin{tikzpicture}
	\begin{pgfonlayer}{nodelayer}
		\node [style=port] (0) at (-4.5, -0.25) {};
		\node [style=port] (1) at (-4.5, 3) {};
		\node [style=duplicate] (2) at (-4.5, 0.75) {$\nabla$};
		\node [style=port] (3) at (-3.5, 3) {};
		\node [style=codifferential] (4) at (-4, 1.75) {{\bf =\!=\!=}};
		\node [style=port] (5) at (-5.5, 3) {};
		\node [style=component] (31) at (-1.5, 0.5) {$\delta$};
		\node [style=port] (32) at (-1.5, 3.5) {};
		\node [style=duplicate] (33) at (-1.5, 1.25) {$\nabla$};
		\node [style=port] (34) at (-0.5, 3.5) {};
		\node [style=codifferential] (35) at (-1, 2.25) {{\bf =\!=\!=}};
		\node [style=port] (36) at (-2.5, 3.5) {};
		\node [style=function] (37) at (-1.5, -0.25) {$\varepsilon$};
		\node [style=port] (38) at (-1.5, -1) {};
		\node [style=port] (39) at (-3.25, 1) {$=$};
		\node [style=port] (40) at (-3.25, 0.5) {(\ref{comonad})};
		\node [style=port] (41) at (-0.25, 1) {$=$};
		\node [style=function2] (42) at (1.75, 1.5) {$\bigotimes$};
		\node [style=component] (43) at (2.75, 2.5) {$\delta$};
		\node [style=component] (44) at (1, 3) {$\delta$};
		\node [style=port] (45) at (1, 4.75) {};
		\node [style=duplicate] (46) at (1.75, -0.5) {$\nabla$};
		\node [style=port] (47) at (2.5, 0.5) {};
		\node [style=port] (48) at (1, 0.5) {};
		\node [style=port] (49) at (0.5, 0.5) {};
		\node [style=port] (50) at (3, 0.5) {};
		\node [style=port] (51) at (3, -1.75) {};
		\node [style=port] (52) at (0.5, -1.75) {};
		\node [style=component] (53) at (1.75, -1.25) {$\varepsilon$};
		\node [style=port] (54) at (1.75, 0.5) {};
		\node [style=port] (55) at (1.75, 0.5) {};
		\node [style=port] (56) at (-0.25, 0.5) {(\ref{nablam1})};
		\node [style=port] (57) at (2.25, 4.75) {};
		\node [style=port] (58) at (3.25, 4.75) {};
		\node [style=codifferential] (59) at (2.75, 3.5) {{\bf =\!=\!=}};
		\node [style=port] (60) at (1.75, -2.25) {};
		\node [style=port] (63) at (5.5, 5) {};
		\node [style=function2] (70) at (6.25, 1) {$\bigotimes$};
		\node [style=component] (81) at (5.5, 2) {$\delta$};
		\node [style=differential] (84) at (7, 2) {{\bf =\!=\!=}};
		\node [style=component] (85) at (6.25, 3) {$\delta$};
		\node [style=differential] (86) at (7.5, 3.25) {{\bf =\!=\!=}};
		\node [style=duplicate] (87) at (6.75, 4) {$\Delta$};
		\node [style=port] (88) at (8, 5) {};
		\node [style=port] (89) at (6.75, 5) {};
		\node [style=port] (90) at (4, 1) {$=$};
		\node [style=port] (91) at (4, 0.5) {\textbf{[d.4]}};
		\node [style=component] (93) at (10, 4.5) {$\delta$};
		\node [style=differential] (95) at (12.25, 4) {{\bf =\!=\!=}};
		\node [style=port] (97) at (11.5, 5.75) {};
		\node [style=component] (105) at (11, 3.75) {$\delta$};
		\node [style=port] (106) at (10, 5.75) {};
		\node [style=duplicate] (109) at (11.5, 4.75) {$\Delta$};
		\node [style=port] (111) at (12.75, 5.75) {};
		\node [style=duplicate] (117) at (10, 3.5) {$\Delta$};
		\node [style=differential] (118) at (12, 1.75) {$\bigotimes$};
		\node [style=component] (119) at (10.5, 2.75) {$\varepsilon$};
		\node [style=function2] (120) at (9.75, 1.75) {$\bigotimes$};
		\node [style=codifferential] (121) at (11, 0.75) {{\bf =\!=\!=}};
		\node [style=port] (122) at (8.5, 1) {$=$};
		\node [style=port] (123) at (8.5, 0.5) {\textbf{[d.m]}};
		\node [style=duplicate] (135) at (6.25, -1) {$\nabla$};
		\node [style=port] (136) at (7, 0) {};
		\node [style=port] (137) at (5.5, 0) {};
		\node [style=port] (138) at (5, 0) {};
		\node [style=port] (139) at (7.5, 0) {};
		\node [style=port] (140) at (7.5, -2.25) {};
		\node [style=port] (141) at (5, -2.25) {};
		\node [style=component] (142) at (6.25, -1.75) {$\varepsilon$};
		\node [style=port] (143) at (6.25, 0) {};
		\node [style=port] (144) at (6.25, 0) {};
		\node [style=port] (145) at (6.25, -2.75) {};
		\node [style=duplicate] (146) at (11, -1.25) {$\nabla$};
		\node [style=port] (147) at (11.75, -0.25) {};
		\node [style=port] (148) at (10.25, -0.25) {};
		\node [style=port] (149) at (9.75, -0.25) {};
		\node [style=port] (150) at (12.25, -0.25) {};
		\node [style=port] (151) at (12.25, -2.5) {};
		\node [style=port] (152) at (9.75, -2.5) {};
		\node [style=component] (153) at (11, -2) {$\varepsilon$};
		\node [style=port] (154) at (11, -0.25) {};
		\node [style=port] (155) at (11, -0.25) {};
		\node [style=port] (156) at (11, -3) {};
		\node [style=differential] (158) at (17, 4) {{\bf =\!=\!=}};
		\node [style=port] (159) at (16.25, 5.75) {};
		\node [style=component] (160) at (15.75, 3.75) {$\delta$};
		\node [style=port] (161) at (14.5, 5.75) {};
		\node [style=duplicate] (162) at (16.25, 4.75) {$\Delta$};
		\node [style=port] (163) at (17.5, 5.75) {};
		\node [style=duplicate] (164) at (14.5, 4.75) {$\Delta$};
		\node [style=differential] (165) at (16.75, 1.75) {$\bigotimes$};
		\node [style=component] (166) at (15, 2.75) {$\varepsilon$};
		\node [style=function2] (167) at (14.5, 1.75) {$\bigotimes$};
		\node [style=codifferential] (168) at (15.75, 0.75) {{\bf =\!=\!=}};
		\node [style=duplicate] (169) at (15.75, -1.25) {$\nabla$};
		\node [style=port] (170) at (16.5, -0.25) {};
		\node [style=port] (171) at (15, -0.25) {};
		\node [style=port] (172) at (14.5, -0.25) {};
		\node [style=port] (173) at (17, -0.25) {};
		\node [style=port] (174) at (17, -2.5) {};
		\node [style=port] (175) at (14.5, -2.5) {};
		\node [style=component] (176) at (15.75, -2) {$\varepsilon$};
		\node [style=port] (177) at (15.75, -0.25) {};
		\node [style=port] (178) at (15.75, -0.25) {};
		\node [style=port] (179) at (15.75, -3) {};
		\node [style=component] (180) at (14, 3.75) {$\delta$};
		\node [style=component] (181) at (15, 3.75) {$\delta$};
		\node [style=port] (182) at (13.25, 1) {$=$};
		\node [style=port] (183) at (13.25, 0.5) {(\ref{coalgeq})};
		\node [style=differential] (184) at (21.5, 4) {{\bf =\!=\!=}};
		\node [style=port] (185) at (20.75, 5.75) {};
		\node [style=component] (186) at (19.5, 3) {$\delta$};
		\node [style=port] (187) at (19, 5.75) {};
		\node [style=duplicate] (188) at (20.75, 4.75) {$\Delta$};
		\node [style=port] (189) at (22, 5.75) {};
		\node [style=duplicate] (190) at (19, 4.75) {$\Delta$};
		\node [style=differential] (191) at (21.25, 1.75) {$\bigotimes$};
		\node [style=function2] (193) at (19, 1.75) {$\bigotimes$};
		\node [style=codifferential] (194) at (20.25, 0.75) {{\bf =\!=\!=}};
		\node [style=duplicate] (195) at (20.25, -1.25) {$\nabla$};
		\node [style=port] (196) at (21, -0.25) {};
		\node [style=port] (197) at (19.5, -0.25) {};
		\node [style=port] (198) at (19, -0.25) {};
		\node [style=port] (199) at (21.5, -0.25) {};
		\node [style=port] (200) at (21.5, -2.5) {};
		\node [style=port] (201) at (19, -2.5) {};
		\node [style=component] (202) at (20.25, -2) {$\varepsilon$};
		\node [style=port] (203) at (20.25, -0.25) {};
		\node [style=port] (204) at (20.25, -0.25) {};
		\node [style=port] (205) at (20.25, -3) {};
		\node [style=component] (206) at (18.5, 3) {$\delta$};
		\node [style=port] (208) at (18, 1) {$=$};
		\node [style=port] (209) at (18, 0.5) {(\ref{comonad})};
	\end{pgfonlayer}
	\begin{pgfonlayer}{edgelayer}
		\draw [style=wire, bend right=15, looseness=1.25] (2) to (4);
		\draw [style=wire, bend left=15] (2) to (5);
		\draw [style=wire] (0) to (2);
		\draw [style=wire, bend right=15] (4) to (3);
		\draw [style=wire, bend left=15] (4) to (1);
		\draw [style=wire, bend right=15, looseness=1.25] (33) to (35);
		\draw [style=wire, bend left=15] (33) to (36);
		\draw [style=wire] (31) to (33);
		\draw [style=wire, bend right=15] (35) to (34);
		\draw [style=wire, bend left=15] (35) to (32);
		\draw [style=wire] (31) to (37);
		\draw [style=wire] (37) to (38);
		\draw [style=wire, in=0, out=-90, looseness=1.25] (43) to (42);
		\draw [style=wire] (45) to (44);
		\draw [style=wire, in=165, out=-90, looseness=1.25] (44) to (42);
		\draw [style=wire, bend left] (46) to (48);
		\draw [style=wire, bend right] (46) to (47);
		\draw [style=wire] (52) to (51);
		\draw [style=wire] (52) to (49);
		\draw [style=wire] (49) to (50);
		\draw [style=wire] (51) to (50);
		\draw [style=wire] (46) to (53);
		\draw [style=wire] (42) to (55);
		\draw [style=wire, bend right=15] (59) to (58);
		\draw [style=wire, bend left=15] (59) to (57);
		\draw [style=wire] (59) to (43);
		\draw [style=wire] (53) to (60);
		\draw [style=wire, in=-90, out=180, looseness=1.50] (70) to (81);
		\draw [style=wire] (63) to (81);
		\draw [style=wire, bend right] (87) to (85);
		\draw [style=wire] (89) to (87);
		\draw [style=wire, in=-90, out=44] (86) to (88);
		\draw [style=wire, in=150, out=-30, looseness=1.25] (87) to (86);
		\draw [style=wire, in=30, out=-90] (86) to (84);
		\draw [style=wire, in=150, out=-90] (85) to (84);
		\draw [style=wire, in=0, out=-90, looseness=1.25] (84) to (70);
		\draw [style=wire] (106) to (93);
		\draw [style=wire, bend right] (109) to (105);
		\draw [style=wire] (97) to (109);
		\draw [style=wire, in=-90, out=44] (95) to (111);
		\draw [style=wire, in=150, out=-30, looseness=1.25] (109) to (95);
		\draw [style=wire, in=180, out=180] (117) to (120);
		\draw [style=wire, in=150, out=-90] (120) to (121);
		\draw [style=wire, in=30, out=-90, looseness=1.25] (118) to (121);
		\draw [style=wire, in=90, out=0, looseness=1.25] (117) to (119);
		\draw [style=wire, in=180, out=-75, looseness=1.25] (119) to (118);
		\draw [style=wire] (93) to (117);
		\draw [style=wire, in=0, out=-90, looseness=1.25] (105) to (120);
		\draw [style=wire, in=0, out=-90, looseness=0.75] (95) to (118);
		\draw [style=wire, bend left] (135) to (137);
		\draw [style=wire, bend right] (135) to (136);
		\draw [style=wire] (141) to (140);
		\draw [style=wire] (141) to (138);
		\draw [style=wire] (138) to (139);
		\draw [style=wire] (140) to (139);
		\draw [style=wire] (135) to (142);
		\draw [style=wire] (142) to (145);
		\draw [style=wire] (70) to (144);
		\draw [style=wire, bend left] (146) to (148);
		\draw [style=wire, bend right] (146) to (147);
		\draw [style=wire] (152) to (151);
		\draw [style=wire] (152) to (149);
		\draw [style=wire] (149) to (150);
		\draw [style=wire] (151) to (150);
		\draw [style=wire] (146) to (153);
		\draw [style=wire] (153) to (156);
		\draw [style=wire] (121) to (155);
		\draw [style=wire, bend right] (162) to (160);
		\draw [style=wire] (159) to (162);
		\draw [style=wire, in=-90, out=44] (158) to (163);
		\draw [style=wire, in=150, out=-30, looseness=1.25] (162) to (158);
		\draw [style=wire, in=150, out=-90] (167) to (168);
		\draw [style=wire, in=30, out=-90, looseness=1.25] (165) to (168);
		\draw [style=wire, in=180, out=-75, looseness=1.25] (166) to (165);
		\draw [style=wire, in=0, out=-90, looseness=1.25] (160) to (167);
		\draw [style=wire, in=0, out=-90, looseness=0.75] (158) to (165);
		\draw [style=wire, bend left] (169) to (171);
		\draw [style=wire, bend right] (169) to (170);
		\draw [style=wire] (175) to (174);
		\draw [style=wire] (175) to (172);
		\draw [style=wire] (172) to (173);
		\draw [style=wire] (174) to (173);
		\draw [style=wire] (169) to (176);
		\draw [style=wire] (176) to (179);
		\draw [style=wire] (168) to (178);
		\draw [style=wire] (161) to (164);
		\draw [style=wire] (181) to (166);
		\draw [style=wire, in=90, out=-30] (164) to (181);
		\draw [style=wire, in=90, out=-150] (164) to (180);
		\draw [style=wire, in=-180, out=-105] (180) to (167);
		\draw [style=wire, bend right] (188) to (186);
		\draw [style=wire] (185) to (188);
		\draw [style=wire, in=-90, out=44] (184) to (189);
		\draw [style=wire, in=150, out=-30, looseness=1.25] (188) to (184);
		\draw [style=wire, in=150, out=-90] (193) to (194);
		\draw [style=wire, in=30, out=-90, looseness=1.25] (191) to (194);
		\draw [style=wire, in=0, out=-90, looseness=1.25] (186) to (193);
		\draw [style=wire, in=0, out=-90, looseness=0.75] (184) to (191);
		\draw [style=wire, bend left] (195) to (197);
		\draw [style=wire, bend right] (195) to (196);
		\draw [style=wire] (201) to (200);
		\draw [style=wire] (201) to (198);
		\draw [style=wire] (198) to (199);
		\draw [style=wire] (200) to (199);
		\draw [style=wire] (195) to (202);
		\draw [style=wire] (202) to (205);
		\draw [style=wire] (194) to (204);
		\draw [style=wire] (187) to (190);
		\draw [style=wire, in=90, out=-150] (190) to (206);
		\draw [style=wire, in=-180, out=-105] (206) to (193);
		\draw [style=wire, in=180, out=-30, looseness=1.25] (190) to (191);
	\end{pgfonlayer}
\end{tikzpicture}
}%
\]
\[\resizebox{!}{5cm}{%
\begin{tikzpicture}
	\begin{pgfonlayer}{nodelayer}
		\node [style=differential] (0) at (19.5, 4.25) {{\bf =\!=\!=}};
		\node [style=port] (1) at (18.75, 6) {};
		\node [style=component] (2) at (17.5, 3.25) {$\delta$};
		\node [style=port] (3) at (17, 6) {};
		\node [style=duplicate] (4) at (18.75, 5) {$\Delta$};
		\node [style=port] (5) at (20, 6) {};
		\node [style=duplicate] (6) at (17, 5) {$\Delta$};
		\node [style=function2] (8) at (17, 2) {$\bigotimes$};
		\node [style=codifferential] (9) at (18.5, -2.25) {{\bf =\!=\!=}};
		\node [style=duplicate] (10) at (17, 0) {$\nabla$};
		\node [style=port] (11) at (17.75, 1) {};
		\node [style=port] (12) at (16.25, 1) {};
		\node [style=port] (13) at (15.75, 1) {};
		\node [style=port] (14) at (18.25, 1) {};
		\node [style=port] (15) at (18.25, -1.25) {};
		\node [style=port] (16) at (15.75, -1.25) {};
		\node [style=component] (17) at (17, -0.75) {$\varepsilon$};
		\node [style=port] (18) at (17, 1) {};
		\node [style=port] (19) at (17, 1) {};
		\node [style=port] (20) at (18.5, -3) {};
		\node [style=component] (21) at (16.5, 3.25) {$\delta$};
		\node [style=port] (24) at (14.75, 1) {$=$};
		\node [style=port] (25) at (14.75, 0.5) {Nat of $\mathsf{d}$};
		\node [style=port] (26) at (17, -1.25) {};
		\node [style=duplicate] (27) at (19, 2.75) {$\nabla$};
		\node [style=component] (28) at (19, 1.75) {$\varepsilon$};
		\node [style=differential] (29) at (23.75, 3.25) {{\bf =\!=\!=}};
		\node [style=port] (30) at (23, 5) {};
		\node [style=port] (32) at (21.25, 5) {};
		\node [style=duplicate] (33) at (23, 4) {$\Delta$};
		\node [style=port] (34) at (24.25, 5) {};
		\node [style=duplicate] (35) at (21.25, 4) {$\Delta$};
		\node [style=codifferential] (37) at (22, -0.5) {{\bf =\!=\!=}};
		\node [style=port] (48) at (22, -1.25) {};
		\node [style=duplicate] (51) at (23.25, 1.75) {$\nabla$};
		\node [style=component] (52) at (23.25, 0.75) {$\varepsilon$};
		\node [style=duplicate] (53) at (21, 1.75) {$\nabla$};
		\node [style=port] (56) at (20, 1) {$=$};
		\node [style=port] (57) at (20, 0.5) {(\ref{comonad})};
		\node [style=differential] (58) at (28, 3.25) {{\bf =\!=\!=}};
		\node [style=port] (59) at (27.25, 5) {};
		\node [style=port] (60) at (25.5, 5) {};
		\node [style=duplicate] (61) at (27.25, 4) {$\Delta$};
		\node [style=port] (62) at (28.5, 5) {};
		\node [style=duplicate] (63) at (25.5, 4) {$\Delta$};
		\node [style=codifferential] (64) at (26.25, -0.5) {{\bf =\!=\!=}};
		\node [style=port] (65) at (26.25, -1.25) {};
		\node [style=component] (67) at (28, 2.25) {$\varepsilon$};
		\node [style=duplicate] (68) at (25.25, 1.75) {$\nabla$};
		\node [style=component] (69) at (27, 2.25) {$e$};
		\node [style=differential] (70) at (32.5, 3.25) {{\bf =\!=\!=}};
		\node [style=port] (71) at (31.75, 5) {};
		\node [style=port] (72) at (30, 5) {};
		\node [style=duplicate] (73) at (31.75, 4) {$\Delta$};
		\node [style=port] (74) at (33, 5) {};
		\node [style=duplicate] (75) at (30, 4) {$\Delta$};
		\node [style=codifferential] (76) at (30.75, -0.5) {{\bf =\!=\!=}};
		\node [style=port] (77) at (30.75, -1.25) {};
		\node [style=component] (78) at (32.5, 2.25) {$e$};
		\node [style=duplicate] (79) at (29.75, 1.75) {$\nabla$};
		\node [style=component] (80) at (31.5, 2.25) {$\varepsilon$};
		\node [style=port] (81) at (28.75, 1.5) {$+$};
		\node [style=port] (82) at (24.25, 1) {$=$};
		\node [style=port] (83) at (24.25, 0.5) {(\ref{nablam1})};
		\node [style=port] (85) at (35.75, 5) {};
		\node [style=port] (86) at (34.25, 5) {};
		\node [style=duplicate] (87) at (35.75, 4) {$\Delta$};
		\node [style=port] (88) at (37, 5) {};
		\node [style=duplicate] (89) at (34.25, 4) {$\Delta$};
		\node [style=codifferential] (90) at (35, -0.5) {{\bf =\!=\!=}};
		\node [style=port] (91) at (35, -1.25) {};
		\node [style=duplicate] (93) at (34, 1.75) {$\nabla$};
		\node [style=component] (94) at (34.75, 3.25) {$e$};
		\node [style=port] (106) at (37.25, 1.5) {$+$};
		\node [style=component] (107) at (36.25, 3.25) {$e$};
		\node [style=port] (108) at (33, 1.5) {$=$};
		\node [style=port] (109) at (33, 1) {\textbf{[d.3]}};
		\node [style=port] (110) at (33, 0.5) {+ \textbf{[d.1]}};
		\node [style=port] (111) at (39.75, 3) {};
		\node [style=port] (112) at (41.75, 3) {};
		\node [style=codifferential] (113) at (40.75, 0.75) {{\bf =\!=\!=}};
		\node [style=port] (114) at (40.75, -0.25) {};
		\node [style=port] (115) at (40.75, 3) {};
		\node [style=duplicate] (116) at (40.25, 1.75) {$\nabla$};
		\node [style=port] (117) at (39, 1.5) {$=$};
		\node [style=port] (118) at (39, 1) {(\ref{coalgeq})};
		\node [style=port] (120) at (38, 1.5) {$0$};
	\end{pgfonlayer}
	\begin{pgfonlayer}{edgelayer}
		\draw [style=wire, bend right] (4) to (2);
		\draw [style=wire] (1) to (4);
		\draw [style=wire, in=-90, out=44] (0) to (5);
		\draw [style=wire, in=150, out=-30, looseness=1.25] (4) to (0);
		\draw [style=wire, in=0, out=-90, looseness=1.25] (2) to (8);
		\draw [style=wire, bend left] (10) to (12);
		\draw [style=wire, bend right] (10) to (11);
		\draw [style=wire] (16) to (15);
		\draw [style=wire] (16) to (13);
		\draw [style=wire] (13) to (14);
		\draw [style=wire] (15) to (14);
		\draw [style=wire] (10) to (17);
		\draw [style=wire] (3) to (6);
		\draw [style=wire, in=90, out=-150] (6) to (21);
		\draw [style=wire, in=-180, out=-105] (21) to (8);
		\draw [style=wire] (8) to (19);
		\draw [style=wire] (17) to (26);
		\draw [style=wire, in=120, out=-90] (26) to (9);
		\draw [style=wire] (9) to (20);
		\draw [style=wire] (27) to (28);
		\draw [style=wire, in=165, out=-15] (6) to (27);
		\draw [style=wire, in=15, out=-90, looseness=0.75] (0) to (27);
		\draw [style=wire, in=60, out=-90, looseness=0.50] (28) to (9);
		\draw [style=wire] (30) to (33);
		\draw [style=wire, in=-90, out=44] (29) to (34);
		\draw [style=wire, in=150, out=-30, looseness=1.25] (33) to (29);
		\draw [style=wire] (32) to (35);
		\draw [style=wire] (37) to (48);
		\draw [style=wire] (51) to (52);
		\draw [style=wire, in=165, out=-15] (35) to (51);
		\draw [style=wire, in=15, out=-90, looseness=0.75] (29) to (51);
		\draw [style=wire, in=60, out=-90] (52) to (37);
		\draw [style=wire, in=15, out=-135, looseness=1.25] (33) to (53);
		\draw [style=wire, in=135, out=-135] (35) to (53);
		\draw [style=wire, in=135, out=-90] (53) to (37);
		\draw [style=wire] (59) to (61);
		\draw [style=wire, in=-90, out=44] (58) to (62);
		\draw [style=wire, in=150, out=-30, looseness=1.25] (61) to (58);
		\draw [style=wire] (60) to (63);
		\draw [style=wire] (64) to (65);
		\draw [style=wire, in=60, out=-90] (67) to (64);
		\draw [style=wire, in=15, out=-135, looseness=1.25] (61) to (68);
		\draw [style=wire, in=135, out=-135] (63) to (68);
		\draw [style=wire, in=135, out=-90] (68) to (64);
		\draw [style=wire] (58) to (67);
		\draw [style=wire, in=90, out=-15, looseness=0.75] (63) to (69);
		\draw [style=wire] (71) to (73);
		\draw [style=wire, in=-90, out=44] (70) to (74);
		\draw [style=wire, in=150, out=-30, looseness=1.25] (73) to (70);
		\draw [style=wire] (72) to (75);
		\draw [style=wire] (76) to (77);
		\draw [style=wire, in=15, out=-135, looseness=1.25] (73) to (79);
		\draw [style=wire, in=135, out=-135] (75) to (79);
		\draw [style=wire, in=135, out=-90] (79) to (76);
		\draw [style=wire] (70) to (78);
		\draw [style=wire, in=90, out=-15, looseness=0.75] (75) to (80);
		\draw [style=wire, in=45, out=-90] (80) to (76);
		\draw [style=wire] (85) to (87);
		\draw [style=wire, style=wire] (86) to (89);
		\draw [style=wire] (90) to (91);
		\draw [style=wire, in=15, out=-135, looseness=1.25] (87) to (93);
		\draw [style=wire, in=135, out=-135] (89) to (93);
		\draw [style=wire, in=135, out=-90] (93) to (90);
		\draw [style=wire, in=90, out=-30, looseness=0.75] (89) to (94);
		\draw [style=wire, in=45, out=-90, looseness=0.50] (88) to (90);
		\draw [style=wire, in=105, out=-30] (87) to (107);
		\draw [style=wire, bend left=15, looseness=1.25] (113) to (116);
		\draw [style=wire, bend right=15] (113) to (112);
		\draw [style=wire] (114) to (113);
		\draw [style=wire, bend left=15] (116) to (111);
		\draw [style=wire, bend right=15] (116) to (115);
	\end{pgfonlayer}
\end{tikzpicture}
}%
\]
\end{proof}  

\begin{corollary} For the induced additive bialgebra modality of an additive linear category, the following are equivalent for a natural transformation $\mathsf{d}: \oc A \otimes A \to A$: 
\begin{enumerate}[{\em (i)}]
\item $\mathsf{d}$ is a deriving transformation; 
\item $\mathsf{d}$ satisfies the product rule {\bf [d.2]}, the linear rule {\bf [d.3]}, and the chain rule {\bf [d.4]};
\item $\mathsf{d}$ satisfies the linear rule {\bf [d.3]}, the chain rule {\bf [d.4]}, and the $\nabla$-rule {\bf [d.$\nabla$]};
\item $\mathsf{d}$ satisfies the linear rule {\bf [d.3]}, the chain rule {\bf [d.4]}, and the monoidal rule {\bf [d.m]}. 
\end{enumerate}
\end{corollary} 

Finally, this gives the following theorem: 

\begin{theorem} 
For the monoidal coalgebra modality of an additive linear category, all deriving transformations satisfy the monoidal rule ${\bf [d.m]}$ and are induced by a codereliction (for the induced additive bialgebra modality). 
\end{theorem} 

%%%%%%%%%%%%%%%%%%%%%%%%%%%%%%%%%%%%%%%%%%%%%%%%%%%%%%%%%%%%%%%%%%%%%%%
%%%%%%%%%%%%%%%%%%%%%%%%%%%%%%%%%%%%%%%%%%%%%%%%%%%%%%%%%%%%%%%%%%%%%%%

\section{Seely Isomorphisms and the Biproduct Completion}\label{Seelysec}

%%%%%%%%%%%%%%%%%%%%%%%%%%%%%%%%%%%%%%%%%%%%%%%%%%%%%%%%%%%%%%%%%%%%%%%

In this section we discuss additive bialgebra modalities in the presence of biproducts, which are equivalently described by the Seely isomorphisms and additive monoidal storage categories.  

\begin{definition}\label{Seelydef} In a symmetric monoidal category with finite products $\times$ and terminal object $\mathsf{T}$, a coalgebra modality has \textbf{Seely isomorphisms} \cite{seely1987linear,bierman1995categorical,blute2015cartesian} if the map $\chi_{\mathsf{T}}: \oc \mathsf{T} \to K$ and natural transformation ${\chi: \oc(A \times B) \to \oc A \otimes \oc B}$ defined respectively as:
\begin{equation}\label{}\begin{gathered} \xymatrixcolsep{2.5pc}\xymatrix{  \oc(\mathsf{T}) \ar[r]^-{e} & K &&\oc(A \times B) \ar[r]^-{\Delta} & \oc(A \times B) \otimes \oc(A \times B) \ar[rr]^-{\oc(\pi_0) \otimes \oc(\pi_1)} && \oc A \otimes \oc B
  } \end{gathered}\end{equation}
are isomorphisms, so $\oc(\mathsf{T}) \cong K$ and $\oc(A \times B) \cong \oc A \otimes \oc B$.  A \textbf{monoidal storage category}  \cite{blute2015cartesian} is a symmetric monoidal category with finite products and a coalgebra modality which has Seely isomorphisms. 
\end{definition}

It is worth pointing out that monoidal storage categories were called \textbf{new Seely categories} in \cite{bierman1995categorical,mellies2003categorical}. As explained in \cite{blute2015cartesian}, every coalgebra modality which has Seely isomorphisms is a monoidal coalgebra modality, where $m_\otimes$ is defined as 
$$\xymatrixcolsep{2.5pc}\xymatrix{\oc A \otimes \oc B \ar[r]^-{\chi^{-1}} & \oc(A \times B) \ar[r]^-{\delta} &   \oc \oc(A \times B) \ar[r]^-{\oc(\chi)} & \oc (\oc A \otimes \oc B) \ar[r]^-{\oc(\varepsilon \otimes \varepsilon)} & \oc(A \otimes B) 
  }$$
  and $m_K$ is defined as
  $$\xymatrixcolsep{2.5pc}\xymatrix{K \ar[r]^-{\chi^{-1}_{\mathsf{T}}} & \oc(\mathsf{T}) \ar[r]^-{\delta} & \oc \oc(\mathsf{T}) \ar[r]^-{\oc(\chi_\mathsf{T})} & \oc(K)
  }$$
  Conversly, in the presence of finite products, every monoidal coalgebra modality has Seely isomorphisms \cite{bierman1995categorical} where the inverse of $\chi$ is 
  $$ \xymatrixcolsep{2.5pc}\xymatrix{\oc A \otimes \oc B \ar[r]^-{\delta \otimes \delta} & \oc \oc A \otimes \oc \oc B \ar[r]^-{m_\otimes} & \oc(\oc A \otimes \oc B) \ar[rr]^-{\oc \left(\left \langle \varepsilon \otimes e, e \otimes \varepsilon \right \rangle \right)} && \oc(A \times B)  
  } $$
 while the inverse of $\chi_{\mathsf{T}}$ is
    $$  \xymatrixcolsep{2.5pc}\xymatrix{ K \ar[r]^-{m_K} & \oc(K) \ar[r]^-{\oc(\mathsf{t})} & \oc(\mathsf{T})
  } $$
where $\mathsf{t}: K \to \mathsf{T}$ is the unique map to the terminal object. Therefore we obtain the following:

\begin{theorem}\label{linstorthm}  \cite[Theorem 3.1.6]{blute2015cartesian} Every monoidal storage category is a linear category and conversely, every linear category with finite products is a monoidal storage category.   
\end{theorem} 

We now turn our attention to monoidal storage categories with additive structure:

\begin{definition} An \textbf{additive monoidal storage category} is a monoidal storage category which is also an additive symmetric monoidal category. 
\end{definition}

Notice, this implies that additive monoidal storage categories have finite biproducts $\times$ and a zero object $0$. As noted in \cite{blute2006differential}, the coalgebra modality of an additive monoidal storage category is an additive bialgebra modality where the multiplication and unit are defined respectively as: 
\[\xymatrixcolsep{2.5pc}\xymatrix{\oc A \otimes \oc A \ar[r]^-{\chi^{-1}} & \oc (A \times A) \ar[r]^-{\oc (\nabla_\times)} & \oc(A) & K \ar[r]^-{\chi^{-1}_{0}} & \oc 0 \ar[r]^-{\oc(0)} & \oc A
  }\]
  where $\nabla_\times$ is the codiagonal map of the biproduct. Conversly, every additive bialgebra modality satisfies the Seely isomorphisms where $\chi^{-1}$ and $\chi^{-1}_0$ are defined respectively as: 
 \[  \xymatrixcolsep{2.5pc}\xymatrix{\oc A \otimes \oc B  \ar[rr]^-{\oc(\iota_0) \otimes \oc(\iota_1)} && \oc(A \times B) \otimes \oc(A \times B) \ar[r]^-{\nabla} & \oc(A \times B) & K \ar[r]^-{u} & \oc 0 
  } \]
 where $\iota_0$ and $\iota_1$ are the injection maps of the biproduct. It is easy to check that this indeed gives the Seely isomorphisms, and therefore we have that: 

\begin{theorem} The following are equivalent: 
\begin{enumerate}[{\em (i)}]
\item An additive monoidal storage category;
\item An additive linear category with finite biproducts;
\item An additive symmetric monoidal category with finite biproducts and an additive bialgebra modality. 
\end{enumerate}
\end{theorem} 

Every additive symmetric monoidal category with an additive bialgebra modality induces an additive monoidal storage category via the biproduct completion. We first recall the biproduct completion for an additive category \cite{mac2013categories}. Let $\mathbb{X}$ be an additive category. Define the {biproduct completion} of $\mathbb{X}$, $\mathsf{B}[\mathbb{X}]$, as the category whose objects are list of objects of $\mathbb{X}$: $(A_1, \hdots, A_n)$, including the empty list $()$, and whose maps are matrices of maps of $\mathbb{X}$, including the empty matrix:
\[\xymatrixcolsep{2.5pc}\xymatrix{(A_1, \hdots, A_n) \ar[r]^-{[f_{i,j}]} & (B_1, \hdots, B_m)  
  }\]
  where $f_{i,j}: A_i \to B_j$. The composition in $\mathsf{B}[\mathbb{X}]$ is the standard matrix multiplication: 
\[[f_{i,j}][g_{l,k}]=[\sum f_{i,k} g_{k, j}]\]
while the identity is the standard identity matrix:
\[\xymatrixcolsep{2.5pc}\xymatrix{(A_1, \hdots, A_n) \ar[r]^-{[\delta_{i,j}]} & (A_1, \hdots, A_n)  
  } \]
  where $\delta_{i,j}=0$ if $i \neq j$, and $\delta_{i,i}=1$. It is easy to see that $\mathsf{B}[\mathbb{X}]$ does in fact have biproducts: 
\begin{lemma} $\mathsf{B}[\mathbb{X}]$ is a well-defined category with finite biproducts. 
\end{lemma}

If $\mathbb{X}$ is an additive symmetric monoidal category, then so is $\mathsf{B}[\mathbb{X}]$. The monoidal unit is the same as in $\mathbb{X}$, the tensor product of objects is:
\[(A_1, \hdots, A_n) \otimes (B_1, \hdots, B_m)= (A_1 \otimes B_1, \hdots, A_1 \otimes B_m, \hdots, A_n \otimes B_n)\]
while the tensor product of maps is the standard Kronecker product of matrices. 

\begin{lemma} If $\mathbb{X}$ is an additive symmetric monoidal category, then so is $\mathsf{B}[\mathbb{X}]$. 
\end{lemma}

If $\mathbb{X}$ admits an additive bialgebra modality, then $\mathsf{B}[\mathbb{X}]$ is an additive monoidal storage category where the Seely isomorphisms are strict, i.e., equalities, and in particular it is an additive linear category. We give the additive bialgebra modality of $\mathsf{B}[\mathbb{X}]$, and leave it to the reader to check that it is in fact an additive bialgebra modality. The functor $\oc: \mathsf{B}[\mathbb{X}] \to \mathsf{B}[\mathbb{X}]$ is defined on objects as:
\[\oc(A_1, \hdots, A_n)=\oc A_1 \otimes \hdots \otimes \oc A_n\]
and on a map $[f_{i,j}]: (A_1, \hdots, A_n) \to (B_1, \hdots, B_m)$, $\oc([f_{i,j}])$ is represented in the graphical calculus as: 
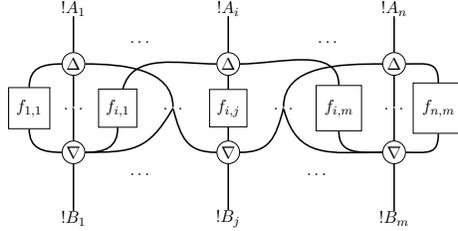
\[\resizebox{!}{3cm}{%
\begin{tikzpicture}
	\begin{pgfonlayer}{nodelayer}
		\node [style=port] (0) at (1.25, -1.5) {$\oc B_1$};
		\node [style=duplicate] (1) at (1.25, 0) {$\nabla$};
		\node [style=duplicate] (2) at (1.25, 2) {$\Delta$};
		\node [style=port] (3) at (1.25, 3.25) {$\oc A_1$};
		\node [style=port] (4) at (4.75, 3.25) {$\oc A_i$};
		\node [style=duplicate] (5) at (4.75, 2) {$\Delta$};
		\node [style=port] (6) at (8.5, 3.25) {$\oc A_n$};
		\node [style=duplicate] (7) at (8.5, 2) {$\Delta$};
		\node [style=port] (8) at (2.75, 2.5) {$\hdots$};
		\node [style=port] (9) at (7, 2.5) {$\hdots$};
		\node [style=port] (10) at (4.75, -1.5) {$\oc B_j$};
		\node [style=duplicate] (11) at (4.75, 0) {$\nabla$};
		\node [style=port] (12) at (8.5, -1.5) {$\oc B_m$};
		\node [style=duplicate] (13) at (8.5, 0) {$\nabla$};
		\node [style=port] (14) at (2.75, -0.5) {$\hdots$};
		\node [style=port] (15) at (6.75, -0.5) {$\hdots$};
		\node [style=function2] (16) at (0.25, 1) {$f_{1,1}$};
		\node [style=function2] (17) at (7.25, 1) {$f_{i,m}$};
		\node [style=function2] (18) at (2.25, 1) {$f_{i,1}$};
		\node [style=function2] (19) at (4.75, 1) {$f_{i,j}$};
		\node [style=function2] (20) at (9.5, 1) {$f_{n,m}$};
		\node [style=port] (21) at (3.5, 1) {$\hdots$};
		\node [style=port] (22) at (1.25, 1) {$\hdots$};
		\node [style=port] (23) at (8.5, 1) {$\hdots$};
		\node [style=port] (24) at (6, 1) {$\hdots$};
	\end{pgfonlayer}
	\begin{pgfonlayer}{edgelayer}
		\draw [style=wire] (0) to (1);
		\draw [style=wire] (3) to (2);
		\draw [style=wire] (4) to (5);
		\draw [style=wire] (6) to (7);
		\draw [style=wire] (10) to (11);
		\draw [style=wire] (12) to (13);
		\draw [style=wire, in=90, out=180, looseness=1.25] (2) to (16);
		\draw [style=wire, in=180, out=-90, looseness=1.25] (16) to (1);
		\draw [style=wire] (2) to (22);
		\draw [style=wire] (22) to (1);
		\draw [style=wire, in=75, out=180] (5) to (18);
		\draw [style=wire] (5) to (19);
		\draw [style=wire] (19) to (11);
		\draw [style=wire, in=90, out=0] (5) to (17);
		\draw [style=wire, in=90, out=0] (2) to (21);
		\draw [style=wire] (7) to (23);
		\draw [style=wire] (23) to (13);
		\draw [style=wire, in=0, out=-90, looseness=1.50] (20) to (13);
		\draw [style=wire, in=90, out=0, looseness=1.50] (7) to (20);
		\draw [style=wire, in=180, out=-90, looseness=1.25] (17) to (13);
		\draw [style=wire, in=0, out=-90, looseness=1.25] (18) to (1);
		\draw [style=wire, in=-105, out=0, looseness=1.25] (11) to (24);
		\draw [style=wire, in=180, out=-75, looseness=1.25] (21) to (11);
		\draw [style=wire, in=90, out=180] (7) to (24);
		\draw [style=wire, in=0, out=-120, looseness=1.25] (21) to (1);
		\draw [style=wire, in=180, out=-75, looseness=1.25] (24) to (13);
	\end{pgfonlayer}
\end{tikzpicture}
  }% 
\]
The bialgebra structure is given by the standard tensor product of bialgebras, the comonad comultiplication $\oc(A_1, \hdots, A_n) \to \oc\oc(A_1, \hdots, A_n)$ is represented in the graphical calculus as:
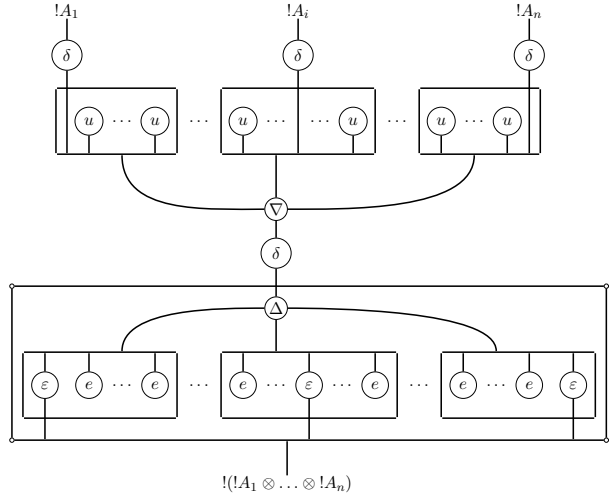
\[\resizebox{!}{6.5cm}{%
\begin{tikzpicture}
	\begin{pgfonlayer}{nodelayer}
		\node [style=component] (0) at (2.25, 3) {$u$};
		\node [style=port] (1) at (1.75, 5.5) {$\oc A_1$};
		\node [style=port] (2) at (4.25, 2.25) {};
		\node [style=duplicate] (3) at (0.5, -0.75) {};
		\node [style=duplicate] (4) at (6.5, -1.25) {$\Delta$};
		\node [style=duplicate] (5) at (6.5, 1) {$\nabla$};
		\node [style=port] (6) at (2.25, 2.25) {};
		\node [style=port] (7) at (4.25, 3.75) {};
		\node [style=component] (8) at (1.75, 4.5) {$\delta$};
		\node [style=port] (9) at (1.5, 2.25) {};
		\node [style=duplicate] (10) at (14, -0.75) {};
		\node [style=component] (11) at (6.5, 0) {$\delta$};
		\node [style=port] (12) at (1.75, 2.25) {};
		\node [style=port] (13) at (1.5, 3.75) {};
		\node [style=port] (14) at (3, 2.25) {};
		\node [style=port] (15) at (3, 3) {$\hdots$};
		\node [style=component] (16) at (3.75, 3) {$u$};
		\node [style=port] (17) at (3.75, 2.25) {};
		\node [style=component] (18) at (11.75, 3) {$u$};
		\node [style=component] (19) at (10.25, 3) {$u$};
		\node [style=port] (20) at (11.75, 2.25) {};
		\node [style=component] (21) at (12.25, 4.5) {$\delta$};
		\node [style=port] (22) at (10.25, 2.25) {};
		\node [style=port] (23) at (9.75, 2.25) {};
		\node [style=port] (24) at (12.25, 5.5) {$\oc A_n$};
		\node [style=port] (25) at (12.25, 2.25) {};
		\node [style=port] (26) at (11, 3) {$\hdots$};
		\node [style=port] (27) at (12.5, 3.75) {};
		\node [style=port] (28) at (11, 2.25) {};
		\node [style=port] (29) at (12.5, 2.25) {};
		\node [style=port] (30) at (9.75, 3.75) {};
		\node [style=component] (31) at (5.75, 3) {$u$};
		\node [style=component] (32) at (8.25, 3) {$u$};
		\node [style=port] (33) at (5.75, 2.25) {};
		\node [style=component] (34) at (7, 4.5) {$\delta$};
		\node [style=port] (35) at (8.25, 2.25) {};
		\node [style=port] (36) at (8.75, 2.25) {};
		\node [style=port] (37) at (7, 5.5) {$\oc A_i$};
		\node [style=port] (38) at (7, 2.25) {};
		\node [style=port] (39) at (7.5, 3) {$\hdots$};
		\node [style=port] (40) at (5.25, 3.75) {};
		\node [style=port] (41) at (6.5, 2.25) {};
		\node [style=port] (42) at (5.25, 2.25) {};
		\node [style=port] (43) at (8.75, 3.75) {};
		\node [style=port] (44) at (4.75, 3) {$\hdots$};
		\node [style=port] (45) at (9.25, 3) {$\hdots$};
		\node [style=port] (46) at (6.5, 3) {$\hdots$};
		\node [style=port] (47) at (5.25, -3.75) {};
		\node [style=port] (48) at (13.75, -3.75) {};
		\node [style=port] (49) at (7.25, -4.25) {};
		\node [style=port] (50) at (5.25, -2.25) {};
		\node [style=component] (51) at (8.75, -3) {$e$};
		\node [style=port] (52) at (13.25, -4.25) {};
		\node [style=port] (53) at (3, -2.25) {};
		\node [style=port] (54) at (9.75, -3) {$\hdots$};
		\node [style=component] (55) at (12.25, -3) {$e$};
		\node [style=port] (56) at (13.75, -2.25) {};
		\node [style=port] (57) at (11.5, -3) {$\hdots$};
		\node [style=port] (58) at (7.25, -2.25) {};
		\node [style=component] (59) at (10.75, -3) {$e$};
		\node [style=port] (60) at (8.75, -2.25) {};
		\node [style=component] (61) at (5.75, -3) {$e$};
		\node [style=component] (62) at (7.25, -3) {$\varepsilon$};
		\node [style=port] (63) at (3, -3) {$\hdots$};
		\node [style=port] (64) at (6.5, -3) {$\hdots$};
		\node [style=port] (65) at (0.75, -2.25) {};
		\node [style=port] (66) at (11.5, -2.25) {};
		\node [style=port] (67) at (4.25, -3.75) {};
		\node [style=component] (68) at (2.25, -3) {$e$};
		\node [style=port] (69) at (13.25, -2.25) {};
		\node [style=port] (70) at (4.75, -3) {$\hdots$};
		\node [style=component] (71) at (1.25, -3) {$\varepsilon$};
		\node [style=port] (72) at (8, -3) {$\hdots$};
		\node [style=port] (73) at (4.25, -2.25) {};
		\node [style=port] (74) at (1.25, -4.25) {};
		\node [style=port] (75) at (9.25, -3.75) {};
		\node [style=port] (76) at (9.25, -2.25) {};
		\node [style=component] (77) at (3.75, -3) {$e$};
		\node [style=port] (78) at (2.25, -2.25) {};
		\node [style=port] (79) at (10.25, -3.75) {};
		\node [style=port] (80) at (6.5, -2.25) {};
		\node [style=port] (81) at (3.75, -2.25) {};
		\node [style=port] (82) at (10.75, -2.25) {};
		\node [style=component] (83) at (13.25, -3) {$\varepsilon$};
		\node [style=port] (84) at (12.25, -2.25) {};
		\node [style=port] (85) at (5.75, -2.25) {};
		\node [style=port] (86) at (10.25, -2.25) {};
		\node [style=port] (87) at (1.25, -2.25) {};
		\node [style=port] (88) at (0.75, -3.75) {};
		\node [style=duplicate] (89) at (14, -4.25) {};
		\node [style=duplicate] (90) at (0.5, -4.25) {};
		\node [style=port] (91) at (6.75, -5.25) {$\oc(\oc A_1 \otimes \hdots \otimes \oc A_n)$};
		\node [style=port] (92) at (6.75, -4.25) {};
	\end{pgfonlayer}
	\begin{pgfonlayer}{edgelayer}
		\draw [style=wire] (1) to (8);
		\draw [style=wire] (2) to (7);
		\draw [style=wire] (6) to (0);
		\draw [style=wire] (9) to (13);
		\draw [style=wire] (2) to (9);
		\draw [style=wire] (7) to (13);
		\draw [style=wire] (5) to (11);
		\draw [style=wire] (11) to (4);
		\draw [style=wire] (3) to (10);
		\draw [style=wire] (8) to (12);
		\draw [style=wire] (17) to (16);
		\draw [style=wire] (24) to (21);
		\draw [style=wire] (23) to (30);
		\draw [style=wire] (20) to (18);
		\draw [style=wire] (29) to (27);
		\draw [style=wire] (23) to (29);
		\draw [style=wire] (30) to (27);
		\draw [style=wire] (21) to (25);
		\draw [style=wire] (22) to (19);
		\draw [style=wire] (37) to (34);
		\draw [style=wire] (36) to (43);
		\draw [style=wire] (33) to (31);
		\draw [style=wire] (42) to (40);
		\draw [style=wire] (36) to (42);
		\draw [style=wire] (43) to (40);
		\draw [style=wire] (34) to (38);
		\draw [style=wire] (35) to (32);
		\draw [style=wire, in=180, out=-90] (14) to (5);
		\draw [style=wire] (41) to (5);
		\draw [style=wire, in=0, out=-90, looseness=0.75] (28) to (5);
		\draw [style=wire] (74) to (71);
		\draw [style=wire] (73) to (67);
		\draw [style=wire] (78) to (68);
		\draw [style=wire] (65) to (88);
		\draw [style=wire] (73) to (65);
		\draw [style=wire] (67) to (88);
		\draw [style=wire] (71) to (87);
		\draw [style=wire] (81) to (77);
		\draw [style=wire] (52) to (83);
		\draw [style=wire] (86) to (79);
		\draw [style=wire] (84) to (55);
		\draw [style=wire] (56) to (48);
		\draw [style=wire] (86) to (56);
		\draw [style=wire] (79) to (48);
		\draw [style=wire] (83) to (69);
		\draw [style=wire] (82) to (59);
		\draw [style=wire] (49) to (62);
		\draw [style=wire] (76) to (75);
		\draw [style=wire] (85) to (61);
		\draw [style=wire] (50) to (47);
		\draw [style=wire] (76) to (50);
		\draw [style=wire] (75) to (47);
		\draw [style=wire] (62) to (58);
		\draw [style=wire] (60) to (51);
		\draw [style=wire] (90) to (89);
		\draw [style=wire] (3) to (90);
		\draw [style=wire] (10) to (89);
		\draw [style=wire, in=90, out=180, looseness=0.75] (4) to (53);
		\draw [style=wire] (4) to (80);
		\draw [style=wire, in=90, out=0, looseness=0.50] (4) to (66);
		\draw [style=wire] (92) to (91);
	\end{pgfonlayer}
\end{tikzpicture}
}%
\]
while the comonad counit is the following matrix:
\[\begin{bmatrix}\varepsilon_{A_1} \otimes e \otimes \hdots \otimes e, & 
\hdots & , e \otimes \hdots \otimes \varepsilon_{A_i} \otimes \hdots \otimes e, &
\hdots & , e \otimes e \otimes \hdots \otimes \varepsilon_{A_n}
\end{bmatrix} : \oc A_1 \otimes \hdots \otimes \oc A_n \longrightarrow (A_1, \hdots, A_n) \]

\begin{proposition}\label{biprodaddmod} If $\mathbb{X}$ has an additive bialgebra modality, then $\mathsf{B}[\mathbb{X}]$ is an additive monoidal storage category. 
\end{proposition}

Note that, as discussed in the introduction, Proposition \ref{biprodaddmod} together with Theorem \ref{linstorthm} provides a rather indirect verification of Theorem \ref{addlinaddbialg}.

If the additive bialgebra modality of $\mathbb{X}$ comes equipped with a codereliction $\eta$ then the additive bialgebra modality of $\mathsf{B}[\mathbb{X}]$ comes equipped with a codereliction defined as follows: 
$$ \begin{bmatrix}\eta_{A_1} \otimes u \otimes \hdots \otimes u \\ 
\hdots \\
u \otimes \hdots \otimes \eta_{A_i} \otimes \hdots \otimes u \\
\hdots \\
u \otimes u \otimes \hdots \otimes \eta_{A_n}
\end{bmatrix} : (A_1, \hdots, A_n) \longrightarrow  \oc A_1 \otimes \hdots \otimes \oc A_n$$

\begin{proposition} If $\mathbb{X}$ is a differential category with an additive bialgebra modality, then $\mathsf{B}[\mathbb{X}]$ is a differential category which is an additive monoidal storage category. 
\end{proposition} 

%%%%%%%%%%%%%%%%%%%%%%%%%%%%%%%%%%%%%%%%%%%%%%%%%%%%%%%%%%%%%%%%%%%%%%%
%%%%%%%%%%%%%%%%%%%%%%%%%%%%%%%%%%%%%%%%%%%%%%%%%%%%%%%%%%%%%%%%%%%%%%%

\section{Constructing Non-Additive Bialgebra Modalities}\label{nonaddbialg}

%%%%%%%%%%%%%%%%%%%%%%%%%%%%%%%%%%%%%%%%%%%%%%%%%%%%%%%%%%%%%%%%%%%%%%%

In this section we give a construction of non-additive bialgebra modalities induced by an additive algebra modalities. Given an additive bialgebra modality $(\oc, \delta, \varepsilon, \Delta, e, \nabla, u)$ on an additive symmetric monoidal category $\mathbb{X}$, for each object $B$ consider the functor $\oc^B: \mathbb{X} \to \mathbb{X}$ defined on objects as $\oc^B A= \oc B \otimes \oc A$, and on a map $f: A \to C$ as $\oc^B (f) = 1 \otimes \oc(f): \oc B \otimes \oc A \to \oc B \otimes \oc C$. Consider the natural transformations $\delta^B: \oc^B(A) \to \oc^B( \oc^B(A))$ and $\varepsilon^B: \oc^B(A) \to A$ defined as follows: 
  \[\delta^B :=  \xymatrixcolsep{5pc}\xymatrix{\oc B \otimes \oc A \ar[r]^-{\Delta \otimes 1} & \oc B \otimes \oc B \otimes \oc A \ar[r]^-{1 \otimes \delta \otimes \delta} & \oc B \otimes \oc \oc B \otimes \oc \oc A \ar[r]^-{1 \otimes m_\otimes} & \oc B \otimes \oc (\oc B \otimes \oc A)  
  } \]
  \[ \varepsilon^B := \xymatrixcolsep{5pc}\xymatrix{\oc B \otimes \oc A \ar[r]^-{e \otimes \varepsilon} & A  
  } \]
which drawn in the graphical calculus gives:
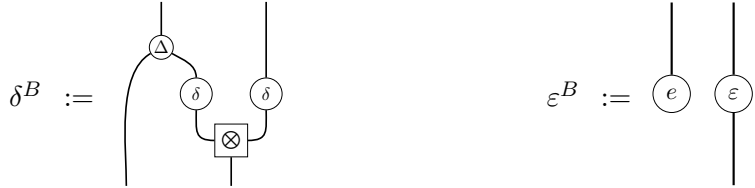
\begin{align*}
\begin{array}[c]{c}
\delta^B  
   \end{array} := 
 \begin{array}[c]{c}
\resizebox{!}{2.5cm}{%
\begin{tikzpicture}
	\begin{pgfonlayer}{nodelayer}
		\node [style=duplicate] (0) at (2.75, 3) {$\Delta$};
		\node [style=port] (1) at (2.75, 4) {};
		\node [style=port] (2) at (5, 4) {};
		\node [style=component] (3) at (3.5, 2) {$\delta$};
		\node [style=port] (4) at (2, 0) {};
		\node [style=function2] (5) at (4.25, 1) {$\bigotimes$};
		\node [style=port] (6) at (4.25, 0) {};
		\node [style=component] (7) at (5, 2) {$\delta$};
	\end{pgfonlayer}
	\begin{pgfonlayer}{edgelayer}
		\draw [style=wire] (5) to (6);
		\draw [style=wire] (1) to (0);
		\draw [style=wire, in=90, out=-30] (0) to (3);
		\draw [style=wire, in=180, out=-90, looseness=1.50] (3) to (5);
		\draw [style=wire, in=0, out=-90, looseness=1.50] (7) to (5);
		\draw [style=wire] (2) to (7);
		\draw [style=wire, in=90, out=-150, looseness=0.75] (0) to (4);
	\end{pgfonlayer}
\end{tikzpicture}
}%
   \end{array}  &&   \begin{array}[c]{c}
\varepsilon^B 
   \end{array}  :=   \begin{array}[c]{c}
 \resizebox{!}{2.5cm}{%
 \begin{tikzpicture}
	\begin{pgfonlayer}{nodelayer}
		\node [style=component] (10) at (9.75, 2.25) {$\varepsilon$};
		\node [style=port] (11) at (9.75, 3.75) {};
		\node [style=port] (12) at (9.75, 0.75) {};
		\node [style=component] (13) at (8.75, 2.25) {$e$};
		\node [style=port] (14) at (8.75, 3.75) {};
	\end{pgfonlayer}
	\begin{pgfonlayer}{edgelayer}
		\draw [style=wire] (11) to (10);
		\draw [style=wire] (10) to (12);
		\draw [style=wire] (14) to (13);
	\end{pgfonlayer}
\end{tikzpicture}
}%
   \end{array} \end{align*}

\begin{lemma} $(\oc^B, \delta^B, \varepsilon^B)$ is a comonad. 
\end{lemma}  
\begin{proof} We must show the following three identities: \\
\noindent 1. $\delta^B \delta^B=\delta^B \oc^B(\delta^B)$: Here we use that $\delta$ is a monoidal transformation, the naturality of $m_\otimes$, the co-associativity of $\Delta$, the associativity of $m_\otimes$, the co-associativity of the comonad, and that $\Delta$ is a $\oc$-coalgebra morphism:
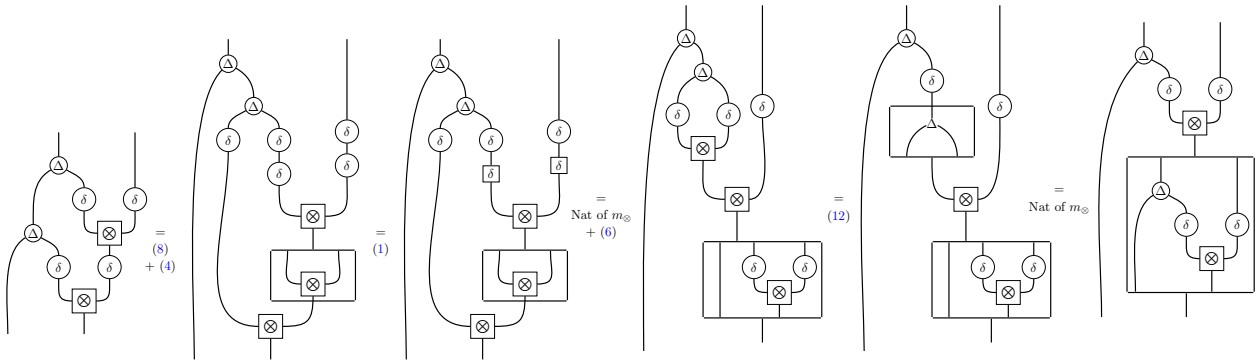
\[\resizebox{!}{4.75cm}{%
%
}%
\]
2. $\delta^B\varepsilon^B=1$: Here we use the counit law of the comultiplication, that $\varepsilon$ is a monoidal transformation, and the triangle identities of the comonad: 
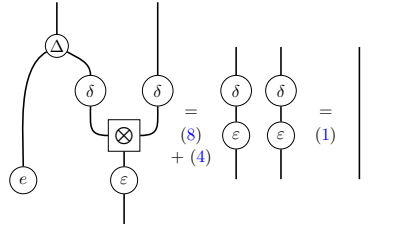
\[  \resizebox{!}{3cm}{%
%
}%
  \]
3. $\delta^B\oc^B(\varepsilon^B)=1$: Here we use the naturality of $m_\otimes$, that $e$ is a monoidal transformation, the comonad triangle identities, the unit law of $m_\otimes$, and the counit law for the comultiplication:
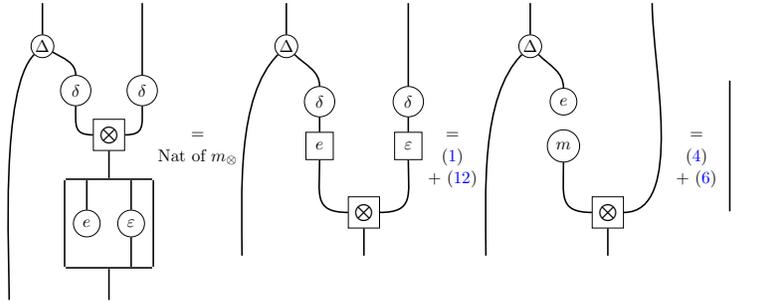
\[   \resizebox{!}{4cm}{%
%
}%
\]
\end{proof} 

The bialgebra structure of $\oc^B A$ is given by the standard tensor product of bialgebras, that is, the comultiplication $\Delta^B$ and multiplication $\nabla^B$ are defined respectively as: 
  \[ \Delta^B := \xymatrixcolsep{5pc}\xymatrix{\oc B \otimes \oc A \ar[r]^-{\Delta \otimes \Delta} & \oc B \otimes \oc B \otimes \oc A \otimes \oc A \ar[r]^-{1 \otimes \sigma \otimes 1} & \oc B \otimes \oc A \otimes \oc B \otimes \oc A  
  } \] 
    \[ \nabla^B := \xymatrixcolsep{5pc}\xymatrix{ \oc B \otimes \oc A \otimes \oc B \otimes \oc A\ar[r]^-{1 \otimes \sigma \otimes 1} & \oc B \otimes \oc B \otimes \oc A \otimes \oc A \ar[r]^-{\nabla \otimes \nabla} &   \oc B \otimes \oc A 
  } \] 
while the counit $e^B$ and the unit $u^B$ are:
  \[  e^B:= \xymatrixcolsep{5pc}\xymatrix{\oc B \otimes \oc A \ar[r]^-{e \otimes e} &K   
  } ~~~~~ u^B:=  \xymatrixcolsep{5pc}\xymatrix{K \ar[r]^-{u \otimes u} & \oc B \otimes \oc A} \]  
  which drawn in the graphical calculus is: 
  \begin{align*}
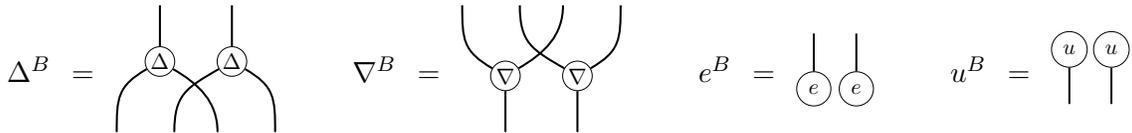
\begin{array}[c]{c}
\Delta^B 
   \end{array}=
   \begin{array}[c]{c}
\resizebox{!}{1.75 cm}{%
\begin{tikzpicture}
	\begin{pgfonlayer}{nodelayer}
		\node [style=port] (0) at (6.25, 5.25) {};
		\node [style=duplicate] (1) at (6.25, 4.25) {$\Delta$};
		\node [style=port] (2) at (5.5, 3) {};
		\node [style=port] (3) at (6.5, 3) {};
		\node [style=port] (4) at (7.25, 3) {};
		\node [style=port] (5) at (8.25, 3) {};
		\node [style=duplicate] (6) at (7.5, 4.25) {$\Delta$};
		\node [style=port] (7) at (7.5, 5.25) {};
	\end{pgfonlayer}
	\begin{pgfonlayer}{edgelayer}
		\draw [style=wire] (0) to (1);
		\draw [style=wire] (7) to (6);
		\draw [style=wire, in=90, out=-150, looseness=1.25] (6) to (3);
		\draw [style=wire, in=90, out=-30] (1) to (4);
		\draw [style=wire, in=-150, out=90, looseness=1.25] (2) to (1);
		\draw [style=wire, in=90, out=-30, looseness=1.25] (6) to (5);
	\end{pgfonlayer}
\end{tikzpicture}
}%
   \end{array} && \begin{array}[c]{c}
\nabla^B
   \end{array}=
   \begin{array}[c]{c}
\resizebox{!}{1.75cm}{%
\begin{tikzpicture}
	\begin{pgfonlayer}{nodelayer}
		\node [style=port] (0) at (6.25, 3) {};
		\node [style=duplicate] (1) at (6.25, 4) {$\nabla$};
		\node [style=port] (2) at (5.5, 5.25) {};
		\node [style=port] (3) at (6.5, 5.25) {};
		\node [style=port] (4) at (7.25, 5.25) {};
		\node [style=port] (5) at (8.25, 5.25) {};
		\node [style=duplicate] (6) at (7.5, 4) {$\nabla$};
		\node [style=port] (7) at (7.5, 3) {};
	\end{pgfonlayer}
	\begin{pgfonlayer}{edgelayer}
		\draw [style=wire] (0) to (1);
		\draw [style=wire] (7) to (6);
		\draw [style=wire, in=-90, out=150, looseness=1.25] (6) to (3);
		\draw [style=wire, in=-90, out=30] (1) to (4);
		\draw [style=wire, in=150, out=-90, looseness=1.25] (2) to (1);
		\draw [style=wire, in=-90, out=30, looseness=1.25] (6) to (5);
	\end{pgfonlayer}
\end{tikzpicture}
}%
   \end{array} && 
   \begin{array}[c]{c}
e^B
   \end{array}=
   \begin{array}[c]{c}
\resizebox{!}{1cm}{%
\begin{tikzpicture}
	\begin{pgfonlayer}{nodelayer}
		\node [style=port] (0) at (6.25, 4) {};
		\node [style=component] (1) at (6.25, 3) {$e$};
		\node [style=component] (6) at (7, 3) {$e$};
		\node [style=port] (7) at (7, 4) {};
	\end{pgfonlayer}
	\begin{pgfonlayer}{edgelayer}
		\draw [style=wire] (0) to (1);
		\draw [style=wire] (7) to (6);
	\end{pgfonlayer}
\end{tikzpicture}
}%
   \end{array} &&
   \begin{array}[c]{c}
u^B 
   \end{array}=
   \begin{array}[c]{c}
\resizebox{!}{1cm}{%
\begin{tikzpicture}
	\begin{pgfonlayer}{nodelayer}
		\node [style=port] (0) at (6.25, 3) {};
		\node [style=component] (1) at (6.25, 4) {$u$};
		\node [style=component] (6) at (7, 4) {$u$};
		\node [style=port] (7) at (7, 3) {};
	\end{pgfonlayer}
	\begin{pgfonlayer}{edgelayer}
		\draw [style=wire] (0) to (1);
		\draw [style=wire] (7) to (6);
	\end{pgfonlayer}
\end{tikzpicture}
}%
   \end{array}
\end{align*}

\begin{proposition} $(\oc^B, \delta^B, \varepsilon^B, \Delta^B, e^B, \nabla^B, u^B)$ is a bialgebra modality. 
\end{proposition} 
\begin{proof} We first show that $\delta^B$ preserves the comultiplication, which follows from the fact that $\delta$, $\Delta$, and $m_\otimes$ are all comonoid morphisms: 
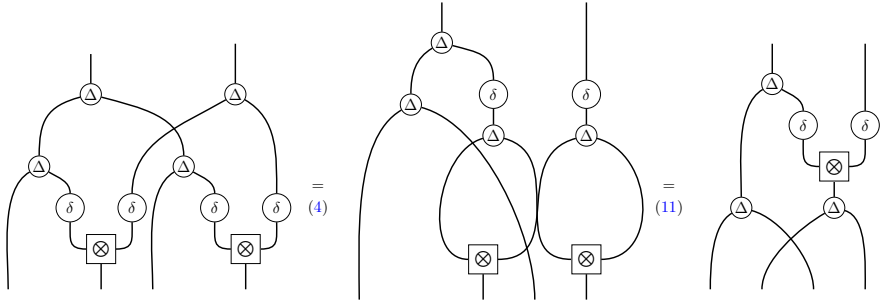
\[\resizebox{!}{4cm}{%
%
}%
\]
Next we show the compatibility relation between $\varepsilon^B$ and the multiplication. Here we use the compatibility between the counit and the multiplication, and the compatibility between $\varepsilon$ and the multiplication: 
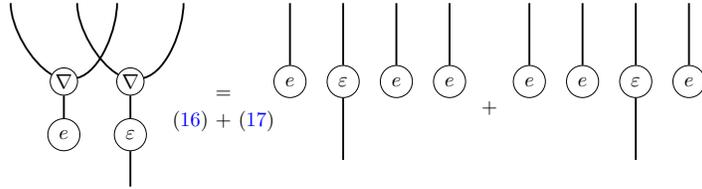
\[   \resizebox{!}{2.5cm}{%
%
}%
\]
\end{proof} 

In general, this bialgebra modality is not additive (unless $\oc(B) \cong K$). In particular, for the zero map we have that: 
\[\oc^B (0)= 1 \otimes \oc (0) = (1 \otimes e)(1 \otimes u) \neq (e \otimes e)(u \otimes u) = e^B u^B\]

Every codereliction $\eta$ on the additive bialgebra modality induces a codereliction $\eta^B: A \to \oc^B(A)$ defined as follows: 
\begin{align*}
\begin{array}[c]{c}
\xymatrixcolsep{5pc}\xymatrix{A \ar[r]^-{u \otimes \eta} & \oc B \otimes \oc A  
  } 
   \end{array}&&
   \begin{array}[c]{c}
\eta^B
   \end{array} :=    \begin{array}[c]{c}
   \resizebox{!}{2cm}{%
    \begin{tikzpicture}
	\begin{pgfonlayer}{nodelayer}
		\node [style=port] (0) at (-1.75, 2) {};
		\node [style=port] (1) at (-0.75, 5) {};
		\node [style=port] (2) at (-0.75, 2) {};
		\node [style={circle, draw}] (3) at (-1.75, 3.5) {$u$};
		\node [style={circle, draw}] (4) at (-0.75, 3.5) {$\eta$};
	\end{pgfonlayer}
	\begin{pgfonlayer}{edgelayer}
		\draw [style=wire] (2) to (4);
		\draw [style=wire] (4) to (1);
		\draw [style=wire] (0) to (3);
	\end{pgfonlayer}
\end{tikzpicture}
}%
   \end{array}
\end{align*}
  
\begin{proposition} $\eta^B$ is a codereliction for $(\oc^B, \delta^B, \varepsilon^B, \Delta^B, e^B, \nabla^B, u^B)$. 
\end{proposition}
\begin{proof}We must show \textbf{[dC.2]}, \textbf{[dC.3]}, and \textbf{[dC.4]}: \\Ê\\
\textbf{[dC.2]}: Here we use the bialgebra identity between the unit and the comultiplication, and that $\eta$ satisfies the product rule \textbf{[dC.2]}: 
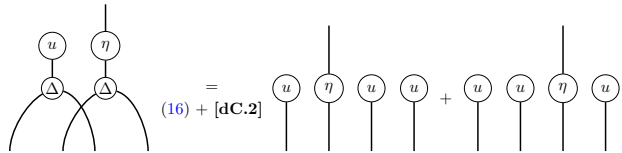
\[  \resizebox{!}{2cm}{%
%
}%
\]
\textbf{[dC.3]}: Here we use the bialgebra identity between the unit and counit, and that $\eta$ satisfies the linear rule \textbf{[dC.3]}: 
\[\resizebox{!}{2cm}{%
%
}%
 \]
\textbf{[dC.4]}: Here we use coassociativity of the comultiplication, one of the identities of Proposition \ref{monoidalnabla}, that $\delta$ preserves the comultiplication, and that $\eta$ satisfies the chain rule \textbf{[dC.4]} and the monoidal rule \textbf{[dC.m]}:
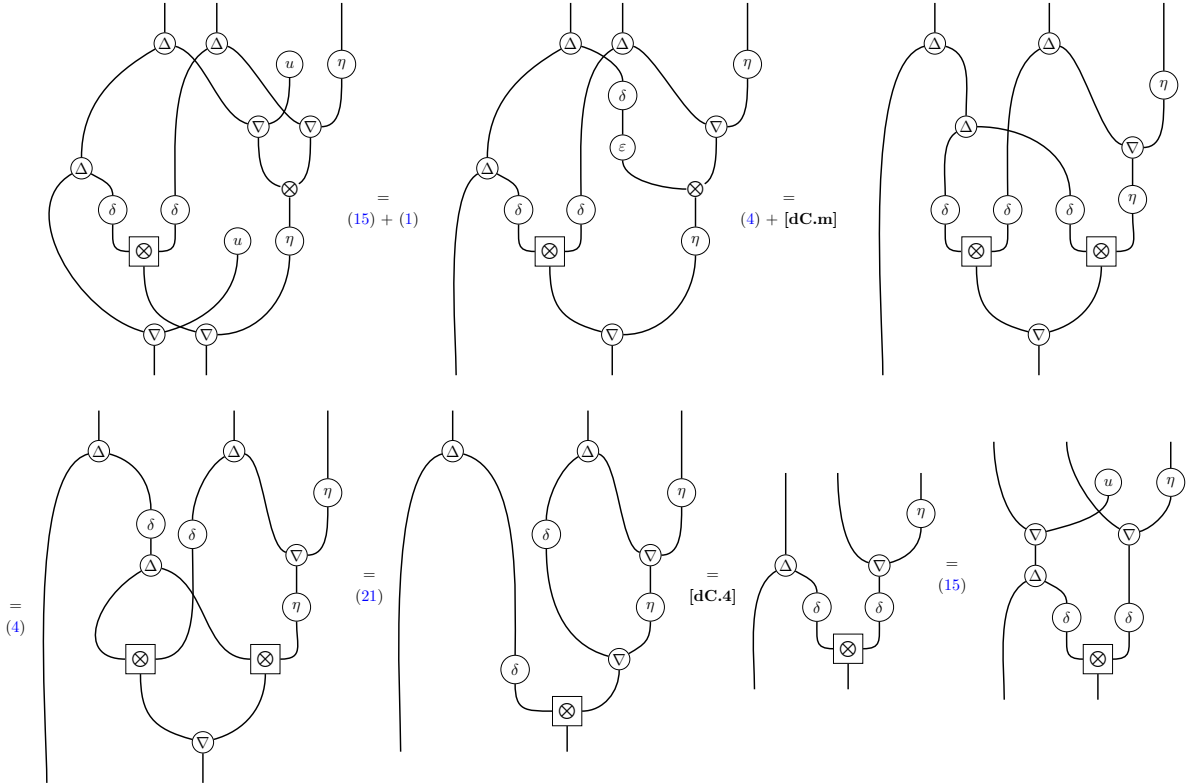
\[\resizebox{!}{5cm}{%
%
}%
\]
\end{proof} 

The induced deriving transformation $\mathsf{d}^B$:
  \[  \xymatrixcolsep{5pc}\xymatrix{\oc B \otimes \oc A \otimes A \ar[r]^-{1 \otimes \mathsf{d}} & \oc B \otimes \oc A  
  } \]   
intuitively, this should be thought of as the partial derivative with respect to $A$. 

%%%%%%%%%%%%%%%%%%%%%%%%%%%%%%%%%%%%%%%%%%%%%%%%%%%%%%%%%%%%%%%%%%%%%%%

\section{Separating Examples}\label{examples}

Here we present an overview of separating examples between the various structures defined throughout this paper. To help understand what the examples illustrate we present a Venn diagram which classifies the examples 
we shall give below. 
 \[  \begin{array}[c]{c}\resizebox{!}{12cm}
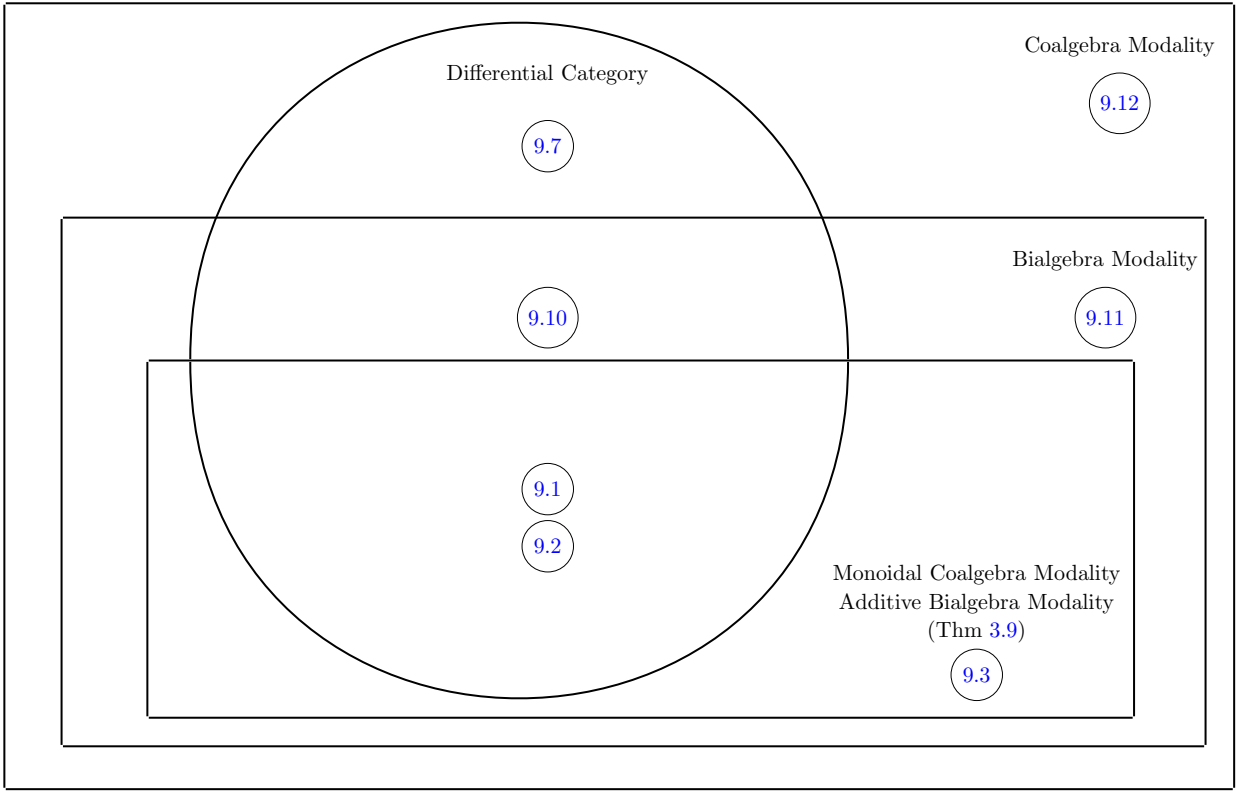
{%
\begin{tikzpicture}
	\begin{pgfonlayer}{nodelayer}
		\node [style=port] (0) at (-8.5, -5) {};
		\node [style={circle, draw}] (1) at (7.5, -3.75) {\ref{ex3}};
		\node [style=port] (2) at (-9.5, -5.75) {};
		\node [style=port] (3) at (11.5, 4.25) {};
		\node [style=port] (4) at (12, 8) {};
		\node [style=port] (5) at (-9.5, 8) {};
		\node [style=port] (6) at (-8.5, 4.25) {};
		\node [style=port] (7) at (11.5, -5) {};
		\node [style=port] (8) at (10, 7.25) {Coalgebra Modality};
		\node [style={circle, draw}] (9) at (10, 6.25) {\ref{ex7}};
		\node [style=port] (10) at (12, -5.75) {};
		\node [style=port] (11) at (7.5, -2) {Monoidal Coalgebra Modality};
		\node [style=port] (12) at (-7, 1.75) {};
		\node [style=port] (13) at (10.25, 1.75) {};
		\node [style=port] (14) at (10.25, -4.5) {};
		\node [style=port] (15) at (-7, -4.5) {};
		\node [style=port] (16) at (7.5, -2.5) {Additive Bialgebra Modality};
		\node [style=port] (17) at (9.75, 3.5) {Bialgebra Modality};
		\node [style={circle, draw}] (18) at (9.75, 2.5) {\ref{ex6}};
		\node [style=port] (19) at (5.25, 1.75) {};
		\node [style=port] (20) at (0, 6.75) {Differential Category};
		\node [style={circle, draw}] (21) at (0, 5.5) {\ref{ex4}};
		\node [style={circle, draw}] (22) at (0, -0.5) {\ref{ex1}};
		\node [style={circle, draw}] (23) at (0, 2.5) {\ref{ex5}};
		\node [style=port] (24) at (-6.25, 1.75) {};
		\node [style={circle, draw}] (25) at (0, -1.5) {\ref{ex2}};
		\node [style={}] (25) at (7.5, -3) {(Thm \ref{addlinaddbialg})};
	\end{pgfonlayer}
	\begin{pgfonlayer}{edgelayer}
		\draw [style=wire] (2) to (5);
		\draw [style=wire] (2) to (10);
		\draw [style=wire] (10) to (4);
		\draw [style=wire] (5) to (4);
		\draw [style=wire] (0) to (7);
		\draw [style=wire] (6) to (3);
		\draw [style=wire] (0) to (6);
		\draw [style=wire] (7) to (3);
		\draw [style=wire] (15) to (14);
		\draw [style=wire] (15) to (12);
		\draw [style=wire] (14) to (13);
		\draw [style=wire, bend left=90, looseness=1.75] (24) to (19);
		\draw [style=wire, bend right=90, looseness=1.75] (24) to (19);
		\draw [style=wire] (12) to (13);
	\end{pgfonlayer}
\end{tikzpicture}
}% 
\end{array}\]

A well-known example of a differential category (whose coalgebra modality also happens to be monoidal) comes from the free symmetric algebra construction which actually gives a co-differential category. We briefly recall this example (see \cite{blute2006differential} for more details):

\begin{example}\label{ex1} \normalfont Let $R$ be a commutative ring. For an $R$-module $M$, define $\mathsf{Sym}(M)$, called the free symmetric algebra over $M$, as follows (see Section 8, Chapter XVI in \cite{lang2002algebra} for more details): 
$$\mathsf{Sym}(M)= \bigoplus^{\infty}_{n=0} \mathsf{Sym}^n(M)= R \oplus M \oplus \mathsf{Sym}^2(M) \oplus ... $$ 
where $\mathsf{Sym}^n(M)$ is simply the quotient of $M^{\otimes^n}$ by the tensor symmetry equalities: 
$$a_1 \otimes ... \otimes a_i \otimes ... \otimes a_n = a_{\sigma(1)} \otimes ... \otimes a_{\sigma(i)} \otimes ... \otimes a_{\sigma(n)}$$
$\mathsf{Sym}(M)$ is a commutative algebra where the multiplication $\nabla: \mathsf{Sym}(M) \otimes \mathsf{Sym}(M) \to \mathsf{Sym}(M)$ is the concatenation of words $\nabla(v_1 \otimes ... \otimes v_n, w_1 \otimes ... \otimes w_m)= v_1 \otimes ... \otimes v_n \otimes w_1 \otimes ... \otimes w_m$ which we then extend by linearity, and the unit $u: R \to \mathsf{Sym}(M)$ is the injection map of $R$ into $\mathsf{Sym}(M)$. Furthermore, $\mathsf{Sym}(M)$ is the free commutative algebra over $M$, that is, we obtain an adjunction: 
$$\xymatrixcolsep{2.5pc}\xymatrix{\mathsf{MOD}_R \ar@<+1.1ex>[r]^-{\mathsf{Sym}} & \mathsf{CALG}_{R} \ar@<+1ex>[l]_-{\bot}^-{U}
  }$$
The unit $\eta: M \to \mathsf{Sym}(M)$ is the injection map of $M$ into $\mathsf{Sym}(M)$ and for an algebra $A$, the counit $\epsilon: \mathsf{Sym}(A) \to A$ is defined on pure tensors as $\epsilon(a_1 \otimes ... \otimes a_n) = a_1...a_n$, which we then extend by linearity. The induced monad $(\mathsf{Sym}, \mu, \eta)$ is an algebra modality (the dual of a coalgebra modality) the multiplication of the monad is an algebra morphism (as it's a map in the category of algebras). Furthermore, this algebra modality satisfies the Seely isomorphism \cite{lang2002algebra}, that is:
$$\mathsf{Sym}(M \oplus N) \cong \mathsf{Sym}(M) \otimes \mathsf{Sym}(N) \quad \mathsf{Sym}(0)\cong R$$
which implies that the free symmetric algebra adjunction induces an additive linear category (or equivalently an additive Seely category) structure on $\mathsf{MOD}_R$. Furthermore, it comes equipped with a deriving transformation, making $\mathsf{MOD}^{op}_R$ into a co-differential category. The deriving transformation $\mathsf{d}: \mathsf{Sym}(M) \to ~\mathsf{Sym}(M) \otimes M$ on pure tensors is defined as follows:
$$\mathsf{d}(a_1 \otimes ... \otimes a_n)= \sum_{i=1}^{n} (a_1 \otimes ... \otimes a_{i-1} \otimes a_{i+1} \otimes ... \otimes a_n) \otimes a_i$$ 
which we then extend by linearity (if this map looks backwards, recall that $\mathsf{MOD}_R$ is a co-differential category). 
\end{example}

It is important to note that this differential category structure on $\mathsf{MOD}^{op}_R$ can be generalized to the category of modules over any ring $R$. In fact, this example can be generalized further. Indeed, the free symmetric algebra construction on appropriate additive symmetric monoidal categories induces a differential category structure, such as on the category of sets and relations (see \cite{blute2006differential} for more details).

\begin{example}\label{ex2} \normalfont Convenient vector spaces provides another example (given by R. Blute, T. Ehrhard and C. Tasson) of a co-differential category with a monoidal algebra modality \cite{blute2010convenient}. 
\end{example} 

Interestingly, the free differential algebra construction -- which one might suppose would give rise rather naturally to a modality with a differential -- gives an example of an additive bialgebra modality which {\em does not\/} 
admit a deriving transformation:

\begin{example} \label{ex3}
\normalfont Let $R$ be a commutative ring. A (commutative) \textbf{differential algebra} (of weight $0$) over $R$ (see \cite{guo2008differential}) is a pair $(A, \mathsf{D})$ consisting of a commutative $R$-algebra $A$ and a linear map $\mathsf{D}: A \to A$ such that $\mathsf{D}$ satisfies the Leibniz rule:
$$\mathsf{D}(ab)=\mathsf{D}(a)b+a\mathsf{D}(b) \quad \forall a,b \in A$$
where the multiplication of the $R$-algebra $A$ has been written as juxtaposition. 
A \textbf{map of differential algebras} $f: (A, \mathsf{D}) \to (C, \mathsf{D}')$ is an $R$-algebra morphism $f: A \to C$ such that $f\mathsf{D}' = \mathsf{D}f$. 

The forgetful functor from the category of differential algebras over $R$, $\mathsf{CDA}_R$ to modules over $R$ has a left adjoint:
$$\xymatrixcolsep{2.5pc}\xymatrix{\mathsf{MOD}_{R}  \ar@<+1.1ex>[r]^-{\mathsf{DIFF}} & \mathsf{CDA}_{R} \ar@<+1ex>[l]_-{\bot}^-{U} }$$
which induces an algebra modality, $\mathsf{Diff}$, on the category of modules:  we shall now give an explicit description of this modality and, furthermore, show that it is an additive bialgebra modality which does {\em not\/} admit a deriving transformation.
  
  Let $M$ be an $R$-module, then the free commutative differential $R$-algebra $\mathsf{Diff}(M)$ is defined as follows: 
  \[\mathsf{Diff}(M)=\mathsf{Sym}(\bigoplus\limits_{n=0}^\infty M)\]
 where the unit and multiplication are just that of the symmetric algebra, $\mathsf{u}: R \to \mathsf{Sym}(\bigoplus\limits_{n=0}^\infty M)$ and $\nabla: \mathsf{Sym}(\bigoplus\limits_{n=0}^\infty M) \otimes \mathsf{Sym}(\bigoplus\limits_{n=0}^\infty M) \to \mathsf{Sym}(\bigoplus\limits_{n=0}^\infty M)$.  The differential is obtained by ``shifting'' the 
 infinite sum (on which the symmetric algebras is built) up one  
 $$\phi := \langle \iota_{n+1} \rangle^\infty_{n=0}: \bigoplus\limits_{n=0}^\infty M \to \bigoplus\limits_{n=0}^\infty M$$  
 and defining the map 
  $\mathsf{D}: \mathsf{Diff}(M) \to \mathsf{Diff}(M)$ as: 
  \[  \xymatrixcolsep{3pc}\xymatrix{\mathsf{Sym}(\bigoplus\limits_{n=0}^\infty M) \ar[r]^-{\mathsf{d}} & \mathsf{Sym}(\bigoplus\limits_{n=0}^\infty M) \otimes \bigoplus\limits_{n=0}^\infty M \ar[r]^-{1 \otimes \phi} &\mathsf{Sym}(\bigoplus\limits_{n=0}^\infty M) \otimes \bigoplus\limits_{n=0}^\infty M  \ar[r]^-{\mathsf{d}^\circ} & \mathsf{Sym}(\bigoplus\limits_{n=0}^\infty M)
  } \]
where $\mathsf{d}$ is the deriving transformation of the free symmetric algebra modality and $\mathsf{d}^\circ := (1 \otimes \eta)\nabla$ (where $\eta$ is the unit of the free symmetric algebra monad). By the associativity and unit laws of the multiplication, we note the following identities $\mathsf{d}^\circ$ satisfies:
\begin{eqnarray*}
(\mathsf{u} \otimes 1)\mathsf{d}^\circ &=& \eta \\
(\nabla \otimes 1) \mathsf{d}^\circ &=& (1 \otimes 1 \otimes \mathsf{d}^\circ) \nabla
\end{eqnarray*}
We then have:
\begin{eqnarray*}
\nabla \mathsf{D} &= & \nabla \mathsf{d} (1 \otimes \phi) \mathsf{d}^\circ \\
&=& (1 \otimes \mathsf{d})(\nabla\otimes 1)(1 \otimes \phi) \mathsf{d}^\circ + (\mathsf{d} \otimes 1)(1 \otimes \sigma)(\nabla \otimes 1)(1 \otimes \phi) \mathsf{d}^\circ \\
&=& (1 \otimes \mathsf{d})(1 \otimes 1 \otimes \phi)(\nabla \otimes 1) \mathsf{d}^\circ + (\mathsf{d} \otimes 1)(1 \otimes \phi \otimes 1)(1 \otimes \sigma)(\nabla \otimes 1)\mathsf{d}^\circ \\
&=& (1 \otimes \mathsf{d})(1 \otimes 1 \otimes \phi)(1 \otimes 1 \otimes \mathsf{d}^\circ) \nabla + (\mathsf{d} \otimes 1)(1 \otimes \phi \otimes 1)(\mathsf{d}^\circ \otimes 1) \nabla \\
&=& (1 \otimes \mathsf{D}) \nabla + (\mathsf{D} \otimes 1)\nabla
\end{eqnarray*}
showing that it is a differential on the algebra.

To show it is a monad set the unit to be the natural transformation $\alpha: M \to \mathsf{Diff}(M)$ defined as 
    \[ \alpha := \xymatrixcolsep{5pc}\xymatrix{M \ar[r]^-{\iota_0}  & \bigoplus\limits_{n=0}^\infty M \ar[r]^-{\eta} & \mathsf{Sym}(\bigoplus\limits_{n=0}^\infty M) } \]
where $\eta$ is the unit of the free symmetric algebra monad, and define the multiplication as the natural transformation $\nu: \mathsf{Diff}(\mathsf{Diff}(M)) \to \mathsf{Diff}(M)$ to be
  \[  \xymatrixcolsep{5pc}\xymatrix{\mathsf{Sym}(\bigoplus\limits_{n=0}^\infty (\mathsf{Sym}(\bigoplus\limits_{n=0}^\infty M))) \ar[r]^-{\mathsf{Sym}(\psi)} & \mathsf{Sym}(\mathsf{Sym}(\bigoplus\limits_{n=0}^\infty M)) \ar[r]^-{\mu}  & \mathsf{Sym}(\bigoplus\limits_{n=0}^\infty M)
  } \]
where $\mu$ is the multiplication of the free symmetric algebra monad and $\psi:= \langle \mathsf{D}^n \rangle_{n=0}^\infty: \bigoplus\limits_{n=0}^\infty A \to A$ is the map with $\iota_k \psi = \mathsf{D}^k$ where $A$ is a 
differential algebra: notice that $\psi$ is natural for differential algebra maps.  We then have:
  
\begin{lemma} $(\mathsf{Diff}, \nu, \alpha)$ is a monad on $\mathsf{MOD}_R$. 
\end{lemma}
\begin{proof} We must verify the three monad identities:
\begin{description}
\item{$\mathsf{Diff}(\nu) \nu=\nu \nu$:} Here we use that $\psi$ is natural with respect to differential algebra morphisms, that $\nu$ is a differential algebra morphism, the monad associativity of $\mu$, and the naturality of $\mu$: 
\begin{eqnarray*}
\mathsf{Diff}(\nu) \nu &= & \mathsf{Sym}(\bigoplus\limits^\infty_{n=0} \nu) \nu 
=  \mathsf{Sym}(\bigoplus\limits^\infty_{n=0} \nu) \mathsf{Sym}(\psi) \mu 
= \mathsf{Sym}((\bigoplus\limits^\infty_{n=0} \nu )\psi) \mu \\
& =  & \mathsf{Sym}(\psi \nu) \mu = \mathsf{Sym}(\psi) \mathsf{Sym}(\nu) \mu = \mathsf{Sym}(\psi) \mathsf{Sym}(\mathsf{Sym}(\psi) \mu) \mu \\
&= &  \mathsf{Sym}(\psi) \mathsf{Sym}(\mathsf{Sym}(\psi)) \mu \mu = \mathsf{Sym}(\psi) \mu \mathsf{Sym}(\psi)  \mu = \nu \nu
\end{eqnarray*}
\item{$\alpha \nu=1$:} Here we use the naturality of $\eta$, the definition of $\psi$, and the monad triangle identity of $\mu$ and $\eta$: 
$$\alpha \nu = \iota_0 \eta \mathsf{Sym}(\psi) \mu = \iota_0 \psi \eta  \mu = 1.$$
\item{$\mathsf{Diff}(\alpha) \nu=1$:} Here we have:
\begin{eqnarray*}
\mathsf{Diff}(\alpha) \nu & = & \mathsf{Sym}(\bigoplus\limits^\infty_{n=0} \alpha) \nu =  \mathsf{Sym}(\bigoplus\limits^\infty_{n=0} \alpha) \mathsf{Sym}(\psi) \mu \\
&=  & \mathsf{Sym}((\bigoplus\limits^\infty_{n=0} \alpha) \psi) \mu = \mathsf{Sym}(\eta) \mu =  1
\end{eqnarray*}
where the equality $(\bigoplus\limits^\infty_{n=0} \alpha) \psi = \eta$ is used in the penultimate step which we must now establish.
Notice, using the linear rule \textbf{[d.3]}, the following identity holds: 
\begin{align*}
\eta\mathsf{D} = \eta \mathsf{d} (1 \otimes \phi) \mathsf{d}^\circ = (\mathsf{u} \otimes 1) (1 \otimes \phi) \mathsf{d}^\circ  = \phi (\mathsf{u} \otimes 1)\mathsf{d}^\circ= \phi \eta
\end{align*}
This allows us to observe that 
$$ \iota_k  (\bigoplus\limits^\infty_{n=0} \alpha) \psi = \alpha \iota_k  \psi  =\alpha \mathsf{D}^k  = \iota_0 \eta \mathsf{D}^k = \iota_0 \phi^k \eta = \iota_k \eta$$
so that $(\bigoplus\limits^\infty_{n=0} \alpha) \psi= \eta$ as desired.
\end{description}
\end{proof} 

\begin{proposition} 
$(\mathsf{Diff}, \nu, \alpha, \nabla, \mathsf{u}, \Delta, \mathsf{e})$ is an additive bialgebra modality. 
\end{proposition}

\begin{proof} First observe that $\nu$ is an algebra morphism as both $\mathsf{Sym}(\psi)$ and $\mu$ are algebra morphisms. 

To establish that the modality is an additive  bialgebra modality we exhibit the Seely isomorphisms. Since $\mathsf{Sym}$ has Seely isomorphisms it follows that: 

$$\mathsf{Diff}(0) =  \mathsf{Sym}(\bigoplus\limits^\infty_{n=0} 0) = \mathsf{Sym}(0) \cong R$$
 $$\mathsf{Sym}(\bigoplus\limits^\infty_{n=0}(M \oplus N))\cong \mathsf{Sym}(\bigoplus\limits^\infty_{n=0}M \oplus \bigoplus\limits^\infty_{n=0}N)  \cong \mathsf{Sym}(\bigoplus\limits^\infty_{n=0}M) \otimes \mathsf{Sym}(\bigoplus\limits^\infty_{n=0}N)$$
\end{proof}

We now set about proving that this modality does not admit a deriving transformation.  This is accomplished by proving that if there is a deriving transformation, then the ring $R$ over which the modules are taken 
must be trivial: that is in $R$ we must have $1=0$.  More precisely we prove that $1_{M \otimes M} = - \sigma: M \otimes M\to M \otimes M$, however, by substituting $R+R$ for $M$ this gives the matrix equality 
$$\left( \begin{array}[c]{cccc} 1 & 0 & 0 & 0\\ 0 &1 & 0 & 0 \\ 0 & 0 & 1 &  0 \\  0 & 0 & 0 & 1 \end{array} \right) = \left( \begin{array}[c]{cccc} -1 & 0 & 0 & 0 \\ 0 & 0 & -1 & 0 \\ 0 & -1 & 0 &  0 \\  0 & 0 & 0 & -1 \end{array} \right)$$
but this then immediately gives.
\begin{eqnarray*}
0 & =  & \left( \begin{array}[c]{cccc} 0 & 1 & 0 & 0 \end{array} \right) \left( \begin{array}[c]{c} 0 \\ 0 \\ -1 \\ 0 \end{array} \right) 
              = \left( \begin{array}[c]{cccc} 0 & 1 & 0 & 0\end{array} \right) \left( \begin{array}[c]{cccc} 1 & 0 & 0 & 0\\ 0 &1 & 0 & 0 \\ 0 & 0 & 1 &  0 \\  0 & 0 & 0 & 1 \end{array} \right) \left( \begin{array}[c]{c} 0 \\ 0 \\ -1 \\ 0\end{array} \right) \\
& = & \left( \begin{array}[c]{cccc} 0 & 1 & 0 & 0\end{array} \right) \left( \begin{array}[c]{cccc}  -1 & 0 & 0 & 0 \\ 0 & 0 & -1 & 0 \\ 0 & -1 & 0 &  0 \\  0 & 0 & 0 & -1 \end{array} \right) \left( \begin{array}[c]{c} 0 \\ 0 \\ -1 \\ 0\end{array} \right) = 1
\end{eqnarray*}

\begin{theorem} 
For any category of modules over a non-trivial commutative ring the free differential algebra modality, $\mathsf{Diff}$, does not admit a deriving transformation.
\end{theorem} 

\begin{proof}
Suppose then that there is a natural transformation $\mathsf{b}: \mathsf{Diff}(M) \to \mathsf{Diff}(M) \otimes M$ which is a deriving transformation. Then for each $R$-module $M$, the Leibniz rule 
implies that $\mathsf{b}$ is a $\mathsf{Sym}(\bigoplus\limits^\infty_{n=0}M)$-derivation, and therefore by the universality of the deriving transformation of $\mathsf{Sym}$ (see \cite{blute2011kahler,blute2015derivations} for more details), there exists a unique $\mathsf{Sym}(\bigoplus\limits^\infty_{n=0}M)$-module homomorphism $f^\sharp$ making 
  \[  \xymatrixcolsep{5pc}\xymatrix{\mathsf{Sym}(\bigoplus\limits^\infty_{n=0}M) \ar[dr]_-{\mathsf{b}} \ar[r]^-{\mathsf{d}} & \mathsf{Sym}(\bigoplus\limits^\infty_{n=0}M) \otimes\bigoplus\limits^\infty_{n=0}M \ar@{..>}[d]^-{f^\sharp} \\
  & \mathsf{Sym}(\bigoplus\limits^\infty_{n=0}M) \otimes M
  } \]
 commute. However, as $f^\sharp$ is a morphism between free modules
it is determined by a map 
$$f: \bigoplus\limits^\infty_{n=0}M \to \mathsf{Sym}(\bigoplus\limits^\infty_{n=0}M) \otimes M~~\mbox{where}~~f^\sharp = (1 \otimes f)(\nabla \otimes 1).$$
Now both $\mathsf{d}$ and $\mathsf{b}$ are determined by derelictions, respectively $\widehat{\mathsf{d}}:= \mathsf{d} (e \otimes 1)$ and $\widehat{\mathsf{b}} = \mathsf{b}(e \otimes 1)$, as they are on additive algebra modalities.
Thus, setting $\widehat{f} = f(e \otimes 1)$ we have:
$$\widehat{\mathsf{b}} = \mathsf{b} (e \otimes 1) = \mathsf{d} f^\sharp (e \otimes 1) = \mathsf{d}  (1 \otimes f)(\nabla \otimes 1) (e \otimes 1) = \mathsf{d}  (e \otimes f (e \otimes 1))
      = \widehat{\mathsf{d}}f (e \otimes 1) = \widehat{\mathsf{d}} \widehat{f}.$$
But then $f^\sharp = 1 \otimes \widehat{f}$ as $\mathsf{d} (1 \otimes \widehat{f}) = \Delta (1 \otimes (\widehat{\mathsf{d}}~\widehat{f})) =  \Delta (1 \otimes \widehat{\mathsf{b}}) = \mathsf{b}$.  Thus we have now shown that 
\[  \xymatrixcolsep{5pc}\xymatrix{\mathsf{Sym}(\bigoplus\limits^\infty_{n=0}M) \ar[dr]_-{\mathsf{b}} \ar[r]^-{\mathsf{d}} & 
\mathsf{Sym}(\bigoplus\limits^\infty_{n=0}M) \otimes\bigoplus\limits^\infty_{n=0}M \ar@{..>}[d]^-{1 \otimes \widehat{f}} \\
  & \mathsf{Sym}(\bigoplus\limits^\infty_{n=0}M) \otimes M
  } \]
  commutes.  Furthermore, $\widehat{f}$ is a natural transformation as if $g:M \to N$ then, using that $\mathsf{d}$ is natural and $\mathsf{b}$ is assumed to be a natural, both $\left(\bigoplus\limits^\infty_{n=0}g \right) \widehat{f}$ and $\widehat{f} g$ provide the unique mediating maps between $\mathsf{d}$ and $\mathsf{Diff}(g) \mathsf{b}$.

Consider the map $\pi_1: \bigoplus\limits^\infty_{n=0} A_n \to A_1$ defined as $\pi_1 = \langle \delta_{k, 1} \rangle$, that is, the unique map which makes the following diagram commute for each injection map: 
  \[  \xymatrixcolsep{5pc}\xymatrix{A_k \ar[r]^-{\iota_k} \ar[dr]_-{\delta_{k,1}} &\bigoplus\limits_{n=0}^\infty A_k \ar[d]^-{\pi_1} \\
  & A_1 
  } \]
where $\delta_{1,1}= 1$ and $\delta_{k,1}=0$ for $k \neq 1$. Because $\mathsf{b}$ satisfies the chain rule, $\nu \mathsf{b} = \mathsf{b} (\nu \otimes \mathsf{b}) (\nabla \otimes 1)$, 
the following equality of natural transformations, obtained by sandwiching the chain rule between the same maps, must also hold:
$$(\alpha \otimes \alpha) \nabla \iota_1 \eta \nu \mathsf{b} ((\varepsilon \pi_1) \otimes 1)= (\alpha \otimes \alpha) \nabla \iota_1 \eta \mathsf{b}(\nu \otimes \mathsf{b})(\nabla \otimes 1) (\varepsilon \otimes 1)(\pi_1 \otimes 1)$$ 
However, we will show that this forces $1 = -\sigma: M \otimes M \to M \otimes M$ which forces the module category to be trivial. Preliminary to this we note the following useful identities: 
\begin{eqnarray*}
\alpha & = & \iota_0 \eta \\
\alpha \mathsf{D} & = & \iota_0 \eta \mathsf{D} = \iota_0 \phi \eta = \iota_1 \eta \\
\iota_0 \widehat{f} & =  & \iota_0 \eta \widehat{\mathsf{d}} \widehat{f} = \alpha \widehat{\mathsf{b}} = 1 \\
(\eta \otimes \eta) \nabla \mathsf{d} &=  & (\eta \otimes \eta) (1 \otimes \mathsf{d})(\nabla \otimes 1) + (\eta \otimes \eta) (\mathsf{d} \otimes 1)(1 \otimes \sigma)(\nabla \otimes 1) \\
&= & (\eta \otimes \mathsf{u} \otimes 1) (\nabla \otimes 1) + (\mathsf{u} \otimes 1 \otimes \eta) (1 \otimes \sigma)(\nabla \otimes 1) \\
&= & \eta \otimes 1 + \sigma(\eta \otimes 1)
\end{eqnarray*}
Now explicitly calculating out the first map above we have: 
\begin{eqnarray*}
\lefteqn{ (\alpha \otimes \alpha) \nabla \iota_1 \eta \nu \mathsf{b} (\varepsilon \pi_1 \otimes 1)} \\
&= &(\alpha \otimes \alpha) \nabla \iota_1 \eta  \mathsf{Sym}(\psi) \mu \mathsf{d} (\varepsilon \pi_1 \otimes f) \\
&= & (\alpha \otimes \alpha)  \nabla \iota_1 \psi \eta \mu \mathsf{d} (\varepsilon\pi_1 \otimes f) \\
&= & (\alpha \otimes \alpha)  \nabla \mathsf{D} \mathsf{d} (\varepsilon \pi_1 \otimes f) \\
&= & ( (\alpha \otimes \alpha \mathsf{D}) \nabla + (\alpha \mathsf{D} \otimes \alpha) \nabla))\mathsf{d} (\varepsilon \pi_1 \otimes f) \\
&= & ((\iota_0 \eta \otimes \iota_1\eta) \nabla  + (\iota_1\eta  \otimes \iota_0\eta)  \nabla) \mathsf{d} (\varepsilon \pi_1 \otimes f) \\
&= & (\iota_0 \otimes \iota_1 + \iota_1 \otimes \iota_0)(\eta \otimes \eta) \nabla \mathsf{d} (\varepsilon \pi_1 \otimes f) \\
&=  & (\iota_0 \otimes \iota_1 + \iota_1 \otimes \iota_0)((\eta \otimes 1) + \sigma (\eta \otimes 1)) (\varepsilon \pi_1 \otimes f) \\
&= &  (\iota_0 \otimes \iota_1 + \iota_1 \otimes \iota_0)((\pi_1 \otimes f) + \sigma (\pi_1 \otimes f)) \\
&= & (\iota_0 \otimes \iota_1)(\pi_1 \otimes f) + (\iota_1 \otimes \iota_0)(\pi_1 \otimes f) \\ & & ~~+ \sigma (\iota_1 \otimes \iota_0)(\pi_1 \otimes f) + \sigma (\iota_0 \otimes \iota_1)(\pi_1 \otimes f)\\
&= & 0 + 1 \otimes 1 + \sigma + 0 \\
&= & 1 \otimes 1 + \sigma 
\end{eqnarray*}
while for the second map: 
\begin{eqnarray*}
\lefteqn{ (\alpha \otimes \alpha) \nabla \iota_1 \eta \mathsf{b}(\nu \otimes \mathsf{b})(\nabla \otimes 1) (\varepsilon \pi_1 \otimes 1) } \\
& = & (\alpha \otimes \alpha) \nabla \iota_1 \eta \mathsf{d}(\nu \otimes \widehat{f}\mathsf{b})(\nabla \otimes 1) (\varepsilon \pi_1 \otimes 1) \\
& = & (\alpha \otimes \alpha) \nabla \iota_1 (u \otimes 1)(\nu \otimes \widehat{f}\mathsf{b})(\nabla \otimes 1) (\varepsilon \pi_1 \otimes 1) \\
& = & (\alpha \otimes \alpha) \nabla \iota_1 (u\nu \otimes \widehat{f}\mathsf{b})(\nabla \otimes 1) (\varepsilon \pi_1 \otimes 1) \\
& = & (\alpha \otimes \alpha) \nabla \iota_1 (u \otimes \widehat{f} \mathsf{b})(\nabla \otimes 1) (\varepsilon \pi_1 \otimes 1) \\
& = & (\alpha \otimes \alpha) \nabla \iota_1 \widehat{f} \mathsf{b} ((u \otimes 1)\nabla \otimes 1) (\varepsilon \pi_1 \otimes 1) \\
& = & (\alpha \otimes \alpha) \nabla \iota_1 \widehat{f} \mathsf{b} (\varepsilon \pi_1 \otimes 1) \\
& = & (\alpha \otimes \alpha) \nabla \iota_1  \left(\bigoplus\limits^\infty_{n=0} \mathsf{b} (\varepsilon \pi_1 \otimes 1) \right) \widehat{f} \\
& = & (\alpha \otimes \alpha) \nabla \mathsf{b} (\varepsilon \pi_1 \otimes 1) \iota_1 \widehat{f} \\
& = & (\iota_0 \otimes \iota_0)(\eta \otimes \eta) \nabla \mathsf{d} (\varepsilon \pi_1 \otimes \widehat{f}) \iota_1 \widehat{f}  \\
& = & (\iota_0 \otimes \iota_0)(\eta \otimes 1 + \sigma( \eta \otimes 1)) (\varepsilon \pi_1 \otimes \widehat{f}f) \iota_1 \widehat{f}  \\
&= & (\iota_0 \otimes \iota_0)(\pi_1 \otimes f) \iota_1  f + \sigma (\iota_0 \otimes \iota_0) (\pi_1 \otimes \widehat{f})\iota_1  \widehat{f}  \\
&= & 0 + 0 = 0 
\end{eqnarray*}
Therefore, if $\mathsf{b}$ satisfies the chain rule, this would imply that $1_{M \otimes M} = - \sigma$ for every $R$-module $M$ which only happens when $R$ is the trivial ring as discussed above. Therefore, $\mathsf{Diff}$ does not have a deriving transformation when the ring is non-trivial. 
\end{proof}
\end{example}

By way of contrast, Rota-Baxter algebras -- an algebraic abstraction of integration -- whose algebra modality is not additive, always give a differential category: 

\begin{example}\label{ex4} \normalfont Let $R$ be a commutative ring. A (commutative) \textbf{Rota-Baxter algebra (of weight $0$)} \cite{guo2012introduction} over $R$ is a pair $(A, \mathsf{P})$ consisting of a commutative $R$-algebra $A$ and an $R$-linear map $\mathsf{P}: A \to A$ such that $\mathsf{P}$ satisfies the  Rota-Baxter equation, that is, the following equality holds:
$$\mathsf{P}(a)\mathsf{P}(b)=\mathsf{P}(a\mathsf{P}(b))+\mathsf{P}(\mathsf{P}(a)b) \quad \forall ab, \in A$$
 The map $\mathsf{P}$ is called a \textbf{Rota-Baxter operator} (we refer the reader to \cite{guo2012introduction} for more details on Rota-Baxter algebras). It turns out that there is a left adjoint to the forgetful functor between the category of Rota-Baxter algebras, $\mathsf{CRBA}_{R}$, and the category of commutative algebras over $R$, $\mathsf{CALG}_{R}$. We quickly review the construction of the free Rota-Baxter algebra over an algebra (for more details see chapter $3$ of \cite{guo2012introduction}): let $M$ be an $R$-module and consider the \textbf{shuffle algebra}, $\mathsf{Sh}(M)$, over $M$ which is defined as follows:
$\mathsf{Sh}(M)=R \oplus M \oplus (M \otimes M) \oplus (M \otimes M \otimes M) \oplus \hdots = \bigoplus_{n \in \mathbb{N}} M^{\otimes^n}$
where $M^{\otimes^0}=R$ and where the multiplication $\shuffle$, called the \textbf{shuffle product} \cite{guo2012introduction}, is defined inductively on pure tensors $w= a \otimes w^\prime$ and $v= b \otimes v^\prime$ as follows:
$$w \shuffle v= a \otimes (w^\prime \shuffle v) + b \otimes (w \shuffle v^\prime)$$
which we then extend by linearity (notice that the unit for the shuffle product is $1_R$). Denote the multiplication and unit maps of the shuffle algebra by $\blacktriangledown: \mathsf{Sh}(M) \otimes \mathsf{Sh}(M) \to \mathsf{Sh}(M)$ and $v: R \to \mathsf{Sh}(M)$ respectively. The free commutative Rota-Baxter over a commutative $R$-algebra $A$, $\mathsf{RB}(A)$, is then the tensor product of shuffle algebra and $A$ itself: $\mathsf{RB}(A)=\mathsf{Sh}(A) \otimes A$. The Rota-Baxter operator $\mathsf{P}: \mathsf{RB}(A) \to \mathsf{RB}(A)$ is defined on pure tensors as follows: 
$$w \otimes b = \begin{cases} (w\cdot b) \otimes 1_A & \text{ if } w \in R \\
(w \otimes b) \otimes 1_A & \text{ otherwise}
\end{cases}$$
which we then extend by linearity. The induced functor $\mathsf{RB}: \mathsf{CALG}_{R} \to \mathsf{CRBA}_{R}$ is the left adjoint to the forgetful functor $U:\mathsf{CRBA}_{R}  \to \mathsf{CALG}_{R}$:
$$\xymatrixcolsep{2.5pc}\xymatrix{\mathsf{CALG}_{R} \ar@<+1.1ex>[r]^-{\mathsf{RB}} & \mathsf{CRBA}_{R} \ar@<+1ex>[l]_-{\bot}^-{U}
  }$$
(for more details on this adjunction and monad see \cite{zhang2016monads}) where for an algebra $A$, the unit of the adjunction is defined as $v \otimes 1: A \to \mathsf{Sh}(A) \otimes A$, while for a Rota-Baxter algebra $(B, \mathsf{Q})$, the counit $\omega: \mathsf{Sh}(B) \otimes B \to B$ is defined on pure tensors as 
$$\omega_B((b_1 \otimes ... \otimes b_n) \otimes b )= \mathsf{Q}(\hdots \mathsf{Q}(\mathsf{Q}(b_1)b_2) \hdots b_n)b$$
which we then extend by linearity. 
To obtain an algebra modality on $\mathsf{MOD}_R$, we compose the free Rota-Baxter algebra adjunction and the free symmetric algebra adjunction: 
$$\xymatrixcolsep{2.5pc}\xymatrix{\mathsf{MOD}_R \ar@<+1.1ex>[r]^-{\mathsf{Sym}} & \mathsf{CALG}_{R} \ar@<+1ex>[l]_-{\bot}^-{U} \ar@<+1.1ex>[r]^-{\mathsf{RB}} & \mathsf{CRBA}_{R} \ar@<+1ex>[l]_-{\bot}^-{U}
  }$$
 The monad induced by the resulting adjunction between $\mathsf{MOD}_R$ and $ \mathsf{CRBA}_{R}$ is clearly an algebra modality by construction again. After some simplifications, the unit and multiplication of the monad are represented in string diagrams respectively as:
 $$   \begin{array}[c]{c}\resizebox{!}{1.5cm}
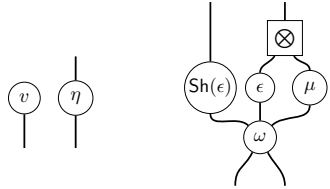
{%
\begin{tikzpicture}
	\begin{pgfonlayer}{nodelayer}
		\node [style=circle] (0) at (1.25, 0.5) {};
		\node [style={circle, draw}] (1) at (1.25, -0.5) {$\eta$};
		\node [style={circle, draw}] (2) at (0.25, -0.5) {$v$};
		\node [style=port] (3) at (0.25, -1.5) {};
		\node [style=port] (4) at (1.25, -1.5) {};
	\end{pgfonlayer}
	\begin{pgfonlayer}{edgelayer}
		\draw [style=wire] (1) to (0);
		\draw [style=wire] (3) to (2);
		\draw [style=wire] (4) to (1);
	\end{pgfonlayer}
\end{tikzpicture}
  }% 
\end{array}~~~~~~ \begin{array}[c]{c}\resizebox{!}{2.5cm}{%
 \begin{tikzpicture}
	\begin{pgfonlayer}{nodelayer}
		\node [style=port] (0) at (-0.5, -1.75) {};
		\node [style={circle, draw}] (1) at (0, -0.75) {$\omega$};
		\node [style={circle, draw}] (2) at (1, 0.25) {$\mu$};
		\node [style={regular polygon,regular polygon sides=4, draw, inner sep=1pt,minimum size=1pt}] (3) at (0.5, 1.25) {$\bigotimes$};
		\node [style=port] (4) at (-1, 2) {};
		\node [style=port] (5) at (0.5, -1.75) {};
		\node [style={circle, draw}] (6) at (0, 0.25) {$\epsilon$};
		\node [style=port] (7) at (0.5, 2) {};
		\node [style={circle, draw, inner sep=1pt,minimum size=1pt}] (8) at (-1, 0.25) {$\mathsf{Sh}(\epsilon)$};
	\end{pgfonlayer}
	\begin{pgfonlayer}{edgelayer}
		\draw [style=wire] (4) to (8);
		\draw [style=wire, in=90, out=-120, looseness=1.00] (3) to (6);
		\draw [style=wire, in=90, out=-63, looseness=1.00] (3) to (2);
		\draw [style=wire] (7) to (3);
		\draw [style=wire] (6) to (1);
		\draw [style=wire, in=45, out=-90, looseness=1.00] (2) to (1);
		\draw [style=wire, in=135, out=-90, looseness=1.00] (8) to (1);
		\draw [style=wire, in=90, out=-117, looseness=1.25] (1) to (0);
		\draw [style=wire, in=90, out=-63, looseness=1.00] (1) to (5);
	\end{pgfonlayer}
\end{tikzpicture}
  }% 
\end{array}$$
While the unit and multiplication of the algebra structure are represented in string diagrams respectively as: 
$$  \begin{array}[c]{c}\resizebox{!}{1cm}
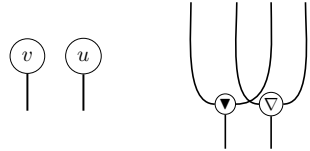
{%
\begin{tikzpicture}
	\begin{pgfonlayer}{nodelayer}
		\node [style={circle, draw}] (1) at (1.25, -0.5) {$u$};
		\node [style={circle, draw}] (2) at (0.25, -0.5) {$v$};
		\node [style=port] (3) at (0.25, -1.5) {};
		\node [style=port] (4) at (1.25, -1.5) {};
	\end{pgfonlayer}
	\begin{pgfonlayer}{edgelayer}
		\draw [style=wire] (3) to (2);
		\draw [style=wire] (4) to (1);
	\end{pgfonlayer}
\end{tikzpicture}
  }% 
\end{array}~~~~~   \begin{array}[c]{c}\resizebox{!}{2cm}{%
\begin{tikzpicture}
	\begin{pgfonlayer}{nodelayer}
		\node [style=port] (0) at (-0.25, -0.25) {};
		\node [style=duplicate] (1) at (-0.25, 0.75) {$\blacktriangledown$};
		\node [style=duplicate] (2) at (0.75, 0.75) {$\nabla$};
		\node [style=port] (3) at (0.75, -0.25) {};
		\node [style=port] (4) at (0, 3) {};
		\node [style=port] (5) at (0.75, 3) {};
		\node [style=port] (6) at (1.5, 3) {};
		\node [style=port] (7) at (-1, 3) {};
	\end{pgfonlayer}
	\begin{pgfonlayer}{edgelayer}
		\draw [style=wire] (1) to (0);
		\draw [style=wire] (2) to (3);
		\draw [style=wire, in=0, out=-90, looseness=1.00] (5) to (1);
		\draw [style=wire, in=0, out=-90, looseness=0.75] (6) to (2);
		\draw [style=wire, in=180, out=-90, looseness=0.75] (4) to (2);
		\draw [style=wire, in=180, out=-90, looseness=0.75] (7) to (1);
	\end{pgfonlayer}
\end{tikzpicture}
  }% 
\end{array}$$ 
However this algebra modality, $\mathsf{RB}(\mathsf{Sym}(M))$, does not have the Seely isomorphism as $\mathsf{Sh}$ is not a strong monoidal functor (i.e. $\mathsf{Sh}(A \otimes B) \ncong \mathsf{Sh}(A) \otimes \mathsf{Sh}(B)$):
\begin{align*}
\mathsf{RB}(\mathsf{Sym}(M \oplus N)) &\cong \mathsf{RB}(\mathsf{Sym}(M) \otimes \mathsf{Sym}(N))\\
& = \mathsf{Sh}(\mathsf{Sym}(M) \otimes \mathsf{Sym}(N)) \otimes \mathsf{Sym}(M) \otimes \mathsf{Sym}(N)\\
&\ncong \mathsf{Sh}(\mathsf{Sym}(M)) \otimes \mathsf{Sh}(\mathsf{Sym}(N)) \otimes \mathsf{Sym}(M) \otimes \mathsf{Sym}(N)\\
&\cong  \mathsf{Sh}(\mathsf{Sym}(M)) \otimes \mathsf{Sym}(M) \otimes \mathsf{Sh}(\mathsf{Sym}(N)) \otimes \mathsf{Sym}(N)\\
&= \mathsf{RB}(\mathsf{Sym}(M)) \otimes \mathsf{RB}(\mathsf{Sym}(N))
\end{align*}
Therefore, this algebra modality is not a bialgebra modality or a comonoidal algebra modality. We should mention that while it is true that the shuffle algebra is a bialgebra, its comultiplication is not cocommutative \cite{guo2012introduction}. However, we may still use the free Rota-Baxter adjunction to obtain a differential category structure on $\mathsf{MOD}_R$. The deriving transformation is defined as 
\begin{align*}
\begin{array}[c]{c}
 \xymatrixcolsep{5pc}\xymatrix{ \mathsf{Sh}(\mathsf{Sym}(M)) \otimes \mathsf{Sym}(M) \ar[r]^-{1 \otimes \mathsf{d}} & \mathsf{Sh}(\mathsf{Sym}(M)) \otimes \mathsf{Sym}(M) \otimes M} 
   \end{array} && 
   \begin{array}[c]{c}\resizebox{!}{1.75cm}{%
\begin{tikzpicture}
	\begin{pgfonlayer}{nodelayer}
		\node [style=port] (0) at (1.25, 3) {};
		\node [style=differential] (1) at (1.25, 2) {{\bf =\!=\!=}};
		\node [style=port] (2) at (2, 0.75) {};
		\node [style=port] (3) at (0.5, 0.75) {};
		\node [style=port] (4) at (0, 0.75) {};
		\node [style=port] (5) at (0, 3) {};
	\end{pgfonlayer}
	\begin{pgfonlayer}{edgelayer}
		\draw [style=wire, bend left, looseness=1.00] (1) to (2);
		\draw [style=wire] (0) to (1);
		\draw [style=wire, bend right, looseness=1.00] (1) to (3);
		\draw [style=wire] (4) to (5);
	\end{pgfonlayer}
\end{tikzpicture}
  }% 
   \end{array}
\end{align*}
where recall that $\mathsf{d}$ is the deriving transformation of $\mathsf{Sym}$ (if this looks upside-down, recall that we are working in a co-differential category). It may seem trivial that this is a deriving transformation, but in fact proving the chain rule is quite non-trivial! We will need the following lemma to prove the chain rule:
 \begin{lemma}\label{RBderiving} Let $R$ be a commutative ring. 
\begin{enumerate}[{\em (i)}]
\item For every commutative $R$-algebra $A$, $(1 \otimes \lozenge)\omega=(\omega \otimes 1 \otimes 1)\lozenge$
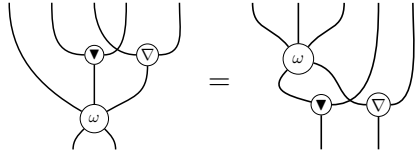
$$\begin{array}[c]{c}\resizebox{!}{2cm}{%
\begin{tikzpicture}
	\begin{pgfonlayer}{nodelayer}
		\node [style=port] (0) at (-1, 3) {};
		\node [style=port] (1) at (-2, 3) {};
		\node [style=port] (2) at (0, 3) {};
		\node [style=port] (3) at (0.75, 3) {};
		\node [style=port] (4) at (2, 3) {};
		\node [style=duplicate] (5) at (0, 1.75) {$\blacktriangledown$};
		\node [style={circle, draw}] (6) at (0, 0.25) {$\omega$};
		\node [style=port] (7) at (0.5, -0.5) {};
		\node [style=port] (8) at (-0.5, -0.5) {};
		\node [style=duplicate] (9) at (1.25, 1.75) {$\nabla$};
	\end{pgfonlayer}
	\begin{pgfonlayer}{edgelayer}
		\draw [style=wire, in=90, out=-135, looseness=1.00] (6) to (8);
		\draw [style=wire] (5) to (6);
		\draw [style=wire, in=90, out=-45, looseness=1.00] (6) to (7);
		\draw [style=wire, in=30, out=-90, looseness=1.25] (9) to (6);
		\draw [style=wire, in=150, out=-90, looseness=1.00] (1) to (6);
		\draw [style=wire, in=180, out=-90, looseness=1.25] (0) to (5);
		\draw [style=wire, in=0, out=-90, looseness=1.00] (3) to (5);
		\draw [style=wire, in=180, out=-90, looseness=1.00] (2) to (9);
		\draw [style=wire, in=0, out=-90, looseness=1.25] (4) to (9);
	\end{pgfonlayer}
\end{tikzpicture}}
   \end{array}=
   \begin{array}[c]{c}\resizebox{!}{2cm}{%
\begin{tikzpicture}
	\begin{pgfonlayer}{nodelayer}
		\node [style=port] (0) at (-0.5, -0.25) {};
		\node [style=duplicate] (1) at (-0.5, 0.75) {$\blacktriangledown$};
		\node [style=duplicate] (2) at (0.75, 0.75) {$\nabla$};
		\node [style=port] (3) at (0.75, -0.25) {};
		\node [style={circle, draw}] (4) at (-1, 1.75) {$\omega$};
		\node [style=port] (5) at (-1, 3) {};
		\node [style=port] (6) at (-2, 3) {};
		\node [style=port] (7) at (0, 3) {};
		\node [style=port] (8) at (0.75, 3) {};
		\node [style=port] (9) at (1.5, 3) {};
	\end{pgfonlayer}
	\begin{pgfonlayer}{edgelayer}
		\draw [style=wire, in=165, out=-135, looseness=1.75] (4) to (1);
		\draw [style=wire, in=180, out=-15, looseness=1.25] (4) to (2);
		\draw [style=wire] (1) to (0);
		\draw [style=wire] (2) to (3);
		\draw [style=wire] (5) to (4);
		\draw [style=wire, in=129, out=-90, looseness=1.25] (6) to (4);
		\draw [style=wire, in=51, out=-90, looseness=1.25] (7) to (4);
		\draw [style=wire, in=0, out=-90, looseness=1.00] (8) to (1);
		\draw [style=wire, in=0, out=-90, looseness=0.75] (9) to (2);
	\end{pgfonlayer}
\end{tikzpicture}}
   \end{array}$$
   where $\lozenge$ is the multiplication on $\mathsf{RB}(A)$. 
   \item For every $R$-module $M$, $\omega(1 \otimes \mathsf{d})=(1 \otimes 1 \otimes \mathsf{d})(\omega \otimes 1)$
$$\begin{array}[c]{c}\resizebox{!}{2cm}
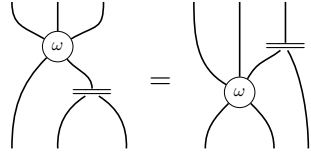
{%
\begin{tikzpicture}
	\begin{pgfonlayer}{nodelayer}
		\node [style=port] (0) at (-1, 3) {};
		\node [style=port] (1) at (-2, 3) {};
		\node [style=port] (2) at (0, 3) {};
		\node [style=port] (3) at (-2, -0.25) {};
		\node [style=codifferential] (4) at (-0.25, 1) {{\bf =\!=\!=}};
		\node [style={circle, draw}] (5) at (-1, 2) {$\omega$};
		\node [style=port] (6) at (0.5, -0.25) {};
		\node [style=port] (7) at (-1, -0.25) {};
	\end{pgfonlayer}
	\begin{pgfonlayer}{edgelayer}
		\draw [style=wire, bend left, looseness=1.00] (4) to (6);
		\draw [style=wire, bend right, looseness=1.00] (4) to (7);
		\draw [style=wire, in=90, out=-135, looseness=1.00] (5) to (3);
		\draw [style=wire, in=90, out=-53, looseness=1.00] (5) to (4);
		\draw [style=wire] (0) to (5);
		\draw [style=wire, in=30, out=-90, looseness=1.50] (2) to (5);
		\draw [style=wire, in=150, out=-90, looseness=1.25] (1) to (5);
	\end{pgfonlayer}
\end{tikzpicture}}
   \end{array}=
   \begin{array}[c]{c}\resizebox{!}{2cm}{%
\begin{tikzpicture}
	\begin{pgfonlayer}{nodelayer}
		\node [style=port] (0) at (-1, 3) {};
		\node [style=port] (1) at (-2, 3) {};
		\node [style=port] (2) at (0, 3) {};
		\node [style=codifferential] (3) at (0, 2) {{\bf =\!=\!=}};
		\node [style=port] (4) at (-0.25, -0.25) {};
		\node [style={circle, draw}] (5) at (-1, 1) {$\omega$};
		\node [style=port] (6) at (-1.75, -0.25) {};
		\node [style=port] (7) at (0.5, -0.25) {};
	\end{pgfonlayer}
	\begin{pgfonlayer}{edgelayer}
		\draw [style=wire, in=90, out=-60, looseness=0.75] (3) to (7);
		\draw [style=wire, in=90, out=-135, looseness=1.00] (5) to (6);
		\draw [style=wire, in=60, out=-143, looseness=1.25] (3) to (5);
		\draw [style=wire, in=90, out=-45, looseness=1.00] (5) to (4);
		\draw [style=wire] (0) to (5);
		\draw [style=wire] (2) to (3);
		\draw [style=wire, in=150, out=-90, looseness=1.00] (1) to (5);
	\end{pgfonlayer}
\end{tikzpicture}}
   \end{array}$$
\end{enumerate}
\end{lemma}
\begin{proof} $i)$: Let $A$ be a commutative $R$-algebra. It suffices to prove this equality on pure tensors. Consider the following pure tensor of $\mathsf{Sh}(\mathsf{Sh}(A) \otimes A) \otimes \mathsf{Sh}(A) \otimes A$: 
$$([A_1 \otimes a_1] \otimes \hdots \otimes [A_n \otimes a_n]) \otimes A \otimes a \in \mathsf{Sh}(\mathsf{Sh}(A) \otimes A) \otimes \mathsf{Sh}(A) \otimes A$$
where $A, A_1 \hdots A_n \in \mathsf{Sh}(A)$ and $a, a_1 \hdots, a_n \in A$. By definition of $\omega_A$, we obtain that:
\begin{align*}
&\omega_A(([A_1 \otimes a_1] \otimes \hdots \otimes [A_n \otimes a_n]) \otimes A \otimes a)\\
&= \mathsf{P}(\mathsf{P}(\hdots \mathsf{P}(\mathsf{P}(A_1 \otimes a_1) \lozenge (A_2 \otimes a_2))\hdots \lozenge (A_n \otimes a_n)))\lozenge (A \otimes a) \\
&= (\mathsf{P}(\hdots \mathsf{P}(\mathsf{P}(A_1 \otimes a_1) \lozenge (A_2 \otimes a_2))\hdots \lozenge (A_n \otimes a_n))) \otimes 1) \lozenge (A \otimes a) \\
&= (\mathsf{P}(\hdots \mathsf{P}(\mathsf{P}(A_1 \otimes a_1) \lozenge (A_2 \otimes a_2))\hdots \lozenge (A_n \otimes a_n)) \shuffle A) \otimes a \end{align*}
Notice that $a$ is unaffected by $\omega$. Now let $B \otimes b \in \mathsf{Sh}(A) \otimes A$. Then we have the following equality by associativity of the shuffle product:
\begin{align*}
&=\omega_A(([A_1 \otimes a_1] \otimes \hdots \otimes [A_n \otimes a_n]) \otimes ((A \otimes a) \lozenge (B \otimes b))) \\
&= \omega_A(([A_1 \otimes a_1] \otimes \hdots \otimes [A_n \otimes a_n]) \otimes ((A \shuffle B) \otimes (ab))) \\
&= (\mathsf{P}(\hdots \mathsf{P}(\mathsf{P}(A_1 \otimes a_1) \lozenge (A_2 \otimes a_2))\hdots \lozenge (A_n \otimes a_n)) \shuffle (A\shuffle B)) \otimes (ab) \\
&= ((\mathsf{P}(\hdots \mathsf{P}(\mathsf{P}(A_1 \otimes a_1) \lozenge (A_2 \otimes a_2))\hdots \lozenge (A_n \otimes a_n)) \shuffle A )\shuffle B) \otimes (ab)  \\
&= ((\mathsf{P}(\hdots \mathsf{P}(\mathsf{P}(A_1 \otimes a_1) \lozenge (A_2 \otimes a_2))\hdots \lozenge (A_n \otimes a_n)) \shuffle A) \otimes a) \lozenge (B \otimes b) \\
&= \omega_A(([A_1 \otimes a_1] \otimes \hdots \otimes [A_n \otimes a_n]) \otimes A \otimes a) \lozenge (B \otimes b)
\end{align*}
$ii)$: Let $M$ be an $R$-module. It suffices to prove this equality on pure tensors. Consider the following pure tensor of $\mathsf{Sh}(\mathsf{Sh}(\mathsf{Sym}(M)) \otimes \mathsf{Sym}(M)) \otimes \mathsf{Sh}(\mathsf{Sym}(M)) \otimes \mathsf{Sym}(M)$:
 $$([W_1 \otimes w_1] \otimes \hdots \otimes [W_n \otimes w_n]) \otimes W_0 \otimes w_0 \in \mathsf{Sh}(\mathsf{Sh}(\mathsf{Sym}(M)) \otimes \mathsf{Sym}(M)) \otimes \mathsf{Sh}(\mathsf{Sym}(M)) \otimes \mathsf{Sym}(M)$$
where $W, W_1, \hdots W_n \in \mathsf{Sh}(\mathsf{Sym}(M))$ and $w, w_1 \hdots, w_n \in \mathsf{Sym}(M)$. By definition of $\omega_{\mathsf{Sym}(M)}$, we obtain that:
\begin{align*}
&\omega_W(([W_1 \otimes w_1] \otimes \hdots \otimes [W_n \otimes w_n]) \otimes W \otimes w)\\
&= \mathsf{P}(\mathsf{P}(\hdots \mathsf{P}(\mathsf{P}(W_1 \otimes w_1) \lozenge (W_2 \otimes w_2))\hdots \lozenge (W_n \otimes w_n)))\lozenge (W \otimes w) \\
&= (\mathsf{P}(\hdots \mathsf{P}(\mathsf{P}(W_1 \otimes w_1) \lozenge (W_2 \otimes w_2))\hdots \lozenge (W_n \otimes w_n))) \otimes 1) \lozenge (W \otimes w) \\
&= (\mathsf{P}(\hdots \mathsf{P}(\mathsf{P}(W_1 \otimes w_1) \lozenge (W_2 \otimes w_2))\hdots \lozenge (W_n \otimes w_n)) \shuffle W) \otimes w \end{align*}
Notice that $w$ is unaffected by $\omega$. Let $w= m_1 \otimes \hdots \otimes m_k$, $m_i \in M$, then we have that:
$$\mathsf{d}_M(w)= \sum^k_{i=1} (m_1 \otimes \hdots m_{i-1} \otimes m_{i+1} \otimes \hdots m_k) \otimes m_i$$
Then we obtain the following equality: 
\begin{align*}
&(1 \otimes \mathsf{d}_M)(\omega_{\mathsf{Sym}(M)}(([W_1 \otimes w_1] \otimes \hdots \otimes [W_n \otimes w_n]) \otimes W \otimes w))\\
&=(1 \otimes \mathsf{d}_M)((\mathsf{P}(\hdots \mathsf{P}(\mathsf{P}(W_1 \otimes w_1) \lozenge (W_2 \otimes w_2))\hdots \lozenge (W_n \otimes w_n)) \shuffle W) \otimes w) \\
&= (\mathsf{P}(\hdots \mathsf{P}(\mathsf{P}(W_1 \otimes w_1) \lozenge (W_2 \otimes w_2))\hdots \lozenge (W_n \otimes w_n)) \shuffle W) \otimes \mathsf{d}_M(w) \\
&= (\mathsf{P}(\hdots \mathsf{P}(\mathsf{P}(W_1 \otimes w_1) \lozenge (W_2 \otimes w_2))\hdots \lozenge (W_n \otimes w_n)) \shuffle W)\\
&~~~\otimes (\sum^k_{i=1} (m_1 \otimes \hdots m_{i-1} \otimes m_{i+1} \otimes \hdots m_k) \otimes m_i) \\
&= \sum^k_{i=1} (\mathsf{P}(\hdots \mathsf{P}(\mathsf{P}(W_1 \otimes w_1) \lozenge (W_2 \otimes w_2))\hdots \lozenge (W_n \otimes w_n)) \shuffle W)\\
&~~~\otimes (m_1 \otimes \hdots m_{i-1} \otimes m_{i+1} \otimes \hdots m_k) \otimes m_i \\
&= \sum^k_{i=1} (\omega_{\mathsf{Sym}(M)} \otimes 1)(([W_1 \otimes w_1] \otimes \hdots \otimes [W_n \otimes w_n])\otimes W \\
&~~~\otimes (m_1 \otimes \hdots m_{i-1} \otimes m_{i+1} \otimes \hdots m_k)) \otimes m_i \\
&= (\omega_{\mathsf{Sym}(M)} \otimes 1)(([W_1 \otimes w_1] \otimes \hdots \otimes [W_n \otimes w_n])\otimes W \\
&~~~\otimes (\sum^k_{i=1} (m_1 \otimes \hdots m_{i-1} \otimes m_{i+1} \otimes \hdots m_k) \otimes m_i))\\
&= (\omega_{\mathsf{Sym}(M)} \otimes 1)(([W_1 \otimes w_1] \otimes \hdots \otimes [W_n \otimes w_n]) \otimes W \otimes  \mathsf{d}_M(w))
\end{align*}
\end{proof} 

\begin{proposition} For the free Rota-Baxter monad, $1 \otimes \mathsf{d}$ is a deriving transformation. 
\end{proposition} 
\begin{proof} By construction $1 \otimes \mathsf{d}$ is a natural transformation. We need to show {\bf [d.2]} to {\bf [d.5]} (recall that {\bf [d.1]} is redundant). The linear rule {\bf [d.3]}, the Leibniz rule {\bf [d.2]} and the interchange rule {\bf [d.5]} are straightforward but the chain rule {\bf [d.4]} requires some work:

\noindent {\bf [d.2]}: Here we use the Leibniz rule {\bf [d.2]}: 
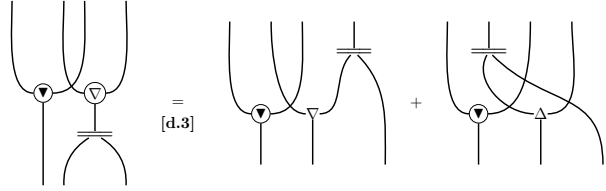
\[\resizebox{!}{2.5cm}{%
\begin{tikzpicture}
	\begin{pgfonlayer}{nodelayer}
		\node [style=port] (8) at (10.5, 0.5) {$=$};
		\node [style=port] (9) at (10.5, 0) {\textbf{[d.3]}};
		\node [style=port] (29) at (7.25, -1.5) {};
		\node [style=duplicate] (30) at (7.25, 0.75) {$\blacktriangledown$};
		\node [style=duplicate] (31) at (8.5, 0.75) {$\nabla$};
		\node [style=port] (32) at (7.75, 3) {};
		\node [style=port] (33) at (8.25, 3) {};
		\node [style=port] (34) at (9.25, 3) {};
		\node [style=port] (35) at (6.5, 3) {};
		\node [style=port] (36) at (7.75, -1.5) {};
		\node [style=differential] (37) at (8.5, -0.25) {{\bf =\!=\!=}};
		\node [style=port] (38) at (9.25, -1.5) {};
		\node [style=port] (39) at (12.5, -1) {};
		\node [style=duplicate] (40) at (12.5, 0.25) {$\blacktriangledown$};
		\node [style=port] (41) at (13.5, 2.5) {};
		\node [style=port] (42) at (11.75, 2.5) {};
		\node [style=differential] (43) at (14.75, 1.75) {{\bf =\!=\!=}};
		\node [style=port] (44) at (14.75, 2.5) {};
		\node [style=port] (45) at (15.5, -1) {};
		\node [style=port] (46) at (13.75, 0.25) {$\nabla$};
		\node [style=port] (47) at (13.75, -1) {};
		\node [style=port] (48) at (12.75, 2.5) {};
		\node [style=port] (49) at (17.75, -1) {};
		\node [style=duplicate] (50) at (17.75, 0.25) {$\blacktriangledown$};
		\node [style=port] (51) at (19, 2.5) {};
		\node [style=port] (52) at (17, 2.5) {};
		\node [style=differential] (53) at (18, 1.75) {{\bf =\!=\!=}};
		\node [style=port] (54) at (18, 2.5) {};
		\node [style=port] (55) at (20, 2.5) {};
		\node [style=port] (56) at (19.25, -1) {};
		\node [style=port] (57) at (19.25, 0.25) {$\Delta$};
		\node [style=port] (58) at (20.75, -1) {};
		\node [style=port] (59) at (16.25, 0.5) {$+$};
	\end{pgfonlayer}
	\begin{pgfonlayer}{edgelayer}
		\draw [style=wire] (30) to (29);
		\draw [style=wire, in=0, out=-90] (33) to (30);
		\draw [style=wire, in=0, out=-90, looseness=0.75] (34) to (31);
		\draw [style=wire, in=180, out=-90, looseness=0.75] (32) to (31);
		\draw [style=wire, in=180, out=-90, looseness=0.75] (35) to (30);
		\draw [style=wire, bend left] (37) to (38);
		\draw [style=wire, bend right] (37) to (36);
		\draw [style=wire] (31) to (37);
		\draw [style=wire] (40) to (39);
		\draw [style=wire, in=0, out=-90] (41) to (40);
		\draw [style=wire, in=180, out=-90, looseness=0.75] (42) to (40);
		\draw [style=wire, in=90, out=-45] (43) to (45);
		\draw [style=wire, in=-91, out=180, looseness=0.75] (46) to (48);
		\draw [style=wire, in=-135, out=0] (46) to (43);
		\draw [style=wire] (47) to (46);
		\draw [style=wire] (43) to (44);
		\draw [style=wire] (50) to (49);
		\draw [style=wire, in=0, out=-90] (51) to (50);
		\draw [style=wire, in=180, out=-90, looseness=0.75] (52) to (50);
		\draw [style=wire, in=90, out=-45, looseness=1.25] (53) to (58);
		\draw [style=wire, in=-120, out=180] (57) to (53);
		\draw [style=wire] (56) to (57);
		\draw [style=wire] (53) to (54);
		\draw [style=wire, in=-90, out=0] (57) to (55);
	\end{pgfonlayer}
\end{tikzpicture}
}%
\]
 {\bf [d.3]}: Here we use the linear rule {\bf [d.3]}:
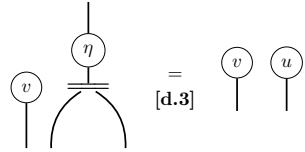
\[\resizebox{!}{2cm}{%
\begin{tikzpicture}
	\begin{pgfonlayer}{nodelayer}
		\node [style=port] (8) at (3.25, -3) {$=$};
		\node [style=port] (9) at (3.25, -3.5) {\textbf{[d.3]}};
		\node [style=component] (18) at (1.5, -2.5) {$\eta$};
		\node [style=port] (19) at (0.25, -4.5) {};
		\node [style=port] (20) at (0.75, -4.5) {};
		\node [style=differential] (21) at (1.5, -3.25) {{\bf =\!=\!=}};
		\node [style=port] (22) at (2.25, -4.5) {};
		\node [style=component] (23) at (0.25, -3.25) {$v$};
		\node [style=port] (24) at (1.5, -1.5) {};
		\node [style=component] (25) at (5.5, -2.75) {$u$};
		\node [style=component] (26) at (4.5, -2.75) {$v$};
		\node [style=port] (27) at (4.5, -3.75) {};
		\node [style=port] (28) at (5.5, -3.75) {};
	\end{pgfonlayer}
	\begin{pgfonlayer}{edgelayer}
		\draw [style=wire, bend left] (21) to (22);
		\draw [style=wire] (18) to (21);
		\draw [style=wire, bend right] (21) to (20);
		\draw [style=wire] (19) to (23);
		\draw [style=wire] (24) to (18);
		\draw [style=wire] (27) to (26);
		\draw [style=wire] (28) to (25);
	\end{pgfonlayer}
\end{tikzpicture}
}%
\]
{\bf [d.4]}: Here we use the counit law, Lemma \ref{RBderiving}, the chain rule {\bf [d.4]}, the monoidal rule {\bf [d.m]}, that $\epsilon$ is an algebra morphism and the triangle identities of the free symmetric algebra adjunction: 
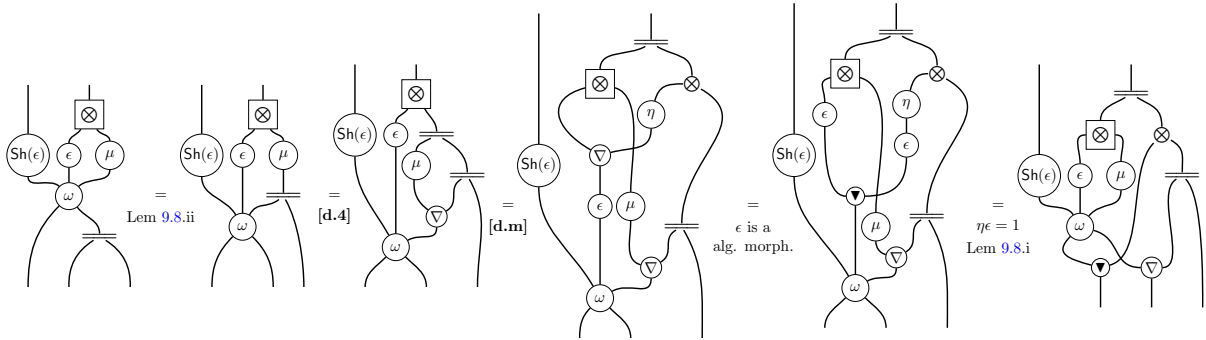
\[\resizebox{!}{4.5cm}{%
\begin{tikzpicture}
	\begin{pgfonlayer}{nodelayer}
		\node [style=port] (0) at (-0.25, -0.75) {};
		\node [style=component] (1) at (0.75, 1.5) {$\omega$};
		\node [style={circle, draw, inner sep=1pt,minimum size=1pt}] (2) at (-0.25, 2.5) {$\mathsf{Sh}(\epsilon)$};
		\node [style=component] (3) at (0.75, 2.5) {$\epsilon$};
		\node [style=component] (4) at (1.75, 2.5) {$\mu$};
		\node [style=function2] (5) at (1.25, 3.5) {$\bigotimes$};
		\node [style=port] (6) at (-0.25, 4.25) {};
		\node [style=port] (7) at (1.25, 4.25) {};
		\node [style=port] (8) at (0.75, -0.75) {};
		\node [style=port] (9) at (2.25, -0.75) {};
		\node [style=codifferential] (10) at (1.5, 0.5) {{\bf =\!=\!=}};
		\node [style=port] (11) at (4.25, -0.75) {};
		\node [style=component] (12) at (5, 0.75) {$\omega$};
		\node [style={circle, draw, inner sep=1pt,minimum size=1pt}] (13) at (4, 2.5) {$\mathsf{Sh}(\epsilon)$};
		\node [style=component] (14) at (5, 2.5) {$\epsilon$};
		\node [style=component] (15) at (6, 2.5) {$\mu$};
		\node [style=function2] (16) at (5.5, 3.5) {$\bigotimes$};
		\node [style=port] (17) at (4, 4.25) {};
		\node [style=port] (18) at (5.5, 4.25) {};
		\node [style=port] (19) at (6.5, -0.75) {};
		\node [style=codifferential] (20) at (6, 1.5) {{\bf =\!=\!=}};
		\node [style=port] (21) at (5.75, -0.75) {};
		\node [style=port] (22) at (3, 1.5) {$=$};
		\node [style=port] (23) at (3, 1) {Lem \ref{RBderiving}.ii};
		\node [style=port] (24) at (8, -0.75) {};
		\node [style=component] (25) at (8.75, 0.25) {$\omega$};
		\node [style={circle, draw, inner sep=1pt,minimum size=1pt}] (26) at (7.75, 3) {$\mathsf{Sh}(\epsilon)$};
		\node [style=component] (27) at (8.75, 3) {$\epsilon$};
		\node [style=function2] (28) at (9.25, 4) {$\bigotimes$};
		\node [style=port] (29) at (7.75, 4.75) {};
		\node [style=port] (30) at (9.25, 4.75) {};
		\node [style=port] (31) at (10.75, -0.75) {};
		\node [style=port] (32) at (9.5, -0.75) {};
		\node [style=duplicate] (33) at (9.75, 1) {$\nabla$};
		\node [style=component] (34) at (9.25, 2.25) {$\mu$};
		\node [style=differential] (35) at (9.75, 3) {{\bf =\!=\!=}};
		\node [style=differential] (36) at (10.5, 2) {{\bf =\!=\!=}};
		\node [style=port] (37) at (7.25, 1.5) {$=$};
		\node [style=port] (38) at (7.25, 1) {\textbf{[d.4]}};
		\node [style=port] (39) at (13.25, -2) {};
		\node [style=component] (40) at (13.75, -1) {$\omega$};
		\node [style={circle, draw, inner sep=1pt,minimum size=1pt}] (41) at (12.25, 2.25) {$\mathsf{Sh}(\epsilon)$};
		\node [style=component] (42) at (13.75, 1.25) {$\epsilon$};
		\node [style=port] (43) at (12.25, 6) {};
		\node [style=port] (44) at (16.25, -2) {};
		\node [style=port] (45) at (14.5, -2) {};
		\node [style=duplicate] (46) at (15, -0.25) {$\nabla$};
		\node [style=component] (47) at (14.5, 1.25) {$\mu$};
		\node [style=differential] (48) at (15.75, 0.75) {{\bf =\!=\!=}};
		\node [style=component] (49) at (15, 3.5) {$\eta$};
		\node [style=differential] (50) at (16, 4.25) {$\bigotimes$};
		\node [style=port] (51) at (15, 6) {};
		\node [style=function2] (52) at (13.75, 4.25) {$\bigotimes$};
		\node [style=duplicate] (53) at (13.75, 2.5) {$\nabla$};
		\node [style=codifferential] (54) at (15, 5.25) {{\bf =\!=\!=}};
		\node [style=port] (55) at (11.5, 1.25) {};
		\node [style=port] (56) at (11.5, 0.75) {\textbf{[d.m]}};
		\node [style=port] (57) at (11.5, 1.25) {$=$};
		\node [style=port] (58) at (19.25, -1.75) {};
		\node [style=component] (59) at (20, -0.75) {$\omega$};
		\node [style={circle, draw, inner sep=1pt,minimum size=1pt}] (60) at (18.5, 2.5) {$\mathsf{Sh}(\epsilon)$};
		\node [style=duplicate] (61) at (20, 1.5) {$\blacktriangledown$};
		\node [style=port] (62) at (18.5, 6.25) {};
		\node [style=port] (63) at (22.25, -1.75) {};
		\node [style=port] (64) at (20.75, -1.75) {};
		\node [style=duplicate] (65) at (21, 0) {$\nabla$};
		\node [style=component] (66) at (20.5, 0.75) {$\mu$};
		\node [style=differential] (67) at (21.75, 1) {{\bf =\!=\!=}};
		\node [style=component] (68) at (21.25, 3.75) {$\eta$};
		\node [style=differential] (69) at (22, 4.5) {$\bigotimes$};
		\node [style=port] (70) at (21, 6.25) {};
		\node [style=function2] (71) at (19.75, 4.5) {$\bigotimes$};
		\node [style=codifferential] (72) at (21, 5.5) {{\bf =\!=\!=}};
		\node [style=component] (73) at (19.25, 3.5) {$\epsilon$};
		\node [style=component] (74) at (21.25, 2.75) {$\epsilon$};
		\node [style=port] (75) at (17.5, 1.25) {$=$};
		\node [style=port] (76) at (17.5, 0.75) {$\epsilon$ is a};
		\node [style=port] (77) at (17.5, 0.25) {alg. morph.};
		\node [style=port] (93) at (23.5, 1.25) {};
		\node [style=port] (94) at (23.5, 0.75) {$\eta \epsilon =1$};
		\node [style=port] (95) at (23.5, 1.25) {$=$};
		\node [style={circle, draw}] (96) at (25.5, 0.75) {$\omega$};
		\node [style={circle, draw, inner sep=1pt,minimum size=1pt}] (97) at (24.5, 2) {$\mathsf{Sh}(\epsilon)$};
		\node [style=duplicate] (98) at (26, -0.25) {$\blacktriangledown$};
		\node [style=port] (99) at (24.5, 4.75) {};
		\node [style=port] (100) at (28.5, -1.25) {};
		\node [style=duplicate] (101) at (27.25, -0.25) {$\nabla$};
		\node [style={circle, draw}] (102) at (26.5, 2) {$\mu$};
		\node [style=differential] (103) at (28, 2) {{\bf =\!=\!=}};
		\node [style=differential] (104) at (27.5, 3) {$\bigotimes$};
		\node [style=port] (105) at (26.75, 4.75) {};
		\node [style={regular polygon,regular polygon sides=4, draw, inner sep=1pt,minimum size=1pt}] (106) at (26, 3) {$\bigotimes$};
		\node [style=codifferential] (107) at (26.75, 4) {{\bf =\!=\!=}};
		\node [style={circle, draw}] (108) at (25.5, 2) {$\epsilon$};
		\node [style=port] (109) at (27.25, -1.25) {};
		\node [style=port] (110) at (26, -1.25) {};
		\node [style=port] (111) at (23.5, 0.25) {Lem \ref{RBderiving}.i};
	\end{pgfonlayer}
	\begin{pgfonlayer}{edgelayer}
		\draw [style=wire] (6) to (2);
		\draw [style=wire, in=90, out=-120] (5) to (3);
		\draw [style=wire, in=90, out=-63] (5) to (4);
		\draw [style=wire] (7) to (5);
		\draw [style=wire] (3) to (1);
		\draw [style=wire, in=45, out=-90] (4) to (1);
		\draw [style=wire, in=135, out=-90] (2) to (1);
		\draw [style=wire, bend left] (10) to (9);
		\draw [style=wire, bend right] (10) to (8);
		\draw [style=wire, in=90, out=-135] (1) to (0);
		\draw [style=wire, in=90, out=-53] (1) to (10);
		\draw [style=wire] (17) to (13);
		\draw [style=wire, in=90, out=-120] (16) to (14);
		\draw [style=wire, in=90, out=-63] (16) to (15);
		\draw [style=wire] (18) to (16);
		\draw [style=wire, in=90, out=-60, looseness=0.75] (20) to (19);
		\draw [style=wire, in=90, out=-135] (12) to (11);
		\draw [style=wire] (15) to (20);
		\draw [style=wire, in=120, out=-90] (13) to (12);
		\draw [style=wire] (14) to (12);
		\draw [style=wire, in=60, out=-143, looseness=1.25] (20) to (12);
		\draw [style=wire, in=90, out=-45] (12) to (21);
		\draw [style=wire] (29) to (26);
		\draw [style=wire, in=90, out=-120] (28) to (27);
		\draw [style=wire] (30) to (28);
		\draw [style=wire, in=90, out=-135] (25) to (24);
		\draw [style=wire, in=120, out=-90] (26) to (25);
		\draw [style=wire] (27) to (25);
		\draw [style=wire, in=90, out=-45] (25) to (32);
		\draw [style=wire, bend left] (33) to (34);
		\draw [style=wire, in=-150, out=30, looseness=1.25] (33) to (36);
		\draw [style=wire, in=-30, out=90] (36) to (35);
		\draw [style=wire, in=-150, out=90] (34) to (35);
		\draw [style=wire, in=90, out=-60] (28) to (35);
		\draw [style=wire, in=30, out=-90, looseness=1.25] (33) to (25);
		\draw [style=wire, in=90, out=-60] (36) to (31);
		\draw [style=wire] (43) to (41);
		\draw [style=wire, in=90, out=-135] (40) to (39);
		\draw [style=wire, in=120, out=-90] (41) to (40);
		\draw [style=wire] (42) to (40);
		\draw [style=wire, in=90, out=-45] (40) to (45);
		\draw [style=wire, in=-90, out=165] (46) to (47);
		\draw [style=wire, in=-150, out=30, looseness=1.25] (46) to (48);
		\draw [style=wire, in=30, out=-90, looseness=1.25] (46) to (40);
		\draw [style=wire, in=90, out=-60] (48) to (44);
		\draw [style=wire] (54) to (51);
		\draw [style=wire, in=-165, out=150, looseness=1.75] (53) to (52);
		\draw [style=wire, in=-150, out=90] (52) to (54);
		\draw [style=wire, in=-30, out=90, looseness=1.25] (50) to (54);
		\draw [style=wire, in=-90, out=15, looseness=1.25] (53) to (49);
		\draw [style=wire, in=180, out=75, looseness=1.25] (49) to (50);
		\draw [style=wire] (53) to (42);
		\draw [style=wire, in=90, out=-15, looseness=0.75] (52) to (47);
		\draw [style=wire, in=90, out=-30, looseness=1.25] (50) to (48);
		\draw [style=wire] (62) to (60);
		\draw [style=wire, in=90, out=-135] (59) to (58);
		\draw [style=wire, in=120, out=-90] (60) to (59);
		\draw [style=wire] (61) to (59);
		\draw [style=wire, in=90, out=-45] (59) to (64);
		\draw [style=wire, in=-90, out=165] (65) to (66);
		\draw [style=wire, in=-150, out=30, looseness=1.25] (65) to (67);
		\draw [style=wire, in=30, out=-90, looseness=1.25] (65) to (59);
		\draw [style=wire, in=90, out=-60] (67) to (63);
		\draw [style=wire] (72) to (70);
		\draw [style=wire, in=-150, out=90] (71) to (72);
		\draw [style=wire, in=-30, out=90, looseness=1.25] (69) to (72);
		\draw [style=wire, in=180, out=75, looseness=1.25] (68) to (69);
		\draw [style=wire, in=90, out=-15, looseness=0.75] (71) to (66);
		\draw [style=wire, in=90, out=-30, looseness=1.25] (69) to (67);
		\draw [style=wire, in=90, out=-165, looseness=1.50] (71) to (73);
		\draw [style=wire, in=180, out=-90, looseness=0.75] (73) to (61);
		\draw [style=wire, in=0, out=-90, looseness=1.25] (74) to (61);
		\draw [style=wire] (68) to (74);
		\draw [style=wire] (99) to (97);
		\draw [style=wire, in=120, out=-90] (97) to (96);
		\draw [style=wire, in=90, out=-45, looseness=0.50] (103) to (100);
		\draw [style=wire] (107) to (105);
		\draw [style=wire, in=-150, out=90] (106) to (107);
		\draw [style=wire, in=-30, out=90, looseness=1.25] (104) to (107);
		\draw [style=wire, in=90, out=0, looseness=1.25] (106) to (102);
		\draw [style=wire, in=90, out=-30, looseness=1.25] (104) to (103);
		\draw [style=wire, in=90, out=180, looseness=1.25] (106) to (108);
		\draw [style=wire] (108) to (96);
		\draw [style=wire, in=51, out=-90, looseness=1.25] (102) to (96);
		\draw [style=wire, in=165, out=-135, looseness=1.75] (96) to (98);
		\draw [style=wire, in=180, out=-15, looseness=1.25] (96) to (101);
		\draw [style=wire] (98) to (110);
		\draw [style=wire] (101) to (109);
		\draw [style=wire, in=15, out=-135, looseness=0.75] (104) to (98);
		\draw [style=wire, in=0, out=-150, looseness=0.75] (103) to (101);
	\end{pgfonlayer}
\end{tikzpicture}
}%
\]
{\bf [d.5]}: Here we use the interchange rule {\bf [d.5]}: 
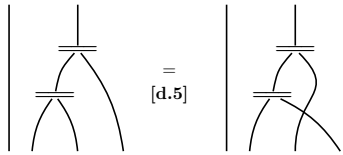
\[\resizebox{!}{2cm}{%
\begin{tikzpicture}
	\begin{pgfonlayer}{nodelayer}
		\node [style=port] (0) at (0.25, -0.25) {};
		\node [style=port] (1) at (0.25, 3) {};
		\node [style=differential] (2) at (1.25, 1) {{\bf =\!=\!=}};
		\node [style=port] (3) at (0.75, -0.25) {};
		\node [style=port] (4) at (1.75, 3) {};
		\node [style=port] (5) at (1.75, -0.25) {};
		\node [style=port] (6) at (2.75, -0.25) {};
		\node [style=differential] (7) at (1.75, 2) {{\bf =\!=\!=}};
		\node [style=port] (8) at (3.75, 1.5) {$=$};
		\node [style=port] (9) at (3.75, 1) {\textbf{[d.5]}};
		\node [style=port] (10) at (5, -0.25) {};
		\node [style=port] (11) at (5, 3) {};
		\node [style=differential] (12) at (6.5, 2) {{\bf =\!=\!=}};
		\node [style=differential] (13) at (6, 1) {{\bf =\!=\!=}};
		\node [style=port] (14) at (5.5, -0.25) {};
		\node [style=port] (15) at (7.5, -0.25) {};
		\node [style=port] (16) at (6.5, -0.25) {};
		\node [style=port] (17) at (6.5, 3) {};
	\end{pgfonlayer}
	\begin{pgfonlayer}{edgelayer}
		\draw [style=wire] (0) to (1);
		\draw [style=wire, bend right=15, looseness=1.25] (7) to (2);
		\draw [style=wire, bend left=15] (7) to (6);
		\draw [style=wire] (4) to (7);
		\draw [style=wire, bend right=15] (2) to (3);
		\draw [style=wire, bend left=15] (2) to (5);
		\draw [style=wire] (10) to (11);
		\draw [style=wire, bend right=15, looseness=1.25] (12) to (13);
		\draw [style=wire] (17) to (12);
		\draw [style=wire, bend right=15] (13) to (14);
		\draw [style=wire, in=-45, out=90, looseness=1.50] (16) to (12);
		\draw [style=wire, bend right=15] (15) to (13);
	\end{pgfonlayer}
\end{tikzpicture}
}%
\]
\end{proof} 
\end{example}

To obtain examples of non-additive bialgebra modalities, we need simply apply the construction for Section \ref{nonaddbialg} to our examples of additive bialgebra modalities. 

\begin{example}\label{ex5} \normalfont The free symmetric algebra modality induces a non-additive bialgebra modality on $\mathsf{MOD}_R$ which has a deriving transformation. 
\end{example}

\begin{example}\label{ex6} \normalfont The free differential algebra modality induces a non-additive bialgebra modality on $\mathsf{MOD}_R$ which does not have a deriving transformation (since if it did, $\mathsf{DIFF}$ would have one). 
\end{example}

Finally to obtain a coalgebra modality which is not a bialgebra modality and does not have a deriving transformation, we look towards differential Rota-Baxter algebras. 

\begin{example}\label{ex7} \normalfont Let $R$ be a commutative ring. A (commutative) \textbf{differential Rota-Baxter algebra (of weight $0$)} \cite{guo2008differential} over $R$ is a triple $(A, \mathsf{D}, \mathsf{P})$ consisting of a differential algebra $(A, \mathsf{D})$ and a Rota-Baxter algebra $(A, \mathsf{P})$ such that $\mathsf{P}\mathsf{D}=1_A$. It turns out that the free Rota-Baxter algebra over a differential algebra is also its free differential Rota-Baxter algebra, and therefore inducing the following adjunction between the category of differential algebras, $\mathsf{CDA}_{R}$, and the category of differential Rota-Baxter algebras, $\mathsf{CDRBA}_{R}$: 
$$\xymatrixcolsep{2.5pc}\xymatrix{\mathsf{CDA}_{R} \ar@<+1.1ex>[r]^-{\mathsf{RB}} & \mathsf{CDRBA}_{R} \ar@<+1ex>[l]_-{\bot}^-{U}
  }$$
  The full construction can be found in \cite{guo2008differential}. Once again, to obtain an algebra modality we compose this adjunction with the free differential algebra adjunction: 
  $$\xymatrixcolsep{2.5pc}\xymatrix{\mathsf{MOD}_{R} \ar@<+1.1ex>[r]^-{\mathsf{DIFF}} & \mathsf{CDA}_{R} \ar@<+1ex>[l]_-{\bot}^-{U} \ar@<+1.1ex>[r]^-{\mathsf{RB}} & \mathsf{CDRBA}_{R} \ar@<+1ex>[l]_-{\bot}^-{U}
  }$$
  This algebra modality is not comonoidal for the same reasons as the free Rota-Baxter algebra, and is not a differential category failing the chain rule like the free differential algebra. 
\end{example}

%%%%%%%%%%%%%%%%%%%%%%%%%%%%%%%%%%%%%%%%%%%%%%%%%%%%%%%%%%%%%%%%%%%%%%%
%%%%%%%%%%%%%%%%%%%%%%%%%%%%%%%%%%%%%%%%%%%%%%%%%%%%%%%%%%%%%%%%%%%%%%%

\section{Conclusion}

%%%%%%%%%%%%%%%%%%%%%%%%%%%%%%%%%%%%%%%%%%%%%%%%%%%%%%%%%%%%%%%%%%%%%%%

There is a tendency to assume that the only important coalgebra modalities are those which arise through linear logic: that is those which are monoidal coalgebra modalities 
(or equivalently additive coalgebra modalities -- in the sense of this paper).  Certainly, it is true that it is the monoidal coalgebra modalities that have filled the lion's share of the literature.  Of course, this does not mean that their relatives, the mere coalgebra modalities, are not worthy of scientific attention -- indeed, we believe that they have been wrongly overlooked. Thus, one objective of this paper is to reemphasize the importance of these mere coalgebra modalities and to provide a 
stock of examples.  

Of course, it would have been very much simpler if all these observations had been made in \cite{blute2006differential} -- the original paper.  Unfortunately, they were not.  
Marcelo Fiore's paper \cite{fiore2007differential}  made us realize that there was much more to say.   In particular, the relation between deriving transformations and coderelictions had not been fully developed.  Here we have revisited this relation and, in addition, filled in some of the gaps left in the original paper. That there will be yet more to say we are sure!

Significantly, it is not that the notion of differentiation varies but rather, as the setting varies, significantly different presentations of the differential become possible.  
Philosophically this is certainly how things should be ... and apparently it is how they are!

\bibliographystyle{spmpsci}      % mathematics and physical sciences
\bibliography{DCRbib}   % name your BibTeX data base

%%%%%%%%%%%%%%%%%%%%%%%%%%%%%%%%%%%%%%%%%%%%%%%%%%%%%%%%%%%%%%%%%%%%%%%
%%%%%%%%%%%%%%%%%%%%%%%%%%%%%%%%%%%%%%%%%%%%%%%%%%%%%%%%%%%%%%%%%%%%%%%

% Activate the appendix
% from now on sections are numerated with capital letters
\appendix

%%%%%%%%%%%%%%%%%%%%%%%%%%%%%%%%%%%%%%%%%%%%%%%%%%%%%%%%%%%%%%%%%%%%%%%
%%%%%%%%%%%%%%%%%%%%%%%%%%%%%%%%%%%%%%%%%%%%%%%%%%%%%%%%%%%%%%%%%%%%%%%

\section{Proof of Proposition \ref{monoidaltoaddbialg}}\label{montoadd}

%%%%%%%%%%%%%%%%%%%%%%%%%%%%%%%%%%%%%%%%%%%%%%%%%%%%%%%%%%%%%%%%%%%%%%%

In this appendix we provide the complete proof that the monoidal coalgebra modality of an additive linear category is an additive bialgebra modality. Recall the definitions of $\nabla$ and $u$ in their string diagram representations:

\begin{equation}\label{}\begin{gathered} 

 \end{gathered}\end{equation}

\begin{lemma} $\nabla$ and $u$ are natural transformations. 
\end{lemma}
\begin{proof} By construction, $\nabla$ is a natural transformation since it is the composition of natural transformation. The unit $u$ on the other hand is not automatically a natural transformation by construction, since $m_K$ is a map and not a natural transformations. Let $f: A \to B$, and since $\oc(0)\oc(f) = \oc(0)$, we obtain the following equality: 
\[   \begin{array}[c]{c}\resizebox{!}{2.5cm}{%
\begin{tikzpicture}
	\begin{pgfonlayer}{nodelayer}
		\node [style=component] (0) at (4.25, 2.75) {$u$};
		\node [style=port] (1) at (4.25, 0.25) {};
		\node [style=function2] (2) at (4.25, 1.5) {$f$};
		\node [style=port] (3) at (5.25, 1.5) {$=$};
		\node [style=function] (4) at (6.25, 2.75) {$0$};
		\node [style=port] (5) at (6.25, 0.25) {};
		\node [style=function2] (6) at (6.25, 1.5) {$f$};
		\node [style=component] (7) at (6.25, 4) {$m$};
		\node [style=port] (8) at (7.5, 1.5) {$=$};
		\node [style=component] (9) at (8.75, 3) {$m$};
		\node [style=function] (10) at (8.75, 1.75) {$0$};
		\node [style=port] (11) at (8.75, 0.5) {};
		\node [style=component] (12) at (10.75, 2) {$u$};
		\node [style=port] (14) at (10.75, 0.5) {};
		\node [style=port] (15) at (9.75, 1.5) {$=$};
	\end{pgfonlayer}
	\begin{pgfonlayer}{edgelayer}
		\draw [style=wire] (0) to (2);
		\draw [style=wire] (2) to (1);
		\draw [style=wire] (4) to (6);
		\draw [style=wire] (6) to (5);
		\draw [style=wire] (7) to (4);
		\draw [style=wire] (9) to (10);
		\draw [style=wire] (10) to (11);
		\draw [style=wire] (12) to (14);
	\end{pgfonlayer}
\end{tikzpicture}
  }% 
\end{array}
\]
\end{proof} 

For space and simplification, we define $\phi= e \otimes \varepsilon + \varepsilon \otimes e$. Therefore, $\nabla$ can also be drawn as:
\begin{equation}\label{}\begin{gathered} 
   %
 \end{gathered}\end{equation}
 
\begin{lemma}\label{philem} $\phi$ satisfies the following equalities: 
\[\resizebox{!}{2cm}{%
%
}%
\] 
\end{lemma} 
\begin{proof}
For the first equality we have that:
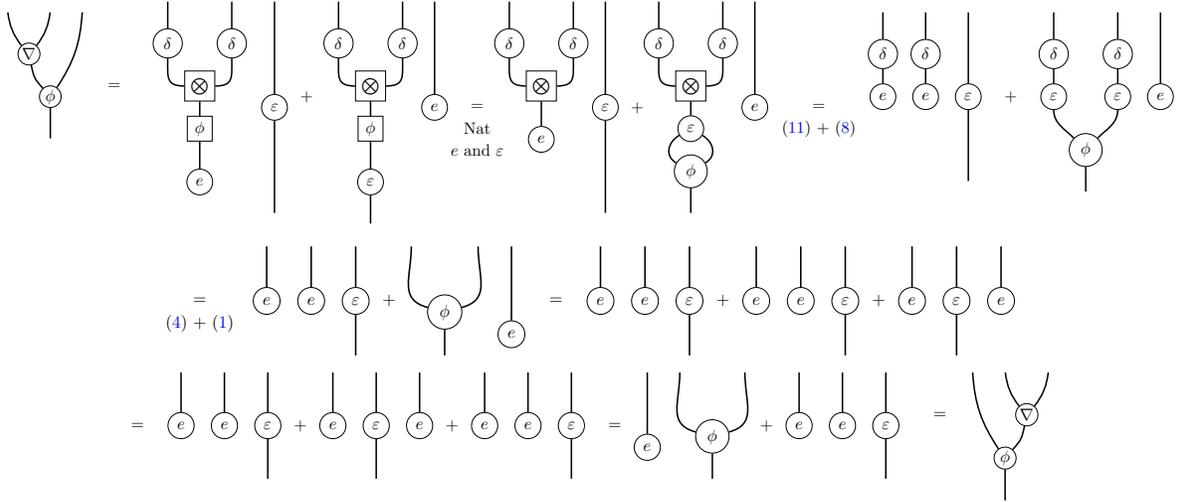
\[\resizebox{!}{3cm}{%
%
}%
\]
For left identity of the second equality we have that: 
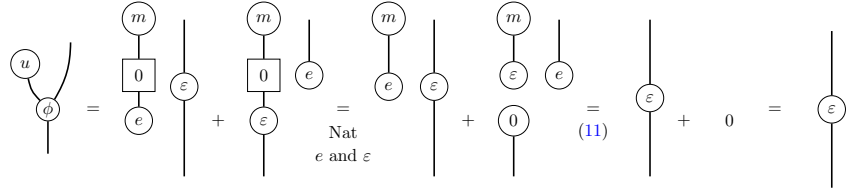
\[\resizebox{!}{2.5cm}{%
%
}%
\]
And similarly for the right identity. Lastly, for the third identity we have that: 
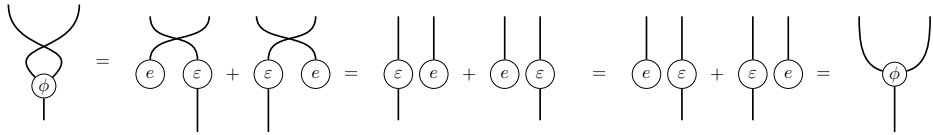
\[\resizebox{!}{1.75cm}{%
%
}%
\]
\end{proof}

 \begin{lemma} For each object $A$, $(\oc A, \nabla, u)$ is a commutative monoid. 
\end{lemma}
\begin{proof} This follows mostly from Lemma \ref{philem}. \\
1. Associativity: 
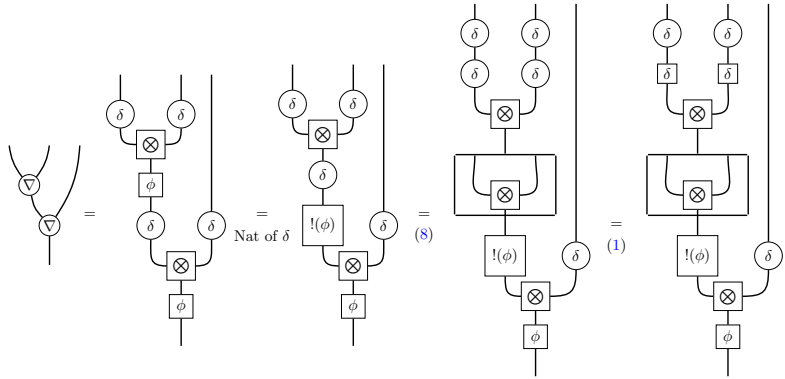
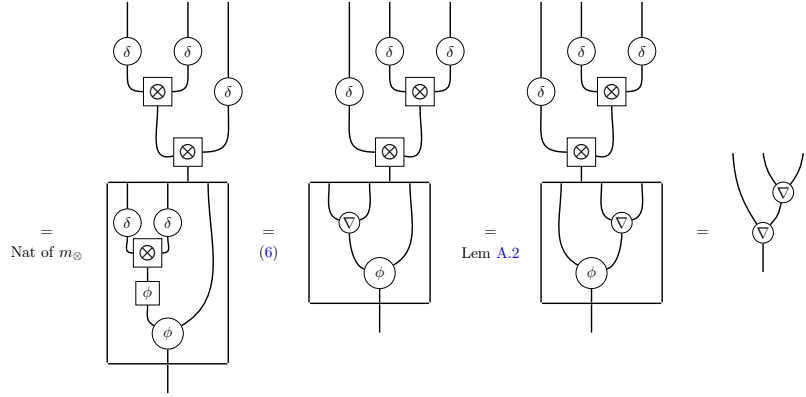
\[\resizebox{!}{5cm}{%
%
}%
\]
2. Unit Laws: (We only show one of them, since calculation for the other is similar) 
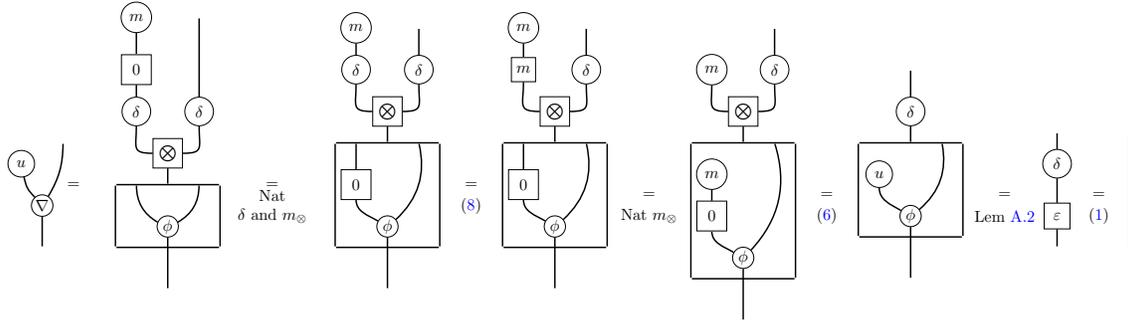
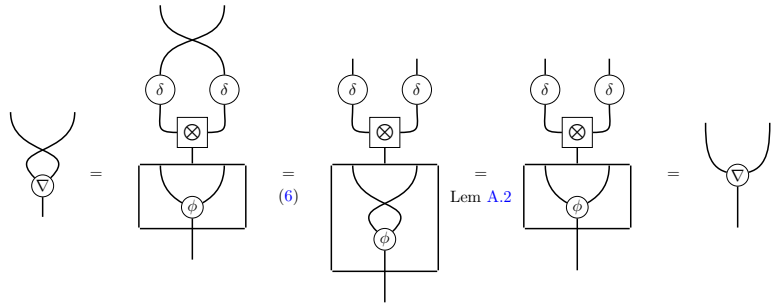
\[\resizebox{!}{4.25cm}{%
%
}%
\]
\end{proof} 

\begin{lemma} For each object $A$, $(\oc A, \nabla, u, \Delta, e)$ is a bialgebra. 
\end{lemma}
\begin{proof} We need to check the four bialgebra compatibility relations: 
\begin{enumerate}
\item Multiplication and comultiplication compatibility: 
\[\resizebox{!}{3.5cm}{%
%
}%
\]
\item Multiplication and counit compatibility:
\[\resizebox{!}{2.5cm}{%
%
}%
\]
   \item Comultiplication and unit compatibility: 
\[\resizebox{!}{2cm}{%
%
}%
\]
\item Counit and unit compatibility: 
\[\resizebox{!}{1.75cm}{%
%
}%
\]
\end{enumerate}
\end{proof}

\begin{proposition} The monoidal coalgebra modality of an additive linear category is an additive bialgebra modality. 
\end{proposition} 
\begin{proof} We first prove that we have a bialgebra modality by proving the compatibility of $\varepsilon$ and $\nabla$: 
\[\resizebox{!}{3cm}{%
%
}%
\]
Finally, we show that the bialgebra modality is in fact additive by proving the compatibility with the additive structure:

\noindent 1. $\oc(0)=eu$:
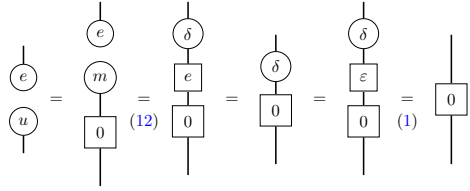
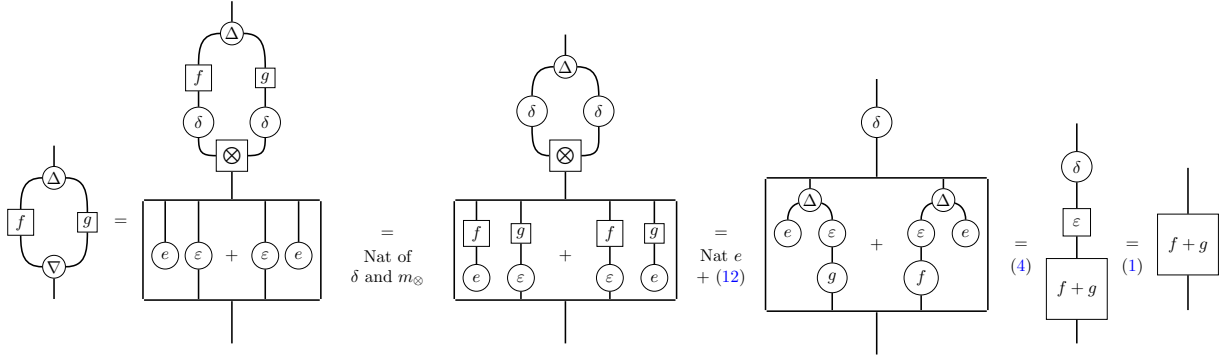
\[\resizebox{!}{2.5cm}{%
%
}%
\]
\end{proof}

%%%%%%%%%%%%%%%%%%%%%%%%%%%%%%%%%%%%%%%%%%%%%%%%%%%%%%%%%%%%%%%%%%%%%%%

\section{Proof of Proposition \ref{monoidalnabla}}\label{Propappendix}

%%%%%%%%%%%%%%%%%%%%%%%%%%%%%%%%%%%%%%%%%%%%%%%%%%%%%%%%%%%%%%%%%%%%%%%

We use the same notation introduced in Appendix \ref{montoadd}. 

\begin{lemma} $\nabla$ and $u$ are $\oc$-coalgebra morphisms.
\end{lemma} 
\begin{proof} We first show that $\nabla$ is a $\oc$-coalgebra morphism:
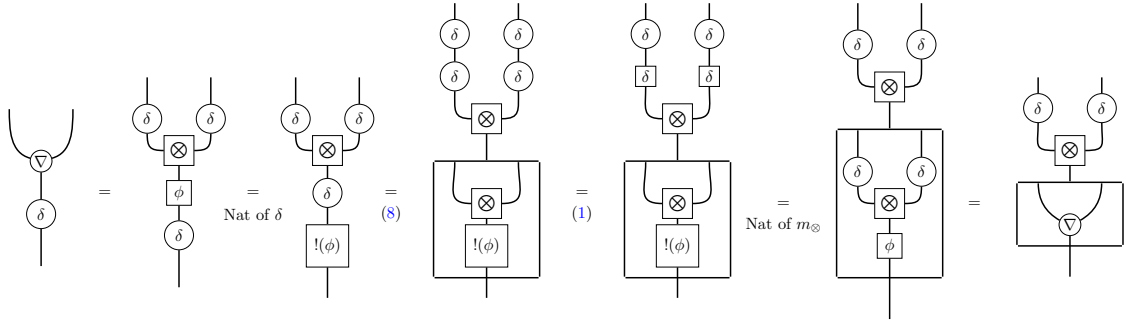
\[\resizebox{!}{4.25cm}{%
%
}%
\]
Next we show that $u$ is also a $\oc$-coalgebra morphism:
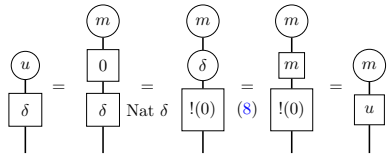
\[\resizebox{!}{2cm}{%
%
}%
\]
\end{proof} 

\begin{lemma} The following equalities hold: 
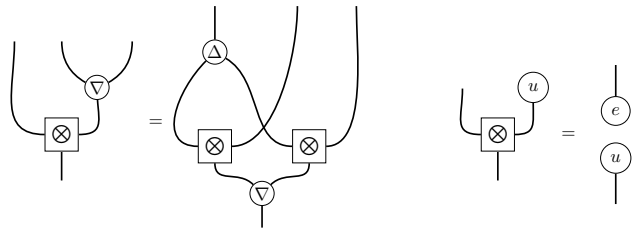
\[ \resizebox{!}{3cm}{%
%
}%
\]
\end{lemma}
\begin{proof} We first prove the equality on the left: 
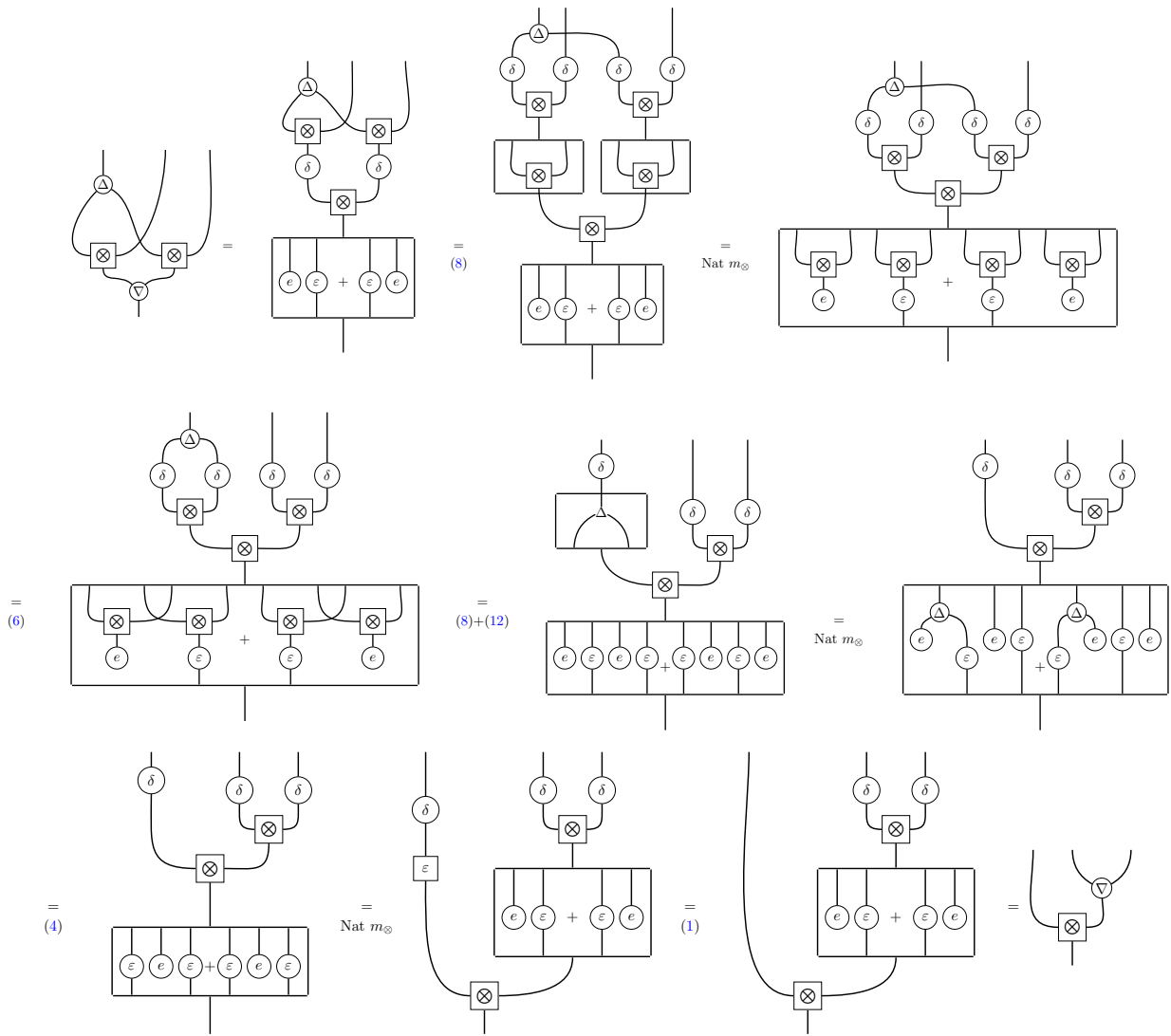
\[\resizebox{!}{5.25cm}{%
%
}%
\]
For the equality on the right, we have that:
\[\resizebox{!}{2.5cm}{%
%
}%
\]
\end{proof} 

%%%%%%%%%%%%%%%%%%%%%%%%%%%%%%%%%%%%%%%%%%%%%%%%%%%%%%%%%%%%%%%%%%%%%%%
%%%%%%%%%%%%%%%%%%%%%%%%%%%%%%%%%%%%%%%%%%%%%%%%%%%%%%%%%%%%%%%%%%%%%%%

\section{Proof of Proposition \ref{addbialgtomonoidal}}\label{addtomon}

%%%%%%%%%%%%%%%%%%%%%%%%%%%%%%%%%%%%%%%%%%%%%%%%%%%%%%%%%%%%%%%%%%%%%%%

In this appendix we prove the converse of Appendix \ref{montoadd}: that an additive bialgebra modality is a monoidal coalgebra modality. Recall the definitions of $m_\otimes$ and $m_K$ in their string diagram representations:

\begin{equation}\label{}\begin{gathered} 
\resizebox{!}{6cm}{%

}%
 \end{gathered}\end{equation}
 
 By construction we have that $m_\otimes$ is natural. 

\begin{lemma} $m_\otimes$ is a natural transformation. 
\end{lemma}

For multiple parts of the following proofs, we will need the following useful identities:

\begin{lemma}\label{apelem} The following equality holds: 
\[\resizebox{!}{4cm}{%
%
}%
\]
\end{lemma}
\begin{proof} By the additive bialgebra modality identities, we have that:
\[\resizebox{!}{4.25cm}{%
%
}%
\]
\end{proof} 

\begin{lemma}\label{apelem2} The following equalities hold:
\[\resizebox{!}{6cm}{%
%
}%
\]
\end{lemma}
\begin{proof} We only prove the equality on the left and the proof for the equality on the right is similar. 
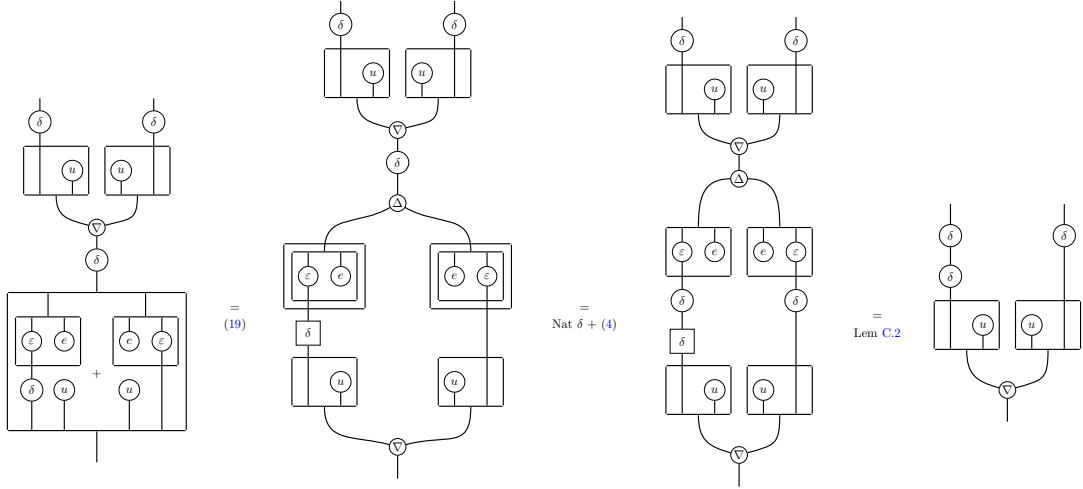
\[\resizebox{!}{6.5cm}{%
%
}%
\]
\end{proof} 

\begin{lemma}\label{apelem3} The following equalities holds: 
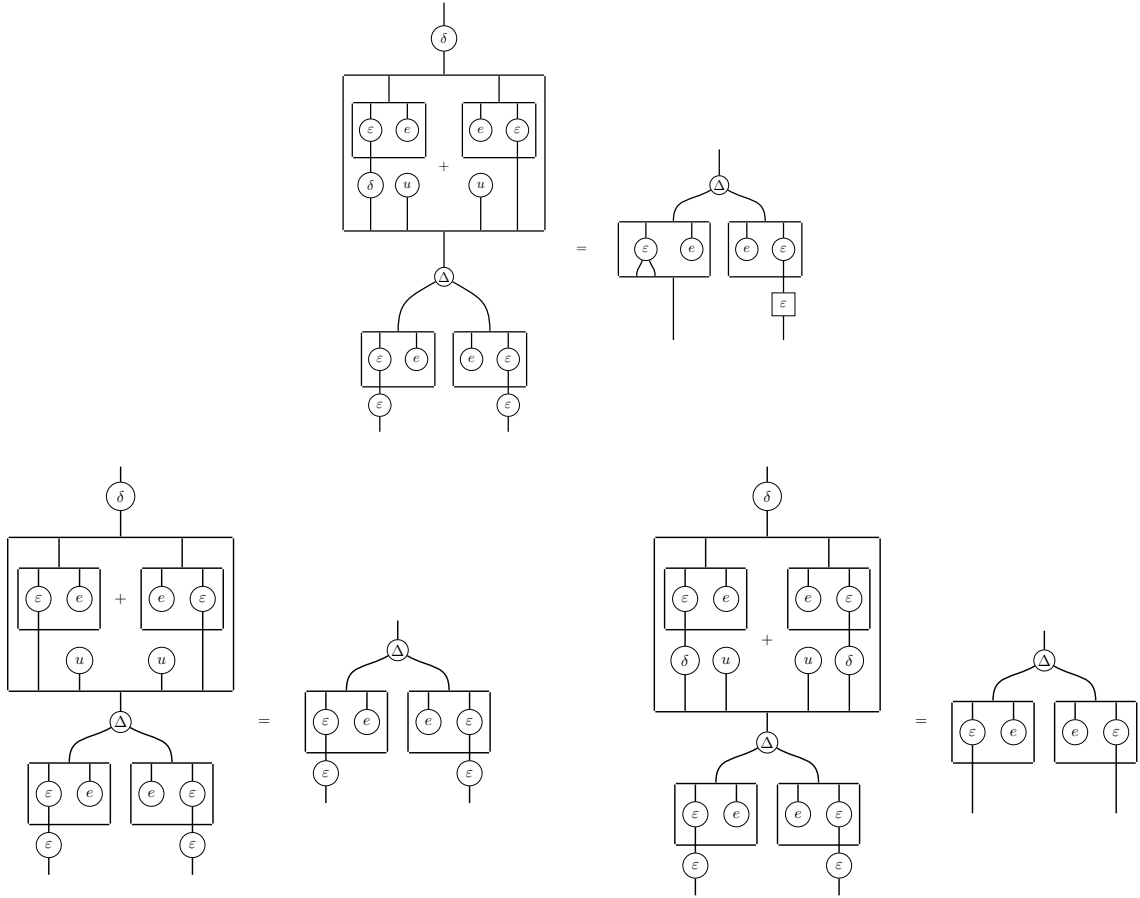
\[\resizebox{!}{5.75cm}{%
%
}%
\]
\end{lemma}
\begin{proof} The proofs are similar to the proof of Lemma \ref{apelem2}, so we leave these as an exercise. 
\end{proof}

\begin{lemma} $(\oc, m_\otimes, m_K)$ is a symmetric monoidal functor. 
\end{lemma}
\begin{proof} We must prove the associativity, unit, and symmetry coherences. 

\noindent 1. Associativity: Consider the following series of equalities: 
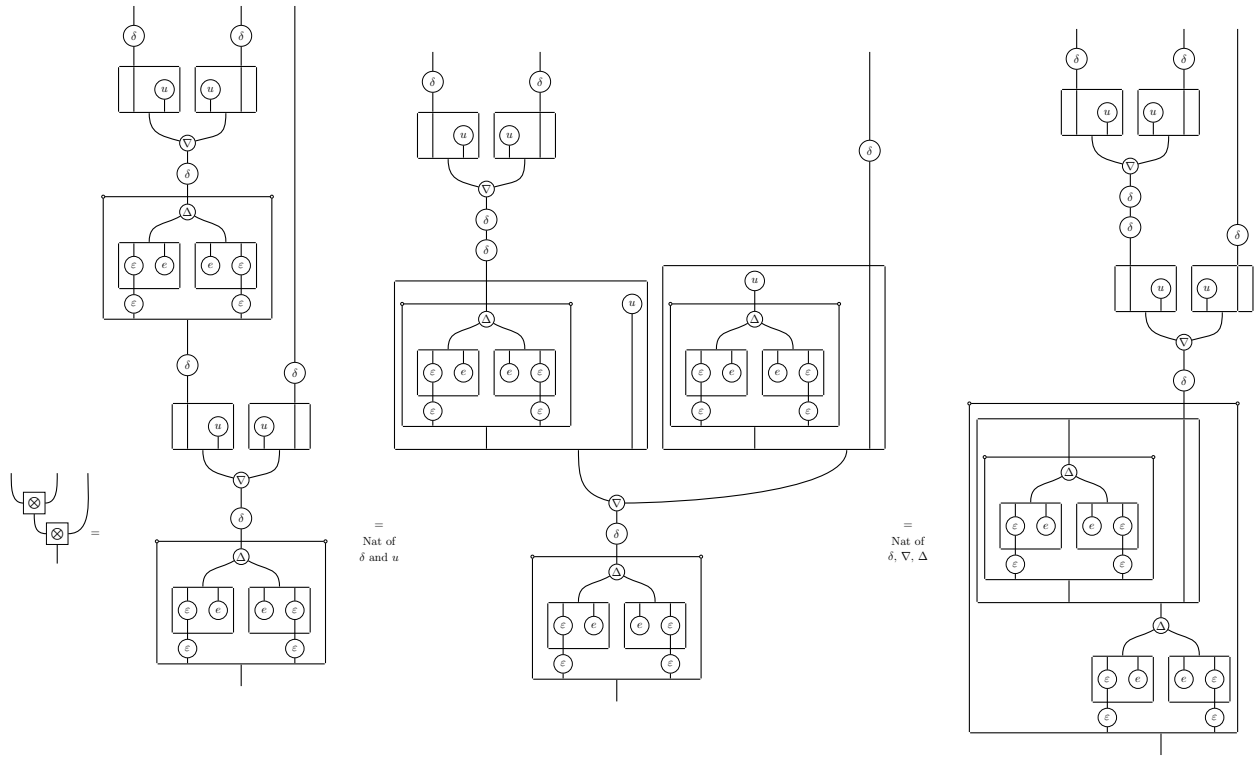
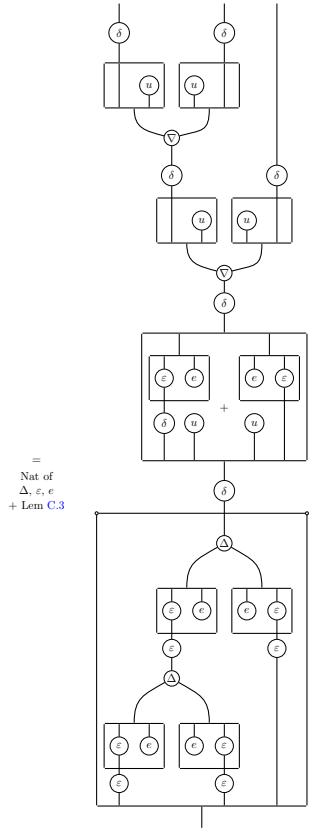
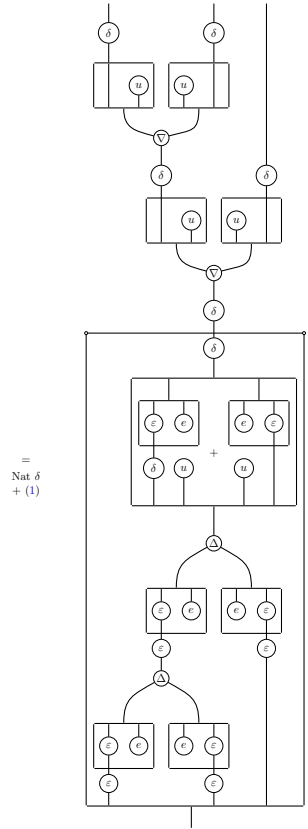
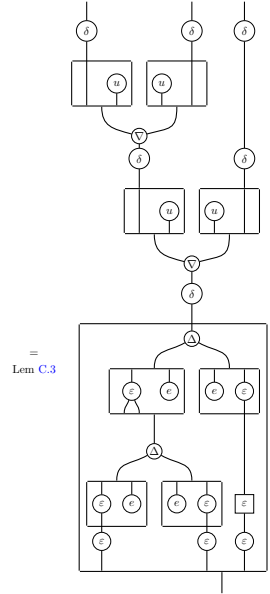
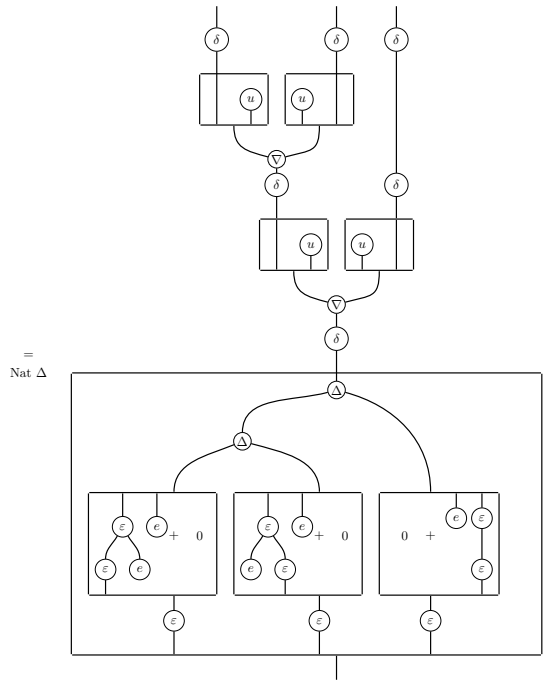
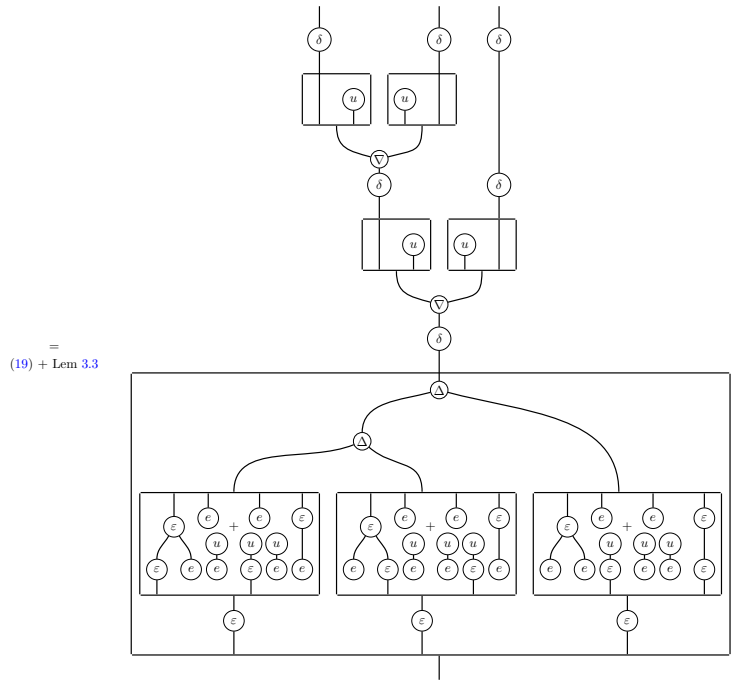
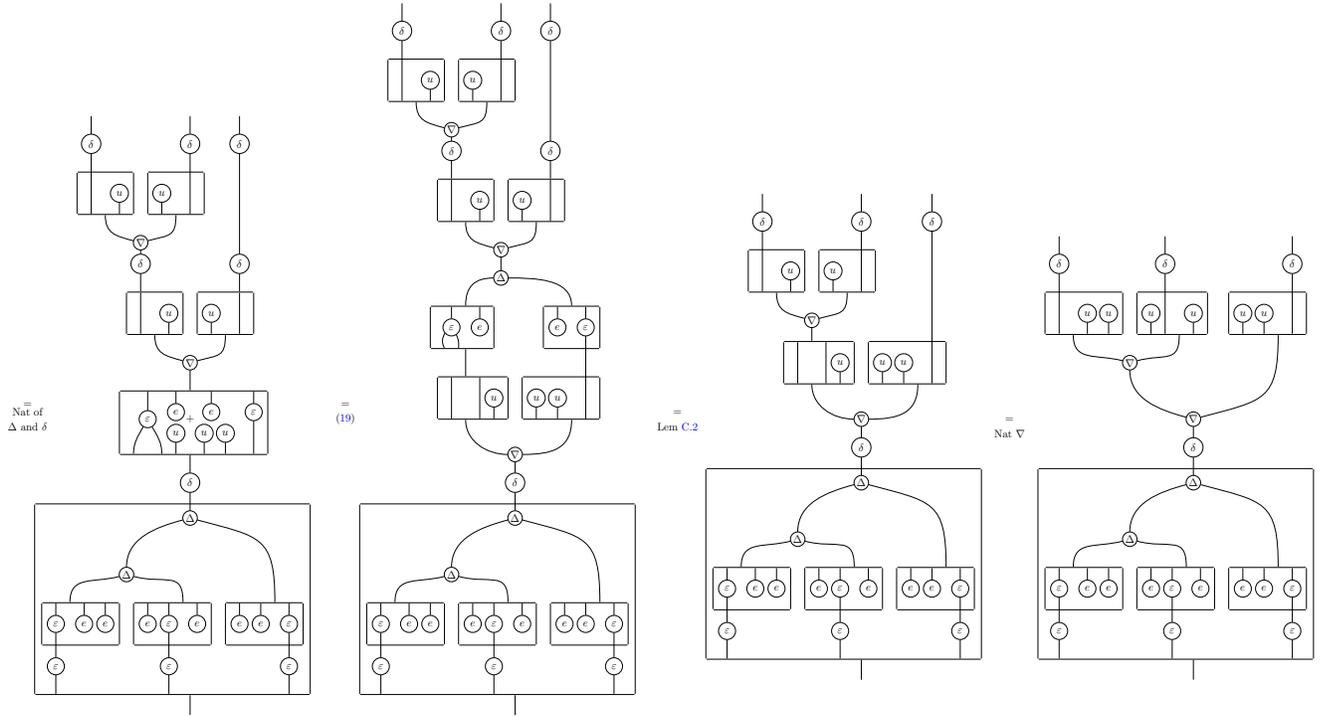
\[\resizebox{!}{10cm}{%
%
}%
\]
The last circuit is symmetric by associativity of the bialgebra, therefore reversing the sequence of equalities by symmetry gives the right associativity of $m_\otimes$. This proves the associativity coherence for $m_\otimes$.  \\

\noindent 2. Unit: We only prove the left unit identity, as the proof for the right unit identity is similar: 
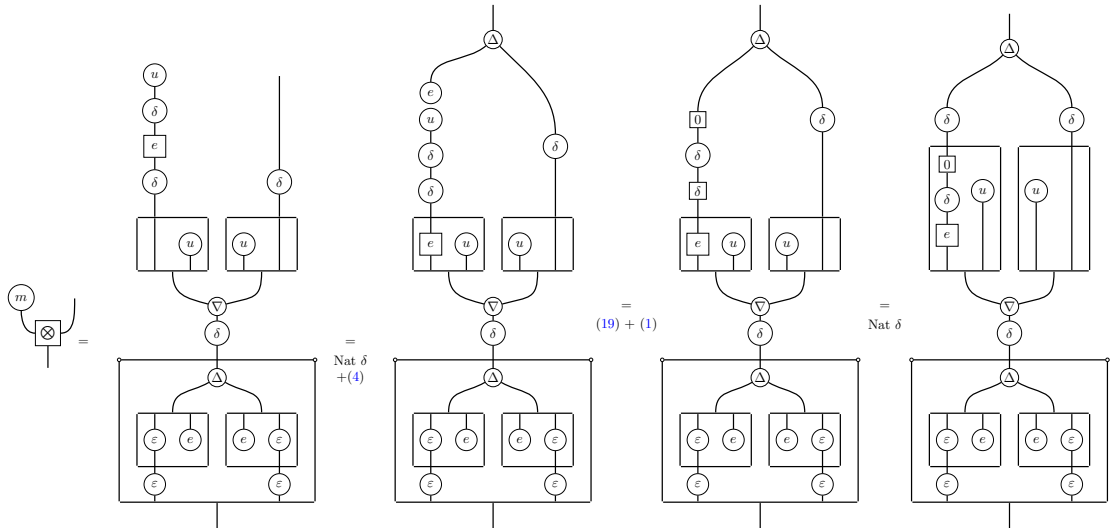
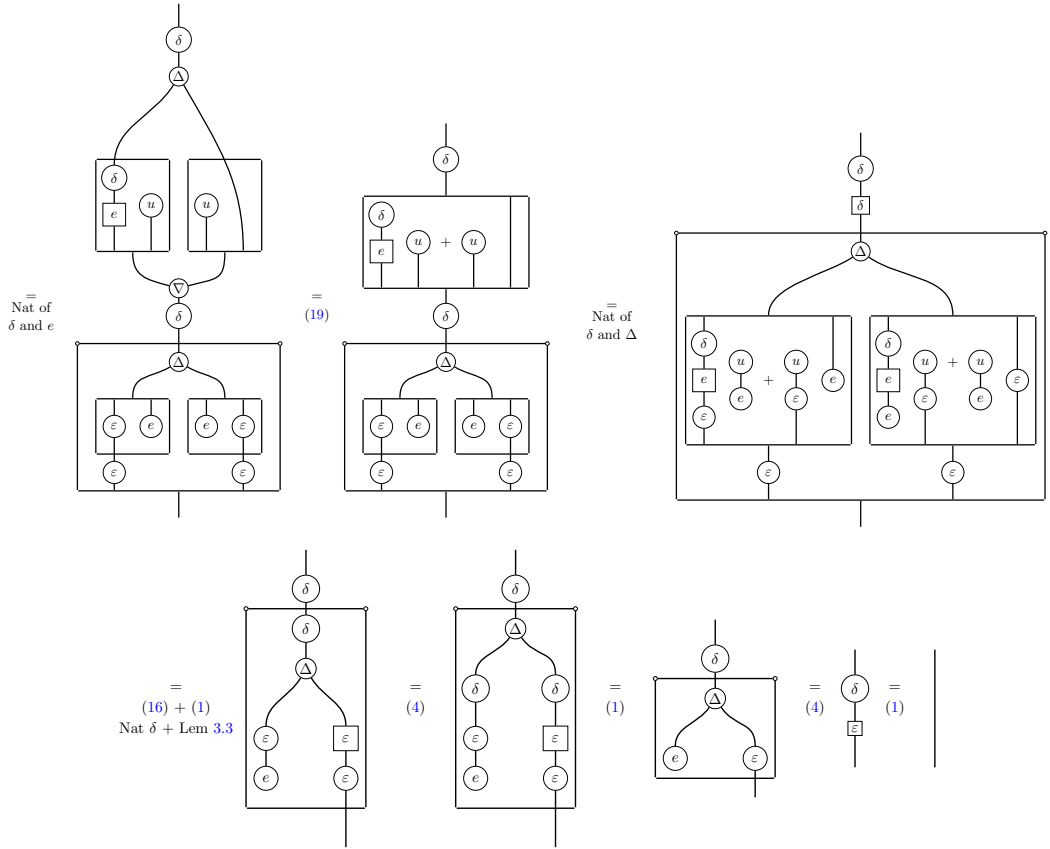
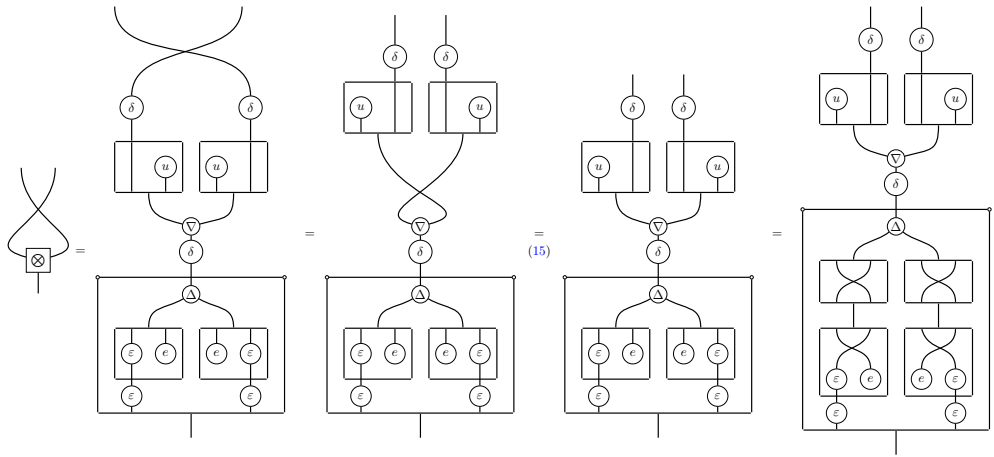
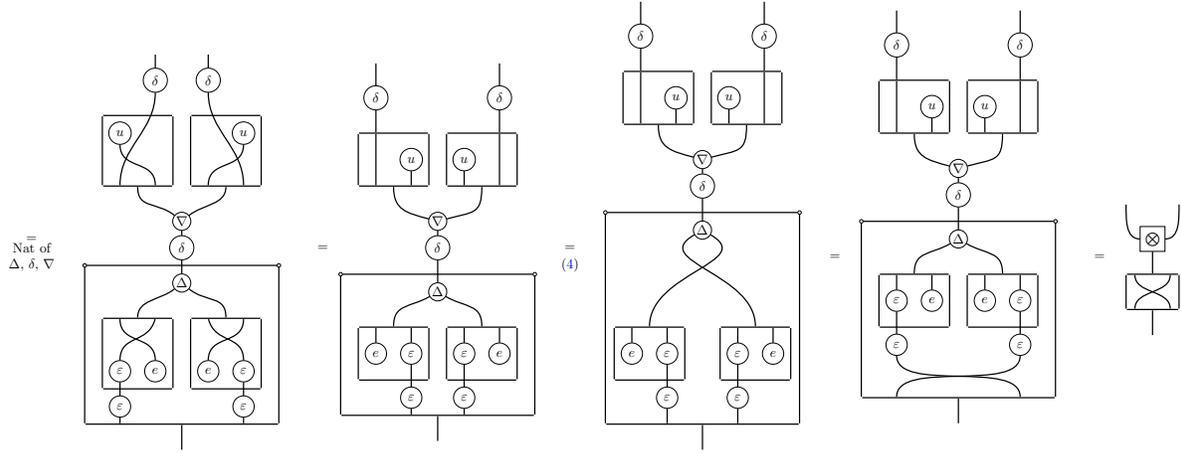
\[\resizebox{!}{7cm}{%
%
}%
\]
\end{proof}

\begin{lemma} $(\oc, \delta, \varepsilon, m_\otimes, m_K)$ is a symmetric monoidal comonad. 
\end{lemma}
\begin{proof} We need to show that $\delta$ and $\varepsilon$ are monoidal natural transformation. \\

\noindent 1. Compatibility between $\delta$ and $m_\otimes$:
\[\resizebox{!}{12cm}{%
%
}%
\]
\end{proof} 

\begin{proposition} Every additive bialgebra modality is a monoidal coalgebra modality. 
\end{proposition} 
\begin{proof} We begin by proving that $\Delta$ and $\varepsilon$ are monoidal transformations. \\

\noindent 1. Compatibility between $\Delta$ and $m_\otimes$
\[\resizebox{!}{6cm}{%
%
}%
\]
Now we prove that $\Delta$ is a $\oc$-coalgebra morphism:
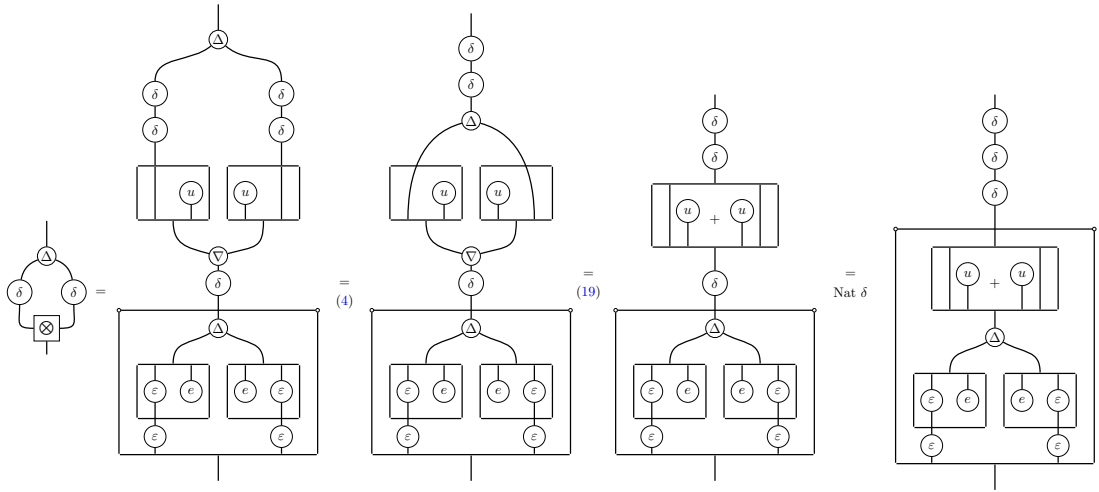
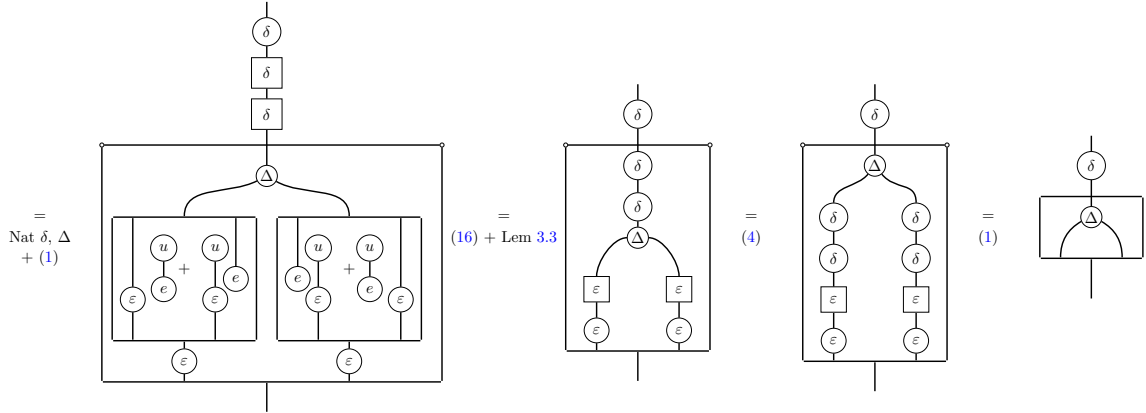
\[\resizebox{!}{6.5cm}{%
%
}%
\]
And finally we prove that $e$ is a $\oc$-coalgebra morphism:
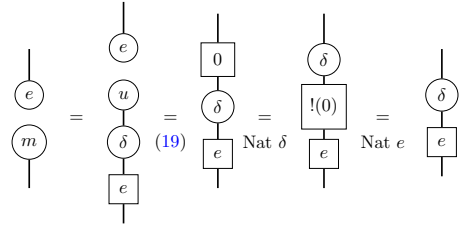
\[\resizebox{!}{3cm}{%
%
}%
\]
\end{proof}

%%%%%%%%%%%%%%%%%%%%%%%%%%%%%%%%%%%%%%%%%%%%%%%%%%%%%%%%%%%%%%%%%%%%%%%
%%%%%%%%%%%%%%%%%%%%%%%%%%%%%%%%%%%%%%%%%%%%%%%%%%%%%%%%%%%%%%%%%%%%%%%

\section{Proof of Theorem \ref{addlinaddbialg}}\label{addmonequiv}

%%%%%%%%%%%%%%%%%%%%%%%%%%%%%%%%%%%%%%%%%%%%%%%%%%%%%%%%%%%%%%%%%%%%%%%

Starting from a monoidal coalgebra modality, we first check that we re-obtain $m_\otimes$: 
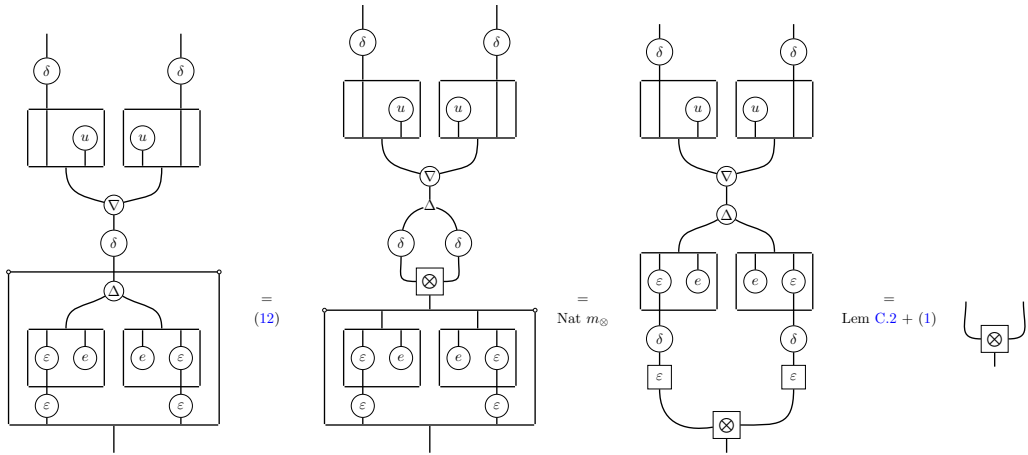
\[\resizebox{!}{6cm}{%
%
}%
\]
Next we check that we get back $m_K$: 
\[\resizebox{!}{3cm}{%
%
}%
\]
Starting with an additive bialgebra modality, we first check that we re-obtain the multiplication:
\[\resizebox{!}{4.25cm}{%
%
}%
\]
Finally we prove that we re-obtain the unit:
\[\resizebox{!}{2.5cm}{%
%
}%
\]

\end{document}